 \documentclass[11pt,reqno,twoside,final]{amsart}

\usepackage[dvips]{graphicx}

\usepackage{amssymb}

\usepackage{mathrsfs}
\usepackage{array}
\usepackage{cite}

\usepackage{pst-plot}

\usepackage[dvips]{hyperref}

\makeatletter

\def\@makechapterhead#1{\global\topskip 7.5pc\relax
  \begingroup
  \fontsize{\@xivpt}{18}\bfseries\centering
    \ifnum\c@secnumdepth>\m@ne
      \leavevmode \hskip-\leftskip
      \rlap{\vbox to\z@{\vss
          \centerline{\Large    \bfseries
                        \@xp{\chaptername}\enspace\thechapter}
          \vskip 3pc}}\hskip\leftskip\fi
     #1\par \endgroup
  \skip@34\p@ \advance\skip@-\normalbaselineskip
  \vskip\skip@ }
\newcommand{\@chapapp}{\chaptername}

\def\section{\@startsection{section}{1}%
  \z@{.7\linespacing\@plus\linespacing}{.5\linespacing}%
  {\Large\normalfont\bfseries\centering}}

\makeatother

\usepackage{fancyhdr}

\def\sn{{\text{sn}}}\def\cn{{\text{cn}}}\def\dn{{\text{dn}}}

\usepackage{multirow}

\begin{document}

\vspace*{-2cm}
\begin{center} { \bf

    \textsf{}\\
\textbf{}\\
       {\Small }

\medskip

  {\Large
Difference operators and difference equations on  lattices, or grids, up to the
elliptic hypergeometric case.
 }}

\bigskip
  Alphonse P. Magnus,    \\
  Institut de  Math\'ematique Pure et Appliqu\'ee,  \\
 Universit\'e catholique de Louvain, \\
 Chemin du Cyclotron,2, \\
  B-1348 Louvain-la-Neuve \\
      (Belgium)\\    

\href{mailto:almagnus@proximus.be}{\tt mailto:almagnus@proximus.be}\ \ \ \ \ \href{https://almagnus.jimdofree.com/}{\tt https://almagnus.jimdofree.com/ }

\medskip


\end{center}

\bigskip

\begin{flushright} {\small
And that has made all the difference.\quad\quad\quad  R. Frost.
}   
\end{flushright}

\medskip

\noindent\textbf{Abstract.} It is shown how to define difference operators and
equations on particular lattices $\{x_n\}$, $2n\in\mathbb{Z}$,
such that
$F(x_n,x_{n\pm 1/2})=0$, where
$F$ s a biquadratic polynomial. We also have the property that
the divided difference operator $(\mathcal{D}f)(x_{n+1/2})=
(f(x_{n+1})-f(x_n))/(x_{n+1}-x_n)$  is a
rational function of degree $\leqslant 2d$ when $f$ is a rational function of
degree $d$. It is then shown that most  known progressions (linear, quadratic, exponential) are covered
by this single theory, up to the case where  the $x_n$s are  values of an elliptic function at a sequence of arguments in
arithmetic progression (\emph{elliptic lattice}). First and second order difference operators and equations are constructed and solved, up to the simplest
elliptic hypergeometric ones. One also shows orthogonality and  biorthogonality properties of rational solutions to some of these difference equations.

Keywords: Elliptic difference operators, Elliptic difference equations, Elliptic hypergeometric expansions, Interpolatory continued fraction, Orthogonal polynomials, Biorthogonal rational functions.

\medskip
2020 Mathematics Subject Classification:

33C05 Classical hypergeometric functions, ${}_2F_1$
33D15 Basic hypergeometric functions in one variable, ${}_r\varphi_s$
33E05 Elliptic functions and integrals
39A06 Linear difference equations
39A13 Difference equations, scaling ($q-$differences)
39A70 Difference operators
41A05 Interpolation
42C05 Orthogonal functions and polynomials



\section{Introduction.}
\chead{\thesection  \ \ \ \   Intro}

\subsection{The subject matter}\   \\

The present study examines particular families of special functions satisfying differential or
difference equations as well as recurrence relations. Known sets of orthogonal polynomials
are followed by new families of biorthogonal rational functions.

An extensive set of apparently unrelated  difference operators (\S \ref{differenceop})  is presented in  a unified theory in section \ref{lattices}. Four equivalent definitions of the relevant
lattices, or grids, $\{x_n\}$ are worked in \S \ref{def1}-\S \ref{def4}.

Properties of the divided difference operator follow in section \ref{divdiff}. For instance, it is shown
in \S \ref{doubling} that the divided difference sends rational functions to rational functions
of double degree or less. One also shows how
this operator acts on appropriate hypergeometric expansion in \S\ref{Dprod}.

Rational interpolatory elliptic expansions follow in section \ref{interp1}, with
particular hypergeometric solutions of simple
linear 1${}^{\text{st}}$ order difference equations in section \ref{lindiffeq}.

Orthogonality and biorthogonality is described  in section \ref{biorths}.
Known special families of orthogonal polynomials are briefly reviewed in \S \ref{LH} and \S \ref{discorthog} , always with
emphasis on the unified analysis, Biorthogonal rational functions,  rational interpolation with free poles\S\ref{biorthinterp} and continued fraction expansion in \S \ref{freepol},
extend the theory of classical orthogonal polynomials to biorthogonal rational functions associated to solutions of particular discrete Riccati equations in \S \ref{Riccatif}.

Difference relations and equations of biorthogonal rational functions extend the related theory of orthogonal polynomials in \S\ref{diffreleq}.

Scalar second order equations are worked in sections \ref{adjointop}, \ref{secondor}, and
 a hypergeometric expansion of a solution of a second order equation is given in \S \ref{elldefL}.




\subsection{ Old and new difference calculus\label{diffsum}}  \   \\

Difference calculus is much older than calculus itself!

The Pythagoreans knew triangular numbers $T_{n+1}-T_n=n+1$, and squares as sums
of odd numbers (\textit{gnomons})  \cite[chap. 4]{Boyer} \cite[\S 3.5]{Kline}   \cite[vol.1, p.161]{BenMenahem}.

Archimedes established  $N^3/3 < 1+2^2+\dotsb + N^2 < (N+1)^3/3$ \cite[\S 17]{BourbakiHist}.

Boole wrote the first treatise on the calculus of finite differences in 1860 \cite{Boole},
even if much had already been done earlier by Newton, Euler, the Bernoullis, or also
by Abel in this identity
\begin{equation}\label{Abelsum}\sum_0^N u_n (v_{n+1}-v_n) = -u_{-1}v_0-\sum_0^N  (u_n-u_{n-1})v_n +u_N v_{N+1},\end{equation}
introduced in 1826 \cite{Abel}, often called "Abel's sum formula", or
partial summation, and which is  a discrete rule of integration by parts   \cite[\S 2.64]{MilneT}.
This is just calculus without infinitesimals, and very elementary!
Why give great names to  such easy identities?  However, the formula
\eqref{Abelsum} appeared incidentally in a much deeper study of convergence of series,
as in recent works \cite{AbelChu1,AbelChu2},
and Boole is mainly concerned by the foundation of the theory of difference equations \cite[chap. VII- chap. X]{Boole}.

\smallskip

\emph{Discrete Calculus} is a field of contemporary applied  mathematics, dealing with graph, or network, theory, algorithmic signal analysis \cite[chap. 1]{DiscreteCalculus}.

\subsection{Difference operators\label{differenceop}}  \   \\

Difference calculus considered initially arithmetic progressions, 
but more instances have been introduced, see Table  \ref{difftable}.

\begin{table}[htbp]
\begin{center}\begin{tabular}{|cc|c|}\hline
 &  $d$  &     $df(x)$ \\
 &$\Delta$ & $f(x+h)-f(x)$ \\ &$\nabla$ & $f(x)-f(x-h)$ \\
  &$\delta$ & $f(x+h/2)-f(x-h/2)$ \\& G & $f(qx)-f(x)$ \\
 &H  &   $f(qx+\omega)-f(x)$\   \\
 &   &         \vphantom{\rule{10pt}{10pt}}                   \\
 \multirow{2}{*}{ ANSUW $\left\{ \vphantom{\rule{25pt}{25pt}}\right.$ }    &W & $f((t+i/2)^2)-f((t-i/2)^2)$, $x=t^2$   \\
 &AW & $f(\cos(\theta+\lambda)) -f(\cos(\theta-\lambda))$, $x=\cos\theta$ \\
&NSU & $f(x(s+1))-f(x(s))$, $x(s)=c_1 q^s +c_2 q^{-s}+c_3$  \\
 &  & \ \ \ \ \ \ \ \ \ \ \ \ \ \ \ \ \ \ or\  $x(s)=c_1 s^2+c_2s+c_3$    \\
 &   &         \vphantom{\rule{10pt}{10pt}}                   \\
 & E & $f(\mathcal{E}(s+1))-f(\mathcal{E}(s))$, where $\mathcal{E}$ is an elliptic function. \\
\hline
\end{tabular}
\caption{Some  difference operators. First column gives the name,
or the context: G is for geometric (why this name? see \S~\ref{geom}) or Heine, Jackson \cite{Jackson},
H for Hahn \cite[eq. (1.3)]{Hahn1949}, W for Wilson \cite[p.34]{AW319}, AW for Askey-Wilson \cite[(5.2)-(5.5)]{AW319} with $q=\exp(2i\lambda)$, NSU, equivalent to AW ($c_1 s^2+\dotsb$ to W), for Nikiforov-Suslov-Uvarov  \cite{NS} \cite[\S 3.4]{NSU}, E for elliptic, Baxter \cite{Baxter}, Rosengren \cite{Rosengren16},
Spiridonov \& Zhedanov \cite{SpZ1,SpZ2,SpZ2007} (to be explained in full detail in the next sections).\label{difftable}} \end{center}\end{table}

Each row is a special case of the next rows.

$\Delta$ is found in Boole 1860 \cite{Boole}, $\nabla$ ({\tt nabla}) in N\"orlund \cite{Norlund},  $\nabla$ and $\delta$ in Steffensen \cite{Steffensen} and  Milne-Thomson \cite{MilneT}.

Thiele \cite{Thiele} uses the notation $\delta(a,b,\dotsc )$ for general divided differences.

We will encounter things as strange as the difference operator $f(x+\sqrt{x})-f(x-\sqrt{x})$, see also \S \ref{ellhist}.

\subsubsection{\label{bisp}}Why just these formulas, why not $f(\log x)$ or $f(x^3)$?

All the particular difference operators of table~\ref{difftable} give rise to
the \emph{same} super class of difference equations solved by special functions
$\phi_m$ satisfying simultaneously

\begin{itemize}
 \item  linear combinations relating several instances $\phi_m(x), \phi_{m\pm 1}(x),\dotsc$, for a given $x$;

 \item a differential or difference equation relating $\phi_m(x), \phi_m(x\pm h),\dotsc$, for a given $m$.
\end{itemize}

These relations have sometimes a neat spectral form, called bispectral form \cite{bispectrencyc},
 as
with Hermite polynomials $\phi_m(x)=H_m(x)$: \cite[\S 22.6,7]{AbrS},  \cite[\S 18.8, 9]{NIST}

  $\phi_{m+1}(x)+2m  \phi_{m-1}(x)=2x  \phi_m(x)$   with eigenvalue $2x$

  $-\phi''_m(x)+2x\phi_m'(x)= 2m\phi_m(x)$   with eigenvalue $2m$

or, to show a difference equation, with Charlier orthogonal polynomials $\phi_m(x)=C_m(x;a)$
\cite[9.14]{Koe}

$a\phi_{m+1}(x)-(m+a)\phi_m(x)+m  \phi_{m-1}(x)=-x  \phi_m(x)$,  \cite[\S 18.22]{NIST},

$a\phi_{m}(x+1)-(x+a)\phi_m(x)+x  \phi_{m}(x-1)=-m  \phi_m(x)$.

See also \cite{Haine,Tsujim}.

The recurrence relations and difference equations to be seen here in \S \ref{approximants}, \S  \ref{freepol} will
normally not be so beautiful!

\subsection{Hypergeometric expansions\label{hyperexp} }  \   \\

The first rows of Table~\ref{difftable}, from $\Delta$ to H, lead
 to hypergeometric expansions as a method to solve difference equations. Later on, clever transformations
of such expansions showed the new difference operators $W$ to $NSU$
\cite{AAR,AW319,Koe,Ko,Segovia,magsnul,NS,NSU}.
Completely new expansions are needed in the E case \cite{ellbib,Rosengren16,spir,
SpZ1,SpZ2,SpZ2007}.

Hypergeometric expansions $F(a,b;c;x)=\sum \dfrac{a(a+1)\dotsb(a+j-1)b(b+1)\dotsb(b+j-1)}{j!c(c+1)\dotsb(c+j-1)}x^j$ (the name goes back to Wallis 1655 \cite[chap. 14]{WW})

Gauss's ratio of hypergeometric functions \cite[ 1913, 1929 \S 59, 1957 \S 24]{Perron}; \cite[ chap. XVIII \S 89]{Wall}.
{\small 
\begin{equation}\label{Gaussf}
f(x)=\dfrac{F(\alpha,\beta;\gamma;x)}{F(\alpha,\beta+1;\gamma+1;x)}=1+\dfrac{a_1 x}{1+ \dfrac{a_2 x}{1+\ddots }}, a_{2m}=\dfrac{-(\beta+m)(\gamma-\alpha+m)}{(\gamma+2m-1)(\gamma+2m)},  a_{2m+1}=\dfrac{-(\alpha+m)(\gamma-\beta+m)}{(\gamma+2m)(\gamma+2m+1)}.\end{equation}
}

This ratio happens to be a solution to the \emph{Riccati} equation

\begin{equation}\label{Gaussr}(1-x)xf'(x)=[(\beta-\alpha)x+ \gamma]f(x)-\gamma f^2(x)+\beta(\alpha/\gamma-1)x.\end{equation} \cite[chap. II, \S 10]{Khov}

\medskip

The present study is NOT based on hypergeometric expansions. Instead
it starts with a very elementary construction,   and hypergeometric
expansions will be recovered at the very end.

\subsection{Orthogonality, Pad\'e approximation,  biorthogonality, and rational interpolation\label{Padinterp}}     \     \\

Many special functions $\phi_m(x)$ that will follow satisfy recurrence relations typical
of \emph{orthogonal polynomials}, with or without the bispectrality property.

Relation to Pad\'e\footnote{Henri Eug\`ene Pad\'e, 1863-1953,  \cite[p. 249]{BrezinskiH}.} approximation:
if $\displaystyle \phi_m(t)=\sum_0^m \phi_m^{(j)} t^j$ is orthogonal to all polynomials of degree \\   $<m$
with respect to a (formal) bilinear form $\langle f,g\rangle = \mathscr{L}(fg)$, one has
$\mathscr{L}(t^r\phi_m(t))=\sum_0^m \phi_m^{(j)} \mathscr{L}(t^{j+r})=0$ for $r=0,1,\dotsc, m-1$, showing that
the product of $\phi_m(x)$ and $f(x)=\sum \mathscr{L}(t^{r})x^{-r-1}=\mathscr{L}((x-t)^{-1})$ has a Laurent expansion
with vanishing coefficients of $x^{-1},\dotsc,x^{-m}$, the numerator polynomial is made of
the nonnegative powers
\cite[chap. 7]{Baker} \cite[\S 5.5]{JonesT}

 Example of Hermite polynomials (actually $i^{-m}H_m(ix)$) in the error function \cite[\S 7.1.14, 15]{AbrS},

 $\displaystyle f(x)=2e^{x^2}\int_x^\infty e^{-t^2}dt =  \dfrac{i}{\sqrt{\pi}}
\int_{-\infty}^\infty \dfrac{e^{-t^2}\, dt}{ix-t}  \sim \dfrac{1}{x}- \dfrac{1/2}{x^3}+ \dfrac{3/4}{x^5}-\dfrac{15/8}{x^7}+ \dfrac{105/16}{x^9} -\dotsb \\    = \frac{1}{x+\frac{1/2}{x+\frac{2/2}{x+\frac{3/2}{x+\ddots}}}}$ \cite[13.2.20]{Cuyt}. Here,  $\displaystyle \mathscr{L}(fg)=\int_{-\infty}^\infty e^{-t^2}f(t)g(t)\, dt$.

Other example,  the exponential integral and Laguerre  (actually $L_m(-x)$) \cite[\S 5]{AbrS},

 $\displaystyle e^xE_1(x)=e^x\int_x^\infty e^{-t}dt/t
= \int_0^\infty \dfrac{e^{-t}dt}{x+t}$ \cite[\S 5.1.28]{AbrS},
$\displaystyle =  \frac{1}{x+\frac{1}{1+\frac{1}{x+\frac{2}{1+\frac{2}{x+\ddots}}}}}$ \cite[\S 5.1.22]{AbrS},   \cite[14.1.16]{Cuyt}
$\displaystyle \sim  \sum_0^\infty \dfrac{(-1)^m m!}{ x^{m+1} } $ \cite[5.1.51]{AbrS}. Here,  $\mathscr{L}(fg)=\int_0^\infty e^{-t}f(t)g(t)dt$.

\begin{table}[htbp]
\begin{center}\begin{tabular}{|c|c|}\hline
$\dfrac{1}{x+1}$  &    $\dfrac{1}{x}-\dfrac{1}{x^2}+\dfrac{1}{x^3}-\dfrac{1}{x^4}+\dfrac{1}{x^5}-\dfrac{1}{x^6}+\dfrac{1}{x^7}-\dfrac{1}{x^8}+\dfrac{1}{x^9}-\dotsb$    \\  \hline
 $\dfrac{x+3}{x^2+4x+2}$& $\dfrac{1}{x}-\dfrac{1}{x^2}+\dfrac{2}{x^3}-\dfrac{6}{x^4}+\dfrac{20}{x^5}-\dfrac{68}{x^6}+\dfrac{232}{x^7}-\dfrac{792}{x^8}+\dfrac{2704}{x^9}-\dotsb$
 \\ \hline
$\dfrac{x^2+8x+11}{x^3+9x^2+18x+6}$ &$\dfrac{1}{x}-\dfrac{1}{x^2}+\dfrac{2}{x^3}-\dfrac{6}{x^4}+\dfrac{24}{x^5}-\dfrac{120}{x^6}+\dfrac{684}{x^7}-\dfrac{4140}{x^8}+\dfrac{25668}{x^9}-\dotsb  $ \\ \hline
$\dfrac{x^3+15x^2+58x+50}{x^4+16x^3+72x^2+96x+24}$ &$\dfrac{1}{x}-\dfrac{1}{x^2}+\dfrac{2}{x^3}-\dfrac{6}{x^4}+\dfrac{24}{x^5}-\dfrac{120}{x^6}+\dfrac{720}{x^7}-\dfrac{5040}{x^8}+\dfrac{39744}{x^9}-\dotsb$   \\ \hline
\end{tabular}
\caption{Pad\'e effect of  Laguerre polynomials. Denominator of degree $m$ leads to matching $2m$ power coefficients. The numerator
polynomials satisfy the same recurrence relations.
\label{Hermite}} \end{center}\end{table}

See  \cite[\S 4.5]{BrezinskiH} for Pad\'e history,  \cite[\S 5.2.3]{BrezinskiH} for orthogonal polynomials , see also \cite{deBruin}
for a wonderful summary of Pad\'e and orthogonal  polynomials activity during the 1970s and 1980s.

We will also encounter the following case of \emph{biorthogonality}:

if the rational function $\phi_m(x)=\dfrac{P_m(x)}{(x-a_0)\dotsb(x-a_m)}$ is orthogonal to
rational functions of the form $\dfrac{Q_n(x)}{(x-b_0)\dotsb(x-b_n)}, n\neq m$, i.e. if $P_m$
is orthogonal to any polynomial of lower degree with respect to $\displaystyle \mathscr{L}\left(\dfrac{f(t)g(t)}{(t-a_0)\dotsb(t-a_m)(t-b_0)\dotsb(t-b_{m-1}) }\right)$ leading to

$\displaystyle P_m(x)x^rf(x)= P_m(x)x^r\mathscr{L}\left(\dfrac{1}{x-t}\right)dt=
\text{\ polynomial\ of\ degree\ }m+r < 2m +\mathscr{L}\left(\dfrac{t^rP_m(t)}{x-t}\right)dt$
and we perform the divided difference of order $2m$ at $x=a_0,\dotsc,a_m,b_0,\dotsc,b_{m-1}$ giving

$\displaystyle \mathscr{L}\left(\dfrac{t^rP_m(t)}{(t-a_0)\dotsb(t-a_m)(t-b_0)\dotsb(t-b_{m-1}) }\right)=0$ from the biorthogonality conditions, meaning that there is a polynomial $N_m$ interpolating  $ P_m(x)f(x)$ at the given $2m+1$ points, and that $N_m$ has degree $\leqslant m$
, as $x^rN_m(x)$ has a vanishing divided difference of order $2m$.

Interpolation by  $N_m/P_m$ is called  multipoint Pad\'e approximation in   \cite{Gon5,Lo1,Lo2},  \cite[\S 5.5.2]{JonesT}

Example of  Psi function  \cite{Norlund}  $f(x)=\Psi(x)-\Psi(1)= 0, 1, 3/2, 11/6, 25/12, $ at 1, 2, 3...

$f(x)=\dfrac{x-1}{1+\dfrac{(x-2)/3}{1 +\dfrac{\ddots }{ 1+\dfrac{(m-1)(x-2m+1)/(9m^2-9m+3)}{1+ \dfrac{m(x-2m)/(9m^2-9m+3)}{ 1+\ddots }  }      }}}$.


Successive rational interpolants are $\dfrac{3x-3}{x+1}, \dfrac{9x^2+21x-30}{2x^2+18x+4},
\dfrac{11x^3+162x^2+37x-210}{ 2x^3+60x^2+166x+12}$. The curious continued fraction above (see also  \cite[\S 244] {Norlund})
will be explained in \S\ref{Exxh}.

See E. Rains \cite{Rainsg} for a very extensive highbrow relation to algebraic geometry.

\section{The lattices.\label{lattices}}
\chead{\thesection   \ \ \ \  The lattices.}

\subsection{   Curves   }\   \\

Simplest difference equations relate two values of the unknown
function $f$, for
 instance $f(x)$ and $f(x+h)$, or the more
symmetric $f(x-h/2), f(x+h/2)$, or also $f(x), f(qx)$ in $q-$difference
equations \cite{GrH,Koe,Ko}.


Some of these difference operators allow to formulate, and solve, curious difference equations
such as: \textsl{find $f$ such that $f(x+\sqrt{x})-f(x-\sqrt{x})=1$, starting from
$f(0)=0$.} Answer: $f(n^2+n)=n, n=0,1,\dotsc$, see the last exercise of \S \ref{exer1}. Or also: \textsl{find $f$ such that
$f(x+x/2)-f(x-x/2)=1, with \\   f(1)=1$}.Answer: $f(3^n)=n+1, n=0,1,\dotsc$
Solution of difference equations asks for values of the solution on particular
lattices, or grids, of various kinds given in \S \ref{latticeformulas}.

\medskip

From W onward, arguments of $f$ are given in parametric form, suggesting
a
CURVE. A Cartesian form of the W case is $y=t^2 \pm it -1/4 \Rightarrow
(y-x+1/4)^2+x=0$, a  parabola.


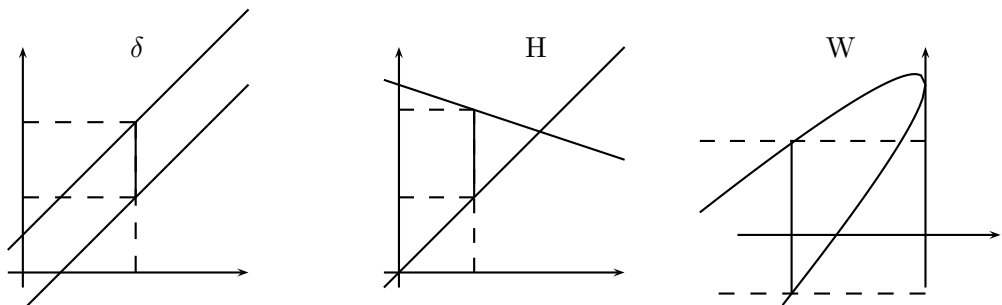
\begin{figure}[htbp]\begin{center}
\psset{xunit=1cm,yunit=1cm}
\begin{pspicture}(-0.2,-0.2)(13,4)
\psline{->}(-0.2,0)(3,0)\psline{->}(0,-0.2)(0,3)
\psline(-0.2,0.3)(3,3.5)\psline(-0.0,-0.5)(3,2.5)
\psline[linestyle=dashed,dash=2mm 2mm](1.5,0)(1.5,2)
\psline[linestyle=dashed,dash=2mm 2mm](1.5,2)(0,2)
\psline[linestyle=dashed,dash=2mm 2mm](1.5,1)(0,1)
\psline(1.5,1)(1.5,2)
\uput[0](1.25,3){$\delta$}
\psline{->}(4.8,0)(8,0)\psline{->}(5,-0.2)(5,3)
\psline(4.8,-0.2)(8,3)\psline(4.8,2.5633)(8,1.5)
\psline[linestyle=dashed,dash=2mm 2mm](6.0,0)(6.0,2.166)
\psline[linestyle=dashed,dash=2mm 2mm](5,2.166)(6.0,2.166)
\psline[linestyle=dashed,dash=2mm 2mm](5,1.0)(6.0,1.0)
\psline(6,1)(6,2.1666)
\uput[0](6.5,3){H}
\psline{->}(9.5,0.5)(13,0.5)\psline{->}(12,-0.2)(12,3)
\psplot{9}{12}{12 x  sub sqrt 0.75 mul x 9.5 sub  add}
\psplot{10.1}{12}{12 x  sub sqrt 0.75 mul neg x 9.5 sub  add}
\psline(10.22,-0.28)(10.22,1.75)
\psline[linestyle=dashed,dash=2mm 2mm](9,1.75)(12,1.75)
\psline[linestyle=dashed,dash=2mm 2mm](9.25,-0.28)(12,-0.28)
\uput[0](10.5,3){W}
\end{pspicture}
\caption{Examples of values involved in the solution of a
difference equation. Why are the axes drawn without variables names?}\label{difffig0}
 \end{center}\end{figure}

We show that the difference operators seen up to now
 involve $f$ at the two
roots $y_1$ and $y_2$ of a quadratic polynomial in $y$. Indeed, we
already saw $(y-x+1/4)^2+x=0$ for the parabola of the third
example of Fig.~\ref{difffig0}. The two easier examples are
$(y-x-h/2)(y-x+h/2)=0$ and $(y-x)(y-qx-\omega)=0$ and do again enter
the quadratic scheme. For the more difficult AW case,
$y_1, y_2= \cos\theta\cos\lambda \pm \sin\theta\sin\lambda \Rightarrow
(y-x\cos\lambda)^2+(x^2-1)\sin^2\lambda=0$, an ellipse if $\lambda$ is real,
a hyperbola if $\lambda$ is pure imaginary, showing how
 conics can be discovered from other
rows of table~\ref{difftable}.

For NSU, where is the curve?
Actually, the NSU formula gives immediately the relevant lattice,
or grid, $x(s), x(s+1), x(s+2),\dotsc$ The next sections follow
with more on lattices and recurrence relations.

\medskip


\subsection{Difference equations and lattices}  \   \\

A first order difference equation relates the
values of the unknown function $f$ at two ordinates on the given curve, say
$y_1$ and $y_2$. If $f(y_1)$ is given (Cauchy problem), the
equation gives $f(y_2)$, and we search another point on the curve of ordinate
$y_2$, the difference equation gives then $f$ at a new point, etc. The solution
of the equation is then described by the values at a sequence of values
$y_1, y_2, y_3,\dotsc$ This rather ridiculous choice of the name $y_k$
for the fundamental lattice, or grid, of arguments of the function $f$
(used by the author in\cite{MagLuminy,magsigma2009})
is avoided by (re)turning to the $x-$axis containing the relevant
$x_0, x_1,\dotsc$ values and the $y-$axis contains now intermediate data \cite[eq. (1.2)]{SpZ2007}. One could consider that the pictures
of fig.~\ref{difffig0} have a horizontal axis $y$ and a vertical
axis $x$.


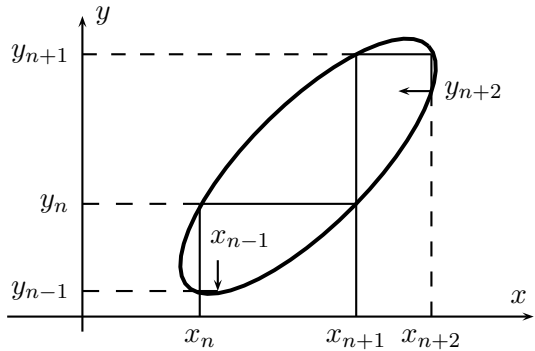
\begin{figure}[htbp]\begin{center}
\psset{xunit=1cm,yunit=1cm}
\begin{pspicture}(2,-0.2)(11,4.05)
\psline{->}(3,0)(10,0)\psline{->}(4,-0.2)(4,4)
\uput[90](9.8,0){$x$}\uput[0](4,4){$y$}
\parametricplot[linewidth=1.5pt]{0}{360}{1.6 t cos mul  0.57 t sin mul  sub 7 add 1.6 t cos mul 0.57 t sin mul add 2 add}
\psline[linestyle=dashed,dash=2mm 2mm](4,1.5)(5.56,1.5)
\psline(5.56,1.5)(7.64,1.5)\psline(5.56,0)(5.56,1.5)
\psline[linestyle=dashed,dash=2mm 2mm](4,3.49)(7.64,3.49)\psline(7.64,0)(7.64,3.49)\psline(7.64,3.49)(8.64,3.49)
\psline(8.64,3.49)(8.64,3.0)\psline[linestyle=dashed,dash=2mm 2mm](8.64,3.0)(8.64,0)
\psline{->}(8.64,3.0)(8.2,3.0)\uput[0](8.64,3.0){$y_{n+2}$}
\psline[linestyle=dashed,dash=2mm 2mm](4,0.34)(5.56,0.34)
\psline(5.56,0.34)(5.8,0.34)\psline{->}(5.8,0.75)(5.8,0.34)
\uput[90](6.1,0.72){$x_{n-1}$}
\uput[270](5.56,0){$x_n$}\uput[270](7.64,0){$x_{n+1}$}
\uput[270](8.64,0){$x_{n+2}$}
\uput[180](4,1.5){$y_{n}$}\uput[180](4,0.34){$y_{n-1}$}
\uput[180](4,3.49){$y_{n+1}$}
 \end{pspicture}
\caption{Some points involved in the solution of a
difference equation.}\label{difffig1}
 \end{center}\end{figure}

\subsubsection{The biquadratic polynomial\label{biquad}.}So, let $(x_n,y_n)$ be a current point of a relevant sequence on the curve
$F(x,y)=0$. Then,
$x_{n+1}$ is the second $x-$ root of $F(x,y_n)=0$ (Fig.~\ref{difffig1}).
To be sure that the equation $F(x,y)=0$ has always exactly two roots in $x$ for
a given $y$, we choose $F$ to be a quadratic polynomial in $x$:

\begin{subequations}
\begin{equation}\label{Fx}
   \hspace*{65pt}\ \ \ \ \  F(x,y)=Y_0(y)+Y_1(y) x +Y_2(y) x^2 =0,
\end{equation}
\noindent where $Y_0, Y_1$, and $Y_2$ are given functions.

\

 But difference equations must allow the recovery of $f$ on a whole set
 of points! An initial-value problem for a first order difference equation
 starts with a value for $f(x)$ at some $x=x_0$. We then get $y_0$ as a $y-$root of $F(x_0,y)=0$. The difference equation  relates then
 $f(x_0)$ to $f(x_1)$, where $x_1$ is the second root of \eqref{Fx} at $y=y_0$, found by \eqref{sumprod} as $x_1=-Y_1(y_0)/Y_2(y_0)-x_0$.
 We now need $y_1$,
 the $y-$root of $F(x_1,y)=0$ which is not $y_0$. Here again, the
simplest case is when $F$ is of degree 2 in $y$:
 \begin{equation}\label{Fy}
   F(x,y)=X_0(x)+X_1(x) y +X_2(x) y^2 =0.
\end{equation}
\end{subequations}
   \   \\
Both forms \eqref{Fx} and \eqref{Fy} hold simultaneously
when $F$ is \emph{\textbf{biquadratic}}:
\begin{equation}\label{Fbiq}
 F(x,y)= \sum_{i=0}^2 \sum_{j=0}^2 c_{i,j} x^i y^j.
\end{equation}

Remark that the sum and the product of the two roots  are
the quadratic rational functions
\begin{equation}\label{sumprod}\begin{split}
 x_n+x_{n+1} &=-Y_1(y_n)/Y_2(y_n),\ \ \ x_nx_{n+1}=Y_0(y_n)/Y_2(y_n), \\
y_{n-1}+y_{n} &=-X_1(x_n)/X_2(x_n),\ \ \ y_{n-1}y_{n}=X_0(x_n)/X_2(x_n).
\end{split}\end{equation}

We also have the two factorizations of $F$

At some $x=x_n$, the two $y-$roots of $F(x_n,y)=0$ are $y_{n-1}$ and
$y_n$,

so $F(x_n,y)=X_2(x_n)(y-y_{n-1})(y-y_n)$; also,
$F(x,y_n)=Y_2(y_n)(x-x_n)(x-x_{n+1})$.

\begin{equation}\label{factorF}\begin{split}
 F(x_m,y_n) &= X_2(x_m)y_n^2+X_1(x_m)y_n+X_0(x_m)=X_2(x_m)(y_n-y_{m-1})(y_n-y_m)\\
           &=  Y_2(y_n)x_m^2+Y_1(y_n)x_m+Y_0(y_n)=Y_2(y_n)(x_m-x_{n})(x_m-x_{n+1}).
\end{split}\end{equation}

\smallskip\begin{flushright}{\small
\textsl{ Il e\^ut \'et\'e si facile ... de faire pr\'ec\'eder et suivre chaque proposition
d'un cort\`ege redoutable d'exemples particuliers! }

(It would have been so easy ... to precede and to follow each statement by a
formidable parade of particular examples!)

\'E. Galois, who hated losing time on examples and exercises.
}\end{flushright}
\smallskip

\subsubsection{\textbf{Exercises.}\label{exer1}} Show that $\dfrac{x_{n+1}-x_{n-1}} {y_n-y_{n-1}  }$
is a rational function of $x_n$. Hint: use $x_{n+1}-x_{n-1}= x_n+x_{n+1}
-(x_{n-1}+x_n)$ and use \eqref{sumprod}. This curious ratio will be used in \S\ref{exer2}  and \ref{elllogar}.

\medskip

Remove the indetermination at $n=0$ of $\dfrac{y_n-y_0}{x_n-x_0}$. Answer: multiply
the two terms of the fraction by $x_{n+1}-x_0$. The denominator is now \\
$(x_n-x_0)(x_{n+1}-x_0)= F(x_0,y_n)/Y_2(y_n)= X_2(x_0)(y_n-y_0)(y_n-y_{-1})/Y_2(y_n)$
from \eqref{factorF}, and the result is $\dfrac{ X_2(x_0)(x_{n+1}-x_0)}{Y_2(y_n)(y_n-y_{-1})}$. Note that the indetermination is now at $n=-1$...

\medskip

How can we relate this to the difference operator $f(x+\sqrt{x})-f(x-\sqrt{x})$?

$x$ is a parameter not to be confused with the $x$ of $F(x,y)$, so, there is
a parameter $s_n$ such that  $x_n=s_n-\sqrt{s_n}$, and $x_{n+1}=s_n+\sqrt{s_n}$,
$x_n+x_{n+1}= 2s_n$ and $x_nx_{n+1}=s_n^2-s_n$. By \eqref{sumprod},
$Y_0(y_n)/Y_2(y_n)=Y_1^2(y_n)/(4Y_2^2(y_n))+Y_1(y_n)/(2Y_2(y_n))$. Let $Y_2(y)\equiv 1$,
and the degrees of $Y_0$ and $Y_1$ be 2 and 1, with $Y_1(y)=2\alpha y+2\beta$, we have
$Y_0(y)= (\alpha y+\beta)^2+\alpha y+\beta$, $F(x,y)=(\alpha y+\beta)^2+\alpha y+\beta
+2(\alpha y+\beta)x+x^2=(x+\alpha y+\beta+1/2)^2-x-1/4=0$ at $x=x_n$ when
$-s_n=\alpha y_n+\beta  =-\sqrt{x_n+1/4}-x_n-1/2= -(1/2+\sqrt{x_n+1/4})^2 $,

$x_{n+1}= 2s_n-x_n= x_n+1+2\sqrt{x_n+1/4},  \sqrt{x_{n+1}+1/4}= 1+\sqrt{x_n+1/4}$,

the sequence is $x_n=(n+c)^2-1/4$, $c=\sqrt{x_0+1/4}$, the Wilson (W) lattice, see below.

\subsection{Lattices formulas\label{latticeformulas}} \  \\

All the examples already seen enter the \eqref{Fbiq} form.

Besides $f(x_{n+1})-f(x_n)=0$, the most ridiculously simple difference
equation is $f(x_{n+1})-f(x_n)=1$, of course solved by $f(x_n)=f(x_0)+n$.
The interesting feature here is not given by the values of $f$, but by the
points $x_n$ where $f$ can be computed. Examples follow.

Most people do not use this $(x,y)$ system, but a particular parametric representation
such that $x_n=x(s_0+n)$. The simplest example is $x(s)=sh$.
The formulas for the various curves $F(x,y)=0$ follow now.

\smallskip

\begin{enumerate}

\item  \label{arith} Two parallel lines, first entries of table~\ref{difftable},
$F(x,y)=(y-ax-b)(y-ax-c)$, we start with $x_0$ and
 $y_0=ax_0+b$,  then, $x_1$ is such that $ax_1+c
=y_0=ax_0+b$ so $x_n=x_0+nh, y_n= y_0+nah$ follow, where $h=(b-c)/a$.

\item  \label{geom}  Two lines $F(x,y)=(y-ax-b)(y-cx-d)$, $a\neq c$, then
$F(x,y_n)=0$ is solved by $x_n=(y_n-b)/a$ and $x_{n+1}=(y_n-d)/c$,
so, $cx_{n+1}-ax_n=b-d$, or $x_{n+1}-x^*=q(x_n-x^*)$ with $q=a/c$
and $x^*=(b-d)/(c-a)$,
whence $x_n=x^*+q^n(x_0-x^*), y_n=y^*+q^n(y_0-y^*)$,
where $y^*=ax^*+b= (bc-ad)/(c-a)$, a \emph{geometric} sequence, rows G and H of table~\ref{difftable}, whence the name.

Note that $q$ is the ratio of the slopes
of the two lines.

When $a$ and $c$ are very close, $a=c(1+\varepsilon),
q=1+\varepsilon, x^*=(d-b)/(c\varepsilon), \\
x_n= (1+\varepsilon)^n x_0 +((d-b)/c)[1-(1+\varepsilon)^n]/\varepsilon \to x_0-n(d-b)/c$ when $\varepsilon\to 0$.

\item  \label{parab} Parabola\footnote{Parametric description of a parabola used in mathematical
typography, the {\tt TrueType} system \cite{Devroye,TrueType}.}, line W of table~\ref{difftable},
  $x_n=as_n^2+bs_n+c, y_n=ds_n^2+es_n+f$,

 or
$F(x,y)=x-a\left( \dfrac{ay-dx-\gamma}{ae-bd}\right)^2-b\dfrac{ay-dx-\gamma}{e-bd}-c$,  where
$\gamma=af-cd$. The two $s-$roots of $y=y_n$ are $s_n$ and $-e/d-s_n$; the two $s-$roots
of $x=x_{n+1}$ are $-e/d-s_n$ and $s_{n+1}=-b/a-(-e/d-s_n)$, so
$s_n=s_0+n(ae-bd)/(ad)$, and $x_n$ and $y_n$ are quadratic polynomials in $n$, that is
the second NSU-line of table~\ref{difftable} .

  \item \label{conic} General centered conic $x_n=a\cos(\theta_n-\theta^*)+b,
 y_n=c\cos\theta_n+d$,

 or $F(x,y)= \left( \dfrac{x-b}{a\sin\theta^*}-\dfrac{y-d}{c\tan\theta^*}\right)^2+\left(\dfrac{y-d}{c}\right)^2-1$,
where the $\theta$s need not be real if the $x$ and $y$s are real.
Then $x_{n+1}$ is the second $x-$root of $F(x,y_n)=0$, so,
$x_{n+1}= b+a\sin\theta^*( \cos\theta_n/\tan\theta^* -\sin\theta_n)$, and
$\theta_{n+1}= \theta_n+2\theta^*$ follows, or
$y_n=d+c\cos(\theta_0+2n\theta^*)=d+c(e^{i\theta_0}q^n+e^{-i\theta_0}q^{-n})/2,
 x_n=b+a\cos(\theta_0+(2n-1)\theta^*)=b+a(e^{i\theta_0}q^{n-1/2}+e^{-i\theta_0}q^{1/2-n})/2$,
with $q=\exp(2i\theta^*)$,
basically the NSU formulas of table~\ref{difftable}.

Note that $q$ is now the ratio of the slopes
of the two asymptotes \cite[p. 264]{Segovia} \cite[p. 255]{magsnul}.

\item The general biquadratic polynomial will of course be the subject matter
of all the following sections.

\end{enumerate}






\begin{figure}[htbp]\begin{center}

\psset{xunit=1.0cm,yunit=1.0cm}
\begin{pspicture*}(-7,-4.5)(6,5)
\psline{->}(-6,0)(5,0)\psline{->}(0,-3.25)(0,4)
\psline[linecolor=red](-10.37,7.384)(3,7.384)(3,1.459)(1.651,1.459)(1.651,1.082)(1.51,1.082)(1.51,1.16)(1.931,1.16)(1.931,2.07)(6.365,2.07)(6.365,2.07)
\uput[90](5.4,0){$x$}\uput[0](0,3.9){$y$}
\psline(-1,-0.05)(-1,0.05)\psline(1,-0.05)(1,0.05)\uput[270](-1,0){$-1$}\uput[270](1,0){1}\psline(-2,-0.05)(-2,0.05)\psline(2,-0.05)(2,0.05)\uput[270](-2,0){$-2$}\uput[270](2,0){2}
\psline(-3,-0.05)(-3,0.05)\psline(3,-0.05)(3,0.05)\uput[270](-3,0){$-3$}\uput[270](3,0){3}
\psline(-4,-0.05)(-4,0.05)\psline(4,-0.05)(4,0.05)\uput[270](-4,0){$-4$}\uput[315](4,0){4}
\psline(-5,-0.05)(-5,0.05)\uput[270](-5,0){$-5$}
\psline(-6,-0.05)(-6,0.05)\uput[270](-6,0){$-6$}
\psline(-0.05,-1)(0.05,-1)\psline(-0.05,1)(0.05,1)\uput[180](0,-1){$-1$}\uput[180](0,1){1}\psline(-0.05,-2)(0.05,-2)\psline(-0.05,2)(0.05,2)\uput[180](0,-2){$-2$}\uput[180](0,2){2}
\psline(-0.05,-3)(0.05,-3)\psline(-0.05,3)(0.05,3)\uput[180](0,-3){$-2$}\uput[180](0,3){3}\psline(-0.05,-4)(0.05,-4)\psline(-0.05,4)(0.05,4)\uput[180](0,-4){$-4$}\uput[180](0,4){4}
\psline[linestyle=dashed,dash=2mm 2mm](4.000,-0.5)(4.000,4)
\psline[linestyle=dashed,dash=2mm 2mm](-5.500,-3.5)(-5.500,0.5)
\psline[linestyle=dashed,dash=2mm 2mm](-5,3.617)(5,3.617)
\psline[linestyle=dashed,dash=2mm 2mm](-5,-2.884)(5,-2.884)
\pscurve[linewidth=0.075](-2.500,-2.437)(-2.400,-2.233)(-2.300,-2.038)(-2.200,-1.849)(-2.100,-1.660)(-2.000,-1.397)(-2.000,-1.393)(-2.100,-1.312)(-2.200,-1.313)(-2.300,-1.325)(-2.400,-1.340)(-2.500,-1.357)(-2.734,-1.400)(-3.279,-1.500)(-3.8128,-1.5897)(-5.0840,-1.7660)(-5.8843,-1.8551)(-6.7673,-1.9385)
\pscurve[linewidth=0.075](-1.000,-0.6434)(-0.9000,-0.6793)(-0.8000,-0.6566)(-0.7000,-0.6236)(-0.6000,-0.5848)(-0.5000,-0.5415)(-0.4000,-0.4944)(-0.3000,-0.4438)(-0.2000,-0.3899)(-0.1000,-0.3325)(-0.00001000,-0.2716)(0.09999,-0.2069)(0.2000,-0.1384)(0.3000,-0.06552)(0.4000,0.01193)(0.5000,0.09449)(0.6000,0.1828)(0.7000,0.2778)(0.8000,0.3812)(0.9000,0.4969)(1.000,0.6674)\pscurve[linewidth=0.075](1.500,1.110)(1.600,1.072)(1.700,1.094)(1.800,1.121)(1.900,1.151)(2.000,1.181)(2.100,1.211)(2.200,1.240)(2.300,1.270)(2.400,1.298)(2.500,1.327)(2.770,1.400)(3.165,1.500)(3.597,1.600)(4.074,1.700)(4.602,1.800)(6.365,2.07)
\pscurve[linewidth=0.075](-1.000,-0.6414)(-0.9000,-0.4750)(-0.8000,-0.3694)(-0.7000,-0.2756)(-0.6000,-0.1891)(-0.5000,-0.1081)(-0.4000,-0.03144)(-0.3000,0.04134)(-0.2000,0.1107)(-0.1000,0.1769)(-0.00001000,0.2403)(0.09999,0.3009)(0.2000,0.3590)(0.3000,0.4146)(0.4000,0.4677)(0.5000,0.5184)(0.6000,0.5664)(0.7000,0.6113)(0.8000,0.6521)(0.9000,0.6853)(1.000,0.6692)
\pscurve[linewidth=0.075](1.500,1.112)(1.600,1.358)(1.700,1.558)(1.800,1.768)(1.900,1.995)(2.000,2.242)(2.100,2.515)(2.200,2.816)(2.300,3.153)(2.400,3.532)(2.500,3.960)
\psline[linecolor=blue](0.66141,0.24028)(0,0.24028)(0,-0.27156)(-0.69548,-0.27156)(-0.69548,-0.62200)(-0.99680,-0.62200)(-0.99680,-0.65853)(-0.80672,-0.65853)(-0.80672,-0.37599)(-0.17524,-0.37599)(-0.17524,0.12738)(0.53805,0.12738)(0.53805,0.53701)(0.93095,0.53701)(0.93095,0.69234)(0.98439,0.69234)(0.98439,0.61976)(0.71975,0.61976)(0.71975,0.29750)(0.094233,0.29750)(0.094233,-0.21076)(-0.62568,-0.21076)(-0.62568,-0.59521)
\psline[linecolor=red](-10.370,7.3839)(3,7.3839)(3,1.4593)(1.6509,1.4593)(1.6509,1.0820)(1.5098,1.0820)(1.5098,1.1601)(1.9312,1.1601)(1.9312,2.0697)(6.3650,2.0697)(6.3650,-7.2209)(-3.8128,-7.2209)(-3.8128,-1.5897)(-2.0651,-1.5897)(-2.0651,-1.3180)(-2.2477,-1.3180)(-2.2477,-1.9385)(-6.7673,-1.9385)(-6.7673,18.081)(3.5103,18.081)(3.5103,1.5806)(1.7110,1.5806)(1.7110,1.0967)(1.5017,1.0967)(1.5017,1.1290)(1.8261,1.1290)(1.8261,1.8258)(4.7479,1.8258)(4.7479,-16.606)(-4.5948,-16.606)(-4.5948,-1.7038)(-2.1228,-1.7038)(-2.1228,-1.3110)(-2.1557,-1.3110)(-2.1557,-1.7660)(-5.0840,-1.7660)(-5.0840,-40.409)(4.2735,-40.409)(4.2735,1.7390)(1.7865,1.7390)(1.7865,1.1175)(1.5003,1.1175)(1.5003,1.1052)(1.7426,1.1052)(1.7426,1.6459)(3.8105,1.6459)(3.8105,51.295)(-5.8843,51.295)(-5.8843,-1.8551)(-2.2032,-1.8551)(-2.2032,-1.3137)(-2.0883,-1.3137)
\end{pspicture*}
\caption{An instance of a biquadratic curve $F(x,y)=X_2(x)y^2+X_1(x)y+X_0(x)=0$, with
$P(x)=42.26667(x+2)(x^2-1)(x-1.5), X_2(x)=(x+5.5)(x-4)$,
and where $X_1(x)= -0.7333x^2 + 27.49x - 0.6882 $ interpolates $\sqrt{P}$ at the zeros of $X_2$ (vertical
asymptotes). }\label{difffigx}
 \end{center}\end{figure}
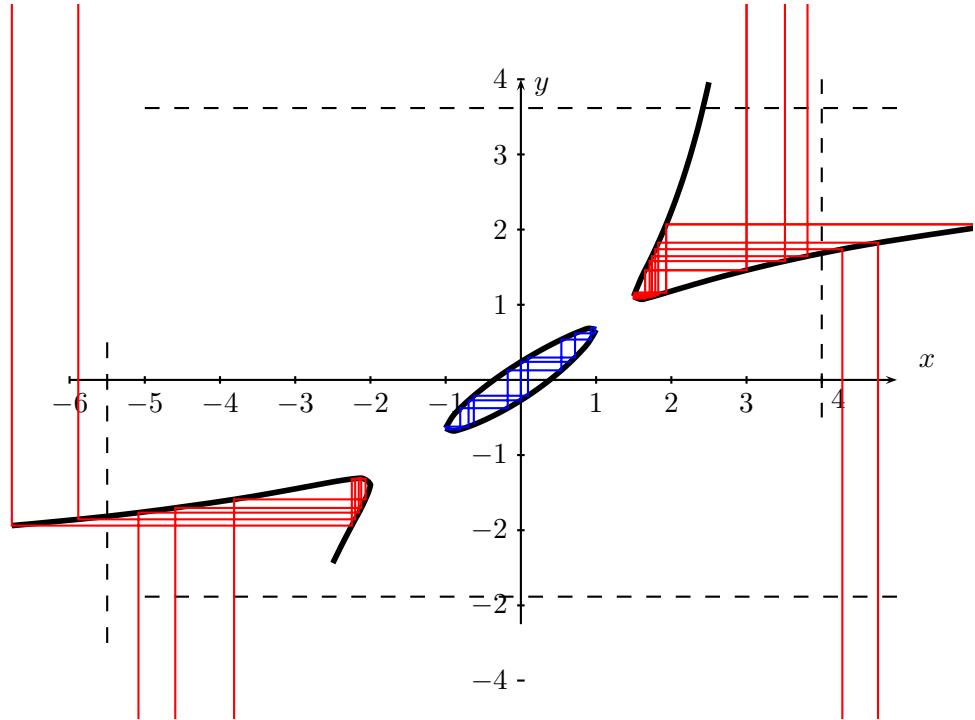

near $(x'_1,y'_0)$, $0=F(x'_1+dx,y'_0+dy)=dx\partial F/\partial x
+dy\partial F/\partial y+o(dy)$,
\begin{equation}\label{tang}
dy/dx=-(2x'_1Y_2(y'_0)+Y_1(y'_0))/(2y'_0X_2(x'_1)+X_1(x'_1))=
-(x'_1-x'_0)Y_2(y'_0)/((y'_0-y'_1)X_2(x'_1))
\end{equation}


\medskip

We come now to equivalent definitions of elliptic lattice, or grid.

\subsection{Definition 1.\label{def1}} \textsl{A sequence $\{\dotsc, x_{-1}, x_0, x_1,\dotsc\}$,
is an elliptic lattice if there exists
a biquadratic polynomial \eqref{Fbiq}
such that $F(x_n,y)=0$ and $F(x_{n+1},y)=0$ have a common root in $y$.}

A very awkward definition! I am probably the only one person using it.

\subsubsection{Exercise.\label{nsquare}} Are $\{an+b\}, \{n^2\}$ special
cases of elliptic sequences? And $\{\sqrt{n}\}$? Impossible to answer with
the present definition, see \S~\ref{nsquarea}.

\subsubsection{Direct formulas.} We have
\begin{equation} \label{ydirect}
y_n \  \text{and\ } y_{n-1} = \dfrac{ -X_1(x_n) \pm \sqrt{P(x_n)}  }{2 X_2(x_n) },\ \ x_n \  \text{and\ } x_{n+1} = \dfrac{ -Y_1(y_n) \pm \sqrt{Q(y_n)}  }{2 Y_2(y_n) },
\end{equation}
where
\begin{equation}\label{P} P=X_1^2-4X_0 X_2,\ \ Q=Y_1^2-4Y_0 Y_2\end{equation} are  polynomials of degree $\leqslant 4$.







Of course, the sequence $\{y_n\}$ is elliptic too, consider $G(x,y)=F(y,x)$.
The construction above is called ``T-algorithm'' in \cite[Theorem~6]{SpZ2007}.

Given $x_0$, we choose $y_0$ as one of the two $y-$roots of $F(x_0,y)=0$,
and $y_1, x_1, y_2$ etc. follow without having to solve any new quadratic
equation, by using \eqref{sumprod} repeatedly.

If we choose the other root at $x_0$, say, $y'_0=y_1$, then $y'_1=y_0, y'_2=y_{-1}$ etc.

\medskip

\emph{Special cases.} We already encountered the usual difference operators
involving $(x,x+h)$ or $(x-h,x)$ or $(x-h/2,x+h/2)$, which must
now be written  $(y+\alpha h ,y+\beta h)$ with various $(\alpha,\beta)$,
corresponding to $F(x,y)=(x-y-\alpha h)(x-y-\beta h)$, so,
$X_2(x)\equiv 1$, $X_1$ of degree 1, $X_0$ of degree 2 with $P=X_1^2-4X_0X_2$ of degree 0.
For the geometric difference operator, $P$ is the square of a first degree polynomial.
For the Askey-Wilson operator \cite{AAR,AW319,Ism2005,Koe,Segovia,magsnul}, $P$ is an arbitrary second degree
polynomial ($=$ a fourth degree polynomial with a double zero at $\infty$).

\medskip

\subsubsection{A glance at the numerical $F(x,y)$ laboratory\label{labo}}.

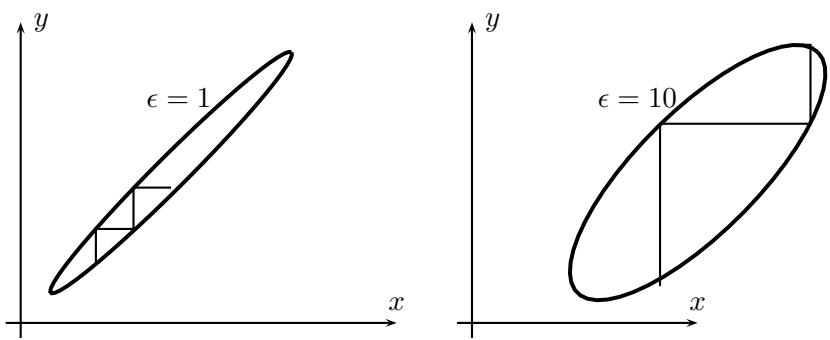
\begin{figure}[htbp]\begin{center}
\psset{xunit=1cm,yunit=1cm}
\begin{pspicture}(0,-0.2)(11,4.05)
\psline{->}(-0.2,0)(5,0)\psline{->}(0,-0.2)(0,4)
\uput[90](5,0){$x$}\uput[0](0,4){$y$}
\psline{->}(5.8,0)(9,0)\psline{->}(6,-0.2)(6,4)
\uput[90](9,0){$x$}\uput[0](6,4){$y$}
\uput[0](1.5,3){$\epsilon=1$}\uput[0](7.5,3){$\epsilon=10$}
\parametricplot[linewidth=1.5pt]{0}{360}{1.6 t cos mul  0.15 t sin mul  sub 2 add 1.6 t cos mul 0.15 t sin mul add 2 add}
  \parametricplot[linewidth=1.5pt]{0}{360}{1.6 t cos mul  0.57 t sin mul  sub 9 add 1.6 t cos mul 0.57 t sin mul add 2 add}
\psline(1,0.8)(1,1.25)(1.5,1.25)(1.5,1.8)(2,1.8)
\psline(8.5,0.49)(8.5,2.65)(10.5,2.65)(10.5,3.7)(10.2,3.7)
 \end{pspicture}
\caption{Influence of $\epsilon$ on the step.}\label{difffig2}
 \end{center}\end{figure}

Choosing the $c_{i,j}$s at random will not yield a readable result.
The most important feature is $P(x)=\epsilon(x-z_1)(x-z_2)(x-z_3)(x-z_4)$
whose square root will be needed in many expressions.
Why this symbol $\epsilon$? This coefficient controls the \emph{step},
the $x_n$s are close together if $\epsilon$ is small, as the square root is
small in \eqref{ydirect}, see fig.~\ref{difffig2}.

How to
reconstruct
the $X$s in \eqref{P}: $P=X_1^2-4X_0X_2$? We could give $X_1$
and recover $X_0$ and $X_2$ by factoring the quartic $X_1^2-P$.

The preferred choice is to give $X_2(x)= c_{2,2}(x-u)(x-v)$ ($u$ and $v$
are the vertical asymptotes in fig. \ref{difffigx}), then $X_1$
interpolates some determinations of the square root of $P$ at $u$
and $v$, and $X_0$ is the (exact) quotient of the division of
$X_1^2-P$ by $4X_2$.

See here a test with the gp-pari software \cite{pari}:

{\small

\begin{verbatim}
GP/PARI CALCULATOR

PARI/GP is free software, covered by the GNU General Public License, and
comes WITHOUT ANY WARRANTY WHATSOEVER.

nd=15;yv=vector(nd+5);xv=vector(nd+5);fv=vector(nd+1);fv1=vector(nd+1);

\\  F(x,y)=X2(x)y2+X1(x)y+X0(x)  = Y2(y)x2+Y1(y)x+Y0(y)   P=(X1)^2-4X0X2
epsi=42.26667;zerP=[-2,-1,1,1.5];
P=epsi*prod(k=1,4,x-zerP[k]);


\\  given X2  u,v = vertical asymptotes
uvert=-5.5;vvert=4.0;X2=(x-uvert)*(x-vvert);
\\  P=X1^2-4X0X2 : X1=sqrt(P) at zeros of X2
sqrP1=-sqrt( subst(P,x,uvert));sqrP2=sqrt( subst(P,x,vvert));
X1=( sqrP2*(x-uvert)-sqrP1*(x-vvert))/(vvert-uvert)-0.7333*X2;
XX0=divrem(P-X1^2,X2,x);X0=-XX0[1]/4;X0check=XX0[2];
print(zerP,P," ",P-X1^2+4*X0*X2," ; check:",X0check);


F=X2*y^2+X1*y+X0;

print("X0= ",X0," ; X1= ",X1," ; X2= ",X2);
Y0=polcoeff(F,0,x);Y1=polcoeff(F,1,x);Y2=polcoeff(F,2,x);
print("Y0= ",Y0," ; Y1= ",Y1," ; Y2= ",Y2);

c=matrix(3,3);for(k=0,2,for(j=0,2,
  c[j+1,k+1]=polcoeff(polcoeff(F,j,x),k,y)));

print1("F=",F);eF=F-sum(j=0,2,sum(k=0,2,c[j+1,k+1]*x^j*y^k));
print(" check F= ",sum(k=0,2,sum(j=0,2,abs(polcoeff(polcoeff(eF,j,x),k,y)))));

\end{verbatim}

results in

\begin{verbatim}


[-2, -1, 1, 1.500]42.27*x^4 + 21.13*x^3 - 169.1*x^2 - 21.13*x + 126.8   0; check: 0*x + 0
X0= -10.43*x^2+ 0.2876*x+ 1.436 ; X1= -0.7333*x^2+27.49*x- 0.6882 ; X2=x^2+ 1.500*x-22.00
Y0= -22.00*y^2-0.6882*y+1.436 ; Y1= 1.500*y^2+27.49*y+ 0.2876 ; Y2=y^2 - 0.7333*y- 10.43
F=(y^2 - 0.7333*y - 10.43)*x^2 + (1.500*y^2 + 27.49*y + 0.2876)*x
+ (-22.00*y^2 - 0.6882*y + 1.436)     check F= 0.E-27


\end{verbatim}

} 

with an almost symmetric configuration $P(x)=\epsilon (x+2)(x^2-1)(x-1.5)$ and
$\epsilon=42.27$ presented in fig.~\ref{difffigx}, we choose
 $u=-5.5, v=4$,
so
$ X_2(x)= x^2 + 1.5x - 22$, then $X_1(x)=$ the linear interpolant to determinations
of the square root of $P$ at
$u$ and $v$ augmented by a multiple of $X_2$, here, we choose
$X_1$ to interpolate the negative square root of $P$ at $u$,
and the positive square root of $P$ at $v$,
 $X_1(x)=-0.73 x^2 +27.49 x +0.2876$ and
$X_0(x)= (X_1^2(x)-P(x))/(4X_2(x))= -10.43 x^2+0.2876 x+1.436$
follows.

$Y0(y)= -22y^2-0.6882y+1.436 ; Y1(y)= 1.5y^2+27.49y+ 0.2876 ; Y2(y)=y^2 - 0.7333y- 10.43$

\noindent We start with $x_0=0$, then $y_0$ is a root of $F(x_0,y)=Y_0(y)=0$, so, \\ $(0.688+\sqrt{0.688^2-4*1.436*(-22)})/(-44)= -0.2716$,
the other root being $y_{-1}=0.2403$. Another run
$(x'_0,y'_0)$, etc. has been performed with $x'_0=3$.

{\Small
\begin{verbatim}
       -1       0       1       2        3      4       5       6       7       8        9       10
x    0.6614  0      -0.6955 -0.9968 -0.8067 -0.1752  0.5380  0.9309  0.9843  0.7197  0.0942 -0.62568
y    0.2403 -0.2716 -0.6220 -0.6585 -0.3760  0.1274  0.5370  0.6923  0.6197  0.2975 -0.2107 -0.59521
x' -10.370   3       1.6509  1.5098  1.9312  6.3650 -3.8128 -2.0651 -2.2477 -6.7673  3.5103  1.7110
y'   7.3839  1.4593  1.0820  1.1601  2.0697 -7.2209 -1.5897 -1.3180 -1.9385 18.081   1.5806  1.0967

asympt hor   3.617 -2.884
asympt vert  4.000 -5.500

\end{verbatim}
}

\medskip

\subsection{Essential and secondary parameters: modulus and step.\label{essential}}\   \\

$F(x,y)=\sum_0^2\sum_0^2 c_{i,j}x^iy^j$ has 9 parameters  in homogeneous
form, leaving 8 degrees of freedom.

The first degree rational transformations

\begin{equation}\label{rat1}
\xi=\dfrac{\alpha x+\beta}{1+\gamma x},\ \ \ \eta=\dfrac{\alpha' y+\beta'}{1+\gamma' y}\end{equation}
still lead to biquadratic polynomials:

$\mathcal{F}(\xi,\eta)= (1+\gamma' y)^{-2}(1+\gamma x)^{-2} \sum_0^2\sum_0^2 \gamma_{i,j}(\alpha' y+\beta')^i(1+\gamma' y)^{2-i}(\alpha x+\beta)^j(1+\gamma x)^{2-j}$

$= (1+\gamma' y)^{-2}(1+\gamma x)^{-2} F(t,x)$.

As each rational function carry 3 coefficients, there remains two
parameters invariant with respect to the transformations of \eqref{rat1},
they are called the \textbf{\emph{essential}} parameters.

The cross ratio of any set of 4 $x-$points, say, $z_1,\dotsc, z_4$
corresponding to $\zeta_1,\dotsc, \zeta_4$ has this invariance
property:

$\lambda =\dfrac{(\zeta_3-\zeta_1)/(\zeta_3-\zeta_2)}{(\zeta_4-\zeta_1)/(\zeta_4-\zeta_2)}
=
\dfrac{(z_3-z_1)/(z_3-z_2)}{(z_4-z_1)/(z_4-z_2)}$.

\noindent We choose the four zeros of $P$ sent to the symmetric Jacobi
configuration \\  $\{\zeta_1,\dotsc, \zeta_4\}=
\{-1/k,-1,1,1/k\}$, then
$\lambda= \dfrac{ (1+1/k)/2}{(2/k)/(1/k+1)}= \dfrac{(1+k)^2}{4k}$
defining the \emph{\textbf{modulus}} $k=\left(\sqrt{\lambda}-\sqrt{\lambda-1}\right)^2$ from the cross-ratio
of the 4 zeros of $P$. Permutations of the 4 zeros lead
to 6 possible values of the cross-ratio
$\lambda, 1/\lambda, 1-\lambda, 1/(1-\lambda), 1-1/\lambda, 1-1/(1-\lambda)=\lambda/(\lambda-1)$, 
the anharmonic group
 ( Maillard\cite[ p.300]{Maillard}).
If $\lambda=\cos^2\theta$, the other values are $\sec^2\theta,
\sin^2\theta, \csc^2\theta, -\tan^2\theta, -\cot^2\theta$.

When $z_1,\dotsc, z_4$ are real numbers in ascending order,

$\lambda= \dfrac{(z_3-z_1)/(z_3-z_2)}{(z_4-z_1)/(z_4-z_2)}
= \dfrac{ 1+(z_2-z_1)/(z_3-z_2) }{1+(z_2-z_1)/(z_4-z_2) }>1$
and $0<k<1$.

With the numerical example of \S \ref{labo}, $[z1,\dotsc,z4]=[-2,-1,1,1.5]$, one finds
$\lambda=15/14=1.0714, k=0.58957; \xi=\dfrac{ x +\gamma}{1+\gamma x}$ with
$\gamma= (z_j-\xi_j)/(z_j\xi_j-1)=-0.12702$.

\begin{equation}\label{zetajm}\zeta_j=\dfrac{\alpha z_j+\beta}{1+\gamma z_j} \Rightarrow
\zeta_j-\zeta_m= \dfrac{\beta-\alpha/\gamma}{1+\gamma z_j}-
\dfrac{\beta-\alpha/\gamma}{1+\gamma z_m}
=\dfrac{(\alpha-\beta\gamma)(z_j-z_m)}{(1+\gamma z_j)(1+\gamma z_m)}\end{equation}

\smallskip

The other essential parameter is the \emph{\textbf{step}},
$x_{n+1}-x_n$ in the uniform case, but should be in general related
to the width of the figure (fig.~\ref{difffig2}). The ratio
$(x_{n+1}-x_n)/(z_2-z_1)$ is independent of dilations of the $x_n$s.
The nonuniformity is almost canceled when one considers that the
computation of $x_{n+1}$ from $x_n$ involves the square root of
$P(x_n)=\epsilon\prod_1^4(x_n-z_k)$ (from \eqref{sumprod}-\eqref{ydirect}), and we
may consider $(x_{n+1}-x_n)/\sqrt{P(x_n)}$.

What works is the definition of the step

\begin{equation}\label{stepdef}
 h=
\dfrac{1}{1-k}
\int_{x_n}^{x_{n+1}} \dfrac{\sqrt{(z_1-z_2)(z_3-z_4)}\;dx}{\sqrt{(x-z_1)(x-z_2)(x-z_3)(x-z_4)}},
\end{equation}
on any set $x_n, x_{n+1}$ of the $x-$lattice. More on the numbering $z_1, z_2$, etc. is
given at the end of the present subsection.

We show that the definition \eqref{stepdef} is invariant under transformations \eqref{rat1}.

$\xi-\zeta_j= \dfrac{\alpha x+\beta}{1+\gamma x}-\dfrac{\alpha z_j+\beta}{1+\gamma z_j}=
(x-z_j)\dfrac{\alpha  -\beta\gamma }{(1+\gamma x)(1+\gamma z_j)}$,

$(\xi-\zeta_1)\dotsb (\xi-\zeta_4)=\dfrac{ (x-z_1)\dotsb(x-z_4)
(\alpha-\beta\gamma)^4 }{(1+\gamma x)^4(1+\gamma z_1)\dotsb(1+\gamma z_4) }$
and, with
$\dfrac{dx}{d\xi}= -\dfrac{(1+\gamma x)^2 }{ \beta\gamma-\alpha}  $, we have \\
$\dfrac{d\xi}{\sqrt{(\xi-\zeta_1)\dotsb (\xi-\zeta_4)}}=
C \dfrac{dx}{\sqrt{(x-z_1)\dotsb(x-z_4)}}$, with

$C=-\dfrac{ \beta\gamma-\alpha} {1+\gamma x)^2}
 \sqrt{   \dfrac{(1+\gamma x)^4(1+\gamma z_1)\dotsb(1+\gamma z_4)}{(\alpha-\beta\gamma)^4}   }
=    \sqrt{   \dfrac    {(z_1-z_2)(z_3-z_4)}{(\zeta_1-\zeta_2)(\zeta_3-\zeta_4)}  }   $,
using \eqref{zetajm}.

This shows that \eqref{stepdef} is independent of rational transformations. The independence
with respect to $n$ is of course much more important and is the main subject of \S \ref{def4}.

With a rational map \eqref{rat1} sending $z_1,\dotsc, z_4$ to
$\zeta_1,\dotsc, \zeta_4=-1/k, -1, 1, 1/k, $ \eqref{stepdef} becomes

\begin{equation}\label{hell}  h=\int_{\xi_{n-1}}^{\xi_n} \dfrac{d\xi}{\sqrt{(1-\xi^2)(1-k^2\xi^2)}},\end{equation}
indeed, $(\zeta_1-\zeta_2)(\zeta_3-\zeta_4)=(1/k-1)^2$, and
$(\xi^2-1)(\xi^2-1/k^2)=(1-\xi^2)(1-k^2\xi^2)k^2$.

With the numerical example of \S \ref{labo}, a crude check of \eqref{stepdef}
with $\sqrt{(z_1-z_2)(z_3-z_4)}/(1-k)=\sqrt{(-2+1)(1-1.5)}/0.4104=1.7229$,
 by
$h\approx 1.7229(x_{n+1}-x_n)/\sqrt{(x-z_1)\dotsb (x-z_4)}$ at

  $x=(x_n+x_{n+1})/2$:

\begin{verbatim}
        n        3     4        5       6         7       8      9     10
       xn     -0.8607 -0.1752 0.5380 0.9309   0.9843   0.7197  0.0942 -0.6257
        h      0.720   0.737  0.689  0.252?  -0.641?? -0.727  -0.735
    Simpson    0.771   0.766  0.785  0.271?  -0.850   -0.768  -0.766
\end{verbatim}

Now we proceed to an accurate computation of $h=\int_{\xi_0}^{\xi_1} dx  /\sqrt{(1-\xi^2)(1-k^2\xi^2)} \\   =\int_{\varphi_0}^{\varphi_1}d\varphi/\sqrt{1-k^2\sin^2\varphi}=F(\varphi_1,k)-F(\varphi_0,k)$, where $\xi=\sin \varphi$.

Landen: $F(\varphi,k)= (\lim\varphi^{(n)}/2^n)(1+k_1)(1+k_2)\dotsb$
where $k_0=k, k_{n+1}=\dfrac{1-\sqrt{1-k_n^2}}{1+\sqrt{1-k_n^2}},
 \varphi^{(n+1)}=\varphi^{(n)}+\arctan(\sqrt{1-k_n^2}\, \tan\varphi^{(n)})$ \cite[15.7]{AbrS}

Here, $k=0.5896, k_1=0.1064, k_2=0.0028, k_3=2.02\; 10^{-6}$.
With $\xi_0=(\alpha x_0+\beta)/(1+\gamma x_0)= -0.12702 =\sin(-0.12736),
\varphi_0^{(1)}/2=-0.11521, \varphi_0^{(2)}/4=-0.11490$, etc. Elliptic integrals
from $\varphi_0$ to $\varphi_1$ and from $\varphi_1$ to $\varphi_2$ yield the
same value of $h=-0.76411$. We have $\xi_n=\sn(nh+g)$, where $g=$ arcsn$(\xi_0)=
F(\varphi_0,k)=-0.11490(1+k_1)(1+k_2)\dotsb = -0.12748$

We also have $nh= \int_{\xi_0}^{\xi_n} \dfrac{d\xi}{\sqrt{(1-\xi^2)(1-k^2\xi^2)}}$,
giving an explicit relation between $n$ and $x_n$, $nh$ being an
\textbf{elliptic integral} of $x_n$, which is its inversion, i.e., an \textbf{elliptic
function} of  $n$: $\xi_n=\sn(nh+g)$. Much more in \S  \ref{def4}.

The simplest $F$ for given $k$ and $h$ will be given in \eqref{snacnadna}.

 R. Baxter uses
the Jacobi canonical setting in the form  $t_n=k^{1/2} \sn(nh+g)$
\cite[\S~15.10]{Baxter} ,so, $nh= \int_{t_0}^{t_n} \dfrac{dt}{\sqrt{(k-t^2)(1-kt^2)}}.$

A third essential parameter depending on the two first ones is the \textbf{period}.
The $n-$period is the (interpolated) index where $(x_n,y_n)$ reproduce $(x_0,y_0)$.
On our example, could it be between 4 and 5?

{\Small
\begin{verbatim}
       -1       0       1       2        3      4      5    ...   8        9    9.1234     10
x    0.6614  0      -0.6955 -0.9968 -0.8067 -0.1752 0.5381  ... 0.7197  0.0942    0    -0.62568
y    0.2403 -0.2716 -0.6220 -0.6585 -0.3760  0.1274 0.5370  ... 0.2975 -0.2107 -0.2716 -0.59521
\end{verbatim}
}

 No, $y_n$ is not close to $y_0$ there,
but between 9 and 10, simple interpolation of $x_n$ and $y_n$ suggest
$n_{\text{period}}=9.1$. The true definition is based on the period of the elliptic
sine, $4K=2\pi(1+k_1)(1+k_2)\dotsb = 6.971$ here, and $n_{\text{period}}=4K/|h|
=-9.1234$ will be explained in \S \ref{def4}.

Simple limit cases are the obvious $h=\int_{x_n}^{x_{n+1}} dx$ on the arithmetic
lattice $x_n=nh$ corresponding to all the $z_i$s $\to\infty$, so that we have
to integrate a mere
constant, and $h=x_{n+1}-x_n$ is a possible choice (there is no more dilation invariance
is such a limit case).

On the geometric lattice $x=q^n$, $P(x)=$ constant $x^2$,and
$h=$ constant times $\int_{x_n}^{x_{n+1}} dx/x$ gives
$h=\log q$ (the log in some basis).

Of course $k=0$ in \eqref{hell} gives $\xi_n=\sin(nh+g)$, basically the formula VI of
Nikiforov et al. \cite{NS,NSU}.

\medskip


\subsection{Definition 2. Symmetry.\label{def2}} \textsl{An elliptic lattice, or grid, is a sequence
satisfying a symmetric biquadratic relation} \cite[Theorem~5]{SpZ2007}
\begin{equation}\label{Ebiq}
 E(x_n,x_{n+1})= \sum_0^2\sum_0^2 e_{i,j}x^i_nx^j_{n+1}=0,\  \   e_{i,j}=e_{j,i}.\end{equation}

As $E$ does not depend on $n$, we also have $E(x_{n-1},x_n)=0$, so,
$E(x_n,x_{n-1})=0$ by the symmetry of $E$.

The two $z-$roots of $E(x_n,z)=0$ are therefore $x_{n\pm 1}$

One has

\begin{equation}\label{sumprodE} x_{n-1}+x_{n+1}= - \dfrac{ e_{1,2}x_n^2+e_{1,1}x_n+e_{0,1}}
{ e_{2,2}x_n^2 +e_{1,2}x_n+e_{0,2} },\ \ \  x_{n-1}x_{n+1}=
\dfrac{ e_{0,2}x_n^2+e_{0,1}x_n+e_{0,0}}
{ e_{2,2}x_n^2 +e_{1,2}x_n+e_{0,2} }.\end{equation}

One also has $E(x,z)
=e_{0,0} +e_{0,1}(x+z) + e_{0,2}(x^2+z^2)
  +e_{1,1}xz
 +e_{1,2}xz(x+z) +e_{2,2}x^2 z^2  $

$   =e_{0,0} +e_{0,1}\Sigma + e_{0,2}\Sigma^2
  +(e_{1,1}-2e_{0,2})\Pi
 +e_{1,2}\Sigma\Pi  +e_{2,2}\Pi^2$.

where  $\Sigma=x+z$ and $\Pi=xz$.

\subsubsection{$1\Rightarrow 2$.}  Let $\{x_n,y_n\}$ be sequences satisfying the definition~1.
A relation involving only $x_n$ and $x_{n+1}$ is obtained by
the elimination of $y_{n}$ through the resultant \cite[\S~8]{bliss},
\cite{discri},
(bigradient in Aitken \cite{Aitken}) of the two
polynomials in $y_{n}$ from \eqref{sumprod}
$P_1(y_{n})= (x_n+x_{n+1})Y_2(y_{n})+Y_1(y_{n})$

and
$P_2(y_{n})= x_n x_{n+1}Y_2(y_{n})-Y_0(y_{n})$.

\[ \begin{vmatrix} \Sigma c_{2,0}+c_{1,0} & \Sigma c_{2,1}+c_{1,1}& \Sigma c_{2,2}+c_{1,2}& 0\\
  0 & \Sigma c_{2,0}+c_{1,0} & \Sigma c_{2,1}+c_{1,1}& \Sigma c_{2,2}+c_{1,2} \\
 \Pi c_{2,0}-c_{0,0} & \Pi c_{2,1}-c_{0,1} &\Pi c_{2,2}-c_{0,0}&0\\
  0 & \Pi c_{2,0}-c_{0,0} & \Pi c_{2,1}-c_{0,1} &\Pi c_{2,2}-c_{0,0}
\end{vmatrix} =0,
\]
were $\Sigma$ and $\Pi$ are here $x_n+x_{n+1}$ and $x_nx_{n+1}$.
 It is not obvious that
the result is biquadratic in $x_n$ and $x_{n+1}$,
that the expression is quadratic in $\Sigma$ and $\Pi$.
The coefficient of $\Sigma^2\Pi^2$ is
$\begin{vmatrix} c_{2,0} & c_{2,1}& c_{2,2}& 0\\
  0 & c_{2,0} & c_{2,1}& c_{2,2}\\
  c_{2,0} & c_{2,1}& c_{2,2}& 0\\
0 &c_{2,0} & c_{2,1}& c_{2,2}\end{vmatrix} =0$, and
the $\Sigma^2\Pi$ and $\Sigma \Pi^2$ terms are still determinants with
at least one pair of identical rows.

The final quadratic formula in $\Sigma,\Pi$ is
\begin{equation} \label{Rconic}
\mathcal{R}(\Sigma,\Pi):=e_{0,0} +e_{0,1}\Sigma + e_{0,2}\Sigma^2
  +(e_{1,1}-2e_{0,2})\Pi
 +e_{1,2}\Sigma\Pi  +e_{2,2}\Pi^2=0,
\end{equation}
which is $E(x,z)=0$ where  $\Sigma=x+z$ and $\Pi=xz$, as above,
\eqref{Rconic}  is the equation of a \emph{conic}\footnote{Called the $\mathcal{R}-$conic.
The  importance of this conic
has been stressed by A.Ronveaux \cite{Ronv}.}.

\subsubsection{$2\Rightarrow 1$.}

We try to arrive at  $\Sigma=-Y_1(y)/Y_2(y), \Pi=Y_0(y)/Y_2(y)$ (we drop the index $n$) considered
as a parametric representation of the curve \eqref{Rconic} with rational
functions  (\emph{ unicursal} curve)   \cite{jjqell}of degree $\leqslant 2$, which  therefore holds for  a conic, as
well known (\cite{BursZhed1,BursZhed,BursZhed09}, \cite[\S~58]{Deaux}). 

Let $(\Sigma^*, \Pi^* )$ be a particular point of the conic $\mathcal{R}(\Sigma,\Pi)=0$, and solve $\Pi=\Pi^*+t(\Sigma-\Sigma^*)$ for $t$:

$0=\mathcal{R}(\Sigma,\Pi)-\mathcal{R}(\Sigma^*,\Pi^*)$

$=  e_{0,1}(\Sigma-\Sigma^*) + e_{0,2}(\Sigma^2-(\Sigma^*)^2)
  +(e_{1,1}-2e_{0,2})(\Pi -\Pi^*)
 +e_{1,2}(\Sigma\Pi- \Sigma^* \Pi^*)  +e_{2,2}(\Pi^2-(\Pi^*)^2),$

$=(\Sigma-\Sigma^*) \{  e_{0,1}+ e_{0,2}(\Sigma+\Sigma^*)+(e_{1,1}-2e_{0,2})t
    +e_{1,2}\Pi^* +e_{1,2}t\Sigma+e_{2,2}t[2\Pi^*+t(\Sigma-\Sigma^*)]\}$, so

$\Sigma=-\dfrac{  e_{0,1}+ e_{0,2}\Sigma^*+(e_{1,1}-2e_{0,2})t
    +e_{1,2}\Pi^* +e_{2,2}t[2\Pi^*-t\Sigma^*]   }
   {  e_{0,2} +e_{1,2}t+e_{2,2}t^2  } \\
= \Sigma^* -\dfrac{  e_{0,1}+ 2e_{0,2}\Sigma^*+(e_{1,1}-2e_{0,2})t
    +e_{1,2}\Pi^* +e_{1,2}t\Sigma^*+2e_{2,2}t\Pi^*   }
   {  e_{0,2} +e_{1,2}t+e_{2,2}t^2  } $,

$\Pi=\Pi^*-t\dfrac{  e_{0,1}+ 2e_{0,2}\Sigma^*+(e_{1,1}-2e_{0,2})t
    +e_{1,2}\Pi^* +e_{1,2}t\Sigma^*+2e_{2,2}t\Pi^*   }
   {  e_{0,2} +e_{1,2}t+e_{2,2}t^2  } $

give a parametric representation $-Y_1(y)/Y_2(y), Y_0(y)/Y_2(y)$ with $y=t$, so build a valid

$F(x,y)=x^2Y_2(y)+xY_1(y)+Y_0(y)$.

End of $2\Rightarrow 1$. \qed

\smallskip


\subsubsection{Direct formula.} We have $x_{n+1}$ in terms of $x_n$
by \eqref{sumprod} $x_{n+1} =-x_n-Y_1(y_n)/Y_2(y_n)$,
 where
$y_{n}$ is replaced by \eqref{ydirect}. This yields $x_{n+1}= R(x_n)+S(x_n)\sqrt{P(x_n)}$,
with the same $P$ as in \eqref{ydirect}, and with rational $R$ and $S$. Note that the second
determination of the square root of $P$ yields $x_{n-1}$, as $E$ is symmetric.
We see  that $R$ and $S$ are rational functions of degree 2,
$R(x)= - \dfrac{ e_{1,2}x^2+e_{1,1}x+e_{0,1}}
{2( e_{2,2}x^2 +e_{1,2}x+e_{0,2}) }$,
$S(x)=  \dfrac{ \ \ \sigma\ \  }
{ e_{2,2}x^2 +e_{1,2}x+e_{0,2} }$,
by \eqref{Ebiq} and \eqref{sumprodE}, where $\sigma$ is a constant, knowing that  $P(x)$ is a constant times
$
(e_{1,2}x^2+e_{1,1}x+e_{0,1})^2-4(e_{0,2}x^2+e_{0,1}x+e_{0,0})(e_{2,2}x^2+e_{1,2}x+e_{0,2})$.


\subsubsection{\label{nsquarea}Answers to \S~\ref{nsquare}.}
$\{an+b+cn^2\}$ is an elliptic lattice, or grid, (a very special one). Indeed, the
sum $a(n-1)+b+c(n-1)^2+a(n+1)+b+c(n+1)^2=2(an+b+cn^2)+2c$
and the product $[a(n-1)+b+c(n-1)^2][a(n+1)+b+c(n+1)^2]=
a^2 n^2 +2abn +2ac (n^3-n)+2bcn^2 +c^2(n^4-2n^2)-a^2+b^2+2bc+c^2
= (an+b+cn ^2)^2 -2c(an+b+cn^2)-a^2+b^2+4bc+c^2$
are rational functions of $an+b+cn^2$. By the way, if $x_n$ is
elliptic, does the same hold for $x_n+an+b$? Normally, no, but
$x_{an+b}$ is, see Definition 4, where $x_s$ may be defined for any real or
complex $s$.

$x_n=n(-1)^n$ is valid: $x_{n-1}+x_{n+1}=2n(-1)^{n-1}=-2x_n, x_{n-1}x_{n+1}=n^2-1=x_n^2-1,  \\   E(x,z)= z^2+2xz+x^2-1$, indeed symmetric.
What is $y$?

$x_n=q^n$ is obviously valid, but $aq_1^n+bq_2^n$? One must have $q_1+1/q_1=q_2+1/q_2$, so, $q_1=q, q_2=1/q$, and
the symmetric $E(x,z)=z^2-(q+1/q)xz+x^2+ab(q-1/q)^2$ follows. Then,
with $E(x,z)=\Sigma^2 -(q+2+1/q)\Pi +ab(q-1/q)^2, \Pi=ab(\sqrt{q}-1/\sqrt{q})^2+y\Sigma
\Rightarrow \Sigma=(q+2+1/q)y, \Pi=ab(q-2+1/q)+(q+2+1/q)y^2,
F(x,y)=x^2-(q+2+1/q)xy+ab(q-2+1/q)+(q+2+1/q)y^2$ is a possible $F$.
Finally, $x=aq^n+b/q^n\Rightarrow y=x_n/2\pm(aq^n-b/q^n)(\sqrt{q}-1/\sqrt{q})/(\sqrt{q}+1/\sqrt{q})$

$\sqrt{n}$ is not elliptic: $\sqrt{n-1}+\sqrt{n+1}$ is not a rational function
of $\sqrt{n}$.

\subsubsection{\label{exer2}Answer to \S~\ref{exer1} exercise.}

One immediately has

  $\chi_n=\dfrac{x_{n+1}-x_{n-1} }{y_n-y_{n-1}  } =
\dfrac {2\sigma\sqrt{P(x_n)} }{  e_{2,2}x^2 +e_{1,2}x+e_{0,2}  }
     \dfrac{X_2(x_n)}{\sqrt{P(x_n)} }   = $ a constant times
$\dfrac  { X_2(x_n)}{  e_{2,2}x_n^2 +e_{1,2}x_n+e_{0,2} } $
from \eqref{ydirect} and the formula for $S_n$ above. How sloppy! we do not
know the signs of the two square roots of $P$.
For a more complete, but tedious, explanation, we follow the suggestion of
\S~\ref{exer1} and computations already done in , so,

$\chi_n=\dfrac{-Y_1(y_n)/Y_2(y_n)+Y_1(y_{n-1})/Y_2(y_{n-1}) }{ y_n-y_{n-1} }$,which is
a symmetric rational function of $y_n$ and $y_{n-1}$, so, of $x_n$. By the simple
fractions expansions $Y_1(y)/Y_2(y)=c_{1,2}/c_{2,2} +Y_1(c)/(Y'2(c)(y-c))+Y_1(d)/(Y'_2(d)/(y-d))$,where $c$ and $d$ are the roots of $Y_2(y)=0$ (horizontal asymptotes).
Remark that $Y'_2(c)=c_{2,2}(c-d)$ and $Y'_2(d)=c_{2,2}(d-c)$ are opposite. Then,

$\chi_n= \dfrac{ Y_1(c)/((y_{n-1}-c)(y_n-c)) -Y_1(d)/((y_{n-1}-d)(y_n-d))}{c_{2,2}(c-d)}
=\dfrac{Y_1(c)/F(x_n,c)-Y_1(d)/F(x_n,d) }{c_{2,2}(c-d)/X_2(x_n)} \\
= \dfrac{1/(x_n+Y_0(c)/Y_1(c)) -1/(x_n+Y_0(d)/Y_1(d))  }{c_{2,2}(c-d)/X_2(x_n)}\\
= \dfrac{(Y_0(d)/Y_1(d)-Y_0(c)/Y_1(d))X_2(x_n) }{ c_{2,2}(c-d)(x_n+Y_0(c)/Y_1(c))(x_n+Y_0(d)/Y_1(d)) } $
is the full answer.

\begin{figure}[htbp]\begin{center}
\psset{xunit=2cm,yunit=2cm}
\begin{pspicture}(-0.2,-0.5)(5.0,2.0)
\psline{->}(0,0)(1.75,0)\psline{->}(0,0)(0,1.75)
\uput[90](1.75,0){$x$}\uput[0](0,1.75){$z$}
\uput[270](0.7,0){$\begin{matrix}x_n\\=x_{n+1}\end{matrix}$}
\uput[270](1.38,0.0){$x_{n+2}=x_{n-1}$}
\uput[180](0,1.140){$x_{n-1}=x_{n+2}$}
\uput[0](1.5,1.5){$z=x$}
\uput[0](1.3,0.9){$E(x,z)=0$}
\psline[linestyle=dashed,dash=2mm 2mm](0.7,0)(0.7,1.14)
\psline[linestyle=dashed,dash=2mm 2mm](1.14,0)(1.14,1.14)
\psline[linestyle=dashed,dash=2mm 2mm](0,0.7)(1.14,0.7)
\psline[linestyle=dashed,dash=2mm 2mm](0,1.14)(1.14,1.14)
\psline{->}(2.5,0)(4.25,0)\psline{->}(2.5,0)(2.5,1.75)
  \parametricplot[linewidth=1.75pt]{110}{240}{0.8 t cos mul  0.29 t sin mul  sub 1.5 add 0.8 t cos mul 0.29 t sin mul add 1.5 add}
 \parametricplot[linewidth=1.75pt]{-60}{60}{0.8 t cos mul  -0.29 t sin mul  sub -0.5 add 0.8 t cos mul -0.29 t sin mul add -0.5 add}
\psline(0,0)(1.5,1.5)
\psline(2.5,0)(4.25,1.75)\uput[0](4.25,1.75){$z=x$}\uput[90](4.5,0){$x$}\uput[0](2.5,1.75){$z$}
\psline[linewidth=1.75pt](2.5,0.5)(4,1.25)
\psline[linewidth=1.75pt](3,0)(4,2)
 \end{pspicture}
\caption{A near fixed point and a true  one.}\label{figfixedp}
 \end{center}\end{figure}
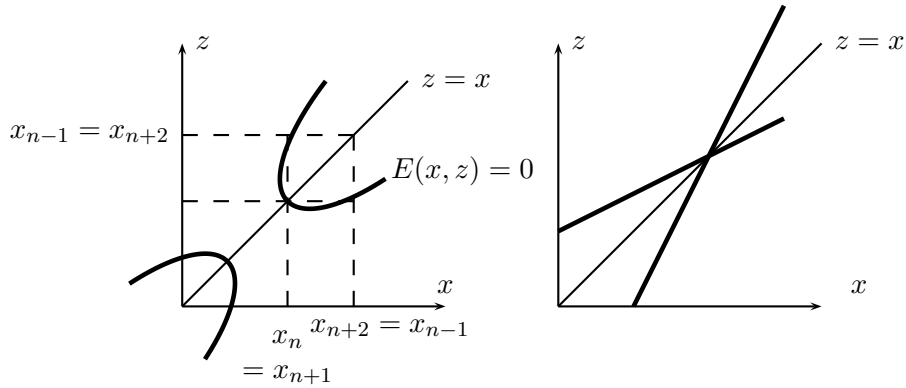

\subsection{Chasing fixed points.\label{fixedp}} \  \\

It seems obvious that there are always 4 fixed points, the 4 roots of $E(x,x)=0$, but
no: we simply have a mirror configuration $x_{n+1}=x_n, x_{n+2}=x_{n-1}$, etc. see
fig.~\ref{figfixedp}.

\textbf{Proposition.} \textsl{There is a fixed point $\dotsb = x_{n-1}=x_n=x_{n+1}=\dotsb
=x_\infty$ only if the equation $P(x)=0$ has a double root, i.e., in the zero-genus case}.

Indeed, as $x_{n\pm 1}=R(x_n)\pm S(x_n)\sqrt{P(x_n)}$, one must have
$P(x_\infty)=
0$, the tangent of the curve $E(x,y)=0$ is vertical at $x=y=x_\infty$,
so, $\partial E(x,y)/\partial y=0$ at this point, but $\partial E/\partial x=0$ as well,
from symmetry of $E$. So, $(x_\infty,x_\infty)$ is a double point of the
curve $E=0$, which looks like two directions $y-x_\infty =$ constant
times $x-x_\infty$ near $(x_\infty,x_\infty)$, see fig. \ref{figfixedp} . \qed

 In the simple geometric case, $E(x,z)= (z-qx)(x-qz)$, the fixed
points are 0 and $\infty$. For ANSUW, the fixed point is $\infty$.

\   \\
\smallskip
What about $y$? For  $x_n, x_{n+1}$ near $x_\infty$, \eqref{sumprod} tells that
$-Y_1(y_n)/Y_2(y_n)$ and $Y_0(y_n)/Y_2(y_n)$ are close to $2x_\infty$ and $x^2_\infty$, so that $y_n$ is close to a zero of $Q=Y_1^2-4Y_0Y_2$, a fixed point
of the $\{y_n\}$ lattice.

\subsection{Definition 3. Continued fraction.\label{def3}}

\subsubsection{Continued fraction.}
Square roots and continued fractions share a long history,
simplest example is $A=\sqrt{a^2+1} = a+\dfrac{1}{a+\sqrt{a^2+1}}
= a+\dfrac{1}{2a+\dfrac{1}{2a+\ddots }}$, going back to Greek mathematics
\cite[\S 1.2]{BrezinskiH},  \cite[\S 2.]{Kru}.

Similar constructions with the square root of a function
$P$ with a known Taylor-Maclaurin expansion have been imagined, such as
$\sqrt{P(x)}=\sqrt{P(0)}+\dfrac{x}{\ \ \  \dfrac{x}{ \sqrt{P(x)}-\sqrt{P(0)}}=\dfrac{ \sqrt{P(x)}+\sqrt{P(0)} }{(P(x)-P(0))/x}\ \ \   }$. If $P$ is a linear or quadratic polynomial,
one  reproduces the same form  and we get a simple periodic continued fraction. Continued fraction expansions built as $f(x)=$ two first terms of
Taylor-Maclaurin expansion $+\dfrac{x^2}{g(x)}$ are called
\emph{associated} continued fractions \cite{Kru,Perron,Wall} etc.

For application to solid-state physics, see \cite{HaydNex}.

When the degree of $P$ is higher,  we encounter again our subject matter:

\subsubsection{The definition.}
\textsl{An elliptic lattice, or grid,
is made of coefficients $x_0, x_1,\dotsc, x_m, \dotsc$ involved in the 
  continued fraction expansion of }  
\begin{equation}\label{fdef3}
f(x)=\dfrac{\sqrt{P(x)}-V(x)}{(x-x_0)(x-v)}     =\dfrac{x-u}{\alpha_0 (x-u)+\beta_0(x-v) -\dfrac{(x-u)(x-v)}{
                \alpha_1 (x-u)+\beta_1(x-v) -\dfrac{(x-u)(x-v)}{ \ddots } } },\end{equation}
 \textsl{containing the square root of a polynomial $P$ of degree}
$\leqslant 4$, \textsl{ where $V$ is a quadratic polynomial interpolating a square root of $P$ at $u$, $v$, and $x_0$.  At the $m^{\text{th}}$ step,}

\begin{equation}\label{def3step}f_m(x)= \dfrac{\sqrt{P(x)}-V_m(x)}{\gamma_m(x-x_m)(x-v)}     =\dfrac{x-u}{ \alpha_m(x-u)+\beta_m(x-v) -(x-v)f_{m+1}(x)},  ,\end{equation}
\textsl{where $\alpha_m(x-u)+\beta_m(x-v)$ is the linear interpolant  to $(x-u)/f_m$ about $x=u$ and $x=v$, and $f_{m+1}(x)$ is the remainder divided by $ -(x-v)$.}

We only need the square root of $P$ as an analytic function  in a neighborhood of $u$, and an independent determination in a neighborhood of $v$. Then, $V_m$ interpolates the chosen determinations, and some (unknown) determination at $x_m$,
 and  $\alpha_m (x-u)+\beta_m(x-v)$ is the linear interpolant to $(x-u)/f_m(x)$. Of course, as $f_m(u)=0$, we already need the derivative in $\beta_m=  1/(f_m'(u)(u-v))=(u-x_0)/(P'(u)/(2\sqrt{P(u)})-V'_m(u)), \alpha_m=1/f(v)$ which is itself a limit: $\alpha_0=(v-x_0)/(P'(v)/(2\sqrt{P(v)})-V'_0(v))$.

We show that $f_{m}$ and $f_{m+1}$ have the same form   in \eqref{def3step}, with $\rho_m$
the fourth root of $V_m^2(x)-P(x)=\delta_m(x-u)(x-v)(x-x_m)(x-\rho_{m})$:

$f_m(x)=\dfrac{1}{\ \ \dfrac{\gamma_m(x-x_m)(x-v)\ \ \  }
           {-V_m(x)+\sqrt{P(x)}  }\ \ \   }
   =\dfrac{V_m^2(x)-P(x)=\delta_m(x-u)(x-v)(x-x_m)(x-\rho_{m})}{\ \ \ \gamma_m(x-x_m)(x-v)[-V_m(x)-\sqrt{P(x)}]   }
$

$= \dfrac{x-u}{\alpha_m(x-u)+\beta_m(x-v)-(x-v)\left(   \dfrac{V_m(x)+\sqrt{P(x)}   }
{ \delta_m(x-v)(x-\rho_{m})/\gamma_m}    + \dfrac{\alpha_m(x-u)}{x-v}+\beta_m ]\right)   }$

So, $V_{m+1}(x)=-V_m(x)  -\delta_m(x-\rho_m) [\alpha_m(x-u)+\beta_m(x-v) ]/\gamma_m$, and  $\gamma_{m+1}=\delta_m/\gamma_m, x_{m+1}=\rho_m$.

More on $\delta_m$, coefficient of $x^4$ of the quartic $V_m^2-P$: as $V_m$ interpolates a square root of $P$ at $u, v$ and $x_m$, let $V_m(x)=X(x)+ z_m(x-u)(x-v)$, where the quadratic (or linear) polynomial $X$ performs the interpolation
at $u$ and $v$, so,
\begin{multline}\label{VXY}V_m^2(x)-P(x)= (x-u)(x-v)[ Y(x)+2z_m X+ z_m^2(x-u)(x-v)], \\   \text{\ where\  } Y  \text{\ is\  the\  polynomial\ } (X^2(x)-P(x))/((x-u)(x-v)).\end{multline}

$\delta_m(x-x_m)(x-x_{m+1})= \gamma_{m}(x-x_{m})\gamma_{m+1}(x-x_{m+1})=Y(x)+2z_m X+ z_m^2(x-u)(x-v)$

This identity appears as $R=r_m^2+s_ms_{m-1}$ in Abel's famous work on the integral of $\rho(x)/\sqrt{R(x}$ \cite[\S 12]{Abelrho}, see also Houzel \cite[\S 4, pp. 73-77]{Houzel}, Zannier \cite{Zannier}.

\begin{figure}[htbp]\begin{center}
 \includegraphics[angle=-90, scale=0.15]{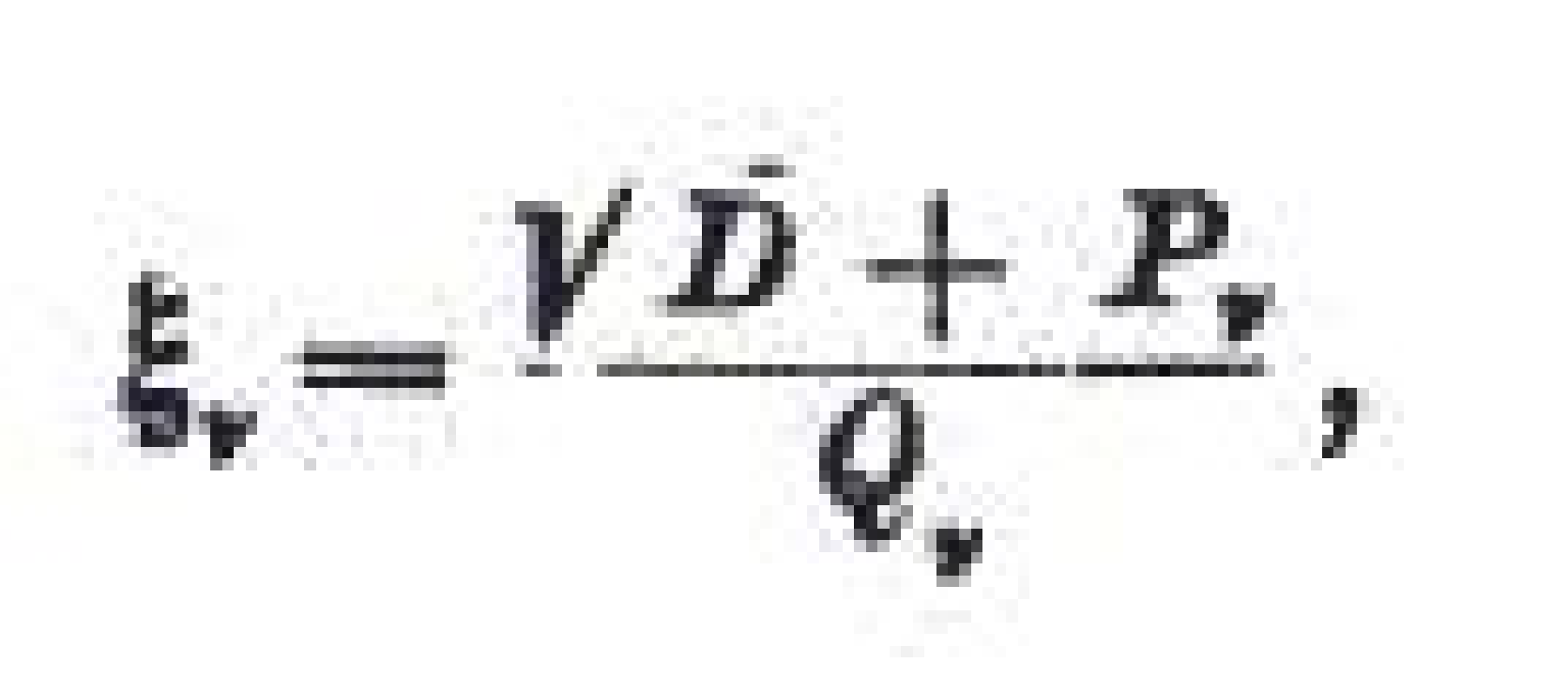} 
\caption{Formula (3) in \S 20 of Perron 1913 \cite{Perron} on
quadratic irrationals, and where it is added that $P_\nu, Q_\nu$, AND $(D-P_\nu^2)/Q_\nu$ are integers.}\label{figper}
 \end{center}\end{figure}

Formulas for quadratic irrationals are strikingly similar, see fig. \ref{figper}.

\subsubsection{$M-$ and $T-$ continued fraction.} The $M-$ or $T-$  continued fraction expansion \cite[\S~6.6, 6.7]{Cuyt} \cite{JonesT1}
\cite[\S~7.3]{JonesT} of  (at least formal) power
expansions  about $x=u$ and $x=v$ is

\begin{equation*}f(x)= f(u) +\dfrac{ z}{\alpha_0 z+\beta_0 -\dfrac{z}{
                \alpha_1 z+\beta_1 -\dfrac{ z}{ \ddots } } },
\end{equation*}

with $z=\dfrac{x-u}{x-v}$,
$f(x)-f(u)=f_0(x)= \dfrac{x-u}{\alpha_0(x-u)+\beta_0(x-v)-(x-v)f_1(x)}$,
where $\alpha_0(x-u)+\beta_0(x-v)$ is the linear interpolant  to $(x-u)/f_0$ about $x=u$ and $x=v$,
 and where the
operation is repeated on $f_1$ etc.
The existence of the $\alpha_n, \beta_n$s depends on the nonvanishing
of relevant determinants \cite[\S~6.6, 6.7]{Cuyt}
\cite[\S~7.3]{JonesT}.

The $T-$fractions were introduced in 1948 by W.J. Thron, and their two-point Pad\'e property
was established in 1976 by J.H. McCabe and J.A. Murphy \cite{Coop,McCabeM,McCabe}, see below in \S \ref{approximants} for
more on the two-point Pad\'e property.

\subsubsection{Proof of $1\Rightarrow 3$.} Let ${x_n, y_n}$ be elliptic grids determined by a biquadratic curve $F(x,y)=X_0(x)+X_1(x)y+X_2(x)y^2=0$, let $V_m(x)=X_1(x)+2 y_m X_2(x), u$ and $v$ be the roots of $X_2(x)=(x-u)(x-v)=0$ (same notation as in \S  \ref{labo}) , and we see
that this $V_m$ interpolates indeed a square root of $P$ at $u, v,x_m, x_{m+1}$:

 $V_m^2(x)-P(x)=X_1^2(x)+4 y_m X_1(x)X_2(x)+4 y_m^2 X_2^2(x)-X_1^2(x)+4X_0(x)X_2(x)      \\   =4X_2(x)\left[\underbrace{X_0(x)+y_m X_1(x)+ y_m^2 X_2(x)}_{\displaystyle F(x,y_m)=Y_2(y_m)(x-x_m)(x-x_{m+1})}\right]$,  so $\delta_m=4Y_2(y_m)$. \qed

\smallskip

Halphen \cite[chap. XIV, (48), p. 603]{Hal} gives the equivalent of our $\alpha, \beta, \gamma, \delta$s in terms  of the Weierstrass $\wp$ function.

\subsubsection{Proof of $3\Rightarrow 1$.}
We have the same form for $f$ and its transformations, but the same symbols may have  different meanings. Now, $P, u, v, x_0$ are given. From $V_0, X, Y$ in \eqref{VXY} at $m=0$, we make $X_0=Y/4, X_1=X, X_2(x)=(x-u)(x-v)$, so we have $F(x,y)$. At the $m^{\text{th}}$ step, again written as

$f_m(x)=\dfrac{-V_m(x)+\sqrt{P(x)}}{\gamma_m(x-x_m)(x-v)}= \dfrac{x-u}{\alpha_m(x-u)+\beta_m(x-v)-(x-v)f_{m+1}(x)} $, where $z_m$ in $V_m(x)=X(x)+z_m(x-u)(x-v)$ is such that
$x_m$ is a root of $V_m^2(x)-P(x)=0$ (together with $u$ and $v$), $x_{m+1}$ is \emph{defined}
here as the fourth root, so that $x_{m}$ and $x_{m+1}$ solve $Y(x)+2z_mX(x)+z_m^2X_2(x)=4F(x,z_m/2)=0$, and  $z_m/2$ is a valid $y_m$. \qed

\medskip

\subsubsection{Approximants\label{approximants}}.

Approximants of $f_0$ are $\dfrac{A_{-1}(x)}{B_{1}(x)}=\dfrac{-1/(x-v)}{0}, \dfrac{A_0(x)}{B_0(x)}=\dfrac{0}{1},  \dfrac{A_1(x)}{B_1(x)}=
\dfrac{x-u}{\alpha_0 (x-u)+\beta_0(x-v)}, \\
\dfrac{A_2(x)}{B_2(x)}= \dfrac{(x-u)(\alpha_1 (x-u)+\beta_1(x-v))}
 {(\alpha_0 (x-u)+\beta_0(x-v))(\alpha_1 (x-u)+\beta_1(x-v))-(x-u)(x-v)}$.


\noindent $A_m$ and $B_m$ satisfy the same recurrence relation, as
well as any linear combination, the most useful being\\
$B_{n+1}(x)f_0(x)-A_{n+1}(x)=(\alpha_n (x-u)+\beta_n(x-v))
(B_{n}(x)f_0(x)-A_{n}(x))-(x-u)(x-v)(B_{n-1}(x)f_0(x)-A_{n-1}(x))$,

$\dfrac{B_{n+1}(x)f_0(x)-A_{n+1}(x)}{B_{n}(x)f_0(x)-A_{n}(x)}=
\alpha_n (x-u)+\beta_n(x-v)
-\dfrac{(x-u)(x-v)}{\ \ \  \dfrac{B_{n}(x)f_0(x)-A_{n}(x)}
      {B_{n-1}(x)f_0(x)-A_{n-1}(x)}\ \ \ }$, or

$\dfrac{B_{n}(x)f_0(x)-A_{n}(x)}
      {B_{n-1}(x)f_0(x)-A_{n-1}(x)}=\dfrac{(x-u)(x-v)}
{ \alpha_n (x-u)+\beta_n(x-v) -  \dfrac{B_{n+1}(x)f_0(x)-A_{n+1}(x)}{B_{n}(x)f_0(x)-A_{n}(x)} }
$, so

$\dfrac{B_{n}(x)f_0(x)-A_{n}(x)}
      {(x-v)(B_{n-1}(x)f_0(x)-A_{n-1}(x))}=f_n(x)$ and

$B_{n}(x)f_0(x)-A_{n}(x) = (x-v)^nf_0(x)\dotsb f_n(x)$
showing that $A_n(x)/B_n(x)-f_0(x)=O((x-u)^{n+1})$ about $x=u$, and
$O((x-v)^n)$ about $x=v$ (\emph{two-point Pad\'e approximation}).
If $u=v$, the approximants of the associated continued fraction are
Pad\'e approximants \cite[(2.3.5), (6.5.1)]{Cuyt}
(\cite[3.10.11]{NIST} where it is called Jacobi)
\cite[\S 61, \S 77 (1913) \S 25, \S 44 (1957)]{Perron}.
When $u=v=\infty$, one must expand in series of $1/x$, and
we have a Jacobi continued fraction \cite{Cuyt,NIST,Perron,Wall}.

\smallskip

$ (x-v)^nf_0(x)\dotsb f_m(x)= B_{m}(x)f_0(x)-A_{m}(x)=
[B_{m}(x)\sqrt{P(x)}-\tilde{A}_{m}(x)]/(\gamma_0(x-x_0)(x-v))$

with $\tilde{A}_m= V_0(x)B_m(x)+\gamma_0(x-x_0)(x-v)A_n(x)$
from $f_0(x)=\dfrac{-V_0(x)+\sqrt{P(x)}}{\gamma_0(x-x_0)(x-v)} $

Outside a system of cuts, $f_m(x)=\dfrac{ -V_m(x)+\sqrt{P(x)}}{\gamma_m(x-x_m)(x-v)}$ is a meromorphic function with a possible pole at $x_m$. With the  conjugate
function $f_m^{(\text{conj})}(x)=\dfrac{ -V_m(x)-\sqrt{P(x)}}{\gamma_m(x-x_m)(x-v)}$,
 the product  $f_m(x)f_m^{\text{conj})}(x)=
\dfrac{\delta_{m+1}(x-u)(x-x_{m+1})}{\delta_m(x-v)(x-x_m)}$, only one of  $f_m, f_m^{(\text{conj})}$
 has a pole at $x_m$. This is an elementary way to cope with a
two-sheeted Riemann surface.

Product of the $f$s:

$(x-v)^{n} f_0(x)\dotsb f_{n-1}(x)\; f_0^{(\text{conj})}(x)\dotsb
f_{n-1}^{(\text{conj})}(x)=(x-u)^{n}\delta_n(x-x_n)/(\delta_0(x-x_0)) \\  =
 (x-v)^{-n+2}\dfrac{ B^2_{n-1}(x)P(x)-\tilde{A}^2_{n-1}(x) }{4\delta_0^2(x-v)^2(x-x_0)^2}$.

\emph{\textbf{Pell}} property: $x_n$ is such that there exists polynomials
$\tilde{A}_{n-1}$ and $B_{n-1}$ such that
\begin{equation}\label{Pell}
B^2_{n-1}(x)P(x)-\tilde{A}^2_{n-1}(x) = 4\delta_0\delta_n(x-u)^n(x-v)^{n}(x-x_0)(x-x_n)
\end{equation}
See \cite{PellLasserre} for a modern treatment of Pell's equation.

\medskip

\subsection{Definition 4. Elliptic functions\label{def4}}\textsl{An elliptic grid is a sequence
$x_n=\mathcal{E}(n h+t_0)$, where $\mathcal{E}$ is any  elliptic function  of order 2 (i.e., bivalent in a fundamental parallelogram of periods)}.

\subsubsection{Proof of $1\Rightarrow 4$.}

How on earth to find  elliptic functions out of $F(x_n,y_n)=0$ and $F(x_{n+1},y_n)=0$ with the biquadratic polynomial $F$ of \S \ref{biquad}?

One may establish that the biquadratic curve $F(x,y)=0$ in \eqref{Fbiq}
has genus 1 and a parametric representation

\[ x=\mathcal{E}_1(t),\ \ y=\mathcal{E}_2(t),      \]
with $\mathcal{E}_1$ and $\mathcal{E}_2$ elliptic functions of order 2.

Indeed, a birational transformation $(x,y) \leftrightarrow (r,s)$ sending the biquadratic
curve \eqref{Fbiq} $F(x,y)=0$ to the canonical cubic  curve, known to have genus 1,  $s^2=R(r)$, where $R$ is a
polynomial of third degree, see Appell \& Goursat \cite[p.292]{Appell}: from \eqref{ydirect}, choose
$w=$ a square root of $P(x)$, so that $y=(-X_1(x)+w)/(2X_2(x)) \leftrightarrow
w= X_1(x)+2y X_2(x)$, and $x=z_1+1/r$, where $z_1$ is one of the four roots
of $P(x)=0$, and take $s=wr^2$.
\begin{equation}\label{Appell}
\left\{ \begin{matrix} x&=&z_1+1/r, &  \  \   & r&=&1/(x-z_1), \\
            y&=&(-X_1(z_1+1/r)+s/r^2)/(2X_2(z_1+1/r)) & \  \   &
            s&=& (X_1(x)+2y X_2(x))/(x-z_1)^2  \end{matrix} \right. \end{equation}

Then, with $P(z_1+1/r)=R(r)/r^4$, $s^2=R(r)$ of third degree,
 the Weierstrass representation holds $r=A+B\wp(at+b), s=C\wp'(at+b)$.

Note that $R(r)= r^4P(z_1+1/r)= P'(z_1)r^3+P''(z_1)r^2/2+P'''(z_1)r/6+\epsilon$.

So, $x=z_1+1/(A+B\wp), y=(-X_1+C\wp'/(A+B\wp)^2)/(2X_2)$, which are our elliptic functions $\mathcal{E}_1$ and $\mathcal{E}_2$! See also \cite{jjqell}. But how to show that the arguments of $x_0, x_1,\dotsc$ make an arithmetic progression?

Let $(r_m,s_m)$ on the curve $s^2=R(r)$ correspond to $(x_m,y_m)=(\mathcal{E}_1(t_m),\mathcal{E}_2(t_m))$ through \eqref{Appell}. The addition theorem of elliptic functions in this cubic curve setting tells that the line joining a point $(r^*,s^*)$ of argument $t^*$ on the cubic curve  to $(r_m,s_m)$ meets the curve at a point of argument $-(t^*+t_m)$ (up to integer combinations of periods)  \cite[chap. XII, \S 226]{Appell},  \cite[\S III.2.]{Silverman}.

One must show that the line joining $(r_{m+1},-s_{m+1})$ to $(r_{m},s_{m})$ must meet the cubic curve at $(r^*,s^*)$ independent of $m$, representing the step.

\begin{figure}[htbp]\begin{center}
\psset{xunit=10cm,yunit=0.5cm}
\begin{pspicture}(-0.2,-5)(1,7.5)
\psline{->}(-0.2,0)(1.05,0)\psline{->}(0,-5)(0,5)
\uput[90](1.05,0){$r$}\uput[0](0,5){$s$}
\psplot{-0.02}{0.285}{ x 0.333  sub x 1 sub mul 0.2857 x  sub mul 443.8 mul  sqrt }
\psplot{0.34}{1}{ x 0.333  sub x 1.001 sub mul 0.2857 x  sub mul 443.8 mul  sqrt }
\psplot{0.1}{0.285}{ x 0.333  sub x 1 sub mul 0.2857 x  sub mul 443.8 mul  sqrt -1 mul}
\psplot{0.34}{1}{ x 0.333  sub x 1.001 sub mul 0.2857 x  sub mul 443.8 mul  sqrt -1 mul}
\pscurve(0.34,0.3)(0.33,0)(0.34,-0.3) \pscurve(0.99681, -0.81683)(1,0)(0.99681, 0.81683) \pscurve(0.285,0.1)(0.286,0)(0.285,-0.1)
\psline(0.5, 2.8151)(0.76657, -4.6456)(0.76657, 4.6456)(0.99681, -0.81683)(0.99681, 0.81683)(0.83802, 4.4763)(0.83802, -4.4763)
\psdot[dotstyle=*](0.5, 2.8151)\uput[90](0.5, 2.8151){$(r_0,s_0)$}
\psdot[dotstyle=*](0.76657, -4.6456)
\psdot[dotstyle=*](0.76657, 4.6456)\uput[90](0.76657, 4.6456){$(r_1,s_1)$}
\psdot[dotstyle=*](0.99681, -0.81683)
\psdot[dotstyle=*](0.99681, 0.81683)\uput[90](1.05, 0.81683){$(r_2,s_2)$}
\psdot[dotstyle=*](0.83802, 4.4763)
\psdot[dotstyle=*](0.83802, -4.4763)\uput[315](0.83802, -4.4763){$(r_3,s_3)$}
\psline[linestyle=dashed,dash=2mm 2mm](0.0018, 6.457)(0.99681, -0.81683)\psdot[dotstyle=*](0.0018, 6.457)
\psline[linestyle=dashed,dash=2mm 2mm](0.0018,6.457)(0.83802, 4.4763)
\end{pspicture}
\caption{Successive grid points $(r,s)$ on the cubic curve.}\label{elladd}
 \end{center}\end{figure}
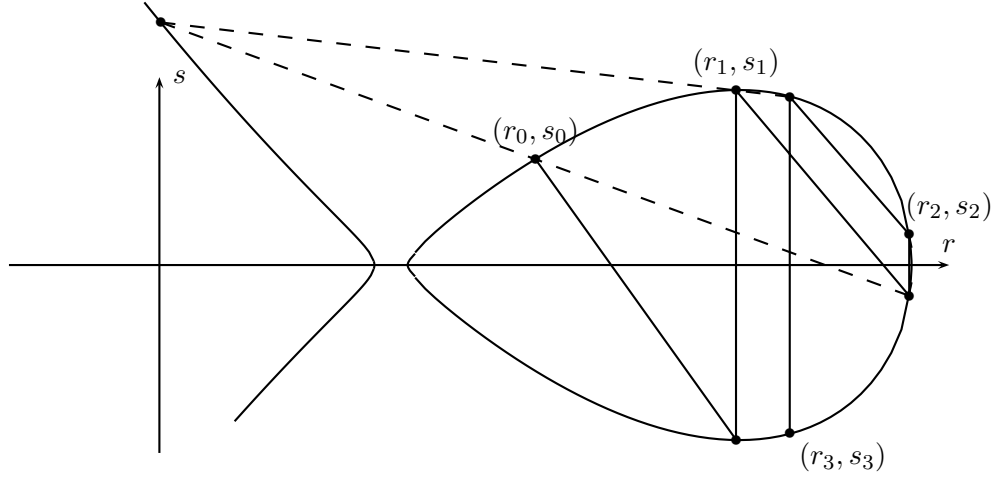

Numerical confirmation has been tried with the biquadratic $F$ of \S  \ref{labo}, $(x_0,y_0)=(0,-0.2716),   \\   (x_1,y_1)=(-0.6955,-0.6220)$ etc. giving by \eqref{Appell} $(r_0,s_0)=(0.5,2.815), (r_1,s_1)=(0.7666,4.646)$, using here $z_1=-2, R(r)=r^4 P(-2+1/r)=-443.8(r-2/7)(r-1/3)(r-1)$.

On fig. \ref{elladd}, the lines joining $(r_{m+1},-s_{m+1})$ to $(r_{m},s_{m})$ seem almost parallel! but the calculations show indeed a single point $(-1.4126, 56.345)$ of course too far. However, looking for the double step,  joining now $(r_{m+2},-s_{m+2})$ to $(r_{m},s_{m})$, we see the meting point $(0.00181,6.4572)$, the dashed lines in fig. \ref{elladd}.

Proof that the line joining $(r_{m+1},-s_{m+1})$ to $(r_{m},s_{m})$
does meet the cubic curve at $(r^*,s^*)$ independent of $m$.

We enter $s=\mu_m (r-r_m)+s_m$ in $s^2-R(r)=P'(z_1)(r-r_{m})(r-r_{m+1})(r-r^*)=0$, where $r^*$ is the unknown.

$r_{m}+r_{m+1}+r^*=(\mu^2_m-P''(z_1)/2)/P'(z_1)$,  
with $\mu_m =\dfrac{-s_{m+1}-s_{m} }{r_{m+1}-r_{m} }   \\  =\dfrac{r_{m+1}^2 (X_1(z_1+1/r_{m+1})+2y_m X_2(z_1+1/r_{m+1}))- r_{m}^2(X_1(z_1+1/r_m)+2y_m X_2(z_1+1/r_m))  }{ r_{m+1}-r_{m}}   \\
=  (r_{m+1}+r_{m}) w(z_1,y_m)  \\  +(\partial w(x,y_m)/ \partial x)(x=z_1)  =
 \dfrac{-4X_2(z_1)   (\partial F(x,y_m)/ \partial x)_{x=z_1}}{w(z_1,y_m)}+(\partial w(x,y_m)/ \partial x)_{x=z_1} $,
where $w(x,y)=X_1(x) +2yX_2(x)=\partial F/\partial y$, and
using from \eqref{ydirect}

 $P(z_1)=X_1^2-4X_0X_2=0\Rightarrow F(z_1,y)=\dfrac{[2X_2(z_1)y+X_1(z_1)]^2}{4X_2(z_1)}=\dfrac{w^2(z_1,y)}{4X_2(z_1)},\\
r_{m+1}+r_{m}=\dfrac{1}{x_{m+1}-z_1}+\dfrac{1}{x_{m}-z_1}=\dfrac{- (\partial F(x,y_m)/ \partial x)_{x=z_1}}{F(z_1,y_m)=w^2(z_1,y)/(4X_2(z_1))}$,    and also
 $y_{m}+y_{m+1} =-X_1(x_{m+1})/X_2(x_{m+1})$ from \eqref{sumprod}, so that $-s_{m+1}= (X_1(x_{m+1})+2y_m X_2(x_{m+1}))/(x_{m+1}-z_1)^2$,

$P'(z_1) r^*+-P''(z_1)/2= -P'(z_1)(r_{m}+r_{m+1})+\mu^2_m   \\   =4X_2(z_1)P'(z_1)\dfrac{ (\partial F(x,y_m)/ \partial x)_{x=z_1}}{w^2(z_1,y_m)}+\left[ 4X_2(z_1) \dfrac{- (\partial F(x,y_m)/ \partial x)_{x=z_1}}{w(z_1,y_m)}+(\partial w(x,y_m)/ \partial x)_{x=z_1} \right]^2$

$=4X_2(z_1) (\partial F(x,y_m)/ \partial x)_{x=z_1}\dfrac{P'(z_1)  +4X_2(z_1) (\partial F(x,y_m)/ \partial x)_{x=z_1}}{w^2(z_1,y_m)}
  \\   -8X_2(z_1) \dfrac{ (\partial F(x,y_m)/ \partial x)(\partial w(x,y_m)/ \partial x)_{x=z_1} }{w(z_1,y_m)}
+(\partial w(x,y_m)/ \partial x)^2_{x=z_1}
=4X_2(z_1) (X'_0+X'_1y+X'_2y^2)_{x=z_1}   \\   \times\dfrac{ 2X'_1+X'_2(2X_2y_m-X_1)/X_2-2(  X'_1+2y_mX'_2 )= X'_2[ -2y -X_1/X_2]=-X'_2w/X_2}{w(z_1,y_m)}  \\
+(X'_1+2yX'_2)^2_{x=z_1}  = -4X'_0(z_1)X'_2(z_1)+(X'_1)^2(z_1)$,

using $P'(z_1)  +4X_2(z_1) (\partial F(x,y_m)/ \partial x)_{x=z_1}=2X_1X'_1-4X_0X'_2-4X'_0X_2+4X_2(X'_2y^2+X'_1y+X'_0)=2X'_1\left(X_1+2yX_2= w\right)+4X'_2[X_2y^2-X_0=(2X_2y-X_1)w/(4X_2)]$, as $X_0(z_1)=X_1^2(z_1)/(4X_2(z_1))$.

This ends the proof that $r, s$ of the step is independent of $m$ on the cubic curve. A very tedious and uninspiring proof, sorry.

\smallskip

\subsubsection{Proof of $4\Rightarrow 1$.}

Let the elliptic function $\mathcal{E}_1$ : $x_n=\mathcal{E}_1(nh)$ and consider the function $\mathcal{E}_2(t)=\mathcal{E}_1(t+h/2)$. These two elliptic functions of order 2 and sharing same periods must be related by a biquadratic equation $F(\mathcal{E}_1,\mathcal{E}_2)=0$ \cite[\S 20.54]{WW}, so, $F(x_n,y_n)=0$, with $y_n=\mathcal{E}_2(nh)$.
This short proof is incomplete, as it does not establish $F(x_{n+1},y_n)=0$...

We have $x_m=\mathcal{E}(mh)$ with $\mathcal{E}$ of order 2. Let $P$ be the quartic polynomial vanishing at the four values $z_1,\dotsc, z_4$ of $\mathcal{E}$ where the derivative $d\mathcal{E}/dt$ vanishes in  a fundamental parallelogram. The ratio $(d\mathcal{E}(t)/dt)^2/P(\mathcal{E}(t))$ has no singularity, so is a constant, so

{$\displaystyle t=\text{\ constant\ }+\text{\ constant\ }\int^{x=\mathcal{E}} \dfrac{du}{\sqrt{P(u)}}  $.    }

First, we find $\alpha, \beta, \gamma$ in $\xi=\dfrac{\alpha x+\beta}{1+\gamma x}$ 
of \eqref{rat1} so to send the four
zeros of $P$ on $\{-k^{-1},-1,1,k^{-1}\}$ as in \S~\ref{essential},

$mh=\int_{\xi_0}^{\xi_m} \dfrac{du}{\sqrt{(1-u^2)(1-k^2u^2)}}  $, i.e., $\xi_m=\xi_0+\sn(mh)$ with the Jacobian ellliptic sine function, and we \emph{define} here $\eta_m=\sn((m+h/2)$.



Consider now these Jacobi functions identities
\begin{equation} \label{Jacadd}
\sn(s+a)+\sn(s-a)=\dfrac{2 \cn a \,\dn a\; \sn s }{1-k^2 \sn^2a\,\sn^2 s},
\sn(s+a)\sn(s-a)=\dfrac{ \sn^2 s-\sn^2 a }{1-k^2 \sn^2a\,\sn^2 s},\end{equation}
 \cite[ell table XIV p.208]{Akhell}, and take $a=h/2, s=(m+1/2)h$, so

$\xi_m+\xi_{m+1}=2\dfrac{2 \cn a \,\dn a\; \eta_m }{1-k^2 \sn^2a\,\eta_m^2 },
\xi_m\xi_{m+1}=\dfrac{\eta_m^2-\sn^2 a }{1-k^2 \sn^2a\,\eta_m^2 }$
reproducing  \eqref{sumprod}  (with $\xi_0=0$),

and we build the biquadratic $\mathcal{F}(\xi,\eta)=\xi^2  \mathcal{Y}_2(\eta)+    \xi \mathcal{Y}_1(\eta)+\mathcal{Y}_0(\eta)$

with $\mathcal{Y}_0(\eta)=\eta^2-\sn^2 a $, $\mathcal{Y}_1(\eta)=-2 \cn a \,\dn a\; \eta$,
$\mathcal{Y}_2(\eta)=1-k^2 \sn^2a\,\eta^2$, or

\begin{equation}\label{snacnadna}
\mathcal{F}(\xi,\eta)= -k^2 \sn^2a \xi^2\eta^2
+\xi^2+\eta^2
-2 \cn a \,\dn a\; \xi\eta
-\sn^2 a=0\end{equation} cf.   Baxter \cite[eq. 15.10.3]{Baxter}, Maillard \cite[eq. (13) p.294]{Maillard}.


Given the symmetric $\mathcal{F}(\xi,\eta)$ of above, starting with
$     \xi_0=0, \eta_0=\pm \sn a, \\   \xi_1=\pm 2 \cn a \,\dn a\; \sn a /(1-k^2 \sn^4 a)=\sn h$ etc.

If we start with the other $\eta-$root $-\sn a$, the step is $-h$,
$\cn a \,\dn a$ is a square root of $(1-\sn^2 a)(1-k^2\sn^2 a)$. What happens if the wrong square
root is used? Answer: add a half period to $h\to h+2K$  \cite[table XII p. 206]{Akhell}.

\medskip

This is the fast and learned proof of $4\Rightarrow 1$.  \qed

%
\medskip
 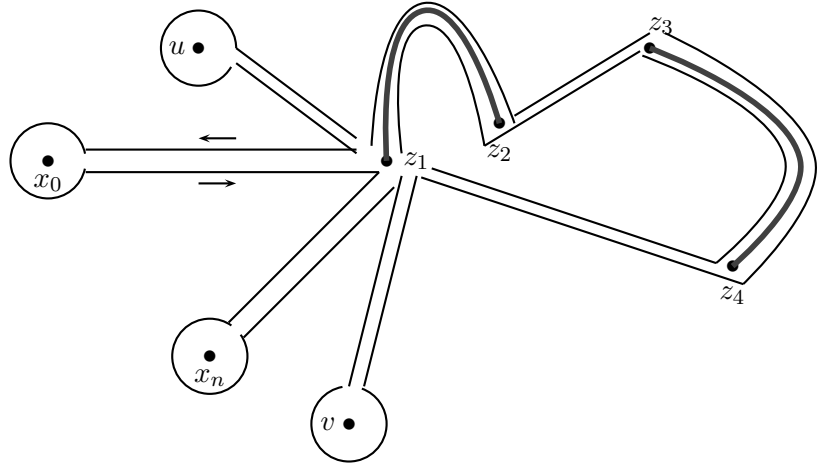
\begin{figure}[htbp]\begin{center}
\psset{unit=1cm}
\begin{pspicture}(-2.5,0)(8,7)
  \psarc(-1,4.5){0.5}{11}{348}
  \psdots[dotsize=0.15](3.5,4.5)     \uput[0](3.55,4.5){$z_1$}
  \psline{<-}(1,4.8)(1.5,4.8)
  \psline{<-}(1.5,4.2)(1.0,4.2)
  \psline(-0.5,4.65)(3.1,4.65)
  \psline(-0.5,4.35)(3.4,4.35)  \psdots[dotsize=0.15](-1,4.5) \uput[270](-1,4.5){$x_0$}
  \psdots[dotsize=0.15](5,5) \uput[270](5,4.9){$z_2$}
  \psdots[dotsize=0.15](7,6) \uput[0](6.8,6.3){$z_3$}
  \psdots[dotsize=0.15](8.1,3.1)
  \psline(3.4,4.35)(1.4,2.35)\psline(3.6,4.15)(1.6,2.15)
 \psarc(1.15,1.9){0.5}{55}{395}\psdots[dotsize=0.15](1.15,1.9)\uput[270](1.15,1.9){$x_n$}
\psline(3.1,4.6)(1.5,5.8) \psline(3.1,4.8)(1.5,6.0)
   \psdots[dotsize=0.15](1.0,6)\psarc(1.0,6){0.5}{0}{330}
   \uput[180](1,6){$u$}
\psline(3.7,4.3)(3,1.5) \psline(3.9,4.3)(3.2,1.5)
   \psdots[dotsize=0.15](3.0,1)\psarc(3.0,1){0.5}{100}{430}
  \pscurve(3.3,4.7)(4,6.6)(5.2,5)\pscurve(3.7,4.6)(4,6.3)(4.8,4.7)\psline(5.2,5.1)(7.0,6.2)\psline(4.8,4.7)(6.87,5.96)
\pscurve(8.24,2.86)(9.2,4.5)(7.2,6.2)\pscurve(7.88,3.14)(8.8,4.5)(6.9,5.9)\psline(8.24,2.86)(3.96,4.28)\psline(7.88,3.14)(4.1,4.4)
   \uput[180](3,1){$v$}
 \psset{linewidth=0.075}
 \pscurve[linecolor=darkgray](3.5,4.5)(4,6.5)(5,5)
 \pscurve[linecolor=darkgray](8.1,3.1)(9,4.5)(7,6)  \uput[270](8.1,3.0){$z_{4}$}
  \end{pspicture}
\caption{Appropriate contour.}\label{ellPell}
 \end{center}\end{figure}

   \smallskip

\noindent


On the other hand,
we look at the integral of $\dfrac{\ \  \log\dfrac{\tilde{A}_{m-1}(x) +B_{m-1}(x)\sqrt{P(x)}}{   \tilde{A}_{m-1}(x)-B_{m-1}(x)\sqrt{P(x)}}dx \ \ }  {2\pi i\sqrt{P(x)} }$ on a big contour, where the polynomials $\tilde{A}_{m-1}, B_{m-1}$
come from the Pell equation \eqref{Pell}.
The integral must vanish, as we integrate $O(x^{-2})$ on the big contour $\to\infty$, and the singularities which are the zeros of $P$ and $\tilde{A}_{m-1}\pm B_{m-1}\sqrt{P}$, i.e. $z_1,\dotsc, z_4, x_0, x_m, u, v$ are in the bounded complex plane. We define the square root of $P$ as continuous far from the singular points just described,
outside two arcs joining $[z_1,z_2]$ and $[z_3,z_4]$.

We now make the contour shrink about these arcs and the singular points, one finds

\begin{equation}  \label{Jacobiprob}
 0=-(m+1) \int_{z_1}^{u} \frac{dx}{\sqrt{P(x)}}-m \int_{z_1}^{v} \frac{dx}{\sqrt{P(x)}}- \int_{z_1}^{x_0} \frac{dx}{\sqrt{P(x)}}
      \pm \int_{z_1}^{x_m}
            \frac{dx}{\sqrt{P(x)}}
+  N\int_{z_2}^{z_{3}}
            \frac{dx}{\sqrt{P(x)}},
\end{equation}
where the first terms account for the zero of multiplicity $m+1,m$ at $u,v$
of $\tilde{A}_{m-1}(x)-B_{m-1}(x)\sqrt{P(x)}$, so that
the logarithm of the ratio on the lower path from $z_0$ to $u$ is the value
on the upper path minus $m+1$ and $m$;
 the $\pm$ signs tell if $x_m$ is a zero of $\tilde{A}_{m-1} +$ or $-B_{m-1}\sqrt{P}$; and where
$N$ is an unknown integer related to the increase of the imaginary part of the logarithm on the two sides. There is only one  integral
$\displaystyle{ \int_{z_2}^{z_{3}} \frac{dx}{\sqrt{P(x)}} }$ to
consider, in the ``gap'' between the arcs.  Integrals on the two sides of the arcs
vanish, as $\sqrt{P}_+ = -\sqrt{P}_-$, so that the logarithm takes
opposite values on the sides.

\smallskip

I learned this technique in Nuttall \& Singh \cite[lemma~5.2]{NuttallS}
allowing polynomials $P$ of degree $2l$ (hyperelliptic case if $l>2$).
The full Abelian case is considered by Nuttall too\cite{Nuttall}.

\medskip

\subsection{A brief history \cite{BursZhed1,BursZhed,Maillard,SpZ2007}\label{ellhist} } Elliptic lattices were developed by Baxter in the solution
of special problems of statistical physics, they appear in works
by Fritz John, in many treatments of a Poncelet problem
\cite{BursZhed1,BursZhed} \cite[\S~6]{SpZ2007}, and go back to
pioneering work by Euler\footnote{and perhaps even to Fermat
[communicated by R.~Askey]!} on the addition formulas of elliptic
functions, that's why the symmetric biquadratic polynomial
\eqref{Ebiq} has been called the Euler polynomial in
\cite[p.~294]{SpZ2007}.

\medskip

Even the name of our subject is not easy to choose: ``elliptic
sequences'' seems perfect, but this name is used by other
sequences related in another way to elliptic functions (sequences
$\{A_n\}$ where $A_{n-1}A_{n+1}/A_n^2$ is our $x_n$, \cite{VdP}),
``elliptic lattice'' may by used for the repetitions of the
periods parallelogram of an elliptic function, ``elliptic grid''
means a convenient mesh for discretizing over ellipses, and
``elliptic difference operator'' is a partial difference operator
extending partial differential operator of elliptic type.
``Elliptic curve'' is best, as seen in the $(r_m,s_m)$ of fig. \ref{elladd}, but this name
suggests nowadays highbrow features of modern number theory.


\medskip

\section{The divided difference operator.\label{divdiff}}
\chead{\thesection  \ \ \ \   Divided difference operator }

We define $\mathcal D$ of a function $f$ defined on the $x_n$s  as the function of $y$
involving the difference of the values of $f$ at the
 two $x-$roots of $F(x,y)=0$ for the given $y$. The value at
this $y$
of $\mathcal{D}f$ will often be noted as the
somewhat awkward $(\mathcal{D} f)(y)$.

Considering the preceding section on lattices, or grids,
we see that $f$ must be known on at least a relevant
sequence $\{x_n\}$, so that $\mathcal{D}f$ is automatically
defined on the companion lattice $\{y_n\}$:

\[   (\mathcal{D} f)(y_n) = \dfrac{f(x_{n+1})-f(x_n)}{x_{n+1}-x_n}.   \]

\begin{table}[htbp]
\begin{center}\begin{tabular}{|c|c|}\hline
$f(x)$ & $(\mathcal{D}f)(y)$ \\  \hline
1 & 0 \\   $x$ & 1 \\  $x^2$ & $-Y_1(y)/Y_2(y)$ \\
$x^3$ & $Y_1^2(y)/Y_2^2(y) -Y_0(y)/Y_2(y)$ \\
$\dfrac{1}{x-A}$ & $-\dfrac{Y_2(y)}{F(A,y)}$
 \\
\hline
\end{tabular}
\caption{Some  instances of the difference operator $\mathcal{D}$.\label{difftable2}} \end{center}\end{table}

The simplest nontrivial example if $f(x)=x^2$, then $\mathcal{D}f$ at $y_n$
is $x_n+x_{n+1} = -Y_1(y_n)/Y_2(y_n)$ from \eqref{sumprod}.


Exercise Find $\mathcal{D}(x-x_\alpha)(x-x_{\alpha+1})$. The result has the factor $y-y_\alpha$.
What is the second factor?The author has no elegant solution.

Exercise. Find $f$ such that $\mathcal{D}f=y$. Answer: not elementary!  See \S \ref{DfR}.

The difference operator applied to a simple rational function is of
special interest. Let $\displaystyle f(x)=\dfrac{1}{x-A}$, then
$\displaystyle \mathcal{D} \dfrac{1}{x-A}$ at  $y_n$ is $\dfrac{1}{x_{n+1}-x_n}\left[
   \dfrac{1}{x_{n+1}-A} -\dfrac{1}{x_n-A} \right] =
 -\dfrac{1}{(x_n-A)(x_{n+1}-A)}  \\
= -\dfrac{Y_2(y)}{Y_0(y)+A Y_1(y)+A^2 Y_2(y)}$,
and let
$\{(x'_n,y'_n), (x'_n,y'_{n+1}),\dotsc\}$ be the  elliptic sequence on the
biquadratic curve $F(x,y)=0$ such that $x'_0=A$, then

\begin{equation} \label{simplerat}
(\mathcal{D} \dfrac{1}{x-x'_0})(y)
=-\dfrac{Y_2(y)}{X_2(A)(y-y'_0)(y-y'_{-1})}
\end{equation}

The divided difference operator useful up to the ANSUW
case can be recovered from the simple property that polynomials
are sent to polynomials of lower degree \cite{Segovia}. But
polynomials are not distinguished rational functions, as $\infty$
is now an ordinary point.

\subsection{\label{doubling}}\textbf{Proposition.}\textsl{ Let $\{x_n\}$ and $\{y_m\}$ be two
sequences related by $F(x_n,y_n)=0$ and $F(x_n,y_{n-1})=0$
for all $n$. Let $\dfrac{f(x_{n+1})-f(x_n)}{x_{n+1}-x_n}$
be a rational function of degree $\leqslant 2\delta$ of the
variable $y_n$ when $f$ is a rational function of degree $\delta$.
Then, $F$ must be biquadratic and the sequences must be elliptic
lattices}.

The proof is easier than the statement.We try $f(x)=(x-A)^{-1}$
as before, so, $(x_n-A)(x_{n+1}-A)$ is a rational function
of degree $\leqslant 2$ of $y_n$, and quadratic in  $A$, say, \\
$\dfrac{N(y_n,A)}{D(y_n)}= \dfrac{Y_0(y_n)+AY_1(y_n)+A^2Y_2(y_n)}{Y_2(y_n)}$,  so,
$x_nx_{n+1}$ and $x_n+x_{n+1}$ are two different rational
functions of $y_n$, and we recover $F(x,y)=Y_0(y)+xY_1(y)+x^2Y_2(y)$
of \eqref{Fx}. \qed

\medskip

\emph{Special cases.} We already encountered the usual difference operators
$(x,x+h)$ or $(x-h,x)$ or $(x-h/2,x+h/2)$ corresponding to
$Y_2(y)\equiv 1$, $Y_1$ of degree 1, $Y_0$ of degree 2 with $Q=Y_1^2-4Y_0Y_2$ of degree 0.
For the geometric difference operator, $Q$ is the square of a first degree polynomial.
For the Askey-Wilson operator \cite{AAR,AW319,Ism2005,Koe,Segovia,magsnul}, $Q$ is an arbitrary second degree
polynomial.

\subsection{Operations on elliptic products} \   \\


We will consider sums of elliptic products $\sum c_m \dfrac{(x-x_0)\dotsb(x-x_{m-1})}{(x-x'_0)\dotsb(x-x'_{m-1})}$.
If, by chance, $c_m$ shows a similar form of ratio of products, we see special
cases of hypergeometric expansions!

\subsubsection{Theorem.\label{Dprod}}
\textsl{Let $\{(x_n,y_n)\}, \{(x'_n,y'_n)\}$ be two elliptic lattices on the same
$F-$curve, then}
\begin{equation} \label{dY}
 \mathcal{D}\dfrac{(x-x_r)\dotsb(x-x_{r+m-1}) }{(x-x'_s)\dotsb (x-x'_{s+m-1})}
=
     C_{m,r,s} Y_2(y) \dfrac{(y-y_r)\dotsb(y-y_{r+m-2}) }{(y-y'_{s-1})\dotsb (y-y'_{s+m-1})},
\end{equation}

\textsl{where $C_{m,r,s}$ is a constant}.

Proof: the operator $\mathcal{D}$  applied to each simple fraction
$1/(x-x'_j)$ yield fractions divided by  \\  $(y-y'_{j-1})(y-y'_j)$ by
\eqref{simplerat},
the same operator acting  on $\dfrac{(x-x_r)\dotsb(x-x_{r+m-1}) }{(x-x'_s)\dotsb (x-x'_{s+m-1})}$ is therefore a  rational function of degree $m+1$ and poles
$y'_{s-1},\dotsc, y'_{s+m-1}$.

We now look at some $y=y_n$ involving the values of $\dfrac{(x-x_r)\dotsb(x-x_{r+m-1}) }{(x-x'_s)\dotsb (x-x'_{s+m-1})}$ at $x_n$
and $x_{n+1}$. We see that
$x_n-x_r$ and $x_{n+1}-x_{r+1}$   both vanish at $y=y_r$, etc.
up to  $x_n-x_{r+m-2}$ and $x_{n+1}-x_{r+m-1}$ both vanishing
at $y=y_{r+m-2}$,
  so that $\dfrac{(x-x_r)\dotsb(x-x_{r+m-1}) }{(x-x'_s)\dotsb (x-x'_{s+m-1})}$ at $x_n$ and at $x_{n+1}$ both vanish at
 $x_n=x_r,\dotsc, x_{r+m-2}$ and
 $x_{n+1}=x_{r+1},\dotsc, x_{r+m-1}$, corresponding to $y=y_r,\dotsc, y_{r+m-2}$ .
  \qed

Of course, $C_{0,r,s}=0$, we now look at $m=1$:



 {\small $\mathcal{D}\dfrac{x-x_r }{x-x'_s}  =\mathcal{D}\dfrac{x'_s-x_r }{x-x'_s} $ at $y_n
=  \dfrac{x_r-x'_s}{(x_n-x'_s)(x_{n+1}-x'_s)}= \dfrac{(x_r-x'_s)Y_2(y_n)}
{F(x'_s,y_n)=X_2(x'_s)(y_n-y'_{s-1})(y_n-y'_{s})}$,} must be $\dfrac{C_{1,r,s}}{(y_n-y'_{s-1})(y_n-y'_s)}$, so, at $n=r, \dfrac{1}{x_{r+1}-x'_s}=\dfrac{ C_{1,r,s}Y_2(y_r)}{(y_r-y'_{s-1})(y_r-y'_s)}$,

\begin{equation}\label{C1}C_{1,r,s}=\dfrac{(y_r-y'_{s-1})(y_r-y'_s)}{(x_{r+1}-x'_s)Y_2(y_r)}=\dfrac{x_r-x'_s}{X_2(x'_s)}.
\end{equation}

One finds $C_{2,r,s}/C_{1,r,s}=
\dfrac{(y_r-y'_s)(y'_{s+1}-y'_{s-1})X_2(x'_s)}
     {(x_r-x'_s)(x'_{s+1}-x'_s)Y_2(y_r) }
= -\dfrac{(y'_{s+1}-y'_{s-1})(x'_{s}-x_{r+1})}
     {(x'_{s+1}-x'_s)(y_r-y'_{s-1}) }
$

\smallskip

The constant $C_{m,r,s}$ is found through particular values of $y_n$, either $y_{r-1}$,
where

$(x_{n+1}-x_r)(x_{n+1}-x_{r+1})\dotsb  (x_{n+1}-x_{r+k-1})=0$, or $y_{r+k-1}$,
where \\ $(x_n-x_r)(x_n-x_{r+1})\dotsb  (x_n-x_{r+k-1})=0$:

\begin{subequations}

\begin{equation}\label{Cxm1}\begin{split}
C_{m,r,s} &= -\dfrac{1}{x_r-x_{r-1}}
\dfrac{(y_{r-1}-y'_{s-1})\dotsb (y_{r-1}-y'_{s+m-1})
         (x_{r-1}-x_r)\dotsb(x_{r-1}-x_{r+m-1}) }
          {Y_2(y_{r-1}) (y_{r-1}-y_r)\dotsb(y_{r-1}-y_{r+m-2})
          (x_{r-1}-x'_s)\dotsb (x_{r-1}-x'_{s+m-1})} \\
&= \dfrac{(y_{r-1}-y'_{s+m+1})\dotsb (y_{r-1}-y'_{s+m-1})   (x_{r}-x'_s)
         (x_{r-1}-x_{r+1})\dotsb(x_{r-1}-x_{r+m-1}) }
          {X_2(x'_s) (y_{r-1}-y_r)\dotsb(y_{r-1}-y_{r+m-2})
         (x_{r-1}-x'_{s+1})\dotsb (x_{r-1}-x'_{s+m-1})}
\end{split}\end{equation}
using $F(x'_s,y_{r-1})=Y_2(y_{r-1}) (x_{r}-x'_s)   (x_{r-1}-x'_s)   =X_2(x'_s)
(y_{r-1}-y'_{s+m-1})(y_{r-1}-y'_s)$
\begin{equation}
  C_{m,r,s}   \\   =\dfrac{
              (y_{r+m-1}-y'_s)\dotsb (y_{r+m-1}-y'_{s+m})
               (x_{r+m}-x_r)\dotsb(x_{r+m}-x_{r+m-1})
             }
          {Y_2(y_{r+m-1}) (y_{r+m-1}-y_r)\dotsb(y_{r+m-1}-y_{r+m-2}) (x_m-x_{r+m-1})
          (x_{r+m}-x'_{s+1})\dotsb (x_{r+m}-x'_{s+m})}
\label{Cxn1}
\end{equation}
(Of course, $C_0=0$).

Or through  residues at $y=y_n=y'_{s-1}$, seeing that $x_{n+1}-x'_s=0$,
replacing $(y-y'_{s-1})/(x-s'_m)$ by $dy/dx=-(x'_s-x'_{s-1})Y_2(y'_{s-1})/((y'_{s-1}-y'_s)X_2(x'_s))$
 of the tangent at $(x'_s,y'_{s-1})$ from \eqref{tang}
\begin{equation}\label{Cxp0}\begin{split}
C_{m,r,s} &= \dfrac{1}{x'_s-x'_{s-1}}\;\dfrac{(x'_s-x_r)\dotsb (x'_s-x_{r+m-1})dy(y'_{s-1}-y'_s)\dotsb(y'_{s-1}-y'_{s+m-1})}
          {dx(x'_s-x'_{s+1})\dotsb (x'_s-x'_{s+m-1})(y'_{s-1}-y_r)\dotsb (y'_{s-1}-y_{r+m-2})Y_2(y'_{s-1})} \\
&= -\dfrac{1}{y'_{s-1}-y'_s}\;\dfrac{(x'_s-x_r)\dotsb (x'_s-x_{r+m-1})(y'_{s-1}-y'_s)\dotsb(y'_{s-1}-y'_{s+m-1})Y_2(y'_{s-1})}
          {(x'_s-x'_{s+1})\dotsb (x'_s-x'_{s+m-1})(y'_{s-1}-y_r)\dotsb (y'_{s-1}-y_{r+m-2})X_2(x'_{s})}
\end{split} \end{equation}
using $Y_2(y'_{s-1}(x'_s-x'_{s-1})dx=X_2(x'_s)(y'_{s-1}-y'_s)dy$.

Finally, through residues at $y=y'_{s+m-1}$,  seeing now that $x_{n}-x'_{s+m-1}=0$,
$(x'_{s+mk}-x'_{s+m-1})C_{m,r,s}= -\dfrac{(x'_{s+m-1}-x_r)\dotsb(x'_{s+m-1}-x_{r+m-1})}
 {x'_{s+m-1}-x'_s)\dotsb (x'_{s+m-1}-x'_{s+m-2}dx}
\dfrac{y'_{s+m-1}-y'_{s-1})\dotsb(y'_{s+m-1}-y'_{s+m-2})dy }
  {(y'_{s+m-1}-y_r)\dotsb(y'_{s+m-1}-y_{r+m-2}) }$.

Following \eqref{tang}, near $(x'_{s+m-1},y'_{s+m-1})$, $0=F(x'_{s+m-1}+dx,y'_{s+m-1}+dy)=dx\partial F/\partial x
+dy\partial F/\partial y+o(dy)$,
$dy/dx=-(2x'_{s+m-1}Y_2(y'_{s+m-1})+Y_1(y'_{s+m-1}))/
(2y'_{s+m-1}X_2(x'_{s+m-1})+X_1(x'_{s+m-1}))   \\   =
-(x'_{s+m}-x'_{s+m-1})Y_2(y'_{s+m-1})/((y'_{s+m-1}-y'_{s+m-2})X_2(x'_{s+m-1}))
$

  \begin{multline}\label{Cxpn}
C_{m,r,s}=\dfrac{(x'_{s+m-1}-x_r)\dotsb(x'_{s+m-1}-x_{r+m-1})}
 {(x'_{s+m}-x'_{s+m-1})(x'_{s+m-1}-x'_s)\dotsb (x'_{s+m-1}-x'_{s+m-2})(y'_{s+m-1}-y'_{s+m-2})X_2(x'_{s+m-1})
}\\
\dfrac{(y'_{s+m-1}-y'_{s-1})\dotsb(y'_{s+m-1}-y'_{s+m-2})(x'_{s+m}-x'_{s+m-1})Y_2(y'_{s+m-1}) }
  {(y'_{s+m-1}-y_r)\dotsb(y'_{s+m-1}-y_{r+m-2}) }
\end{multline}

\end{subequations}

remark from \eqref{Cxm1} and \eqref{Cxp0}
\begin{equation}\label{slopes}
\dfrac{C_{m+1,r,s}}{C_{m,r,s}}=
 \dfrac{ (y_{r-1}-y'_{s+m})(x_{r-1}-x_{r+m})}{(y_{r-1}-y_{r+m-1})(x_{r-1}-x'_{s+m})}=
\dfrac{ (x'_s-x_{r+m})(y'_{s-1}-y'_{s+m})}
          { (x'_s-x'_{s+m})(y'_{s-1}-y_{r+m-1})
}\end{equation}

Check will be done at  section \ref{theta} with theta functions in eq. (\ref{xxptheta}-\ref{yyptheta}).
\medskip

\section{Rational interpolatory elliptic expansions.\label{interp1}}
\chead{\thesection  \ \ \ \   Expansions}

\subsection{Construction of the interpolant} \label{expan1} \   \\


Let $f$ be a function defined on a sequence $x_0, x_1,\dotsc$, which need NOT be elliptic
at the present stage.
Rational interpolants of  $f$ at $x_0$, $x_1,\dotsc,$ with poles at
$x'_0$, $x'_1,\dotsc,$ are successive sums

\begin{equation} \label{expan}
c_0=f(x_0),   \   c_0+c_1 \dfrac{x-x_0}{x-x'_0},  \dotsc, \qquad \sum_0^\infty c_m
\dfrac{(x-x_0)\dotsb(x-x_{m-1}) }{(x-x'_0)\dotsb (x-x'_{m-1})},
\end{equation}
where
$c_0=f(x_0), c_1=(x_1-x'_0)\dfrac{f(x_1)-f(x_0)}{x_1-x_0},\\
c_2=(x_2-x'_0)(x_2-x'_1)\dfrac{f(x_2)-f(x_0)-c_1(x_2-x_0)/(x_2-x'_0)
         }{(x_2-x_0)(x_2-x_1) }
=(x_2-x'_0)(x_2-x'_1)\dfrac{f(x_2)-f(x_0) }{(x_2-x_0)(x_2-x_1)}
\\-c_1\dfrac{x_2-x'_1}{x_2-x_1} =
 \dfrac{x_2-x'_1}{x_2-x_1}\left[  (x_2-x'_0)\dfrac{f(x_2)-f(x_0) }{x_2-x_0} -(x_1-x'_0)\dfrac{f(x_1)-f(x_0)}{x_1-x_0}\right],$

The sum $\sum_0^k c_m (x-x'_{m})\dotsb (x-x'_{k-1})(x-x_0)\dotsb (x-x_{m-1})$
is the polynomial interpolant $I_k$ of degree $k$ to $(x-x'_0)\dotsb (x-x'_{k-1})f(x)$
at $x=x_0,\dotsc, x_k$. In particular, the main coefficient $c_0+\dotsb+c_k$
is the divided difference of order $k$ of $(x-x'_0)\dotsb (x-x'_{k-1})f(x)$
at $x_0,\dotsc, x_k$.

\begin{equation}\label{cdivdiff}
c_0+\dotsb +c_k = [x_0,\dotsc, x_k]\text{ \  of\ } (x-x'_0)\dotsb (x-x'_{k-1})f(x).
\end{equation}

\subsection{The cases of $1/(x-t)$ and $x$  }  \   \\

\emph{Proposition.}
\begin{equation}\label{uxt}
\dfrac{1}{x-t}=\dfrac{1}{x_0-t}+\sum_1^\infty (x'_{n-1}-x_n) \dfrac{(t-x'_0)\dotsb (t-x'_{n-2}) (x-x_0)\dotsb (x-x_{n-1}) }
                  {(x-x'_0)\dotsb (x-x'_{n-1})(t-x_0)\dotsb (t-x_n) }.
\end{equation}
\begin{equation}\label{uxt2}
x=x_0+\sum_{n=1}^\infty (x_n-x'_{n-1}) \dfrac{ (x-x_0)\dotsb (x-x_{n-1}) }
                  {(x-x'_0)\dotsb (x-x'_{n-1}) }.
\end{equation}

See Walsh \cite[\S 8.1, eq. (3)]{Walsh1935}.

Indeed
$c_0=\dfrac{1}{x_0-t}, c_1=\dfrac{x'_0-x_1}{(x_0-t)(x_1-t)}$, $c_0+\dotsb c_n$ is the main coefficient of the
polynomial interpolant of degree $n$ of $(x-x'_0)\dotsb(x-x'_{n-1})/(x-t)$ at $x_0,\dotsb, x_n$, which is
$I_n(x)=(x-x'_0)\dotsb (x-x'_{n-1})\dfrac{1}{x-t}
- \dfrac{(x-x_0)\dotsb (x-x_n) (t-x'_0)\dotsb (t-x'_{n-1})}{(t-x_0)\dotsb (t-x_n)(x-t)}$
(the result must be a rational function with vanishing residue
at $t$), so,
$c_0+\dotsb + c_n= - \dfrac{ (t-x'_0)\dotsb (t-x'_{n-1})}{(t-x_0)\dotsb (t-x_n)}$,
$c_n= \dfrac{ (t-x'_0)\dotsb (t-x'_{n-2})(x'_{n-1}-x_n)}{(t-x_0)\dotsb (t-x_n)}$. \qed


The proof for $x$ is even more elementary: let $S_0=x_0$, $S_N(x)=$ the sum of \eqref{uxt2} up to $n=N$, suppose that $S_N$ interpolates up to $x_N$ (true when $N=0$ and $N=1$). Then, $x-S_N(x)$ is a polynomial of
degree $N+1$ divided by $(x-x'_0)\dotsb(x-x'_{N-1})$, the numerator must be  $(x-x_0)\dotsb(x-x_N)$,

$x-S_{N+1}(x)=x-S_N(x)-(x_{N+1}-x'_N)\dfrac{(x-x_0)\dotsb(x-x_N)}{(x-x'_0)\dotsb(x-x'_{N})}
\\   =(x-x'_N  \ \ -(x_{N+1}-x'_N)=x-x_{N+1})\dfrac{(x-x_0)\dotsb(x-x_N)}{(x-x'_0)\dotsb(x-x'_{N})}$ does indeed vanish at $x=x_{N+1}$.  \qed

The expansion is purely formal, meaning only that a finite sum $n\leqslant N$
achieves interpolation at required points $x_0,\dotsc, x_N$.
The formula is valid for any sequence $x_n,$ etc., not only
for elliptic sequences.

Why spend time on an expansion for $x$? It will be needed in the following,
if $f$ is analytic in a domain containing a contour $C$, $f(x)$ has
has the expansion

\begin{equation}\label{contou}
\int_C \dfrac{f(t)dt/(2\pi i)}{t-x}
=\int_C \left( \dfrac{1}{t-x_0}+
\sum_1^\infty (x_n-x'_{n-1}) \dfrac{(t-x'_0)\dotsb (t-x'_{n-2}) (x-x_0)\dotsb (x-x_{n-1}) }
                     {(x-x'_0)\dotsb (x-x'_{n-1})(t-x_0)\dotsb (t-x_n) }
\right)\dfrac{f(t)dt}{2\pi i}
\end{equation}

Hermite-Walsh formula \cite[eq. (2.5)]{Walsh1932}, \cite[\S 8.1]{Walsh1935},   \cite[Thm 3.6.1]{Davis}.

\medskip

In the $y, y'$ variables,

\begin{equation}\label{uyt}
\dfrac{1}{y-t}=\dfrac{1}{y_0-t}+\sum_1^\infty (y'_{n-1}-y_n) \dfrac{(t-y'_0)\dotsb (t-y'_{n-2}) (y-y_0)\dotsb (y-y_{n-1}) }
                     {(y-y'_0)\dotsb (y-y'_{n-1})(t-y_0)\dotsb (t-y_n) }
\end{equation}
\begin{equation}\label{uyt2}
y=y_0+\sum_{n=1}^\infty (y_n-y'_{n-1}) \dfrac{ (y-y_0)\dotsb (y-y_{n-1}) }
                  {(y-y'_0)\dotsb (y-y'_{n-1}) }.
\end{equation}

\subsubsection{}Recall that  $\mathcal{M}f$ is defined by $( \mathcal{M})(y_n)=(f(x_n)+f(x_{n+1}))/2$,

\begin{equation}\label{My}
 \mathcal{M}\dfrac{(x-x_r)\dotsb(x-x_{r+m-1}) }{(x-x'_s)\dotsb (x-x'_{s+m-1})}
=
     D_{m,r,s}(y)  \dfrac{(y-y_r)\dotsb(y-y_{r+m-2}) }{(y-y'_{s-1})\dotsb (y-y'_{s+m-1})},
\end{equation}
\textsl{ where $D_{m,r,s}$ is a second degree polynomial.}

As in the proof of \ref{Dprod}, we  have a rational function of degree
$m+1$ with the same $m-1$ zeros $y_r,\dotsc,y_{r+m-2}$, whence the factor $D_{m,r,s}(x)$.  \qed

$ D_{0,r,s}$ is such that $\mathcal{M}1=1=
D_{0,r,s}(y)\dfrac{1}{(y-y_{r-1})(y-y'_{s-1})}$, so,

$D_{0,r,s}(y)=(y-y_{r-1})(y-y'_{s-1})$,



$2\mathcal{M}\dfrac{x-x_r }{x-x'_s}$ at $y_n=
\dfrac{x_n-x_r }{x_n-x'_s}+\dfrac{x_{n+1}-x_r }{x_{n+1}-x'_s}
=\dfrac{2x_nx_{n+1} -(x_r+x'_s)(x_n+x_{n+1})+2x_rx'_s}{(x_{n}-x'_s)
( x_{n+1}-x'_s)=F(x'_s,y_n)/Y_2(y_n)}
\\    =
\dfrac{2Y_0(y_n)+(x_r+x'_s)Y_1(y_n)+2x_rx'_sY_2(y_n)}{X_2(x'_s)(y_n-y'_{s-1})(y_n-y'_{s})}$, and
\begin{equation}\label{D1}D_{1,r,s}(y)=\dfrac{Y_0(y)+(x_r+x'_s)Y_1(y)/2+x_rx'_sY_2(y)}{X_2(x'_s)}\end{equation}

The coefficient of $y^2$ is $\dfrac{D''_{1,r,s}}{2}=\dfrac{X_2(0)+(x_r+x'_s)X'_2(0)/2+x_rx'_sX''_2(0)/2}{X_2(x'_s)}  \\   =1+(x_r-x'_s)\dfrac{X'_2(0)+x'_sX''_2(0)}{2X_2(x'_s)}
= 1+(x_r-x'_s)\dfrac{X'_2(x'_s)}{2X_2(x'_s)}
$

and $\rho_{1,r,s}=\dfrac{D''_{1,r,s}/2}{C_{1,r,s}}
= \dfrac{X_2(x'_s)}{x_r-x'_s}+\dfrac{X'_2(x'_s)}{2}
$



\medskip

Interesting values are found at the same point as in \eqref{Cxm1}-\eqref{Cxpn}:
from \eqref{dY}  $\dfrac{D_{n,r,s}(y_\alpha)}{C_{n,r,s} Y_2(y_\alpha)}$ involve sums and
differences of   $\dfrac{(x_{\alpha+1}-x_r)\dotsb(x_{\alpha+1}-x_{r+n-1}) }{(x_{\alpha+1}-x'_s)\dotsb (x_{\alpha+1}-x'_{s+n-1})    }$ and
  $\dfrac{(x_{\alpha}-x_r)\dotsb(x_{\alpha}-x_{r+n-1}) }{(x_{\alpha}-x'_s)\dotsb (x_{\alpha}-x'_{s+n-1})    }$, leaving the ratio
$\dfrac{1/2}{ 1/(x_{\alpha+1}-x_{\alpha})  }$ at $\alpha=r+n-1$ and $s-1$ (where the term
with a division by $x'_{\alpha+1}-x'_s$ is infinitely larger than the other term); and the opposite
values at $\alpha=r-1$ and $s+n-1$  \ \  ($n>0$).

\begin{subequations}
\begin{equation}\label{Dxm1}
D_{m,r,s}(y_{r-1}) = -\dfrac{C_{m,r,s} Y_2(y_{r-1}) (x_r-x_{r-1})}{2},
\end{equation}
\begin{equation}\label{Dxn1}
D_{m,r,s}(y_{r+m-1}) = \dfrac{C_{m,r,s} Y_2(y_{r+m-1}) (x_{r+m}-x_{r+m-1})}{2},
\end{equation}
\begin{equation}\label{Dxp0}
D_{m,r,s}(y'_{s-1}) = \dfrac{C_{m,r,s} Y_2(y'_{s-1}) (x'_{s}-x'_{s-1})}{2},
\end{equation}
\begin{equation}\label{Dxpn}
D_{m,r,s}(y'_{s+m-1}) = -\dfrac{C_{m,r,s} Y_2(y'_{s+m-1}) (x'_{s+m}-x'_{s+m-1})}{2},
\end{equation}

\end{subequations}
when $m>0$.

\bigskip

Let the rational function $P_{m,r,s}(x)=\dfrac{(x-x_r)\dotsb(x-x_{r+m-1}) }{(x-x'_s)\dotsb (x-x'_{s+m-1})}$,
the two left-hand sides of  \eqref{dY} are the difference and the sum
of  $P_{m,r,s}(x)$ at $x=x_n$ and at $x=x_{n+1}$ for some $y=y_n$, so,

 \begin{multline}\label{npn}
P_{m,r,s}(x_{n+1})= (D_{m,r,s}(y_n)+(x_{n+1}-x_n)C_{m,r,s}Y_2(y_n)/2)\dfrac{(y_n-y_r)\dotsb(y_n-y_{r+m-2}) }{(y_n-y'_{s-1})\dotsb (y_n-y'_{s+m-1})},
\\
P_{m,r,s}(x_n)= (D_{m,r,s}(y_n)-(x_{n+1}-x_n)C_{m,r,s}Y_2(y_n)/2)\dfrac{(y_n-y_r)\dotsb(y_n-y_{r+m-2}) }{(y_n-y'_{s-1})\dotsb (y_n-y'_{s+m-1})},
 \end{multline}

the right-hand sides are not rational functions of $y_n$, as
 $Y_2(y_n)(x_{n+1}-x_n)$ is a square root of $Q(y_n)$, from \eqref{ydirect}.

\medskip

The equations \eqref{Dxm1} and \eqref{Dxp0} also show that $D_{m,r,s}$ is the linear interpolant
at $y_{r-1}$ and $y'_{s-1}$ augmented by a constant times
$(y-y_{r-1})(y-y'_{s-1})$

\begin{equation}\label{Dkinterp0}D_{m,r,s}(y)/C_{m,r,s}=R_{r,s}(y)
     +\rho_{m,r,s}(y-y_{r-1})(y-y'_{s-1}),\end{equation}
where
\begin{equation}\label{Dkinterp}R_{r,s}(y)=\dfrac{ (y-y_{r-1})Y_2(y'_{s-1}) (x'_s-x'_{s-1})
+(y-y'_{s-1})Y_2(y_{r-1}) (x_r-x_{r-1})}{2(y'_{s-1}-y_{r-1}) },\end{equation}
 and $\rho_{m,r,s}$ is the coefficient of $y^2$ of $D_{m,r,s}(y)/C_{m,r,s}$.


\begin{equation}\label{rhon}\begin{split}
\rho_{n,r,s}&=\dfrac{(x_{n+r}-x_{n+r-1})Y_2(y_{n+r-1})/2-R_{r,s}(y_{n+r-1})   }
                 {(y_{n+r-1}-y_{r-1})(y_{n+r-1}-y'_{s-1})  } \\
&=\dfrac{-(x'_{n+s}-x'_{n+s-1})Y_2(y'_{n+s-1})/2-R_{r,s}(y'_{n+s-1})   }
                 {(y'_{n+s-1}-y_{r-1})(y'_{n+s-1}-y'_{s-1})  }, \ \ n>0.
\end{split}\end{equation}

and we have, from \eqref{npn} and \eqref{Dkinterp0},

\begin{equation}\label{PP}
\dfrac{P_{m,r,s}(x_{n+1})}{P_{m,r,s}(x_n)}=
\dfrac{R_{r,s}(y_n)   +\rho_{m,r,s}(y_n-y_{r-1})(y_n-y'_{s-1})+(x_{n+1}-x_n)Y_2(y_n)/2}
{R_{r,s}(y_n)  +\rho_{m,r,s}(y_n-y_{r-1})(y_n-y'_{s-1})-(x_{n+1}-x_n)Y_2(y_n)/2}.
\end{equation}

\medskip

\bigskip


\section{Linear first order difference equations.\label{lindiffeq}}
\chead{\thesection  \ \ \ \   Lin diff eq}


We look for linear equations involving $\mathcal{D}f$ and
$(\mathcal{M}f)(y)$, so of the form

\begin{equation} \label{difflin}
A(y)(\mathcal{D} f)(y) = B(y) (\mathcal{M} f)(y)+C(y)
\end{equation}


This already allows elliptic exponentials ($B(y)\equiv A(y)$) or
logarithms ($B(y)\equiv 0$).

We now try to expand a solution to \eqref{difflin} as an interpolatory
series.  This has been done in \cite{magsigma2009} for polynomials $A, B, C$.
Only two special cases allowing hypergeometric expansions are given here

\subsection{"Elliptic logarithm" or elliptic Psi function\label{elllogar}} \   \\

We investigate the solution of
\begin{equation}\label{elllog2}
\mathcal{D}f(y) = \dfrac{Y_2(y)}{y-A}
\end{equation}  on
some $x_n$ lattice, starting with a given value of $f(x_0)$.

On the linear lattice  $F(x,y)=(y-x)(y-x-1)$ where $Y_2(y)\equiv 1$, we have the Psi, or digamma, function $\Psi(x)=d\log\Gamma(x)/dx$ \cite[chap.~6]{AbrS}, \cite[chap.~5]{NIST} etc.

The Psi ($\Psi$) function is the meromorphic function
with poles at $0,-1,\dotsc ,$ satisfying $\Psi(x+1)=\Psi(x)+1/x$
and $\Psi(x)-\log x\to 0$ when $x\to +\infty$.
\cite[chap.~6.3]{AbrS}  \cite[chap. IX]{MilneT}

{\small This function is sometimes found in probability calculus identities, for
instance \cite{x}, let $0<\alpha<1$ and
$\displaystyle M(t)=\sum_n \dfrac{w_n}{1-x_n t}$ which is here $\dfrac{\Gamma(\alpha+1)\Gamma(1-t)}{\Gamma(\alpha-t+1)}
=\sum_1^\infty \Gamma(\alpha+1)\dfrac{(-1)^{n-1}}{(n-1)! (n-t)\Gamma(\alpha+1-n)}
=\dfrac{\alpha}{1-t}-\dfrac{\alpha(\alpha-1)}{2-t}+\dotsb
=\sum_1^\infty \dfrac{\alpha(1-\alpha)(2-\alpha)\dotsb(n-1-\alpha)}{n! (1-t/n)}
$. The expectation is then $\sum x_n w_n= dM(t)/dt$ at $t=0$, so, $-\Psi(1)+\Psi(\alpha+1)$.   }   

We return to \eqref{elllog2}: $f(x_1)=f(x_0)+\dfrac{(x_1-x_0)Y_2(y_0)}{y_0-A}, f(x_2)=f(x_1)+\dfrac{(x_2-x_1)Y_2(y_1)}{y_1-A}
 = f(x_0)+\dfrac{(x_1-x_0)Y_2(y_0)}{y_0-A}+\dfrac{(x_2-x_1)Y_2(y_1)}{y_1-A}$, etc.

\subsubsection{Proposition \label{elllogthm}} \textsl{ Let  $f$ be a solution of the difference equation \eqref{elllog2}
on  the elliptic sequence $\{x_n,y_n\}$ related to the biquadratic polynomial $F(x,y)=x^2Y_2(y)+xY_1(y)+Y_0(y)$,
 and $\{x'_n,y'_n\}$, $n\in\mathbb{Z}$,  another  elliptic sequence
related to the same biquadratic polynomial, such that $y'_{-1}=A$.}

\textsl{Then $f(x)$  has the expansion}
\begin{equation}\label{elllogx} f(x)= \sum_0^\infty c_m\dfrac{(x-x_0)\dotsb(x-x_{m-1}) }{(x-x'_0)\dotsb (x-x'_{m-1})},\end{equation}
\textsl{where $c_0=f(x_0)$ may be chosen at will, and where}
$c_m=
\dfrac{y_{m-1}-y'_{m-1}}{C_{m,0,0} }    $
\textsl{when $m>0$, with the $C_{m,0,0}$ of thm \ref{Dprod}.  }

\subsubsection{\label{simprod}}A simple product (\textbf{hypergeometric}) formula for $c_m$ holds as, from
\eqref{slopes}

$\dfrac{c_{m+1}}{c_m}=\dfrac{y_{m}-y'_{m}}{y_{m-1}-y'_{m-1}}\ \dfrac{C_{m,0,0}}{C_{m+1,0,0}}=
\dfrac{y_{m}-y'_{m}}{y_{m-1}-y'_{m-1}}\   \dfrac{(y_{-1}-y_{m-1})(x_{-1}-x'_{m})}{ (y_{-1}-y'_{m})(x_{-1}-x_{m})}   \\    =
 \dfrac{y_{m}-y'_{m}}{y_{m-1}-y'_{m-1}}\  \dfrac   { (x'_0-x'_{m})(y'_{-1}-y_{m-1})}{ (x'_0-x_{m})(y'_{-1}-y'_{m})} $.

.

Proof:
put the expansion of $f$ in  \eqref{elllog2} and use thm \ref{Dprod} :

$\dfrac{1}{y-A}=(Y_2)^{-1}\mathcal{D}f=\sum_1^\infty c_m C_{m,0,0} \dfrac{(y-y_0)\dotsb(y-y_{m-2})}{(y-y'_{-1})\dotsb(y-y'_{m-1})}$,
multiply by $(y-y'_{-1})(y-y'_0)$:  \\   $(y-y'_{-1})(y-y'_0)/(y-A)=
y-y'_0
=\sum_1^\infty c_m C_{m,0,0} \dfrac{(y-y_0)\dotsb(y-y_{m-2})}{(y-y'_{1})\dotsb(y-y'_{m-1})}$, to compare to the
expansion \eqref{uyt2} with $y'_0$ etc. shifted to $y'_1$ etc.

 $y-y'_0=-y'_0+y_0+\sum_{n=1}^\infty (y_n-y'_{n}) \dfrac{ (y-y_0)\dotsb (y-y_{n-1}) } {(y-y'_1)\dotsb (y-y'_{n}) }$ with $m=n+1$. Check $m=1$:
$c_1=(x_1-x'_0)(f(x_1)-f(x_0))/(x_1-x_0)=(x_1-x'_0)Y_2(y_0)/(y_0-A=y_0-y'_{-1}), c_1C_{1,0,0}=y_0-y'_0$ from \eqref{C1}.
\qed

{\small
\begin{table}[htbp]
\begin{center}\begin{tabular}{|c|ccc|cc|}\hline
 $n$ &       $ x_n$     &      $  y_n$  &        $ f(x_n)$  &        $   f(1)$ &  $  f(-1.75)$  \\   \hline
 0   &     0        &   -0.27156  &        0      &        0      &0\\
 1    &-0.69548   & -0.62200   & -0.92296   &    2.4521    &-1.8068\\
 2   & -0.99680    &-0.65853   & -1.2839    &    -0.33633  &-2.0514\\
 3   & -0.80672    &-0.37599   & -1.0590    &      5.2544  & -2.0803\\
 4   & -0.17524    & 0.12738   & -0.24396   &    -0.94409 & -2.0846\\
 5    & 0.53805    & 0.53701   &  0.78908   &     1.3530  & -2.0859\\
 6    & 0.93095    & 0.69234   &  1.3938     &     1.4958  & -2.0849\\
 7    & 0.98439    & 0.61976   &  1.4773     &     1.5017  & -2.1004\\
 8    & 0.71975    & 0.29750   &  1.0664     &     1.5018  & -1.8982\\
 9    & 0.094233 & -0.21076   &  0.13411    &     1.5018  & -3.7723\\
10  & -0.62568   & -0.59521   & -0.83593    &     1.5018  & -3.7581\\
16   & 0.99403    & 0.64125   &  1.4924      &     1.5018  & -3.7559\\
20   & -0.96796 & -0.67922    & -1.2503     &      1.5018  & -3.7535\\
21   & -0.89935 & -0.47427    & -1.1697     &      1.5018  & -3.7535\\   \hline
\end{tabular}
\caption{Some values of solution to \eqref{elllog2} with $A=y'_{-1}= 7.3839
...$  The two last items of the $n^{\text{th}}$ row are the sums of $n$ terms of the expansion \eqref{elllogx}, convergence has been sketched in \cite{magsigma2009}, see
also \cite{SpiriConvergence2023}.\label{logtable}} \end{center}\end{table}
}   

\bigskip

The interpolatory expansion of $\Psi(x)-\Psi(1)$ above reduces here, with $F(x,y)=(y-x)(y-x+1)$,  $x_n=n+1, x'_n=-n, y_n=n, A=-1,y'_n=-n,  n=0,1,\dotsc$, to
\[  \Psi(x)-\Psi(1)= \sum_1^\infty \dfrac{2(x-1)\dotsb (x-n)}{nx(x+1)\dotsb(x+n-1)} ,  \]

a poorly convergent useless expansion\dots

\subsubsection{\label{DfR}Equation $\mathcal{D}f=R(y) $, with a general  rational function $R$}   Expand
$(y-y'_{-1})(y-y'_0)R(y)/Y_2(y)$ in simple fractions and use (\ref{uyt}-\ref{uyt2}). Higher powers of $1/(y-A)$ and $y$ are not considered here.

\textbf{Exercise.} Extend (\ref{uyt}-\ref{uyt2}) to $1/(y-A)^2$ and $y^2$. The author did not even try.

\subsection{Exponential-like functions.}

Let
\begin{equation}\label{diffexp1}
\mathcal{D}f=a\mathcal{M}f,
\end{equation}

with given $f(x_0)$, $f(x_{n+1})= \dfrac{1+a(x_{n+1}-x_n)/2}{1-a(x_{n+1}-x_n)/2}f(x_n)$

on the arithmetic lattice, still an exponential function $\displaystyle
 f(x)=\left[\dfrac{1+ah/2}{1-ah/2}\right]^{(x-x_0)/h}f(x_0)$

on the $q-$lattice,  $x_{n+1}=qx_n$, $\displaystyle f(x)=
  \dfrac{e_q(ax/2)}{e_q(-ax/2)} f(x_0)$  (see \cite[\S 1.14]{Koe} for $e_q(x)$).

We try an expansion
\begin{equation}\label{expanexpo}f(x)=\sum_0^\infty c_n \dfrac{(x-x_0)\dotsb(x-x_{n-1}) }
  {(x-x'_0)\dotsb (x-x'_{n-1})},
\end{equation}
with given $c_0\neq 0$. Then, $c_1, c_2,\dotsc$ are found from interpolation
at $x=x_0, x_1,\dotsc$

 $c_1=\dfrac{(f(x_1)-f(x_0))(x_1-x'_0) }{x_1-x_0}
  = \dfrac{ a(x_1-x'_0) }{1-a(x_1-x_0)/2} c_0$,

$c2= \dfrac{ \overbrace{\ \ \ f(x_2)\ \ \  }^{\frac{1+a (x_2-x_1)/2}{1-a(x_2-x_1)/2}\left(c_0+c_1\frac{x_1-x_0}{x_1-x'_0}\right)}   -c_0-c_1(x_2-x_0)/(x_2-x'_0) }{(x_2-x_0)(x_2-x_1)/((x_2-x'_0)(x_2-x'_1))}
=\dfrac{a(x_2-x'_1)( 1+a(x_1-x'_0)/2) }{
(1-a\Delta x_0/2)(1-a\Delta x_1/2) }c_0
$

\textbf{Proposition.} If the $x'$s are such that $x'_0-x'_{-1}=2/a$,
 $c_n=\dfrac{x_n-x'_{n-1}}{x_0-x'_{-1}}\eta_0 \dotsb \eta_{n-1} c_0,n=0,1,\dotsc$, where

$\eta_n =\dfrac{  (x_n-x'_{-1}) (1+a (x'_{n}-x'_{n-1})/2)  }
{    (x'_n-x'_{-1})   (1-a(x_{n+1}-x_{n})/2)  }
$

Check $\eta_0= \dfrac{ (x_0-x'_{-1}) (1+a (x'_{0}-x'_{-1})/2=2)  }
{    (x'_0-x'_{-1}=2/a)   (1-a(x_{1}-x_{0})/2)  }
$

For a simpler way to get further coefficients, we put \eqref{expanexpo} in \eqref{diffexp1} to have


$-ac0+\sum_1^\infty c_n  C_{n,0,0} ( Y_2(y)-a D_{n,0,0}(y)/C_{n,0,0})
   \dfrac{(y-y_0)\dotsb(y-y_{n-2}) }{(y-y'_{-1})\dotsb (y-y'_{n-1})}\equiv 0$,
from \eqref{dY}.

$R_n(y)= Y_2(y)-a D_{n,0,0}(y)/C_{n,0,0} $ is a second degree polynomial.

 Let also
$P_n(y)=\dfrac{(y-y_0)\dotsb(y-y_{n-1}) }{(y-y'_{0})\dotsb (y-y'_{n-1})}$.
One has

$-ac_0+\displaystyle\sum_1^\infty c_n C_n R_n(y)\dfrac{P_{n-1}(y)}{(y-y'_{-1})(y-y'_{n-1})}=0
$.

Particular values of $R_n$ are
$R_n(y_{-1})= Y_2(y_{-1}) (1+a(x_0-x_{-1})/2)$,

$R_n(y_{n-1})= Y_2(y_{n-1}) (1-a(x_n-x_{n-1})/2$,

$R_n(y'_{-1})= Y_2(y'_{-1}) (1-a(x'_0-x'_{-1})/2$,

$R_n(y'_{n-1})= Y_2(y'_{n-1}) (1+a(x'_{n}-x'_{n-1})/2$,

from  (\ref{Dxm1}-\ref{Dxpn}) at $r=s=0$.

A big simplification occurs if $R_n(y'_{-1})=0$, or $x'_0-x'_{-1}=2/a$ which  is an equation for $y'_{-1}$, as $Y_2(y'_{-1})(x'_0-x'_{-1})$ is a square
root of $Y_1^2(y'_{-1})-4Y_0(y'_{-1})Y_2(y'_{-1})$ from \eqref{sumprod}

So, let $R_n(y)=(y-y'_{-1})S_n(y)$ with $S_n$ of first degree interpolated
at $y_{n-1}$ and $y'_{n-1}$ as

$S_n(y)= ( S_n(y'_{n-1})(y-y_{n-1})-(S_n(y_{n-1})(y-y'_{n-1}))/(y'_{n-1}-y_{n-1})$,

  we have
$S_n(y'_{n-1})=Y_2(y'_{n-1})(1+a (x'_{n}-x'_{n-1})/2)/(y'_{n-1}-y'_{-1})$,

$S_n(y_{n-1})=Y_2(y_{n-1}) (1-a(x_{n}-x_{n-1})/2)/(y_{n-1}-y'_{-1})$,

and
$-ac_0\displaystyle +\sum_1^\infty \dfrac{c_n C_n}{y'_{n-1}-y_{n-1}}
\dfrac{S_n(y'_{n-1})(y-y_{n-1})-(S_n(y_{n-1})(y-y'_{n-1})}{y-y'_{n-1}}P_{n-1}(y)\equiv 0
$, or

$-ac_0\displaystyle +\sum_1^\infty \dfrac{c_n C_n}{y'_{n-1}-y_{n-1}}
\{S_n(y'_{n-1})P_n(y)-S_n(y_{n-1})P_{n-1}(y)\}\equiv 0
$,

finally, by interpolating at $y=y_0, y_1,\dotsc$, the coefficients of $P_0, P_1;\dotsc$
must vanish:

$
-ac_0 -c_1C_1 S_1(y_0)/(y'_0-y_0)=0$,  \\
$\dfrac{c_n C_n}{y'_{n-1}-y_{n-1}}
S_n(y'_{n-1})-\dfrac{c_{n+1} C_{n+1}}{y'_{n}-y_{n}}
S_{n+1}(y_{n})=0, n=1,2,\dotsc $

See also \cite{kenfack,koorn2022,newtbasis,verdestar} for interesting constructions and  discussions
on factorial expansions.

$\dfrac{c_{n+1}}{c_n}=\dfrac{C_n S_n(y'_{n-1})(y'_n-y_n)}{C_{n+1}S_{n+1}(y_n)(y'_{n-1}-y_{n-1})}
=\dfrac{C_n (y_{n}-y'_{-1}) Y_2(y'_{n-1})(1+a (x'_{n}-x'_{n-1})/2)(y'_n-y_n)}
{C_{n+1}(y'_{n-1}-y'_{-1}) Y_2(y_{n}) (1-a(x_{n+1}-x_{n})(y'_{n-1}-y_{n-1})}
$

$\dfrac{C_n}{C_{n+1}}= \dfrac{(y_{-1}-y_{n-1})(x_{-1}-x'_{n})}{(y_{-1}-y'_{n})(x_{-1}-x_{n})  }
=\dfrac{(x'_0-x'_{n})(y'_{-1}-y_{n-1})}{(x'_0-x_{n})(y'_{-1}-y'_{n})  }$
from \eqref{slopes}

We recall \eqref{factorF} $F(x_m,y_n) = X_2(x_m)(y_n-y_{m-1})(y_n-y_m)=Y_2(y_n)(x_m-x_{n})(x_m-x_{n+1})$,

and we use the second formula together with the little trick of
multiplying the numerator and the denominator by $y'_{n-1}-y_n$,
so that

$(y_{n-1}-y'_{-1})(y_{n}-y'_{-1})=Y_2(y'_{-1})(x_n-x'_{-1})(x_n-x'_0)/X_2(x_n), \\
(y'_n-y_n)(y'_{n-1}-y_n)=(x'_n-x_n)(x'_n-x_{n+1})Y_2(y_n)/X_2(x'_n)$;

$y'_{n-1}-y'_{-1})(y'_{n}-y'_{-1})=Y_2(y'_{-1})(x'_n-x'_{-1})(x'_n-x'_0)/X_2(x'_n),\\
(y'_{n-1}-y_{n-1})(y'_{n-1}-y_n)=(x_n-x'_{n-1})(x_n-x'_n)Y_2(y'_{n-1})/X_2(x_n)$,
from \eqref{factorF}, to have

 \begin{equation}\label{coeffexp}
\dfrac{c_{n+1}}{c_n}=\dfrac{  (x_n-x'_{-1}) (1+a (x'_{n}-x'_{n-1})/2)(x_{n+1}-x'_n)  }
{    (x'_n-x'_{-1})   (1-a(x_{n+1}-x_{n})/2) (x_n-x'_{n-1}) }
\end{equation}

\begin{table}[htbp]
\begin{center}\begin{tabular}{|c|ccc|cc|}\hline
 $n$ &       $ x_n$     &      $  y_n$  &        $ f(x_n)$  &        $   f(1)$ &  $  f(-1.75)$  \\   \hline
 0   &     0        &   -0.27156  &       1      &       1      &1\\
1 &-0.69548& -0.62200 &0.90110 & 1.2627& 0.80640\\
2 &-0.99680 &-0.65853 &0.86138 & 0.91576& 0.77597\\
3 &-0.80672 &-0.37599 &0.88623 & 1.8574 &0.77109\\
4 &-0.17524 &0.12738 &0.97410 & 0.45911 &0.77014\\
5 &0.53805 &0.53701 &1.0839 & 1.1852 &0.76970\\
10 &-0.62568& -0.59521& 0.91057 & 1.1615& 0.77709\\
20 &-0.96796 &-0.67922& 0.86511 & 1.1615 &0.77621\\
\hline
\end{tabular}
\caption{Exponential-like function $ \mathcal{D}f=a \mathcal{Mf}$ with $a=2/(x'_0-x'_{-1})=0.14959...$\label{exptable}} \end{center}\end{table}

\section{Biorthogonal rational functions satisfying differential and difference equations.\label{biorths}}
\chead{\thesection  \ \ \ \   Biortho diff}

\subsection{Orthogonal polynomials: Laguerre and Hahn theories\label{LH}} \ \\

Laguerre wrote a number of studies on continued fraction expansion of solutions of
differential equations. His final survey \cite{Lag1885} written shortly before his death in 1886,
deals with the Jacobi continued fraction of a (possibly formal) power  expansion about $\infty$:
$f(x)=\cfrac{\mu_0}{x-r_0+\cfrac{s_1}{x-r_1+\cfrac{s_2}{x-r_2+\ddots}}} \\  =\dfrac{\mu_0}{x}+\dfrac{\mu_1}{x^2}+\dotsb$
of approximants $ \dfrac{N_0}{P_0}=\dfrac{0}{1},  \dfrac{N_1}{P_1}=\dfrac{\mu_0}{x-r_0},  \dfrac{N_2}{P_2}=\dfrac{\mu_0(x-r_1)}{(x-r_0)(x-r_1)+s_1}, \dotsc$,

where the power expansion of the $n^{\text{th}}$ approximant of the continued fraction has the
first $2n$ terms matching with the expansion of $f$ (\emph{Pad\'e} approximation). The
denominator of the $n^{\text{th}}$ approximant is also the related \emph{orthogonal polynomial}
of degree $n$, see  Brezinski's History \cite[\S 5.2]{BrezinskiH}, Chihara \cite[chap. III, \S 4]{chi}, Perron \cite[\S 67 (\S 30 in the 1957 edition)]{Perron}, Wall \cite[chap. XI]{Wall}.

 Recurrence relation is $P_{m+1}(x)=(x-r_{m})P_{m}(x)+s_mP_{m-1}(x)$ (recurrence relation of \emph{orthogonal polynomials}, often with $-s_m$ instead of $s_m$).

We consider weight functions satisfying a differential equation $w'/w=$ a rational function
on one or several intervals (`` Pearson-like'' differential equation).

Romanovsky and Hildebrandt related Pearson's equations for special probability weights
to orthogonal polynomials\footnote{The reference \cite{Hild1931} was given by A. Ronveaux.} \cite{Hild1931}

\subsubsection{\textbf{Proposition.}\label{Pear1}} \textsl{If the weight function $w$
satisfies $w'/w=V/U$ on a support $S$ made of one or several intervals, or also complex
contour(s), with polynomials $U,V$ such that $Uw=0$ at endpoints (if any), then the Stieltjes-Markov function $f(x)=\int_S w(t)dt/(x-t)$ of the measure
$d\mu(t)=w(t)dt$ satisfies}

\begin{equation}\label{UVW} Uf'=Vf+W,\end{equation} \textsl{where $W$ is a polynomial too}.

Proof: one has $ U(x)f(x)=\int_S U(x)w(t)dt/(x-t)=R_1(x)+\int_S U(t)w(t)dt/(x-t)$, where
$R_1(x)=\int_S( U(x)-U(t))w(t)dt/(x-t)$ is a polynomial. We derivate in $x$ and proceed
with an integration by parts in $t$ to have $(U(x)f(x))'=R'_1(x)+\int_S   \underbrace{(U(t)w(t))'}_{\displaystyle (U'+V)w \,dt/(x-t)}dt/(x-t)=
R'_1(x)+R_2(x)+(U'(x)+V(x))f(x)$, and \eqref{UVW} follows. \qed

The orthogonal polynomials are
related to $w$ on $S$, see also "semi-classical orthogonal polynomials" by Hendriksen and van Rossum \cite{HvR2}.

Laguerre explores functions satisfying the linear differential equation of first order
\eqref{UVW} with polynomials $U,V,W$, and has remarkable results, including linear
differential equations of second order satisfied by the orthogonal polynomials.
For instance $f(x)=\int_0^\infty (x-t)^{-1}e^{-t}dt$ satisfies \\  $xf'(x)+xf(x)=1$, with
$c_j=(j-1)!$, and the orthogonal polynomials are the ... Laguerre polynomials, but much
more difficult cases are treated in \cite{Lag1885}, such as an integral involving
$\exp(1/x^2)$ on a contour \cite[\S 25-28]{Lag1885}, see Belmehdi and Ronveaux \cite{Belm92,BelRonv} for
other examples in modern notation. These difficult examples are NOT classical, they (probably) do not
lead to hypergeometric expansions.

Perron \cite[\S 76 (\S 43 in the 1957 edition)]{Perron} shows the Laguerre theory.

If $f$ happens to have the form $f(x)=\int_S d\mu(t)/(x-t),  x\not\in S$  with a measure $d\mu$ on
a real or complex support $S$ (\emph{Stieltjes-Markov function}), then $P_m$ is orthogonal
with respect to $d\mu$ to any polynomial $p$ of degree $<m: \int_S p(t)P_m(t)d\mu(t)=0$.

$f_m(x)=\cfrac{s_m}{x-r_m+\cfrac{s_{m+1}}{x-r_{m+1}+\ddots}}$ satisfies a $1^{\text{st}}$ order differential equation, linear or Riccati, has a Taylor (or formal) expansion about $x=\infty$.

The Laguerre theory is equivalent to  finding
the continued fraction coefficients $r_m, s_m$s through  the recurrence relations
for polynomials $A_m,\dotsc, D_m$ (the $ABCD$ method \cite{Mag84,Segovia,MagLuminy} ) in

$A_mf'_{m}=B_mf_m^2+C_mf_m+D_m$, with $f_0=f, A_0=U,B_0=0,C_0=V,D_0=W$, \\   and
$f_m(x)=s_m/(x-r_m+f_{m+1}(x)) \Rightarrow$\\
$A_m(x)f'_{m+1}(x)=-A_m(x)-s_mB_m(x)-(x-r_m+f_{m+1}(x))C_m(x)-(x-r_m+f_{m+1}(x))^2D_m(x)/s_m$, so,
\begin{multline}\label{recAD0} s_m=\lim_{x\to\infty}\dfrac{-x^2D_m(x)}{A_m(x)+xC_m(x)},  r_m=\lim_{x\to\infty}\dfrac{A_m(x)+xC_m(x)+x^2D_m(x)/s_m}{C_m(x)+2xD_m(x)/s_m},  \\   A_m=A_0,     B_{m+1}(x)=-D_m(x)/s_m,  C_{m+1}(x)=-C_m(x)-2(x-r_m)D_m(x)/s_m, \\    D_{m+1}=-A_m(x)-s_mB_m(x)-(x-r_m)C_m(x)-(x-r_m)^2D_m/s_m. \end{multline}

The conditions on $s_m$ and $r_m$ follow from degree $D_{m+1}<$ degree of $C_m$.

See also that $C_{m+1}(x)-C_{m}(x)=2(A_m(x)+D_{m+1}(x)+s_mB_m(x))/(x-r_m),\\     C_m^2=C_0^2-4B_mD_m-4A(B_1+\dotsb +B_m)$.

\medskip
Differential relations and equation for denominator $P_m$: \\  $U(N_m/P_m)'-VN_m/P_m-W=
U(N_m/P_m-f)'-V(N_m/P_m-f)=O(x^{\text{\ degree\ }U,V-2m-1})$, so,   \\
$VN_mP_m-U(N'_mP_m-N_mP'_m)+WP_m^2$ are polynomials of bounded degree, $\Theta_m$ (almost the  Laguerre's notation) . Using
$N_mP_{m-1}-N_{m-1}P_m= \sigma_m$,  where  $\sigma_m=(-1)^{m-1}\mu_0 s_1\dotsb s_{m-1}$,  then

\begin{multline*}\Theta_m/P_m^2-\Theta_{m-1}/P_{m-1}^2=V(N_m/P_m-N_{m-1}/P_{m-1})-U(N_m/P_m-N_{m-1}/P_{m-1})'   \\   =\sigma_m[V/(P_mP_{m-1})-U(1/(P_mP_{m-1}))']\end{multline*}, or


$P_{m-1}\{ \Theta_mP_{m-1} -\sigma_m (UP'_m +VP_m/2)\}=P_m\{\Theta_{m-1}P_m+\sigma_m( UP'_{m-1}+VP_{m-1}/2\}$, whence

  $ \Theta_mP_{m-1} -\sigma_m (UP'_m +VP_m/2)=P_m$ times a polynomial called $\Omega_m$, so,
$ \Theta_{m-1}P_m+\sigma_m UP'_{m-1}=(\Omega_m-\sigma_m V/2) P_{m-1}$ ,
 or also the system of first order differential equations

 \begin{equation}\label{diffeqp1}
 \sigma_m U\dfrac{d}{dx}\begin{bmatrix}P_{m-1} \\ P_m\end{bmatrix}=
   \begin{bmatrix}  \Omega_m-\sigma_m V/2&- \Theta_{m-1}  \\  \Theta_m &  -\Omega_m-\sigma_m V/2 \end{bmatrix}  \begin{bmatrix}P_{m-1} \\ P_m\end{bmatrix}\end{equation}
   which can be turned as scalar second order differential equations if needed \cite{Branq}.

From the recurrence relations
of the
products $N_mP_m$ and derivatives:

$N_{m+1}P_{m+1}=(x-r_m)^2N_{m}P_{m}+s_m(x-r_m)(N_{m}P_{m-1}+N_{m-1}P_{m})+s_m^2N_{m-1}P_{m-1}$,

$N_{m+1}P_{m}+N_{m}P_{m+1}=2(x-r_m)N_{m}P_{m}+s_m(N_{m}P_{m-1}+N_{m-1}P_{m})$,

$N'_{m+1}P_m-N_mP'_{m+1}=(x-r_m)N'_mP_m +s_mN'_{m-1}P_m
     -(x-r_m)N_mP'_m -s_mN_mP'_{m-1}$, etc.,
one finds the same recurrence relations for $\Theta_m, \Omega_m$ and $D_m, C_m$, to have,after checking first items,

\begin{equation}\label{diffeqp2}
 \Theta_m=-\sigma_m D_m =-\sigma_{m+1} B_{m+1} ,  \Omega_m=-\sigma_m C_m/2, m=1,2,\dotsc \end{equation}

\medskip

How to find  $A_m$ etc.? By induction followed by mathematical induction.

\subsubsection{Examples with known orthogonal families.\label{Euler1}}

\begin{enumerate}
\item
Knuth \cite[\S 4.5.3,  answers of \S 4.5.3, 16]{Knuth}
is delighted to report how the  24-year old Euler found the continued fraction
expansion of the exponential function through
 \emph{Riccati} differential equations, see  Khovanskii \cite[chap.II, \S 9]{Khov}, Cretney \cite{Cretney} for the
role of Daniel Bernoulli.

{\small $f(x)=e^{1/x}-1=1/x+1/(2x^2)+\dotsb$ satisfies $x^2f'(x)=-f(x)-1, $ so, $-x^2(1+f'_1(x))/(x-1/2+f_1(x))^2=-1/(x-1/2+f_1(x))-1$, and the Riccati equation for $f_1$ is
$x^2f'_1(x)=f_1^2(x)+2xf_1(x)-1/4$,
$s_1=\lim (x^2/4)/(x^2+2x^2)=1/12, r_1=0,B_2=3, C_2(x)=4x,D_2=-1/12
$,  etc.

$ A_m(x)=x^2, B_m(x)=m, C_m(x)=2mx, D_m=-1/(8m-4),r_m=0,s_m=1/(16m^2-4),\ \ m=1,2,\dotsc $
}
 Cuyt \& {\it al.} \textsl{Handb. cont. fractions} \cite[ \S 11.1.2, 11.1.3.]{Cuyt}
  Bessel polynomials $P_{m}(x)=(-1)^m y_m(-2x;2)/(2^m (2m-1)!!)$ NIST \cite[18.34.]{NIST}.

\item {\small
    $f(x)=\int_0^\infty e^{-t}dt/(x-t)=\sum_0^\infty n!/x^{n+1}$ (formal) satisfies $xf'(x)=-xf(x)+1$, so,
    $A_0(x)=x, B_0=0, C_0(x)=-x, D_0=1, s_0=1, r_0=1,
    B_1=-1,C_1(x)=-x+2,D_1=-1, s_1=-1, r_1=3$, etc.

    $A_m(x)=x, B_m=-1, C_m(x)=2m-x, D_m=s_m=-m, r_m=2m+1, m=1,2\dotsc$
    }

   gives the  Laguerre polynomials $P_m(x)=(-1)^mL_m(x)/m!$.

\item {\small
 $f(x)=\displaystyle \int_{-\infty}^\infty \dfrac{ \exp(-t^2)dt}{x-t}\sim \mu_0/ x+\mu_2 /x^3+\dotsb$ (the series is divergent asymptotic, $\mu_p=\Gamma((p+1)/2), p=0,2,4,\dotsc$), $s_0=\mu_0$, $f$
satisfies $f'(x)=-\sum(p+1)\mu_p/ x^{p+2}-\dotsb =-2xf(x)+2\mu_0$, or
$  A_0(x)=1, B_0=0, C_0(x)=-2x, D_0=2\mu_0$. Then,  \eqref{recAD0} gives\\
$s_0=\mu_0, A_1(x)=1, B_1=- 2, C_1(x)=-2x, D_1=-1,  \\
  s_1=-1/2, B_2(x)=- 2, C_2(x)=-2x, D_2=-2,s_2=-1,   \text{\ \ etc.}\\
 A_m(x)=1, B_m(x)=-2, C_m(x)=-2x, D_m=-m,s_m=-m/2, r_m=0, m=1,2,\dotsc$.
}
We then recover the Hermite polynomials,  $P_m(x)= H_m(x)/2^m$.

\smallskip

\item {\small  $f(x)=\dfrac{1}{2}\log\dfrac{x+1}{x-1}=\dfrac{1}{x}+\dfrac{1}{3x^3}+\dotsb,
(x^2-1)f'(x)=-1,   A_m(x)=(x^2-1), B_m=m, C_m(x)=2mx, s_m=-1/(4m^2-1),      D_m=2m-1, m=0,1,\dotsc$  }
 $P_{m}(x)=((2m-1)!!/m!)$ times the Legendre polynomial of degree $m$.

\item  {\small  ${\displaystyle{  f(x)=\int_{-\infty}^{\infty}\dfrac{e^{-t^4}dt}{x-t}=\dfrac{\mu_0}{x}+\dfrac{\mu_2}{x^3}+\dotsb,
\mu_m=\Gamma((m+1)/4)/2, m=0,2,4,\dotsc}}$}  \\   $f'(x)=-4x^3f(x)+4\mu_0 x^2+4\mu_2 $,
{\small
$ B_1(x)=-4x^2-4\mu_2/\mu_0,  C_1(x)=-4x^3-8x\mu_2/\mu_0, D_1(x)=-1-4x^2\mu_2/\mu_0,s_1=-\mu_2/\mu_0, \\
B_2(x)=1/s_1+4x^2, C_2(x)=-4x^3+(8s_1-8+2/s_1) x,
D_2(x)=-1+  (1/s_1  -4s_1) x^2 , s_2=1/(4s_1)-s_1$, etc.

    $A_m(x)=1,  C_m(x)=-4x^3+8s_m x,      D_m(x)=-1+4s_m x^2,  m=1,2,\dotsc  B_m(x)=-4x^2+1/s_{m-1},  \\   4s_m(s_{m-1}+s_m+s_{m+1})=m, m=2,3,\dotsc  $,
 }
  no simple formula, Freud-Painlev\'e polynomials (Clarkson {\it et al.} \cite{ClarksonFreud}, Shohat\cite{Shohat}, Nevai 1984\cite{Nevai1984} , Van Assche\cite{VAssche2018}, Rebocho2025\cite{Rebo2025}, L.Haine 2026\cite{Haine2026}). This
  is a non classical example, with polynomials $A_m$ etc. of degree $>2$.

\item\label{genJac}

Generalized  Jacobi polynomials: let the weight function $w$  be two constants
times \\  $|(x-a)^\alpha (x-b)^\beta (x-c)^\gamma|$ on two intervals $[a,b]$ and $[b,c]$    \cite[Fig.1]{MNR}. We then have $U(x)=(x-a)(x-b)(x-c), U(x)f(x)=\text{\  pol.\  }+
\int_a^c U(t)w(t)dt/(x-t), (U(x)f(x))'= \text{\  pol'.\  }+
\int_a^c (U(t)w(t))'dt/(x-t)$, so, $V(x)/U(x)=\alpha/(x-a)+\beta/(x-b)+\gamma/(x-c)$,
and $B_m,C_m,D_m$ of degrees $1,2,1$ follow, but there are no elementary formulas.  Gammel and Nuttall  \cite{GammelNuttall}
find recurrences (probably) equivalent to \eqref{recAD0}, and conjecture asymptotic
behaviors given by elliptic functions. Laguerre ends his ultimate 1885 paper \cite{Lag1885}
with a similar conjecture, adding ''Je reviendrai du reste sur ce point particulier, en essayant
de compl\'eter  \`a certains \'egards les r\'esultats obtenus par ces illustres g\'eom\`etres.''
(I intend to return to this point, trying to complete in some way the results of these famous mathematicians [Jacobi and Borchardt]). This very last sentence is not found in the \OE uvres
book form copy \cite[vol. II]{Lag1885}.
Later on, the authors of  \cite{Nuttall25} studied the asymptotic properties of $P_m$ through complex  potential theory, without using recurrence relations. See \cite{Erdel} for
inequalities.

\end{enumerate}

The Laguerre theory will be extended here to any difference calculus based on a biquadratic
polynomial $F$, including the differential limit case,  various discretizations,
up to the elliptic lattice. The differential equation for $f$ will be expanded to
a discrete Riccati equation, a very natural step when dealing with
continued fractions \cite{Cretney,Khov}

W. Hahn \cite{Hahn1952,HahnE,Hahnmonh}  looks at the same situation from the more difficult inverse
point of view: what are the orthogonal polynomials satisfying linear differential equations of
fixed order and complexity of coefficients? He finds that the order can be as low as 4 or 2,
recovering the Laguerre theory, whence the name \emph{Laguerre-Hahn} sometimes given to this theory  \cite{BangRama,BelRonv,Branq,Mag84,Segovia,McCabe,Ronv4}.

See \cite{Jbeli,Rebo,Rebo2025} for  recent surveys.

\subsection{The Pearson's like difference  equation for discrete measures\label{secPear2}}

We consider a discrete measure made of weights $w_n$ at points $t_n$ of a particular
lattice related to  a biquadratic polynomial $F$,
as discussed earlier. The Pearson's like differential equation is extended to the following

\textbf{Proposition.} \textsl{Let weights  $w_n$ at points $t_n$ of a particular
lattice such that $F(t_n,u_n)=F(t_{n+1},u_n)=0$ with a biquadratic polynomial
$F(x,y)=X_0(x)+yX_1(x)+y^2X_2(x)$. If the weights satisfy}

\begin{equation}\label{Pear2}
\dfrac{w_{n+1}/(X_2(t_{n+1})(u_{n+1}-u_n))}{w_n/(X_2(t_n)(u_{n}-u_{n-1}))}=\dfrac{A(u_n)+(t_{n+1}-t_n)C(u_n)/2}{A(u_n)-(t_{n+1}-t_n)C(u_n)/2}
\text{\ or\ }
(A\mathcal{D}-C\mathcal{M})\dfrac{w}{X_2\Delta y}=0,
\end{equation}

\textsl{ and $A(u_{n_1-1})-(t_{n_1}-t_{n_1-1})C(u_{n_1-1})/2=A(u_{n_2})-(t_{n_2+1}-t_{n_2})C(u_{n_2})/2=0$,
 if $t_n$ is defined on a finite set $n\in\{n_1,n_2\}$, with polynomials $A$ and $C$, then $f(x)=\sum_n w_n/(x-t_n)$ satisfies the
linear difference equation}

\begin{equation}\label{Pear21} A(y)(\mathcal{D}f)(y)=C(y)(\mathcal{M}f)(y)+D(y),\end{equation}

\textsl{where $D$ is a polynomial too.}

Proof.

At some $y$, with $F(x_m,y)=F(x_{m+1},y)=0$,
$\displaystyle A(y)(\mathcal{D}f)(y)-C(y)(\mathcal{M}f)(y)  \\=
\sum_{n=n_1}^{n_2} w_n  \dfrac{-A(y)-C(y)\left((x_m+x_{m+1})/2-t_n\right)}{\underbrace{(x_m-t_n)(x_{m+1}-t_n)}_{\displaystyle F(t_n,y)/Y_2(y)=(y-u_n)(y-u_{n-1})X_2(t_n)/Y_2(y)}  }
$

is a rational function of $y$ (as $x_m+x_{m+1}=-Y_1(y)/Y_2(y)$ is, remind \eqref{sumprod}) , and we look at the  residues at $u_{n_1-1}, \dotsc,u_{n_2}$: the residue at $y=u_n$ is

{\small\noindent$  w_{n+1}\dfrac{-A(y_n)Y_2(y_n)+C(y)_n\left(t_{n+1}Y_2(y_n)+Y_1(y_n)/2\right)}{(u_n-u_{n+1})X_2(t_{n+1})}
+  w_n \dfrac{-A(y_n)Y_2(y_n)+C(y_n)\left(t_n Y_2(y_n)+Y_1(y_n)/2\right)}{(u_n-u_{n-1})X_2(t_{n})}=0  $      }
from \eqref{Pear2},  seeing that  $t_{n+1}Y_2(y_n)+Y_1(y_n)/2=(t_{n+1}-t_n)Y_2(y_n)/2,  t_{n}Y_2(y_n)+Y_1(y_n)/2   \\   =-(t_{n+1}-t_n)Y_2(y_n)/2$.

Note that $w_{n_1-1}=w_{n_2+1}=0$ follow from the conditions on $u_{n_1-1}, u_{n_2}$.

See also \cite{Manas1,Manas3}.

%



\subsection{ Orthogonal polynomials of a discrete variable\label{discorthog}}  \    \\

Here is a new way to introduce orthogonal polynomials of a discrete variable, the polynomials
are defined through the moments $\mu_0, \mu_1,\dotsc$ of the measure $d\mu$, which are
also the coefficients of the expansion about $\infty$ of the Stieltjes-Markov function
$f(x)=\int_S \mu(t)/(x-t)=\mu_0/x+\mu_1/x^2+\dotsc$ The discrete feature considered
here is that $f$ satisfies a linear difference equation of first order resulting, as seen in
\S \ref{secPear2} just above, from a discrete measure defined on $(t_n,u_n)$ points satisfying $F(t_n,u_n)=0$ with a biquadratic $F$.

  We are able to discuss the behavior
of $f(x)$ near $\infty$ only if the $x_n$ sequence has a fixed point at $\infty$, and
fixed points have been considered in \S \ref{fixedp}. From
\eqref{sumprodE}, $x_n\to\infty$ only if $e_{2,2}=0$, and as
 $P(x)$ is a constant times
$
(e_{1,2}x^2+e_{1,1}x+e_{0,1})^2-4(e_{0,2}x^2+e_{0,1}x+e_{0,0})(e_{2,2}x^2+e_{1,2}x+e_{0,2})$, its degree is 2 if $e_{1,2}=0$ as well.
The $y-$fixed point is sent to $\infty$ too by the unessential transformation $y\to 1/(y-y_\infty)$. One then gets the curve $F(x,y)=0$ as a conic, as in \cite{Segovia}.

\medskip

We will see with biorthogonal rational functions in the next subsection \S \ref{biorthinterp} that there will be no
fixed point condition, and that a general biquadratic $F$ is allowed.

We discuss now $f_m(x)=\cfrac{s_m}{x-r_m+\cfrac{s_{m+1}}{x-r_{m+1}+\cfrac{s_{m+2}}{x-r_{m+2}+\ddots}}}$ in

        \begin{equation}\label{Riccatim} A_m(y_n)(\mathcal{D}f_m)(y_n)=B_m(y_n)f_m(x_n)f_m(x_{n+1})
          +C_m(y_n)(\mathcal{M}f_m)(y_n)+D_m(y_n), \ \ m=0,1,\dotsc,
\end{equation}

\begin{multline}\label{recAD02} s_m=\lim_{y\to \infty}\dfrac{-Y_0(y)D_m(y)}{Y_2(y)A_m(y)-Y_1(y)C_m(y)/2},
 r_m=\lim_{y\to \infty}\dfrac{Y_2(y)A_m(y)-Y_1(y)C_m(y)/2+Y_0(y)D_m(y)/s_m}{Y_2(y)C_m(y)-Y_1(y)D_m(y)/s_m},  \\    A_{m+1}=A_m-QD_m/(2s_mY_2^2), B_{m+1}=-D_m/s_m, \\ C_{m+1}=-C_m+(Y_1+2r_mY_2)D_m/s_m,   \\    D_{m+1}=-A_m-s_mB_m+(Y_1/2+r_mY_2)C_m/Y_2-\left(\underbrace{Y_0+r_mY_1+r_m^2Y_2}_{ F(r_m,y)}\right)D_m/(s_mY_2). \end{multline}

  Proof:      we see in \eqref{Riccatim} how $f_m(x_{n+1})$ is related to $f_m(x_{n})$:
{\small
\begin{equation}\label{Riccatim2}
 f_m(x_{n+1})= \dfrac{    \left[\dfrac{A_m(y_n)}{x_{n+1}-x_n}+\dfrac{C_m(y_n)}{2}\right]f_m(x_n)  +D_m(y_n)  }
              {-B_m(y_n)f_m(x_n)+\dfrac{A_m(y_n)}{x_{n+1}-x_n}-\dfrac{C_m(y_n)}{2}    }
=-\dfrac{U_m(y_n)f_m(x_n)+D_m(y_n)}{  B_m(y_n)f_m(x_n)+V_m(y_n-},
\end{equation}
where
$U_m(y_n)=\dfrac{C_{m}(y_n)}{2}+ \dfrac{A_{m}(y_n)}{x_{n+1}-x_n}, V_m(y_n)=\dfrac{C_{m}(y_n)}{2}- \dfrac{A_{m}(y_n)}{x_{n+1}-x_n}$,
 and we
perform the step $m\to m+1$ by entering $f_m(x)=\dfrac{ s_m  }{x-r_m+f_{m+1}(x) } $

$\dfrac{ s_m  }{x_{n+1}-r_m+f_{m+1}(x_{n+1}) } =\dfrac{U_m(y_n)s_m
+(x_n-r_m+f_{m+1}(x_n))D_m(y_n)   }
{-(x_n-r_m+f_{m+1}(x_n))V_m(y_n)-s_mB_m(y_n) }
$, or

  $  f_{m+1}(x_{n+1}) =-\dfrac{ (x_n-r_m+f_{m+1}(x_n))V_m(y_n)+s_mB_(y_n) }{U_m(y_n)
+(x_n-r_m+f_{m+1}(x_n))D_m(y_n)/s_m   }-x_{n+1}+r_m$, so,

$U_{m+1}(y_n)=- V_m(y_n)+(x_{n+1}-r_m)D_m(y_n)/s_m$,

$V_{m+1}(y_n)=- U_m(y_n)+(x_{n}-r_m)D_m(y_n)/s_m$,

$B_{m+1}(y_n)=- D_m(y_n)/s_m$,

$D_{m+1}(y_n)=-(x_n-r_m)V_m(y_n)+s_mB_m(y_n)+(x_{n+1}-r_m)U_m(y_n)
+(x_n-r_m)(x_{n+1}-r_m)D_m(y_n)/s_m$,  and \eqref{recAD02} follows.
   }   
 \subsubsection{\label{discorthex}}       All the classical  orthogonal polynomials of a discrete variable (Askey \& Wilson \cite{AW319}, Chihara \cite{chi}, Koekoek et al. \cite{Koe98,Koe}, NIST \cite{NIST}, Nikiforov et al. \cite{NS,NSU})
are recovered by the present theory, here are two instances, a very easy, and a most elaborate:

\begin{enumerate}

\item\label{discCheb}
${\displaystyle{ f(x)=\sum_{j=0}^{N-1} \dfrac{1}{x+jh}=\dfrac{N}{x}-\dfrac{(N-1)Nh}{2x^2}+\dotsb \text{\ satisfies\  }   \\
f(x+h)-f(x)=\dfrac{1}{x+Nh}-\dfrac{1}{x}=-\dfrac{Nh}{x(x+Nh)}    }}$.

{\small
With $F(x,y)=(y-x)(y-x-h), Y_0(y)=y(y-h), Y_1(y)=-2y+h, Y_2(y)=1, Q=h^2$, so that $x_n=nh=y_{n-1}\Rightarrow  A_0(y)=(y-h)(y+(N-1)h), B_0=C_0=0,
D_0=-N, s_0=N, r_0=-(N-1)h/2;   A_1(y)=(y-h)(y+(N-1)h)+h^2/2, B_1=1, C_1(y)=2y+(N-2)h, D_1=(N^2-1)h^2/4,s_1=-(N^2-1)h^2/12, r_1=r_0, 
$ etc.

$A_m(y)=(y-h)(y+(N-1)h)+m^2h^2/2, B_m=2m-1, C_m(y)=m(2y+(N-2)h), D_m=m^2(N^2-m^2)h^2/(4(2m-1)), s_m=-D_m/(2m+1), r_m=r_0, \  \   m=1,2,\dotsc$
}

 Discrete Chebyshev p. Abramowitz \cite[\S 22.17]{AbrS},
special case of Hahn polynomials.

\item\label{AWabcd}

Askey-Wilson \cite[\S 14.1]{Koe98,Koe},   NIST \cite[\S 18.28]{NIST}:

the orthogonality measure may have a discrete part with weights $w_n$ on
points \\   $x_n= (aq^n+1/(aq^n))/2,
y_n=x_{n+1/2}:
F(x,y)=(xy-(\sqrt{q}+1/\sqrt{q})/2)^2
-(x^2-1)(y^2-1): Y_1(y)=-(\sqrt{q}+1/\sqrt{q})y, Y_0(y)= y^2+(q-1)^2/(4q), X_2(x)=Y_2(y)\equiv 1, y_{n+1}-y_n=(aq^{n+1}-1/(aq^{n+1}))(\sqrt{q}-1/\sqrt{q})/2$,

and $\dfrac {w_{n+1}}{w_n}=
\dfrac{ (1-a^2q^{2n+2})(1-a^2q^n)
(1-abq^n)(1-acq^n)(1-adq^n)q}
{ (1-a^2q^{2n})(1-q^{n+1})
(1-aq^{n+1}/b)(1-aq^{n+1}/c)(1-aq^{n+1}/d)abcd}
$


With $x_n=(aq^n+1/(aq^n))/2, y_n=x_{n+1/2}$, we see that $w_{n+1}/w_n$ is a rational function
of $q^n$ which is a quadratic algebraic function
of $x_n$ or $y_n$.

From Pearson \S \ref{secPear2},
  $\dfrac{w_{n+1}/(X_2(x_{n+1})(y_{n+1}-y_n))}{w_n/(X_2(x_n)(y_{n}-y_{n-1}))}$ must be the
  ratio of conjugate algebraic functions $ \dfrac{A_0(y_n)+(x_{n+1}-x_n)C_0(y_n)/2}{A_0(y_n)-(x_{n+1}-x_n)C_0(y_n)/2}$, as
$x_{n+1}-x_n=(\sqrt{q}-1/\sqrt{q})
(aq^{n+1/2}-1/(aq^{n+1/2}))/2=
(\sqrt{q}-1/\sqrt{q})\sqrt{y_n^2-1}$.

  $\dfrac{w_{n+1}/(X_2(x_{n+1})(y_{n+1}-y_n))}{w_n/(X_2(x_n)(y_{n}-y_{n-1}))}
  =\dfrac{w_{n+1}/(aq^{n+1}-1/(aq^{n+1}))}{w_n/(aq^{n}-1/(aq^{n}))}
=q\dfrac{w_{n+1}/(1-a^2q^{2n+2})}{w_n/(1-a^2q^{2n})}
\\    =q^2{\displaystyle \prod_{\alpha=a,b,c,d} }
\dfrac{1-\alpha a q^n}{\alpha-a q^{n+1}  }
=\dfrac{q^2-a\mathfrak{a}q^{n+2}+a^2\mathfrak{b}q^{2n+2}-a^3\mathfrak{c}q^{3n+2}+a^5bcdq^{4n+2}   }{ abcd-a\mathfrak{c}q^{n+1}+a^2\mathfrak{b}q^{2n+2}-a^3\mathfrak{a}q^{3n+3}+a^4q^{4n+4}    }
$,

where  $\mathfrak{a}=a+b+c+d, \mathfrak{b}=ab+ac+ad+bc+bd+cd, \mathfrak{c}=abc+abd+acd+bcd$.


We proceed with sum and difference of
terms of the fraction, divided by $q^{2 n+2}$

$2A_0(y_n)=(a^4q^{2n}+1/q^{2n+2})(q^2+abcd) -a\mathfrak{a}
(1+a^2q^{2n+1})q^{-n}
+2a^2\mathfrak{b}
-a\mathfrak{c}(1/q^{n+1}+a^2q^{n})  \\
=4a^2(y_n^2-1/2)(q^2+abcd)/q -2a^2\mathfrak{a}
q^{1/2}y_n
+2a^2\mathfrak{b}
-2a^2\mathfrak{c}q^{-1/2}y_n
$,
\smallskip

$(x_{n+1}-x_n)C_0(y_n)=(-a^4q^{2n}+1/q^{2n+2})(q^2-abcd) -a\mathfrak{a}
(1-a^2q^{2n+1})q^{-n}
+a\mathfrak{b}(1-a^2q^{2n+1})q^{-n-1}  \\
=(a^2q^{n}-1/q^{n+1})a[
-(aq^{n}+1/(aq^{n+1}))(q^2-abcd) +\mathfrak{a}
q
-\mathfrak{c} ]$, so, dividing by
$x_{n+1}-x_n=(q-1)
(a^2q^{n}-1/q^{n+1})/(2a)$,

$C_0(y_n)=\dfrac{2a^2}{q-1}\left(\vphantom{\dfrac{2a^2}{q-1}}  -2  y_n
(q^2-abcd)/\sqrt{q} +\mathfrak{a}q
-\mathfrak{c} \right)
$


$s_0=-D_0/(A_0''/2+(\sqrt{q}+1/\sqrt{q})C_0'/2=D_0(q-1)/(  4a^2(q-abcd) ),  \\
  r_0=(A_0'(0)+(\sqrt{q}+1/\sqrt{q})C_0(0)/2)/(C_0'+(\sqrt{q}+1/\sqrt{q})D_0/s_0)
  =(\mathfrak{a}-\mathfrak{c})/(2(1-abcd)),\\
  A_1(y)=A_0(y)-QD_0/(2s_0)=A_0(y)   +(1/q-1)(y^2-1) 2a^2q^2(q-abcd),
B_1=4a^2(q-abcd)/(1-q),
C_1(y)=4a^2(1-abcd)q^{1/2}y/(1-q)-2a^2( (q-(2-q)abcd)\mathfrak{a}+(1-2q+abcd)\mathfrak{c})/((1-q)(1-abcd)),
D_1(y)=a^2(1-ab)\dotsb/(1-abcd)^2,
s_1=-(1-q)(1-ab)\dotsb/(4(1-abcd)^2(1-abcdq))$,  etc.

\noindent $A_m(y)=A_0(y)+2a^2q(q^{-m}-1)(1-abcdq^{m-2})(y^2-1), B_m=4a^2(q^{2-m}-abcd q^{m-1})/(1-q),
  \\   C_m(y)=\dfrac{4a^2(1-abcdq^{2m-2})y}{q^{m-3/2}(1-q)}-2a^2\dfrac{ (q-(2q^{m-1}-q^{2m-1})abcd)\mathfrak{a}+(1-2q^m+q^{2m-2}abcd)\mathfrak{c}}{(1-q)(1-q^{2m-2}abcd)},
  \\
   D_m=\dfrac{a^2(1-q^m)(1-abcdq^{m-2})(1-abq^{m-1})\dotsb  }{q^{m-1}(1-q)(1-abcdq^{2m-3})(1-abcd^{2m-2})^2}, s_m=-\dfrac{D_mq^{m-1}}{4a^2(1-abcdq^{2m-1})},
   \\    r_m=q^{m-1}\dfrac{ q(1-q^{m-2}(1+q-q^{m+1})abcd)\mathfrak{a}+(1-q^m-q^{m+1}+abcdq^{2m-1})\mathfrak{c}   }{2(1-abcdq^{2m-2})(1-abcdq^{2m})}
$.


\item A non classical example.

Let the orthogonality measure be a unit weight on $N_1, N_1+1,\dotsc, N_2-1$ and $N_3, N_3+1,\dotsc, N_4-1$, so that

$f(x)=\sum_{N_1}^{N_2-1}1/(x-j)+\sum_{N_3}^{N_4-1}1/(x-j)=(N_2-N_1+N_4-N_3)/x +\dotsb
$
With $F(x,y)=(y-x)(y-x-1), x_n=n, y_n=n+1$, one has
$(\mathcal{D})f(y)= \sum 1/(x+1-j)-1/(x-j)=1/(y-N_1)-1/(y-N_2)+1/(y-N_3)-1/(y-N_4)=D_0(y)/A_0(y)$. Later on, $A_m$ etc. have degrees $4,2,3,2$, but, as in \S\ref{Euler1}-(\ref{genJac}), there are no  elementary formulas.
   \end{enumerate}

\subsection{Biorthogonal rational functions,  rational interpolation with free poles\label{biorthinterp}}  \   \\

\subsubsection{\label{biorth}Biorthogonality.} \   \\

Biorthogonality refers to two sequences of objects in vector spaces, possibly of different kind,
but able to interact through a bilinear form $\langle\ ,\ \rangle$:
$\{\mathcal{A}_n\}$ and $\{\mathcal{B}_m\}$ are biorthogonal if $\langle\mathcal{A}_n,\mathcal{B}_m\rangle=0, m\neq n$.

Let
 $\{\mathfrak{a}_n, \mathfrak{b}_m\}$ be sequences of independent elements
in the two vector spaces, and we try triangular combinations $\mathcal{A}_n=\sum_0^n \alpha_j^{(n)} \mathfrak{a}_j, \mathcal{B}_m=\sum_0^m \beta_j^{(m)} \mathfrak{b}_j$.
The required combinations are found from
the triangular factorization of the matrix $[<\mathfrak{a}_r,\mathfrak{ b}_s>]$ (Brezinski \cite{Brezbio}, Davis \cite{Davis} ):

\[
\begin{bmatrix} \alpha_0^{(0)}  &   0         &           &           &          \\
                \alpha_0^{(1)} & \alpha_1^{(1)}  &    0     &            &         \\
                \vdots      &  \vdots      &  \ddots  &     0       &       \\
                \alpha_0^{(n)} & \alpha_1^{(n)}  &  \dotsb  & \alpha_n^{(n)} &  0 \\
                  \vdots    &   \vdots    &           &            & \ddots
\end{bmatrix}
\begin{bmatrix}  \langle\mathfrak{a}_0,\mathfrak{b}_0\rangle  &  \langle\mathfrak{a}_0,\mathfrak{b}_1\rangle  &   \dotsb&                  \\
                 \langle\mathfrak{a}_1,\mathfrak{b}_0\rangle  &  \langle\mathfrak{a}_1,\mathfrak{b}_1\rangle  &   \dotsb&                   \\
                  \vdots     &   \vdots       &                        &     \\
                 \langle\mathfrak{a}_n,\mathfrak{b}_0\rangle  &  \langle\mathfrak{a}_n,\mathfrak{b}_1\rangle  &   \dotsb&                  \\
                     \vdots   &   \vdots      &                     \\
                     \vdots   &   \vdots      &
\end{bmatrix}
\begin{bmatrix} \beta_0^{(0)}  & \beta_0^{(1)}   &  \dotsb    & \beta_0^{(m)}  &\dotsb    \\
                  0          & \beta_1^{(1)}    &  \dotsb    & \beta_1^{(m)}  &\dotsb    \\
                             &     0          &  \ddots    &   \vdots        &       \\
                              &               &    0       & \beta_m^{(m)}    &      \\
                             &                 &           &       0         & \ddots
\end{bmatrix}
\]
is the diagonal matrix of the $\langle\mathcal{A}_n,\mathcal{B}_n\rangle$s. Existence an unicity hold if
the determinants of the  $[\langle\mathfrak{a}_r,\mathfrak{ b}_s\rangle]$s, for $r,s=0,\dotsc n$, for all $n$, do not vanish  \cite{Brezbio,Davis}.

The triangular factors above are the inverses of the corresponding factors of
the Gauss triangular factorization.

\subsubsection{\label{biorthr}Biorthogonal rational functions.} \   \\

We consider a much narrower set made of rational functions,
with $\mathfrak{a}_n=1/(x-a_{n}), \mathfrak{b}_n=1/(x-b_{n})$.

 Here, $\mathcal{A}_n$ and $\mathcal{B}_m$
are rational functions
\[\dfrac{\mathsf{A}_n(x)}{(x-a_0)\dotsb (x-a_n)}\ \  \text{and}\ \   \dfrac{\mathsf{B}_m(x)}{(x-b_0)\dotsb(x-b_m)},\ \  n,m =0,1,\dotsc \]
with $\mathsf{A}_n$ and $\mathsf{B}_m$ polynomials of degrees $n$ and $m$,
and the bilinear form is actually a linear form $\mathscr{L}$ acting on the product, often
given as an integral$\int_Sf(t)g(t)d\mu(t)$, with some
measure $d\mu$ (possibly formal),  $\int_S  \mathcal{A}_n(t)\mathcal{B}_m(t)d\mu(t)=0, m\neq n$. ¨


From the bilinearity conditions $\mathcal{A}_n$ is  orthogonal to any linear combination of the $\mathcal{B}_m$s, $m=0,\dotsc, m-1$, equivalent statements are

    $\mathsf{A}_n(t)$ is orthogonal to all polynomials of degree $<n$ with respect to\\ $d\mu(t)/((t-a_0)\dotsb(t-a_n)(t-b_0)\dotsb (t-b_{n-1}))$.

For $\mathsf{B}_m$:

 $\mathsf{B}_m(t)$ is orthogonal to all polynomials of degree $<m$ with respect to \\ $d\mu(t)/((t-a_0)\dotsb(t-a_{m-1})(t-b_0)\dotsb (t-b_{m}))$.

This cries for $\mathsf{C}_n$ of degree $n$  orthogonal to all polynomials of smaller
degrees with respect to

$d\mu(t)/((t-a_0)\dotsb(t-a_n)(t-b_0)\dotsb (t-b_{n}))$.

Are there biorthogonal functions involving  $\mathsf{C}_n$ ?

\begin{equation}\label{recAm0}\begin{split}
\mathsf{A}_n(t)&= \text{\  a\  constant\  times\  } \dfrac{\mathsf{C}_n(b_n)\mathsf{A}_{n+1}(t)-\mathsf{A}_{n+1}(b_n)\mathsf{C}_n(t)}{t-b_n},\\
\mathsf{B}_m(t)&= \text{\  a\  constant\  times\  } \dfrac{\mathsf{C}_m(a_n)\mathsf{B}_{m+1}(t)-\mathsf{B}_{m+1}(a_n)\mathsf{C}_m(t)}{t-a_m}.
\end{split}\end{equation}

Indeed, we must check that the right-hand side of the first equation is orthogonal to all polynomials $p$ of degree $<n$ with respect to $d\mu(t)/((t-a_0\dotsb(t-a_n)(t-b_0)\dotsb (t-b_{n-1}))$. This works nicely with the $\mathsf{C}_n(t)$ part, thanks to the division by $t-b_n$.
As for $\mathsf{A}_{n+1}(t)$, orthogonal to all polynomials  of degree $\leqslant n$, including the $(t-a_{n+1})p(t)$ polynomials, restoring the expected orthogonality property.

A similar proof holds for $\mathsf{B}_m(t)$.  \qed

These identities give $\dfrac{ \mathsf{C}_n(t)}{(t-a_0)\dotsb(t-a_{n+1})}$ and $\dfrac{ \mathsf{C}_n(t)}{(t-b_0)\dotsb(t-b_{n+1})}$ in terms of   $\mathcal{A}_n, \mathcal{A}_{n+1}, \mathcal{B}_{n}, \mathcal{B}_{n+1}$, but this is perhaps not useful.

A  similar  identity for $\mathsf{C}_n$ will be needed
\begin{equation}\label{recCm0}\mathsf{C}_n(t)= \text{\    const.\  } \dfrac{\mathsf{A}_{n+1}(a_{n+1})\mathsf{C}_{n+1}(t)-\mathsf{C}_{n+1}(a_{n+1})\mathsf{A}_{n+1}(t)}{t-a_{n+1}}.\end{equation}

We  check  the orthogonality of the right-hand side  to all polynomials $p$ of degree $<n$ with respect to $d\mu(t)/((t-a_0)\dotsb(t-a_n)(t-b_0)\dotsb (t-b_{n}))$. This is easy for the $\mathsf{A}_{n+1}(t)$ part thanks to the division by $t-a_{n+1}$. And  $\mathsf{C}_{n+1}(t)$ is orthogonal to all polynomials  of degree $\leqslant n$, including the $(t-b_{n+1})p(t)$ polynomials.

Remark that the identities just found are the \emph{ Christoffel-Darboux} equations for orthogonal polynomials \cite{chi,Ism2005,Kru} etc. 

\medskip

$n=0: \mathcal{A}_0(x)=1/(x-a_0),\mathcal{B}_0(x)=1/(x-b_0), \mathsf{A}_0=\mathsf{B}_0  =\mathsf{C}_0=1$.

$n=1:
\begin{bmatrix} 1 &   0                \\
                \alpha_0^{(1)} & 1
\end{bmatrix}
\begin{bmatrix}  \int_S d\mu(t)/((t-a_0)(t-b_0)) & \int_S d\mu(t)/((t-a_0)(t-b_1)) \\
  \int_S d\mu(t)/((t-a_1)(t-b_0)) & \int_S d\mu(t)/((t-a_1)(t-b_1)) \end{bmatrix}
\begin{bmatrix} 1  & \beta_0^{(1)}      \\
                  0          & 1
\end{bmatrix}   \\   =\begin{bmatrix}  \int_S d\mu(t)/((t-a_0)(t-b_0)) & 0 \\
     0  & \int_S d\mu(t)/((t-a_1)(t-b_1))-\alpha_0^{(1)}\beta_0^{(1)}\int_S d\mu(t)/((t-a_0)(t-b_0))    \end{bmatrix} $,

with $\alpha_0^{(1)} =-\dfrac{ \int_S d\mu(t)/((t-a_1)(t-b_0))}{\int_S d\mu(t)/((t-a_0)(t-b_0))}, \beta_0^{(1)}  =-\dfrac{ \int_S d\mu(t)/((t-a_0)(t-b_1))}{\int_S d\mu(t)/((t-a_0)(t-b_0))}$.

So, $\mathcal{A}_1(x)=\dfrac{\alpha_0^{(1)}}{x-a_0}+\dfrac{1}{x-a_1}=
\dfrac{  \mathsf{A}_1(x)=(1+\alpha_0^{(1)})x  -a_0 -\alpha_0^{(1)}a_1 }{(x-a_0)(x-a_1)}
\\  =
\dfrac{ \displaystyle \int_S d\mu(t)\left(\dfrac{(a_1-a_0)(x-t)}{(t-a_0)(t-a_1)(t-b_0)
}\right)}{(x-a_0)(x-a_1)\int_S d\mu(t)/((t-a_0)(t-b_0))}$

It figures: $x-c$ is orthogonal to the constants with respect to some measure  if $c$ is the integral of $t$  with respect to the same measure. Also, $\mathsf{C}_1(x)$ is a constant times the integral of

$(x-t)d\mu(t)/((t-a_0)(t-a_1)(t-b_0)(t-b_1))$. 

\medskip
\subsubsection{Recurrence relations for$\mathsf{A}_m,\mathsf{B}_m,\mathsf{C}_m$\label{recABC}.} We choose the constants in a convenient way by reorganizing the recurrence relations
(\ref{recAm0}-\ref{recCm0}) as
\begin{subequations}
\begin{equation}\label{recAm1}
\mathsf{A}_{m+1}(x)=\mathsf{C}_m(x)+\rho_{2m}(x-b_m)\mathsf{A}_m(x),
\mathsf{C}_{m+1}(x)=\mathsf{A}_{m+1}(x)+\rho_{2m+1}(x-a_{m+1})\mathsf{C}_m(x),\end{equation}
\begin{equation}\label{recBm1}
\mathsf{B}_{m+1}(x)=\mathsf{C}_m(x)+\rho'_{2m}(x-a_m)\mathsf{B}_m(x),
\mathsf{C}_{m+1}(x)=\mathsf{B}_{m+1}(x)+\rho'_{2m+1}(x-b_{m+1})\mathsf{C}_m(x).
\end{equation}
\end{subequations}
These relations can be contracted in the $R_{II}$ relations of Ismail and Masson \cite{IsMa}

  \begin{equation}\label{recAm2}\begin{split}
\mathsf{A}_{m+1}(x)&=(1+\rho_{2m-1}(x-a_m)+\rho_{2m}(x-b_m))\mathsf{A}_m(x)-\rho_{2m-2}\rho_{2m-1}(x-a_m)(x-b_{m-1})\mathsf{A}_{m-1}(x),\\
\mathsf{B}_{m+1}(x)&=(1+\rho'_{2m-1}(x-b_m)+\rho'_{2m}(x-a_m))\mathsf{B}_m(x)-\rho'_{2m-2}\rho'_{2m-1}(x-b_m)(x-a_{m-1})\mathsf{B}_{m-1}(x),\\
\mathsf{C}_{m+1}(x)&=(1+\rho_{2m+1}(x-a_{m+1})+\rho_{2m}(x-b_m))\mathsf{C}_m(x)-\rho_{2m-1}\rho_{2m}(x-a_m)(x-b_{m})\mathsf{C}_{m-1}(x).
\end{split}\end{equation}

\medskip

\subsubsection{Special example.\label{specmathsf}} Let the measure be $d\mu(t)=\sqrt{1-t^2}\,dt$ on $[-1,1]$,
$a_n=a, \forall n, b_m=b \forall  m, a,b\notin [-1,1]$,\ \  $a=(\alpha+1/\alpha)/2,
b=(\beta+1/\beta)/2, \text{\ with\ } |\alpha|, |\beta|<1$. Then,
$\mathsf{A}_n, \mathsf{B}_m, \mathsf{C}_n$ are \emph{Bernstein-Szeg\H o polynomials} \cite[\S 2.6]{Sz}, \cite[\S 5.3]{graphene} 
: with $x=(z+1/z)/2$,

\begin{subequations}
\begin{equation}\label{recAm11}
\mathsf{A}_n(x)=\text{\ constant\ }[z^{n+1}(1-\alpha/z)^{n+1}(1-\beta/z)^n-z^{-(n+1)}(1-\alpha z)^{n+1}(1-\beta z)^n]/(z-1/z),
\end{equation}\begin{equation}\label{recBm11}
\mathsf{B}_m(x)=\text{\ constant\ }[z^{m+1}(1-\alpha/z)^{m}(1-\beta/z)^{m+1}-z^{-(m+1)}(1-\alpha z)^{m}(1-\beta z)^{m+1}]/(z-1/z),
\end{equation}\begin{equation}\label{recCm11}
\mathsf{C}_n(x)=\text{\ constant\ }[z^{n+1}(1-\alpha/z)^{n+1}(1-\beta/z)^{n+1}-z^{-(n+1)}(1-\alpha z)^{n+1}(1-\beta z)^{n+1}]/(z-1/z).
\end{equation}
\end{subequations}
These formulas are symmetric in $z$ and $1/z$, and are indeed polynomials of degree $n$
in $x$, as combinations of Chebyshev polynomials of second kind $U_r(x)=(z^{r+1}-1/z^{r+1})/(z-1/z), r=0,\dotsc, n$.


The orthogonality of $\mathsf{A}_n$ and any polynomial of degree $<n$ with respect to

\begin{multline*}\sqrt{1-t^2}\,dt/((t-a)^{n+1}(t-b)^n) \\
=i(e^{i\theta}-e^{-i\theta})^2
de^{i\theta}(2\alpha)^{n+1}(2\beta)^n/((1-\alpha
e^{i\theta})^{n+1}(1-\alpha e^{-i\theta})^{n+1}(1-\beta
e^{i\theta})^{n}(1-\beta e^{-i\theta})^{n})\end{multline*} is
established by producing two contour integrals in $z=e^{i\theta}$
that can be shrunken   or dilated infinitely, when one tests the
product of $\mathsf{A}_n$ and $(z^r-z^{-r})/[(z-z^{-1})$  or $
z^{\pm r}]/(z-z^{-1}), r=1,\dotsc, n$.

We redefine the constants in (\ref{recAm11}-\ref{recCm11}) as $C_{\mathsf{A}_n}$ etc.

$\mathsf{A}_0=1=[z(1-\alpha/z)-z^{-1}(1-\alpha z)]/(z-1/z), \mathsf{C}_0=1   \\  =(1-\alpha\beta)^{-1}[z(1-\alpha/z)(1-\beta/z)-z^{-1}(1-\alpha z)(1-\beta z)]/(z-1/z)$,

From \eqref{recAm1}, $x-b=-(1-\beta/z)(1-\beta z)/(2\beta)$ is a factor of $\mathsf{A}_{m+1}(x)-\mathsf{C}_m(x)$, so, of  \\
$C_{\mathsf{A}_{m+1}}   z^{m+2}(1-\alpha/z)^{m+2}(1-\beta/z)^{m+1}        -C_{\mathsf{C}_m}z^{m+1}(1-\alpha/z)^{m+1}(1-\beta/z)^{m+1}=0$
when $z=1/\beta$, and $(1/\beta-\alpha)C_{\mathsf{A}_{m+1}}   = C_{\mathsf{C}_m}$ follows.
Furthermore, $x-a$ is a factor of $\mathsf{C}_{m+1}(x)-\mathsf{A}_{m+1}(x)$,
$C_{\mathsf{A}_{m+1}}=(1-\alpha\beta)C_{\mathsf{C}_{m+1}}$,

$C_{\mathsf{A}_{m+1}}/ C_{\mathsf{A}_{m}}= C_{\mathsf{C}_{m+1}}/ C_{\mathsf{C}_{m}}=\beta/(1-\alpha\beta)^2$.


$\rho_{2m}$ and $\rho_{2m+1}$ are the coefficients of $x$ of
$[\mathsf{A}_{m+1}(x)-\mathsf{C}_m(x)]/\mathsf{A}_{m}(x)=C_{\mathsf{A}_{m+1}}(z+\dotsc)/C_{\mathsf{A}_{m}}$ and
$[ \mathsf{C}_{m+1}(x)-\mathsf{A}_{m+1}(x)]/ \mathsf{C}_{m}(x)=
C_{\mathsf{C}_{m+1}}((1-(\alpha\beta)^{m+2})z+\dotsc )-C_{\mathsf{A}_{m+1}}(z+\dotsc )/C_{\mathsf{C}_{m}}(1-(\alpha\beta)^{m+1}+\dotsc)=
[\beta(1-(\alpha\beta)^{m+2}/(1-\alpha\beta)^2-\beta/(1-\alpha\beta)](z+\dotsc)/(1-(\alpha\beta)^{m+1}=\beta(\alpha\beta-(\alpha\beta)^{m+2})/
((1-\alpha\beta)(1-(\alpha\beta)^{m+1}))$, so,

 \begin{equation}\label{rhorho} \rho_{2m}=2\beta/(1-\alpha\beta)^2, \rho_{2m+1}=2\alpha\beta^2/(1-\alpha\beta)^2. \end{equation}


\medskip

Early items: Rahman in 1981 \cite{Rahman81} considered $\displaystyle \langle f,g\rangle=\mathscr{L}(fg)=\int_{S=\mathbb{Z}} f(t)g(t)d\mu(t)=\sum_{-\infty}^\infty w_jf(j)g(j),$ where $w_j=\dfrac{ \Gamma(j+c+c')\Gamma(j+d+d')}{\Gamma(j+c)\Gamma(j+d)}$ for $j\in\mathbb{Z}$ and with $c'+d'<-1$ (so that $w_j$ is bounded by const. $|j|^{c'+d'}$ for large $|j|$). The rational functions $\mathcal{A}_n$ and $\mathcal{B}_m$ have poles at $-c,-c-1,\dotsc,-c-n$ and $-d,-d-1,\dotsc, -d-m$, and are given as ${}_3F_2$ hypergeometric expansions, there is also a ${}_4F_3$ in the same paper. See also Rahman and Suslov \cite{RahmanSuslov} for $q-$lattices, and
 Spiridonov and    Zhedanov on true elliptic lattices, \cite{SpZ1,SpZ2,SpZComm2000,Spi2024} , also Rosengren \cite{RosengrenRahman}.

\subsection{\label{freepol}Rational interpolation with free poles.} \   \\

We encountered in \S \ref{expan1} rational interpolation through expansions of products
$c_n\frac{  (x-x_0)\dotsb(x-x_{n-1})}{(x-x'_0)\dotsb(x-x'_{n-1})}$. The
values of the sum of the $m$ first terms  at $x=x_0,\dotsc, x_m$ will not change
with the next sums, so that this $m^{\text{th}}$ sum is a rational interpolant
of degree $m$ with $m$ imposed poles $x'_0,\dotsc, x'_{m-1}$. We consider now rational
 interpolants with unknown denominators.

\textbf{Proposition.}
Let $f(z)= \mathscr{L}(1/(x-z))$ or, in the integral setting $f(z)=\int_S d\mu(x)/(x-z)$, and let  $w_{m,n}(t)=1/((t-a_0)\dotsc(t-a_m)(t-b_0)\dotsc(t-b_{n})$,
the \emph{numerators} $\mathsf{A}_m, \mathsf{B}_m, \mathsf{C}_m$ of particular biorthogonal rational functions, fulfilling, as see in \S \ref{biorthr} above:
\begin{enumerate}
   \item $\mathsf{A}_m$ orthogonal to all polynomials of degree $<m$ with respect to $w_{m,m-1}(t)d\mu(t)$,
   \item $\mathsf{B}_m$ orthogonal to all polynomials of degree $<m$ with respect to $w_{m-1,m}(t)d\mu(t)$,
   \item $\mathsf{C}_m$ orthogonal to all polynomials of degree $<m$ with respect to $w_{m,m}(t)d\mu(t)$,
\end{enumerate}
 then,
$\mathsf{A}_m, \mathsf{B}_m, \mathsf{C}_m$ are the \emph{denominators} of rational functions
\begin{enumerate}
   \item $\widetilde{\mathsf{A}}_m/\mathsf{A}_m$ of degrees $m$ and $m$ interpolating $f$ at $a_0,\dotsc, a_m, b_0,\dotsc, b_{m-1}$,
   \item $\widetilde{\mathsf{B}}_m/\mathsf{B}_m$ of degrees $m$ and $m$ interpolating $f$ at $a_0,\dotsc, a_{m-1}, b_0,\dotsc, b_{m}$,
   \item $\widetilde{\mathsf{C}}_m/\mathsf{C}_m$ of degrees $m+1$ and $m$ interpolating $f$ at $a_0,\dotsc, a_m, b_0,\dotsc, b_{m}$
\end{enumerate}

Proof. First, note that $\mathsf{A}_m(z)f(z)= \int_S \mathsf{A}_m(t)d\mu(t)/(t-z)+$ a polynomial of degree $<m$. Indeed, $\int_S \frac{\mathsf{A}_m(z)-\mathsf{A}_m(t)}{t-z}d\mu(t)$ is a polynomial of degree $\leqslant m-1$ in $z$.

We perform now successive polynomial interpolation of $\mathsf{A}_mf$ at $a_0, a_1,\dotsc$, adding
terms of the appropriate Newton expansion $[a_0,\dotsc, a_j]_{\mathsf{A}_mf} (z-a_0)\dotsb(z-a_{j-1})$, where $[a_0,\dotsc, a_j]_{\mathsf{A}_mf}$ is the divided difference of $\mathsf{A}_mf$. The divided difference
of the part involving $\mathsf{A}_m(t)d\mu(t)/(t-z)$ is \\   $\int_S \mathsf{A}_m(t)d\mu(t)/((t-a_0)\dotsb(t-a_j))$,
as seen in \S \ref{expan1}.

 After $a_m$ and interpolation of degree $m$, we have \\
$[a_0,\dotsc, a_m,b_0, \dotsc,b_j]_{\mathsf{A}_mf} (z-a_0)\dotsb(z-a_{m})(z-b_0)\dotsb(z-b_{j-1})$ terms,  the divided difference is now the integral of $\mathsf{A}_m(t)w_{m,j}(t)d\mu(t)=
\mathsf{A}_m(t)w_{m,m-1}(t)p(t)d\mu(t)$, where $p(t)=(t-b_{j+1})\dotsb(t-b_{m-1})$ has degree $m-j-1<m$, for $j=0,\dotsc, m-1   $, so that  these divided differences vanish thanks to
the  orthogonality properties of $\mathsf{A}_m$.  The polynomial of degree $\leqslant m-1$ in $z$ of above is no more involved here. The remainder starts with the first nonvanishing
divided difference times \\  $(z-a_0)\dotsb(z-a_{m})(z-b_0)\dotsb(z-b_{m-1})$.

The proofs for $\mathsf{B}_m$ and $\mathsf{C}_m$ are similar. \qed

\medskip

So, $\dfrac{\widetilde{\mathsf{A}}_0(x)}{\mathsf{A}_0(x)}= \dfrac{f(a_0)}{1};\dfrac{\widetilde{\mathsf{C}}_0(x)}{\mathsf{C}_0(x)}=\dfrac{(b_0-a_0)^{-1}[(x-a_0)f(b_0)-(x-b_0) f(a_0)]}{1};\\   \dfrac{\widetilde{\mathsf{A}}_1(x)}{\mathsf{A}_1(x)}= f(a_0)+\dfrac{(f(a_1)-f(a_0))(x-a_0)}{(1+\alpha_0^{(1)})(x-a_1)+a_1  -a_0}$. Check the interpolation at $x=b_0$:

$\mathsf{A}_1(b_0)=(1+\alpha_0^{(1)})(b_0-a_1)+a_1-a_0 =\left(1-\dfrac{ \int_S d\mu(t)/((t-a_1)(t-b_0))}{\int_S d\mu(t)/((t-a_0)(t-b_0))}\right)(b_0-a_1)+a_1-a_0
= \left(1- \dfrac{(b_0-a_0)(f(b_0)-f(a_1))}{(b_0-a_1)(f(b_0)-f(a_0))}\right)(b_0-a_1)+a_1-a_0= \dfrac{(b_0-a_0)(f(a_1)-f(a_0))      }{f(b_0)-f(a_0)}$.

\medskip
\emph{Recurrence relations} for the $\mathsf{A}$s and the $\mathsf{C}$s.
We extend \S\ref{recABC} to the numerators polynomials as well, with new proofs.

\textbf{Proposition.} Let the biorthogonal rational functions $\mathcal{A}_m$ and $\mathcal{B}_m$ be uniquely defined up to multiplication by constants of the numerators, then
there is a choice of these numerators, and coefficients $\rho_0,\rho_1,\dotsc$  such that
\begin{subequations}
\begin{equation} \label{recAm} \mathsf{A}_{m+1}(x)=\mathsf{C}_m(x)+\rho_{2m}(x-b_m)\mathsf{A}_{m}(x),\widetilde{\mathsf{A}}_{m+1}(x)=\widetilde{\mathsf{C}}_m(x)+\rho_{2m}(x-b_m)\widetilde{\mathsf{A}}_{m}(x),\end{equation}
\begin{equation} \label{recCm}\mathsf{C}_{m+1}(x)=\mathsf{A}_{m+1}(x)+\rho_{2m+1}(x-a_{m+1})\mathsf{C}_{m}(x),\widetilde{\mathsf{C}}_{m+1}(x)=\widetilde{\mathsf{A}}_{m+1}(x)+\rho_{2m+1}(x-a_{m+1})\widetilde{\mathsf{A}}_{m}(x).\end{equation}

\end{subequations}
Indeed, we performed polynomial interpolation of $\mathsf{A}_mf$ at the $2m+2$ points $a_0, b_0,\dotsc, a_m, b_m$ resulting in $\widetilde{\mathsf{A}}_m$ of degree $m$
followed by the first nonvanishing term involving \\  $(x-a_0)(x-b_0)\dotsb (x-a_m)(x-b_{m-1})$.

$\mathsf{A}_m(x)f(x)-\widetilde{\mathsf{A}}_m(x)= \gamma_m (x-a_0)\dotsb(x-a_m)(x-b_0)\dotsb(x-b_{m-1})+\dotsb $,

$\mathsf{C}_m(x)f(x)-\widetilde{\mathsf{C}}_m(x)= \delta_m (x-a_0)\dotsb(x-a_m)(x-b_0)\dotsb(x-b_{m})+\dotsb $,

so that the polynomials $X_m(x)=\mathsf{C}_m(x)+r_{2m}(x-b_m)\mathsf{A}_{m}(x)$ and $\widetilde{X}(x)=\widetilde{\mathsf{C}}_m(x)+r_{2m}(x-b_m)\widetilde{\mathsf{A}}_{m}(x)$ are polynomials of
degrees $m+1$ achieving  interpolation at $a_0,\dotsc, a_m, b_0,\dotsc, b_m$ when $r_{2m}=-\delta_m/\gamma_m$, see also the first equation of \eqref{recAm0}.

$\mathsf{C}_m(x)f(x)-\widetilde{\mathsf{C}}_m(x)= \delta_m (x-a_0)\dotsb(x-a_m)(x-b_0)\dotsb(x-b_{m})+\dotsb $,

$\mathsf{A}_{m+1}(x)f(x)-\widetilde{\mathsf{A}}_{m+1}(x)= \gamma_{m+1} (x-a_0)\dotsb(x-a_{m+1})(x-b_0)\dotsb(x-b_{m})+\dotsb $,

proceeds in the same way, see also  \eqref{recCm0}.

First approximants of $f(x)=f(a_0)+\cfrac{(f(b_0)-(a_0))(x-a_0)/(b_0-a_0)}{1+\cfrac{\rho_0(x-b_0)}{1+\ddots}}$: $\widetilde{\mathsf{A}}_0(x)/\mathsf{A}_0(x)=f(a_0)/1,
\widetilde{\mathsf{C}}_0(x)/\mathsf{C}_0(x)=[f(a_0)+\cfrac{(f(b_0)-(a_0))(x-a_0)/(b_0-a_0)}]/1,
\widetilde{\mathsf{A}}_1(x)/\mathsf{A}_1(x)=f(a_0)+[(f(b_0)-(a_0))(x-a_0)/(b_0-a_0)]/(1+\rho_0(x-b_0))$, etc.


\medskip

We return to \S \ref{specmathsf} and \eqref{rhorho},


$f(x)=\cfrac{\rho_1(x-a)}{1+\cfrac{\rho_0(x-b)}{1+f(x)}}$
 is a root of \\
$ f^2+(1+\rho_0(x-b)-\rho_1(x-a))f-\rho_1(x-a)=0$, so,
$f(x)=\beta[ -x+\alpha  \pm \beta\sqrt{x^2-1}]/(1-\alpha\beta)$ .

See \cite[Example 3.1]{IsMa} for more on the choice of the root, convergence and integral representation.  We also have a very special  $M-$ and $T-$ continued fraction of \S \ref{def3}.

Let $P_{2n}$ be a new definition of $\mathsf{A}_n$ through multiplication by another constant, and the same with $P_{2n+1}$ and $\mathsf{C}_n$, for $n=0,1,...$, so,

$ P_{n+1}(x)=P_n(x)+\rho_n
\left( \left\{\begin{matrix} (x-a_{n/2}) \text{\ if\ } n \text{\ is\ even}\\
   (x-b_{(n-1)/2}) \text{\ if\ } n \text{\ is\ odd} \end{matrix}\right.
\right)P_{n-1}(x), n=0,1,\dotsc, P_{-1}=0$ Ismail \& Masson  \cite[\S 3]{IsMa}, Spiridonov \& Zhedanov \cite{SpZ1,SpZ2},

$f(x)-f(a_0)=\dfrac{\rho_0(x-a_0) }{1+\dfrac{ \rho_1(x-b_0) }{1+\dfrac{\rho_2(x-a_1)}
 {1+\dfrac{\rho_3(x-b_1)}{1+\ddots}}}}
$ Multipoint Pad\'e: Goncar, L\'opez,
 Zhedanov \cite{Gon5,Lo1,Lo2,ZhedPad}

With $N_0, N_1,\dotsc$ the corresponding numerators, and $x_{2j}=a_j,  x_{2j+1}=b_j, j=0,1,\dotsc$

$\dfrac{N_0}{P_0}=\dfrac{f(a_0)}{1}, \dfrac{N_1(x)}{P_1(x)}=\dfrac{f(a_0)+\rho_0(x-x_0)}{1},
   \dfrac{N_2(x)}{P_2(x)}=f(a_0)+\dfrac{\rho_0(x-x_0)}{1+r_1(x-x_1)}$,

with $x_0=a_0, x_1=b_0, \dotsc, \rho_0=\frac{f(x_1)}{x_1-x_0}$,
$ \dfrac{\rho_0(x-x_0) }{1+\dfrac{\rho_1(x-x_1)}{1 +   \dfrac{\rho_2(x-x_2) }{1}   }}=
\dfrac{\rho_0(x-x_0)[1 +   \rho_2(x-x_2)]}{1 +   \rho_1(x-x_1)+ \rho_2(x-x_2)  }$

\subsubsection{\label{singlex}Remark.} From now on, $\{x_n,y_n\}$  and $\{x'_n,y'_n\}$ are successive points on  a given biquadratic curve $F(x,y)=0$. There are no more independent sequences  $\{a_n\}, \{b_n\}$, and the  interpolation set $a_0,\dotsc, a_m, b_0,\dotsc, b_{m-1}$
   of \S \ref{freepol}
is  $x_0,\dotsc, x_{2m}$. Then, $\mathsf{A}_m$ is called $P_{2m}$ here, together with $\mathsf{C}_m=P_{2m+1}$.

\bigskip

\subsubsection{\label{Thielef}}
Thiele \cite[\S 40]{Thiele} builds rational interpolants to a function $f$ at
$x_0, x_1,\dotsc$ as successive approximants of the continued fraction
$f(x)-f(x_0)=\dfrac{x-x_0}{\delta_0 + \dfrac{x-x_1}{\delta_1 + \ddots }}$, where $\delta_0, \delta_1,\dotsc$
are related to $f(x_0), f(x_1),\dotsc$ by a particular algorithm (\emph{reciprocal differences}
algorithm: $\delta_0=(x_1-x_0)/(f(x_1)-f(x_0))$, etc). If the relevant determinants and the $\delta-$coefficients do not vanish,
the $m^{\text{th}}$ approximant of the continued fraction above is the unique
interpolation at $x_0,\dotsc, x_{m}$ of $f$ by a polynomial of degree $n_m=$integer part of $(m+1)/2$ over a polynomial of degree $d_m=$ integer part of $m/2$.

See   \cite{Ce1,Ce2,ACm18}, Peter Wynn's $\rho-$algorithm  \cite[\S 2.5.3]{Brez2020},
see also \cite{Cuyt1,GravesMorris,Gutchn,MCa,We1,We2,Zou} for existence, unicity, and numerical stability  discussions.

 Let $f(x_0)=f(a_0)=0$.  Consider now the continued fraction
\begin{equation} \label{rxmx}f(x)=
    \dfrac{\rho_0(x-x_0) }{1+\dfrac{\rho_1(x-x_1)}{1+\ddots}}
\end{equation}
equivalent to Thiele's fraction ($\delta_0=1/\rho_0, \delta_1=\rho_0/\rho_1, \delta_2=\rho_1/(\rho_0\rho_2)$, etc.)

Reminding that $f(x_0)=0$, the first steps are

$\dfrac{N_0}{P_0}=\dfrac{0}{1}, \dfrac{N_1(x)}{P_1(x)}=\dfrac{\rho_0(x-x_0)}{1},
   \dfrac{N_2(x)}{P_2(x)}=\dfrac{\rho_0(x-x_0)}{1+\rho_1(x-x_1)}$,
with $\rho_0=\frac{f(x_1)}{x_1-x_0}$,

$\rho_1=\frac{1}{x_2-x_1}  \left( \frac{  (x_2-x_0)f(x_1)}{(x_1-x_0)f(x_2)} -1\right)=\frac{(x_2-x_0)   [  f(x_1)/(x_1-x_0)-f(x_2)/(x_2-x_0)]}{(x_2-x_1)f(x_2)}  $
etc.,
and $N_m/P_m$ interpolates at $x_0,\dotsc, x_m$, but $N_m$ and $P_m$ are
polynomials of degrees  $n_m$ and $d_m$, the integer parts of $(m+1)/2$ and $m/2$.

\begin{equation}\label{matrxmx}
   \begin{bmatrix} N_m(x) & N_{m+1}(x)  \\ P_m(x) & P_{m+1}(x)\end{bmatrix} =
  \begin{bmatrix} N_{m-1}(x) & N_{m}(x)  \\ P_{m-1}(x) & P_{m}(x)\end{bmatrix}
  \begin{bmatrix} 0 & \rho_m(x-x_m)  \\ 1 & 1\end{bmatrix}
\end{equation}

A convenient prequel is $N_{-1}=1, P_{-1}=0$.

Casorati determinant:  let $\mathfrak{C}_m(x)
=N_{m-1}(x) P_m(x)-P_{m-1}(x)  N_m(x)$, \\  then
$\mathfrak{C}_{m+1}(x)=-r_m(x-x_m)\mathfrak{C}_m(x)$,
\begin{equation}\label{Caso}\mathfrak{C}_m(x)=(-1)^m \rho_0\dotsb \rho_{m-1}(x-x_0)\dotsb(x-x_{m-1})\end{equation}

$f(x)=
    \dfrac{\rho_0(x-x_0) }{1+\dfrac{\ddots}{1+\dfrac{\rho_{m-1}(x-x_{m-1})}{1+f_{m}(x)}}}
= \dfrac{ N_{m}(x)   +N_{m-1}(x)f_{m}(x))     }{P_{m}(x)+P_{m-1}(x)f_{m}(x)}
$,

where $f_0=f$, and
\begin{equation} \label{rxmxm}f_m(x)=\dfrac{\rho_m(x-x_m)}{1+\dfrac{\rho_{m+1}(x-x_{m+1})}{1+\ddots}}.\end{equation}

A very basic continued fraction identity (Perron etc. ) check the
$m\to m+1$ step:

$ \dfrac{ N_{m}(x)   +N_{m-1}(x)\rho_m(x-x_m)/(1+f_{m+1}(x)))     }{P_{m}(x)+P_{m-1}(x)r_m(x-x_m)/(1+f_{m+1}(x)))     }=
\dfrac{ N_{m+1}(x)   +N_{m}(x)f_{m+1}(x)    }{P_{m+1}(x)+P_{m}(x)f_{m+1}(x)     }$

\bigskip

\emph{Contraction:}

\begin{equation} \label{contrac}
 P_{2m+2}(x)= \left(1  +\rho_{2m}(x-x_{2m})  +\rho_{2m+1}(x-x_{2m+1})  \right)  P_{2m}(x)
                                   -\rho_{2m-1} \rho_{2m}  (x-x_{2m-1})(x-x_{2m})   P_{2m-2}(x),\end{equation}
type II of Ismail \& Masson \cite{biorthKiran,IsMa}.

Koekoek et al \cite{Koe} often give  $xp_n(x) = p_{n+1}(x)+(A_n+C_n) p_n(x)+A_{n-1}C_n p_{n-1}(x)$, we recognize $(A_n,C_n)$ basically as $(\rho_{2m},\rho_{2m+1})$

 The $R_I$ fractions of \cite{IsMa} are NOT
\eqref{rxmx} but continued fractions based on recurrence relations
$X_{m+1}=(x-\alpha_m)X_m-\lambda_m(x-x_m)X_{m-1}$. 

\subsubsection{\label{orthpol}Orthogonal polynomials.}

Let all the $a_n$s and $b_n$s be a single value $x^*$,
we have from \eqref{contrac} $P_{2m+2}(x)/(x-x^*)^{m+1}= (1/(x-x^*)  +\rho_{2m}  +\rho_{2m+1} )  P_{2m}(x)/(x-x^*)^m
         -\rho_{2m-1} \rho_{2m}    P_{2m-2}(x)/(x-x^*)^{m-1}$, so that $P_{2m}(x)/(x-x^*)^m$ is a
polynomial of degree $m$ in $1/(x-x^*)$ satisfying the three-term recurrence relations of  \emph{orthogonal polynomials}
(plain polynomials if $x^*=\infty$). See also \cite{HazeP}.




\subsection{Riccati difference equation for $f$\label{Riccatif}}  \   \\

\subsubsection{}\textbf{Theorem.} Let $\{(x_n,y_n)\}$, $\{(x_{n+1},y_n)\}$ be points
on the biquadratic curve $F(x,y)=Y_0(y)+xY_1(y)+x^2Y_2(y)=0$, and $f$ a function
defined on the $x_n$s, satisfying the \emph{discrete Riccati} difference equation
 at some $y=y_n$,

\begin{subequations}
\begin{equation}\label{Riccati0}
 A(y_n)(\mathcal{D}f)(y_n)=B(y_n)f(x_n)f(x_{n+1})
          +C(y_n)(\mathcal{M}f)(y_n)+D(y_n),
\end{equation}
or, showing how $f(x_{n+1})$ is related to $f(x_n)$:
\begin{equation}\label{Riccati02}
 f(x_{n+1})= \dfrac{    \left[\dfrac{A(y_n)}{x_{n+1}-x_n}+\dfrac{C(y_n)}{2}\right]f(x_n)  +D(y_n)  }
              {-B(y_n)f(x_n)+\dfrac{A(y_n)}{x_{n+1}-x_n}-\dfrac{C(y_n)}{2}    },
\end{equation}
\end{subequations}
and suppose we also have
\begin{equation}\label{singmm0}
\dfrac{A(y_{-1})}{x_{0}-x_{-1}}+\dfrac{C(y_{-1})}{2}=0,\ \ \  D(y_{-1})=0,\end{equation}
showing that $f(x_0)$ does NOT depend on the value of $f(x_{-1})$, that is, that
$x_0$ is a \emph{singular} point of the Riccati equation. The value of $f(x_0)$ happens to be 0.

In \eqref{Riccati0},   $A,\dotsc, D$ are polynomials. If needed, these polynomials are multiplied
by a common factor so that $Y_2$ is a factor of $B, C$, and $D$. Let $d$ be then
the largest degree of $A,\dotsc, D$.

Let us consider now at some $y=y_n$, the discrete Riccati equation

\begin{equation*} A_m(y_n)(\mathcal{D}f_m)(y_n)=B_m(y_n)f_m(x_n)f_m(x_{n+1})
          +C_m(y_n)(\mathcal{M}f_m)(y_n)+D_m(y_n), \ \ m=0,1,\dotsc,
\end{equation*}
with $f_0=f, A_0=A,$ etc., as in \S \ref{Riccatim}, and

$f_m(x)=\dfrac{ \rho_m(x-x_m)   }{1+f_{m+1}(x) } $, as seen in \eqref{rxmxm}.

Then, \begin{enumerate}

\item With  the rational functions $\Upsilon_m(y)=Y_2(y)/(y-y_{m-1})$ if $m$ is even; $1/(y-y_{m-1})$ if $m$ is odd,

\begin{subequations}

 from $f_{m+1}(x_{m+1})=0, f_m(x_{m+1})=\rho_m(x_{m+1}-x_m) $,

\begin{equation}\label{rm}
\rho_m= \dfrac{   D_m(y_m)  }  {A_m(y_m)-(x_{m+1}-x_m) C_m(y_m)/2  }=\dfrac{   D_m(y_m)  }{ A_m(y_m) + (Y_1(y_m)/2 +x_mY_2(y_m)) C_m(y_m)/Y_2(y_m) },
\end{equation}

\begin{equation}\label{recA}\begin{split}
  A_{m+1}&=   \Upsilon_m \left[ \dfrac{Y_1/2 +x_mY_2}{Y_2}    A_m +\dfrac{Q}{4Y_2^2}C_m\right],\end{split}  \end{equation}
\begin{equation}\label{recB}
 B_{m+1}=\Upsilon_m D_m/\rho_m,\end{equation}
\begin{equation}\label{recC}
 C_{m+1}= \Upsilon_m\left[ -   A_m +2D_m/\rho_m -\dfrac{Y_1/2 +x_mY_2}{Y_2} C_m  \right],\end{equation}
\begin{multline}\label{recD} D_{m+1}(y)=\Upsilon_m(y)
      \left[B_m(y)\rho_m\dfrac{F(x_m,y)}{Y_2(y)}
  - A_m(y)-\dfrac{Y_1/2 +x_mY_2}{Y_2(y)}C_m(y)   +D_m(y)/\rho_m\right] \\
=C_{m+1}(y)+\Upsilon_m(y)
      \left[B_m(y)\rho_m\dfrac{F(x_m,y)}{Y_2(y)}  -D_m(y)/\rho_m\right]\end{multline}

\smallskip

\smallskip

\item  The degrees of $A_m,$ etc. are $\leqslant d+ 1$,

\item $Y_2$ is a factor of $B_m, C_m$, and $D_m$,

\item $x_m$ is a singular point of  \eqref{Riccatim}: \begin{equation}\label{singmm1}
\dfrac{A_m(y_{m-1})}{x_{m}-x_{m-1}}+\dfrac{C_m(y_{m-1})}{2}=0,
D_m(y_{m-1})=0,\end{equation}

\end{subequations}

\end{enumerate}

Sometimes $\widetilde{A}_m= A_m + \dfrac{Y_1/2 +x_mY_2}{Y_2} C_m$ is preferred,
then, {\small
\begin{subequations}
 \begin{equation}\label{recAmp1}A_{m+1}=   \Upsilon_m \left[ \dfrac{Y_1/2 +x_mY_2}{Y_2}    A_m +\dfrac{Q}{4Y_2^2}C_m\right]=    \Upsilon_m \left[ \dfrac{Y_1/2 +x_mY_2}{Y_2}   \widetilde{ A}_m -\dfrac{F(x_m,y)}{Y_2}C_m \right],\end{equation}

$\widetilde{A}_{m+1}=\Upsilon_m \left[ \dfrac{Y_1/2 +x_mY_2}{Y_2} \widetilde{ A}_m -\dfrac{F(x_m,y)}{Y_2}C_m   + \dfrac{Y_1/2 +x_{m+1}Y_2}{Y_2}\left[ -  \widetilde{ A}_m +2D_m/\rho_m  \right]\right] $, or
   \begin{equation}\label{recAt} \widetilde{A}_{m+1} =\Upsilon_m \left[ (x_m-x_{m+1})  \widetilde{ A}_m-\dfrac{F(x_m,y)}{Y_2}C_m + \dfrac{Y_1 +2x_{m+1}Y_2}{\rho_mY_2}D_m\right].\end{equation}

\end{subequations}
}
\subsubsection{}\textbf{Proof.} We have \eqref{Riccatim} at $m=0$ to begin with,
if \eqref{Riccatim} holds at some $m$, we adapt \eqref{Riccati02} at this level $m$ showing
 how $f_m(x_{n+1})$ is related to $f_m(x_n)$:
\begin{equation*}
 f_m(x_{n+1})= \dfrac{    \left[\dfrac{A_m(y_n)}{x_{n+1}-x_n}+\dfrac{C_m(y_n)}{2}\right]f_m(x_n)  +D_m(y_n)  }
              {-B_m(y_n)f_m(x_n)+\dfrac{A_m(y_n)}{x_{n+1}-x_n}-\dfrac{C_m(y_n)}{2}    },
\end{equation*}
as in \eqref{Riccatim2}
Let
$U_m(y_n)=\dfrac{C_{m}(y_n)}{2}+ \dfrac{A_{m}(y_n)}{x_{n+1}-x_n}, V_m(y_n)=\dfrac{C_{m}(y_n)}{2}- \dfrac{A_{m}(y_n)}{x_{n+1}-x_n}$,

 and we
perform the step $m\to m+1$ by entering $f_m(x)=\dfrac{ \rho_m(x-x_m)  }{1+f_{m+1}(x) } $
in \eqref{Riccatim2}:

$ \dfrac{ \rho_m(x_{n+1}-x_m)   }{1+f_{m+1}(x_{n+1})  }= - \dfrac{    U_m(y_n)\dfrac{\rho_m(x_n-x_m)  }{1+f_{m+1}(x_n) }   +D_m(y_n)  }
              {B_m(y_n)\dfrac{ \rho_m(x_n-x_m)   }{1+f_{m+1}(x_n) } +V_m(y_n)}
$

and isolate
\[  f_{m+1}(x_{n+1})=-1-\rho_m(x_{n+1}-x_m)  \dfrac{
       B_m(y_n)\rho_m(x_n-x_m)+V_m(y_n)\left[1+f_{m+1}(x_n)\right]  }
    {U_m(y_n) \rho_m (x_n-x_m)  +D_m(y_n)\left[1+f_{m+1}(x_n)\right]
             }
\]

$=-\dfrac{ \left[\rho_m(x_{n+1}-x_m)  V_m(y_n) +D_m(y_n)
    \right]f_{m+1}(x_n)+X
   }{
   D_m(y_n)f_{m+1}(y_n) +\rho_m(x_{n}-x_m)U_m(y_n)+D_m(y_n)
  }$

with $X=B_m(y_n)\rho_m^2(x_n-x_m) (x_{n+1}-x_m) +  \rho_m(x_{n}-x_m)U_m(y_n)+ \rho_m(x_{n+1}-x_m)   V_m(y_n)+D_m(y_n)$

must be $-\dfrac{  U_{m+1}(y_n)f_{m+1}(x_n)+D_{m+1}(y_n) }{B_{m+1}(y_n)f_{m+1}(x_n)+V_{m+1}(y_n) }$, whence

\begin{equation}\label{recUV}\begin{split}
     U_{m+1}(y_n)&=\Upsilon_m(y_n)  \left[(x_{n+1}-x_m)  V_m(y_n) +D_m(y_n)    /\rho_m \right]\\
   D_{m+1}(y_n)&\   \\
=\Upsilon_m(y_n)  &\left[B_m(y_n)\rho_m(x_n-x_m) (x_{n+1}-x_m) +   (x_{n}-x_m)U_m(y_n)+(x_{n+1}-x_m)   V_m(y_n)+D_m(y_n) /\rho_m\right]\\
 B_{m+1}(y_n)&=\Upsilon_m(y_n)   D_m(y_n)/\rho_m\\
V_{m+1}(y_n) &= \Upsilon_m(y_n)  \left[(x_{n}-x_m)U_m(y_n)+D_m(y_n)/\rho_m\right]
\end{split}\end{equation}

We need only now  the $x$s and $y$s to be on a biquadratic curve, by \eqref{sumprod},  \eqref{ydirect},
$x_n+x_{n+1}=-Y_1(y_n)/Y_2(y_n)$ and $(x_{n+1}-x_n)^2=Q(y_n)/Y_2^2(y_n)$ are rational functions of $y_n$. We also have

$F(x_m, y_n)=(x_n-x_m)(x_{n+1}-x_m)Y_2(y_n)=X_2(x_m)(y_n-y_m)(y_n-y_{m-1})$,

We get \eqref{recA}-\eqref{recD} through $A_m(y_n)/(x_{n+1}-x_n)=(U_m(y_n)-V_m(y_n))/2, C_m(y_n)=U_m(y_n)+V_m(y_n)$,

\begin{equation}\label{recAC} \dfrac{C_{m+1}(y_n)}{2}\pm\dfrac{A_{m+1}(y_n)}{x_{n+1}-x_n}=\Upsilon_m(y_n)\left[\left(\left\{\begin{matrix}(x_{n+1}-x_m)\\ (x_n-x_m)\end{matrix}\right)\left[\dfrac{C_m(y_n)}{2}\mp
  \dfrac{A_m(y_n)}{x_{n+1}-x_n}\right]+D_m(y_n)/\rho_m\right)\right],\end{equation}

 $\Upsilon_m$ is not yet known, it  must be   such that $A_{m+1}$ etc. are polynomials.

\medskip

Matrix form:

$ \dfrac{\rho_m}{\Upsilon_m(y_n)} \begin{bmatrix}U_{m+1} &B_{m+1} \\ D_{m+1} &V_{m+1} \end{bmatrix}=
  \underbrace{ \begin{bmatrix}D_{m} &V_{m} \\ V_{m+1}/\Upsilon_m(y_n) &\rho_m(x_n-x_0)B_m+V_{m} \end{bmatrix}  }_{\displaystyle
   \begin{bmatrix}0&1 \\  \rho_m(x_n-x_0) &1 \end{bmatrix}
   \begin{bmatrix}U_{m} &B_{m} \\ D_{m} &V_{m} \end{bmatrix}        }
\begin{bmatrix}1&1 \\  \rho_m(x_{n+1}-x_0) &0 \end{bmatrix}
$,  or

\begin{multline}\label{ABCDmat}
\begin{bmatrix} \dfrac{C_{m+1}(y_n)}{2}+\dfrac{A_{m+1}(y_n)}{x_{n+1}-x_n}&
B_{m+1}(y_n) \\  D_{m+1}(y_n) & \dfrac{C_{m+1}(y_n)}{2}-\dfrac{A_{m+1}(y_n)}{x_{n+1}-x_n}\end{bmatrix}  \\  =  \dfrac{\Upsilon_m(y_n)}{\rho_m}
\begin{bmatrix} 0 & 1 \\ \rho_m (x_{n}-x_m) & 1 \end{bmatrix}
\begin{bmatrix} \dfrac{C_{m}(y_n)}{2}+\dfrac{A_{m}(y_n)}{x_{n+1}-x_n}&
B_m(y_n) \\  D_m(y_n) & \dfrac{C_{m}(y_n)}{2}-\dfrac{A_{m}(y_n)}{x_{n+1}-x_n}\end{bmatrix}
\begin{bmatrix} 1 & 1 \\ \rho_m(x_{n+1}-x_m) & 0 \end{bmatrix}
\end{multline}

Note the determinant of \eqref{ABCDmat}: $ \dfrac{C^2_{m+1}(y_n)}{4}-\dfrac{A^2_{m+1}(y_n)}{(x_{n+1}-x_n)^2}-B_{m+1}(y_n)D_{m+1}(y_n) \\
=\Upsilon^2_m(y_n)  \underbrace{(x_{n}-x_m)(x_{n+1}-x_m)  }_{F(x_m,y_n)/Y_2(y_n)=X_2(x_m)(y_n-y_{m-1})(y_n-y_m)/Y_2(y_n)  }
      \left[
  \dfrac{C^2_{m}(y_n)}{4}-\dfrac{A_m^2(y_n)}{(x_{n+1}-x_n)^2}
-B_m(y_n)D_m(y_n)\right]$, or

\begin{equation}\label{detABCD} \dfrac{C^2_{m}}{4}-\dfrac{A^2_{m}}{Q/Y_2^2}-B_{m}D_{m}
=\dfrac{\Upsilon^2_0\dotsb \Upsilon^2_{m-1}  F(x_0,y)\dotsb F(x_{m-1},y)}{Y_2^m}
\left[\dfrac{C^2_{0}}{4}-\dfrac{A^2_{0}}{Q/Y_2^2}-B_{0}D_{0}\right] \end{equation}

\medskip

We will also need the contraction of two steps in one
{\Small
\begin{multline}\label{ABCDmat2}
\begin{bmatrix} \dfrac{C_{m+2}(y_n)}{2}+\dfrac{A_{m+2}(y_n)}{x_{n+1}-x_n}&
B_{m+2}(y_n) \\  D_{m+2}(y_n) & \dfrac{C_{m+2}(y_n)}{2}-\dfrac{A_{m+2}(y_n)}{x_{n+1}-x_n}\end{bmatrix}    =  \dfrac{\Upsilon_m(y_n) \Upsilon_{m+1}(y_n)}{\rho_m\rho_{m+1} }   
\begin{bmatrix} \rho_m(x_n-x_m) &1 \\ \rho_m (x_{n}-x_m) & (x_n-x_{m+1})\rho_{m+1}+1 \end{bmatrix}   \\    \times
\begin{bmatrix} \dfrac{C_{m}(y_n)}{2}+\dfrac{A_{m}(y_n)}{x_{n+1}-x_n}&
B_m(y_n) \\  D_m(y_n) & \dfrac{C_{m}(y_n)}{2}-\dfrac{A_{m}(y_n)}{x_{n+1}-x_n}\end{bmatrix}
\begin{bmatrix} (x_{n+1}-x_{m+1})\rho_{m+1}+1 & 1 \\ \rho_m(x_{n+1}-x_m) &\rho_m(x_{n+1}-x_m) \end{bmatrix}
\end{multline}    }
\medskip

We proceed now with the proof, supposing \eqref{singmm1} at level $m$

\begin{enumerate}

\item  $y-y_{m-1}$ is a common factor to the right-hand sides of \eqref{recA}-\eqref{recD}
before the multiplication by $\Upsilon$: the part in square brackets  of \eqref{recA2} has
$\widetilde{A}_m$ vanishing at $y=y_{m-1}$, and so does $F(x_m,y)$; obvious for the $B_{m+1}$ line, as $D_m$ does the job; and the same for $C_{m+1}$ using $\widetilde{A}_m$ and $D_m$ in \eqref{recC}; finally, $\widetilde{A}_m$, $F(x_m,y)$ and $D_m$ do the job in \eqref{recD}.

Proof that \eqref{singmm1} still holds at $m+1$:

From \eqref{rm},
$\rho_m= \dfrac{   D_m(y_m)  }  {A_m(y_m)-(x_{m+1}-x_m) C_m(y_m)/2 =\widetilde{A}_m(y_m) }$,
and \eqref{recAt},  $\widetilde{A}_{m+1}(y_m) =\Upsilon_m(y_m) \left[ (x_m-x_{m+1})  \widetilde{ A}_m(y_m)-\dfrac{F(x_m,y_m)}{Y_2}C_m(y_m) \right.   \\    \left.    + \dfrac{Y_1(y_m) +2x_{m+1}Y_2(y_m)=(x_{m+1}-x_m)Y_2}{\rho_mY_2}D_m(y_m)\right]=0$,

and \eqref{recD}, $D_{m+1}(y_m) =\Upsilon_m(y_m) \left[B_m(y_m)\rho_m\dfrac{F(x_m,y_m)}{Y_2(y_m)}
  - \widetilde{A}_m(y_m)   +D_m(y_m)/\rho_m\right]=0$.

A much shorter, perhaps too short proof:

from $0=f_m(x_m)=\rho_m(x_m-x_m=0)/(1+f_{m+1}(x_m))$, $f_{m+1}(x_m)=0/0$ is indefinite,
and $f_{m+1}(x_{m+1})=0$ for any $f_{m+1}(x_m)$,  \eqref{singmm1} follows at $m+1$.

\item

The degrees of $A_m$ etc. $\leqslant d$ when $m$ is even, and $Y_2$ being a factor of $\Upsilon_m$, $A_{m+1}$ etc. are of course polynomials (of degrees  $\leqslant d+1$),
and $Y_2$ is a factor of  $B_{m+1},C_{m+1},D_{m+1}$.

Are $A_{m+2}$ etc. then still polynomials?
We look at \eqref{ABCDmat2} leading to

\begin{subequations}

\begin{equation}\label{recA2}\begin{split}
  A_{m+2}&=   \Upsilon_m\Upsilon_{m+1} \left[\dfrac{Y_0 +(x_m+x_{m+1})Y_1/2+x_mx_{m+1}Y_2}{Y_2}    A_m +\dfrac{Q}{4Y_2^2}[(x_{m+1}-x_m)C_m+2D_m/\rho_m]\right],\end{split}  \end{equation}
\begin{equation}\label{recB2}
 B_{m+2}=-\dfrac{\Upsilon_m\Upsilon_{m+1}}{\rho_{m+1}}\left[A_m-\dfrac{F(x_m,y)\rho_mB_m-(Y_1/2+x_m)C_m}{Y_2}- D_m/\rho_m\right],\end{equation}
\begin{multline}\label{recC2}
 C_{m+2}= \Upsilon_m\Upsilon_{m+1}\left[ [x_{m+1}-x_m-2/\rho_{m+1}]   A_m +2\rho_m\dfrac{F(x_m,y)}{Y_2}B_m/\rho_{m+1}
 \right. \\   \left.
+ \dfrac{[Y_0+(x_m+x_{m+1})Y_1/2+x_mx_{m+1}Y_2]-Y1-2x_mY_2}{Y_2}  C_m  -\dfrac{(Y_1+2x_{m+1}Y_2)-2Y_2/\rho_{m+1}}{Y_2}D_m/\rho_m \right],\end{multline}
\begin{multline}\label{recD2} D_{m+2}=\dfrac{\Upsilon_m\Upsilon_{m+1} }{\rho_{m+1}}
      \left[   ((x_{m+1}-x_m)\rho_{m+1}-1)A_m+\rho_m\dfrac{F(x_m,y)}{Y_2}B_m  \right.  \\   \left.
+\dfrac{ [Y_0+(x_m+x_{m+1})Y_1/2+x_mx_{m+1}Y_2]\rho_{m+1}-Y_1/2-x_mY_2}{Y_2}C_m   \right.   \\    \left.   +\rho_{m+1}\dfrac{\rho_{m+1}F(x_{m+1},y)-Y_1-2x_{m+1}Y_2+Y_2/\rho_{m+1} }{Y_2}D_m/\rho_m\right] \\
=C_{m+2}(y)+\dfrac{\Upsilon_m\Upsilon_{m+1} }{\rho_m\rho_{m+1}}
      \left[\rho_mA_m-\rho_m^2\dfrac{F(x_m,y)}{Y_2}B_m  +\rho_m\dfrac{Y_1/2+x_{m}Y_2}{Y_2}C_m -\dfrac{\rho_{m+1}^2F(x_{m+1},y)+Y_2}{Y_2}D_m\right]\end{multline}

\end{subequations}
and use $\Upsilon_m(y)\Upsilon_{m+1}(y)=\dfrac{Y_2(y)}{(y-y_{m-1})(y-y_m)}$.
\end{enumerate}

end of proof  \qed

\medskip

\subsubsection{Consequences on the interpolant $N_m/P_m$.}
\     \\

$ A(y_n)\left(\mathcal{D}\dfrac{N_m}{P_m}\right)(y_n)-B(y_n)\dfrac{N_m(x_n)N_m(x_{n+1})}{P_m(x_n)P_m(x_{n+1})}
          -C(y_n)\left(\mathcal{M}\dfrac{N_m}{P_m}f\right)(y_n)-D(y_n)  $
is a rational function of $y_n$ vanishing when $n$ is such that $N_m/P_m$ interpolates $f$ at $x_n$
and $x_{n+1}$, so for $n=0,\dotsc,m-1$, and we have after multiplication by  $P_m(x_n)P_m(x_{n+1})$,
\begin{multline}\label{LagTheta}
 A(y_n)\dfrac{N_m(x_{n+1})P_m(x_n)-N_m(x_m)P_m(x_{n+1}}{x_{n+1}-x_n}-B(y_n)N_m(x_n)N_m(x_{n+1})   \\
          -C(y_n)\dfrac{N_m(x_{n+1})P_m(x_n)+N_m(x_m)P_m(x_{n+1}}{2}-D(y_n)P_m(x_n)P_m(x_{n+1}) \\
=
\Theta_m(y_n)(y_n-y_0)\dotsb(y_n-y_{m-1})
\end{multline}
(almost Laguerre's notation! \cite{Lag1885}) where $\Theta_m$  is a polynomial whose degree
is not larger than the degrees of $A,B,$ etc.
Casorati: \\  {\Small $ A(y_n)\dfrac{P_{m+1}(x_{n+1})N_m(x_{n+1})P_m(x_n)-P_{m+1}(x_{n+1})N_m(x_m)P_m(x_{n+1})
  -P_m(x_{n+1})N_{m+1}(x_{n+1})P_{m+1}(x_n)+P_m(x_{n+1})N_{m+1}(x_n)P_{m+1}(x_{n+1})    }{x_{n+1}-x_n}
  -B(y_n)N_m(x_n)N_m(x_{n+1})
          -C(y_n)\dfrac{N_m(x_{n+1})P_m(x_n)+N_m(x_m)P_m(x_{n+1}}{2}-D(y_n)P_m(x_n)P_m(x_{n+1})   \\
=
\Theta_m(y_n)(y_n-y_0)\dotsb(y_n-y_{m-1})$,               }

$f(x)=
    \cfrac{\rho_0(x-x_0)}{1+\cfrac{\ddots}{1+\dfrac{\rho_{m-1}(x-x_{m-1})}{1+f_{m}(x)}} }
= \dfrac{ N_{m}(x)   +N_{m-1}(x)f_{m}(x)    }{P_{m}(x)+P_{m-1}(x)f_{m}(x)}
$,

where $f_0=f$, and $f_m(x)=\dfrac{\rho_m(x-x_m)}{1+\dfrac{\rho_{m+1}(x-x_{m+1})}{1+\ddots}}$.

A very basic continued fraction identity (Perron etc. ) check the
$m\to m+1$ step:

$ \dfrac{ N_{m}(x)   +N_{m-1}(x)\rho_m(x-x_m)/(1+f_{m+1}(x)))     }{P_{m}(x)+P_{m-1}(x)\rho_m(x-x_m)/(1+f_{m+1}(x)))     }=
\dfrac{ N_{m+1}(x)   +N_{m}(x)f_{m+1}(x)    }{P_{m+1}(x)+P_{m}(x)f_{m+1}(x)     }$

\smallskip

\medskip

\subsubsection{Examples of (\ref{rm}-\ref{recD}).\label{Exxh}}

\begin{enumerate}
\item
Function $f(x)=e^{x-x_0}-1$  interpolated at $x_0,x_0+h,x_0+2h,\dotsc$: $F(x,t)=(y-x)(y-x-h), x_n=x_0+nh, y_n=x_0+(n+1)h,
(f(x+h)-f(x))/h=a((f(x)+f(x+h))/2+1)$, with $f(x_0)=0, a=(2/h)\tanh(h/2)$.

$A_0(y)=y-y_{-1}=y-x_0, B_0(y)\equiv 0, C_0(y)=D_0(y)=a(y-y_{-1})$, as we must have \eqref{singmm0},    $A_0$ and $D_0$ must vanish at $y=y_{-1}=x_0$,
$\rho_0=\dfrac{a}{1-ah/2}$.


Then, $\Upsilon_m(y)=1/(y-x_m), A_{2m}(y)=y-x_0-mh, B_{2m}=-2m, C_{2m}(y)=a(y-x_0)-2m(1+ah), D_{2m}(y)=a((1+ah/2)(y-x_0)-2mh)/2, \rho_{2m}=a(1+ah/2)/(4m+2-ah);  \\
 A_{2m+1}(y)=y-x_0-(m+1/2)h-ah^2/4, B_{2m+1}=ah/2-2m-1, C_{2m+1}(y)=-a(y-x0)-(2m+1)(1-ah)+ah/2, D_{2m+1}(y)=a((m+1)(1-ah/2)(y-x_0-(2m+1)h)/(ah/2-2m-1), \rho_{2m+1}=-a(1-ah/2)/(4m+2-ah)$.

 Rational interpolations to the exponential function have been given by  N\"orlund \cite[\S 243]{Norlund},  by A. Iserles \cite[Thm 4]{iser1}
looking for numerical solvers of differential equations. Milne-Thomson \cite[\S 5.9]{MilneT}  uses confluent
interpolation $x_0=x_1=\dotsb =0$ and recovers Pad\'e approximation to the exponential function.

See also
\begin{equation} \label{Thexp}
  (1+\alpha)^x = 1+\cfrac{ \alpha x}
                   { 1 +\cfrac{ \alpha(1-x)}
                                  {2+\cfrac{\alpha(1+x)}
                                            {3 + \cfrac{\alpha(2-x)  }
                                                      {\ddots +\cfrac{\alpha(k\pm x)}
                                                                  { a_k +\ddots}}}}}
\end{equation}
$a_k=2$ or $2k+1$. Already in Thiele (1909)!

\medskip


\item

The Psi ($\Psi$) function.

We interpolate at $x_0\neq 0,-1,\dotsc, x_1=x_0+1, x_2=x_0+2,\dotsc$ the function $f(x)=\Psi(x)-\Psi(x_0)$,
where the Psi ($\Psi$) function has been briefly discussed in \S\ref{elllogar},
so we have $f(x_0)=0, f(x_1)=1/x_0, f(x_2)=1/x_0+1/x_1$, etc.
We interpolate a meromorphic function of poles $0, -1, -2,\dotsc$
The first steps

$f(x)= \cfrac{ (x-x_0)/x_0 }{ 1+  \cfrac{ ((x-x_1)/(2x_0+1)}{1+\ddots}}$.

One finds $\mathcal{D}f(y_n)= \dfrac{f(x_{n+1})-f(x_n)}{x_{n+1}-x_n=1}=  \dfrac{1}{x_n=y_{n-1}=y_n-1}$, so

$A_0(y)=(y-1), B_0=C_0=0, D_0= 1,\rho_0=\dfrac{D_0(y_0)}{A_0(y_0)-C_0(y_0)/2}= \dfrac{1}{y_0-1}, $

we proceed with (\ref{rm}-\ref{recD}) at $m=1, 2,\dotsc$, find exact divisions  with \\   $\Upsilon_m(y)=1/(y-y_0-m+1)$ and get

\noindent for even index $>0$, $A_{2m}(y)=(y-y_0)^2+(x_0-m+1)(y-y_0)-(m-1)x_0-m^2/2, B_{2m}=C_{2m}=-2m(x_0+3m/2-1),
  D_{2m}(y)=m(x_0+m-1)(y-y_0-2m+1)/(2x_0+3m-2),\rho_{2m}=m(x_0+m-1)/[(2x_0+3m-2)((2m+1)x_0+m(3m+1))]$,

for odd index, $A_{2m+1}(y)=(y-y_0)^2+(s-m+1/2)(y-y_0)-(m-1/2)x_0-m(m+1)/2, B_{2m+1}=-(2m+1)x_0-m(3m+1), C_{2m+1}(y)=y-y_0-(2m+1)x_0 -3m(m+1),
  D_{2m+1}(y)=(m+1)^2(x_0+m)(y-y_0-2m)/((2m+1)x_0+m(3m+1)), \\ \rho_{2m+1}=(m+1)(x_0+m)/[(2x_0+3m+1)((2m+1)x_0+m(3m+1))]$.

The example of \S \ref{Padinterp}    is made with $x_0=1$. 


Tsujimoto,  Vinet, and Zhedanov, have built in \cite{Tsujim} a biorthogonal set ${\mathcal{U}_m, \mathcal{V}_n}$ of rational
functions with respect to the Hahn's discrete measure at $0,1,\dotsc,N$. The function $\mathcal{U}_m$
is a ${}_3F_2$ expansion with poles at $\alpha,\dotsc,\alpha+m-1$;  $\mathcal{V}_n$
is another ${}_3F_2$  with poles at $\gamma,\dotsc,\gamma-n+1$, where $\gamma=  N+\alpha-\beta-2$. The numerators  $P_{2m}(x)=\dfrac{(x-\alpha)\dotsb (x-\alpha-m+1)}{(N-\beta)(N-\beta-1)\dotsb(N-\beta-2m+1)}\mathcal{U}_m(x)$  satisfy difference equations and recurrence relations    involving $P_{2m+2}(x), P_{2m}(x), P_{2m-2}(x)$, \cite[\S 6]{Tsujim},  where one recognizes the R-$II$ type recurrence relation \eqref{contrac}, with $\rho_m=$ a
polynomial in $m$ (depending on the evenness of $m$) divided by $(N-\beta-m)(N-\beta-m+1)$.

\item

The $\{\sin(nh+b)\}$ lattice.

Example with $x_n=\sin(2na+b), y_n=\sin((2n+1)a+b), a=h/2, Y_0(y)=y^2-\sin^2(a),
Y_1(y)=-2y\cos a, Y_2(y)\equiv 1, Q(y)=Y_1^2(y)-4Y_0(y)Y_2(y)=4\sin^2 a(1-y^2)$,

and we consider the Psi-like, or elliptic logarithm of \S\ref{elllogar},   satisfying $f(x_0)=0$ and
$\mathcal{D}f(y)=Y_2(y)/(y-\alpha)$.

As already seen, we must have \eqref{singmm0}, so, $A_0(y)=(y-\alpha)(y-y_{-1}),
B_0=C_0=0, D_0(y)=(y-y_{-1})$, we apply (\ref{recA}-\ref{recD}):

{\small
 $\rho_0=D_0(y_0)/A_0(y_0)=1/(y_0-\alpha),  A_1(y)=(y-\alpha)(x_0-y\cos a ),
 B_1(y)=1/\rho_0=y_0-\alpha, C_1(y)=2B_1(y)-(y-\alpha)=2y_0-\alpha-y, D_1(y)=B_1
 -(y-\alpha)=-(y-y_0)$,
 $A_1(y)+(Y_1(y)/2+x_1)C_1(y)=(y-\alpha)(x_0-y\cos a )+(x_1-y\cos a)(2y_0-\alpha-y)=2(y-y_0)(\alpha\cos a-x_1)$ (using $x_0+x_1=2y_0\cos a$),
  $\rho_1=D_1(y_1)/(A_1(y_1)-(x_2-x_1)C_1(y_1)/2)=-1/(2(\alpha\cos a-x_1))$,

    $B_2=C_2=2(\alpha\cos(a) -\sin(2a+b)),  A_2(y) +(Y_2(y)/2+x_2)C_2(y)=
      (y-y_1)(y-\alpha\sin(3a)/\sin a +2y_1)$,
$D_2(y)=B_1\rho_1X_2(x_1)(y-y_1)=(\alpha-y_0) (y-y_1)/(2(\alpha\cos a-x_1)),
\rho_2=D_2(y_2)/((y_2-y_1)(y_2-\alpha\sin(3a)/\sin a+2y_1))=
(\alpha-y_0)/(2(\alpha\cos a-x_1)(y_2-\alpha\sin(3a)/\sin a+2y_1))$
  }

  etc.

After much computer  algebra, the simple (!!) formulas are

{\small
$A_m(y)=( (\sin(h/4))^2+(-1)^m(\cos(h/4))^2 )(y^2-1)+(-1)^m( -\sin((m-1)h/2+b)-\alpha\cos(mh/2) )y   +(-1 )^m \alpha\sin((m-1)h/2+b)+(-1)^m\cos(mh/2)$,

$B_{m}=[(  \sin^2(a/2) +(-1)^m\cos^2(a/2))  \cos((2m-1)a+b)+(-1)^m( \alpha\sin (ma)  -\cos((m-1)a+b)\ )]/\sin a$,

$C_{m}(y)=((-1)^m-1)y/2+   (-1)^m(\cos( (2m-1/2-(-1)^m/2)a+b)
+\alpha\sin(ma)-\cos((m-1)a+b)\ )/\sin a$,

$D_m(y)=(     (\sin((2m'+1)a+b)-sin(-a+b) )/2               -\sin(m' a+b)+\sin(b-a)
+\alpha   (1-\cos((m'+1)a) )                   )(y-y_{m-1})/   (2B_m\sin^2 a)  $      (where $ m'=m-(1+(-1)^m)/2$).

$A_m(y)+(Y_1(y)/2+x_m)C_m(y)=(y-y_{m-1})((1+(-1)^m)y/2  -\sin((2m+1/2)a+b)/(2\cos(a/2))+(-1)^m[-\alpha\sin((m+1)a)
-\cos((2m+1/2)a+b)\cos(a/2)+\cos(ma+b)]/\sin a$,
$\rho_m=
D'_m/B_{m+1} $.
}   

\item

Trying a full elliptic lattice.

\subsubsection{About classical orthogonal and biorthogonal functions on a true elliptic lattice .}The simplest instances of families of classical special functions have explicit formulas for recurrence relations, differential and difference equations etc. So, the Gauss continued fraction \eqref{Gaussf}
has parameters $\rho_m=$ a polynomial in $m$ (depending on the evenness of $m$, so actually a polynomial in $m$ and $(-1)^m$) over $(\gamma+m-1)(\gamma+m)$.
Other simple examples have been given in \S \ref{Padinterp},  \S\ref{Euler1},  \S  \ref{discorthex}.  
    Koekoek et al  \cite[\S 14.1]{Koe}  present the most elaborate recurrence relations of orthogonal  polynomials  in terms of parameters $A_m, C_m$ easily recognized in most cases as our $\rho_{2m+1},\rho_{2m}$ of  \eqref{contrac}, and have $\rho_m=$ a polynomial of degree 4 in $q^{m/2}$ over $(1-abcdq^{m-2})(1-abcdq^{m-1})$, as seen in \ref{AWabcd} of \S  \ref{discorthex}.   

For biorthogonal rational functions,
Rahman \cite{Rahman81} has  $\rho_m=am+b+c(-1)^m$,      Masson \cite{Masson} has $\rho_m=$ a polynomial of degree 8 (amazing!) in $q^{m/2}$ over $q^m(1-sq^{m-2})(1-sq^{m-1})$.

On fully elliptic grids,
Spiridonov and  Zhedanov   gave in  \cite{SpZ1,SpZ2,SpZComm2000} $\rho_m$ as
a product of 7 factors $[m+a_j+b_j(-1)^m]$ (theta functions) over $[2m+s][2m+1+s][m+c+d(-1)^m]$ \cite[pp. 384-385]{SpZ1}. See \cite{Spi2024} for further developments.



\medskip


The degrees of $A_m,\dotsc, D_m$ must be at least $d=3$  on  a genera1 elliptic grid,
as $Y_2$, normally of degree 2, is a factor of $B_m$ etc., and $y-y_{m-1}$ is a factor of $D_m$ from \eqref{singmm1}. Even if $Y_2$ has degree $<2$, the polynomial $Q$, needed in \eqref{recA}, has  degree $\geqslant 3$  in the genus $>0$ case, from \S \ref{fixedp}, as $Q$ must have 4 simple zeros  (a  polynomial $Q$ of degree $<3$ has a multiple zero at $\infty$).

\medskip

\subsubsection{Example} with $x_n=\sn(2na+b), y_n=\sn((2n+1)a+b), Y_0(y)=y^2-\sn^2(a),
Y_1(y)=-2\cn(a)\dn(a)y, Y_2(y)=1-k^2\sn^2(a)y^2, Q(y)=Y_1^2(y)-4Y_0(y)Y_2(y)=4\sn^2(a)(1-y^2)(1-k^2y^2)$,

and we consider the Psi-like, or elliptic logarithm of \S\ref{elllogar},   satisfying $f(x_0)=0$ and
$\mathcal{D}f(y)=Y_2(y)/(y-\alpha)$.

As already seen, we must have \eqref{singmm0}, so, $A_0(y)=(y-\alpha)(y-y_{-1}),
B_0=C_0=0, D_0(y)=(y-y_{-1})Y_2(y)$, we apply (\ref{rm}-\ref{recD}):


 $\rho_0=D_0(y_0)/A_0(y_0)=Y_2(y_0)/(y_0-\alpha),  A_1(y)=(y-\alpha)(Y_1(y)/2+x_0Y_2(y))=(y-\alpha)(x_0-\cn(a)\dn(a)y-x_0k^2\sn^2(a)y^2),
 B_1(y)=Y_2^2(y)/\rho_0, C_1(y)=2B_1(y)-(y-\alpha)Y_2(y), D_1(y)=C_1(y)-B_1(y)
 =((1-k^2\sn^2 ay^2)/\rho_0-y+\alpha)Y_2(y)=[  (1-k^2\sn^2 ay^2)(y_0-\alpha)-(1-k^2\sn^2 ay_0^2)(y-\alpha)]Y_2(y)/Y_2(y_0)=(-1  -k^2\sn^2 ay_0y+\alpha( y+y_0))(y-y_0)Y_2(y)/Y_2(y_0)$

$A_m$,etc. have degree 4 when $m$ is odd $\geqslant 3$, degree 3 when $m$ is even $\geqslant 2$.

Unfortunately, no explicit formulas  have been found here.

Polynomials of degree 3 are completely determined by their values at the four zeros
$\mathfrak{y}_1,\dotsc, \mathfrak{y}_4$ of $Q$, so, from \eqref{recA},
{$\displaystyle A_m(\mathfrak{y}_j)=A_0(\mathfrak{y}_j)\prod_{r=0}^{m-1}
   \Upsilon_r \rho_r (x_r-\mathfrak{x}_j)$, where $\mathfrak{x}_j$ is the (double)
   $x-$root of $F(x,\mathfrak{y}_j)=0$, using
 $Y_1(\mathfrak{y}_j)=-2\mathfrak{x}_j Y_2(\mathfrak{y}_j) $}, $j=1,\dotsc, 4$.

   Such products are related to $\chi_m(y)$ functions introduced in next section.

\end{enumerate}

\subsection{Difference relations and equations\label{diffreleq}}   \   \\

\subsubsection{\label{diff24}}\textbf{Theorem.} Let $\{(x_n,y_n)\}$, $\{(x_{n+1},y_n)\}$ be points
on the biquadratic curve \\ $F(x,y)=Y_0(y)+xY_1(y)+x^2Y_2(y)=0$, and $f$ a function
defined on the $x_n$s, satisfying the \emph{discrete Riccati} difference equation
\eqref{Riccati0}
 $ A(y_n)(\mathcal{D}f)(y_n)=B(y_n)f(x_n)f(x_{n+1})
          +C(y_n)(\mathcal{M}f)(y_n)+D(y_n)$, with $f(x_0)=0$, and we consider these rational interpolants to $f$

$\dfrac{ N_{m}(x)     }{P_{m}(x)}= \dfrac{r_0(x-x_0) }{1+\dfrac{\ddots}{1+r_{m-1}(x-x_{m-1}) }}$, for $m=0,1,\dotsc$


Then, $N_m$ and $P_m$ satisfy linear difference equations of order $\leqslant 4$ with coefficients depending on $A_m,\dotsc, D_m$ defined in \S \ref{Riccatim},

$P_m$ satisfies a linear difference equation of order $\leqslant 2$ if $B(y)\equiv 0$.

\textbf{Proof.}

We recall 
$U_m(y_n)=\dfrac{C_{m}(y_n)}{2}+ \dfrac{A_{m}(y_n)}{x_{n+1}-x_n}, V_m(y_n)=\dfrac{C_{m}(y_n)}{2}- \dfrac{A_{m}(y_n)}{x_{n+1}-x_n}$,
and rewrite  \eqref{recA}-\eqref{recD} and \eqref{recUV}  as

{\small
\noindent$\begin{bmatrix}B_{m+1}(y_n)\\ U_{m+1}(y_n)\\ V_{m+1}(y_n)\\  D_{m+1}(y_n)\end{bmatrix}=
   \dfrac{\Upsilon_m(y_n)}{\rho_m}\begin{bmatrix}0 &0 &0 &1 \\
  0 & 0 & \rho_m(x_{n+1}-x_m)  & 1   \\ 0 & \rho_m(x_n-x_m) & 0 &1  \\
 \rho_m^2(x_n-x_m)(x_{n+1}-x_m) & \rho_m(x_n-x_m) & \rho_m(x_{n+1}-x_m) & 1 \end{bmatrix}
\begin{bmatrix}B_{m}(y_n)\\ U_{m}(y_n)\\V_{m}(y_n)\\  D_{m}(y_n)\end{bmatrix}$
}

which is exactly  the recurrence relation of products $u(x_{n+1})v(x_n)$, where $u_{m+1}(x)=u_m(x)+
\rho_m(x-x_m)u_{m-1}(x)$ (resp. $v$),the recurrence relations of $N_m$ and $P_m$,
for example  $u_{m+1}(x_{n+1})v_m(x_n)=u_{m}(x_{n+1})v_m(x_n)+\rho_m(x_{n+1}-x_m)u_{m-1}(x_{n+1})v_m(x_n)$ etc. and we look at various instances of

$u_{m-1}(x_{n+1})v_{m-1}(x_n), u_{m}(x_{n+1})v_{m-1}(x_n), u_{m-1}(x_{n+1}v_{m}(x_n),u_{m}(x_{n+1})v_{m}(x_n)$ at $m=0$:

$\begin{matrix}  N_{-1}N_{-1}=1 &  N_{0}N_{-1}=0&N_{-1}N_{0}=0&N_{0}N_{0}=0 \\
               P_{-1}N_{-1}=0 &  P_{0}N_{-1}=1&P_{-1}N_{0}=0&P_{0}N_{0}=0 \\
               N_{-1}P_{-1}=0 &  N_{0}P_{-1}=0&N_{-1}P_{0}=1&N_{0}P_{0}=0 \\
              P_{-1}P_{-1}=0 &  P_{0}P_{-1}=0&P_{-1}P_{0}=0&P_{0}P_{0}=1 \end{matrix}$,

and

{\small \begin{equation}\label{ABCDPN}\begin{split}
  &\begin{bmatrix}B_{m}(y_n)\\ U_{m}(y_n)\\V_{m}(y_n)\\  D_{m}(y_n)\end{bmatrix}=
  \dfrac{ \Upsilon_0(y_n)\dotsb \Upsilon_{m-1}(y_n)  }{\rho_0\dotsb\rho_{m-1}}
\\
 & \times   \begin{bmatrix}
      N_{m-1}(x_{n+1})N_{m-1}(x_n)&  P_{m-1}(x_{n+1})N_{m-1}(x_n)& N_{m-1}(x_{n+1})P_{m-1}(x_n)& P_{m-1}(x_{n+1})P_{m-1}(x_n)\\
    N_{m}(x_{n+1})N_{m-1}(x_n)&  P_{m}(x_{n+1})N_{m-1}(x_n)& N_{m}(x_{n+1})P_{m-1}(x_n)& P_{m}(x_{n+1})P_{m-1}(x_n)\\
 N_{m-1}(x_{n+1})N_{m}(x_n)&  P_{m-1}(x_{n+1})N_{m}(x_n)& N_{m-1}(x_{n+1})P_{m}(x_n)& P_{m-1}(x_{n+1})P_{m}(x_n)\\
 N_{m}(x_{n+1})N_{m}(x_n)&  P_{m}(x_{n+1})N_{m}(x_n)& N_{m}(x_{n+1})P_{m}(x_n)& P_{m}(x_{n+1})P_{m}(x_n)
   \end{bmatrix}\begin{bmatrix}B_{}(y_n)\\ U_{}(y_n)\\V_{}(y_n)\\  D_{}(y_n)\end{bmatrix}
\end{split}\end{equation}
}

We recognize that $D_m$ is basically the Laguerre's $\Theta_m$ of \eqref{LagTheta}.
These identities are not meant to give $B_m$ etc.,which we already have
in \eqref{recA}-\eqref{recD},  in terms of $P_m, N_m$,etc. but exactly
the opposite is true, we will get difference relations and equations for $P_m,
N_m$ precisely. But before that, another way to have \eqref{ABCDPN}:

\begin{equation}    \begin{split}
   0&= A(y_n)(\mathcal{D}f)(y_n)-B(y_n)f(x_n)f(x_{n+1})
          -C(y_n)(\mathcal{M}f)(y_n)-D(y_n)\\
  &= \left[\dfrac{A_{}(y_n)}{x_{n+1}-x_n}-\dfrac{C_{}(y_n)}{2}\right]
\dfrac{ N_{m}(x_{n+1})   +N_{m-1}(x_{n+1})f_{m}(x_{n+1})  }{P_{m}(x_{n+1})+P_{m-1}(x_{n+1})f_{m}(x_{n+1})} \\
&\ \  -\left[\dfrac{A_{}(y_n)}{x_{n+1}-x_n}+\dfrac{C_{}(y_n)}{2}\right]
\dfrac{ N_{m}(x_{n})   +N_{m-1}(x_{n})f_{m}(x_{n})     }{P_{m}(x_{n})+P_{m-1}(x_{n})f_{m}(x_{n})}\\
 -B(y_n)&\dfrac{ N_{m}(x_{n})   +N_{m-1}(x_{n})f_{m}(x_{n})     }{P_{m}(x_{n})+P_{m-1}(x_{n})f_{m}(x_{n})} \dfrac{ N_{m}(x_{n+1})   +N_{m-1}(x_{n+1})f_{m}(x_{n+1})     }{P_{m}(x_{n+1})+P_{m-1}(x_{n+1})f_{m}(x_{n+1})} -D(y_n)
\end{split}\end{equation}

Multiply the denominators,

{\small
\noindent$0=-V(y_n)[N_{m-1}(x_{n+1})P_m(x_n)   f_{m}(x_{n+1})+N_{m}(x_{n+1})P_{m-1}(x_{n})f_{m}(x_{n})
+N_{m-1}(x_{n+1})P_{m-1}(x_{n})f_{m}(x_{n})f_{m}(x_{n+1})      \\
+N_{m}(x_{n+1})P_m(x_n)  ] \\
-U(y_n)[N_{m}(x_{n})P_{m-1}(x_{n+1})f_{m}(x_{n+1})
+N_{m-1}(x_{n})P_m(x_{n+1}) f_{m}(x_{n})
 + N_{m-1}(x_{n})P_{m-1}(x_{n+1})f_{m}(x_{n})f_{m}(x_{n+1})
 +  N_{m}(x_{n})P_m(x_{n+1})]\\
-B(y_n)[     N_{m-1}(x_{n+1})N_{m}(x_{n})f_{m}(x_{n+1})
+    N_{m-1}(x_{n})N_{m}(x_{n+1})f_{m}(x_{n})
+    N_{m-1}(x_{n+1})N_{m-1}(x_{n})f_{m}(x_{n})f_{m}(x_{n+1})
+    N_{m}(x_{n+1})N_{m}(x_{n})] \\
-D(y_n)[     P_{m-1}(x_{n+1})P_{m}(x_{n})f_{m}(x_{n+1})
+    P_{m-1}(x_{n})P_{m}(x_{n+1})f_{m}(x_{n})
+    P_{m-1}(x_{n+1})P_{m-1}(x_{n})f_{m}(x_{n})f_{m}(x_{n+1})
+    P_{m}(x_{n+1})P_{m}(x_{n})]
$
}

must be $\left(-V_m(y_n)=\dfrac{A_m(y_n)}{x_{n+1}-x_n}-\dfrac{C_m(y_n)}{2} \right)    f_m(x_{n+1})-     \left(U_m(y_n)=\dfrac{A_m(y_n)}{x_{n+1}-x_n}+\dfrac{C_m(y_n)}{2}\right)f_m(x_n)  \\
              -B_m(y_n)f_m(x_n)  f_m(x_{n+1})-        D_m(y_n)=0$,

\medskip
and one recovers (the second part of) \eqref{ABCDPN}.

\medskip


In \eqref{ABCDPN}, take $B_m(y_n)P_m(x_{n+1}-U_m(y_n)P_{m-1}(x_{n+1})$: we get
$\dfrac{\Upsilon_0(y_n)\dotsb \Upsilon_{m-1}(y_n) }{\rho_0\dotsb\rho_{m-1}}$

$\times \{ (N_{m-1}(x_{n+1})P_m(x_{n+1})-N_m(x_{n+1})P_{m-1}(_{n+1}))N_{m-1}(x_n)B(y_n)
  \\  +(N_{m-1}(x_{n+1})P_m(x_{n+1})-N_m(x_{n+1})P_{m-1}(x_{n+1}))P_{m-1}(x_n)V(y_n)\}$
remains, where one recognizes the Casorati product \eqref{Caso}

\begin{equation}\label{struct}\begin{split}
&\begin{bmatrix}B_{m}(y_n)P_m(x_{n+1})-U_m(y_n)P_{m-1}(x_{n+1})
\\ B_{m}(y_n)P_m(x_{n})-V_m(y_n)P_{m-1}(x_{n})
\\U_{m}(y_n)N_m(x_n)-D_m(y_n)N_{m-1}(x_n)\\ V_m(y_n)N_{m}(x_{n+1})- D_{m}(y_n)N_{m-1}(x_{n+1})\end{bmatrix}=
  (-1)^m \Upsilon_0(y_n)\dotsb \Upsilon_{m-1}(y_n)   \\
\times &
\begin{bmatrix} (x_{n+1}-x_0)\dotsb(x_{n+1}-x_{m-1})\{B(y_n)N_{m-1}(x_n)+V(y_n)P_{m-1}(x_n)\}  \\
 (x_{n}-x_0)\dotsb(x_{n}-x_{m-1})\{B(y_n)N_{m-1}(x_{n+1})+U(y_n)P_{m-1}(x_{n+1})\}  \\
 (x_{n}-x_0)\dotsb(x_{n}-x_{m-1})\{-V(y_n)N_{m}(x_{n+1})-D(y_n)P_{m}(x_{n+1})\}  \\
(x_{n+1}-x_0)\dotsb(x_{n+1}-x_{m-1})\{-U(y_n)N_{m}(x_n)-D(y_n)P_{m}(x_n)\}
\end{bmatrix}
\end{split}\end{equation}

Such relations between $\phi_m$, differences of $\phi_m$ and values of $\phi_{m\pm 1}$ are sometimes called differential or difference  relations,
sometimes raising and lowering operators relations,  sometimes  structure relations \cite{kenfack,Koornstru1,Koornstru2,Marstruc,Mboustru}.


Let $\chi_m(y_n)=\Upsilon_0(y_n)\dotsb \Upsilon_{m-1}(y_n)  (x_{n}-x_0)\dotsb(x_{n}-x_{m-1})$

$= \dfrac{  (Y_2(y_n))^{m-d_m}  (x_{n}-x_0)\dotsb(x_{n}-x_{m-1})   }
         {   (y_{n}-y_{-1})\dotsb(y_{n}-y_{m-2})   }$,  where $d_m=$ the integer part of $m/2$, and

$\chi^+_m(y_{n})=\Upsilon_0(y_n)\dotsb \Upsilon_{m-1}(y_n)   (x_{n+1}-x_0)\dotsb(x_{n+1}-x_{m-1})$

$= \dfrac{  (Y_2(y_n))^{m-d_m}  (x_{n+1}-x_0)\dotsb(x_{n+1}-x_{m-1})   }
         {   (y_{n}-y_{-1})\dotsb(y_{n}-y_{m-2})   }$.

\eqref{struct}  is a first order vector difference equation relating
$[P_{m-1},P_m,N_{m-1},N_m]$ at $x_n$ and $x_{n+1}$:

\begin{equation}\label{eqdif4}\begin{split}
&\begin{bmatrix}-U_m(y_n) & B_m(y_n) & 0 & 0 \\
 -\chi_m(x_n) U_0(y_n) &0&-\chi_m(x_n)  B(y_n)&0 \\
0&\chi_m(x_n)D(y_n) & 0 & \chi_m(x_n)V(y_n)  \\
0 & 0 & -D_m(y_n) & V_m(y_n) \end{bmatrix}
\begin{bmatrix} P_{m-1}(x_{n+1})\\P_m(x_{n+1})\\N_{m-1}(x_{n+1})\\N_m(x_{n+1})\end{bmatrix}            \\
=&\begin{bmatrix} \chi^+_m(x_{n+1})V(y_n)&0&  \chi^+_m(x_{n+1})B(y_n)&0\\
V_m(y_n)&-B_m(y_n) & 0 & 0\\
0 & 0 & D_m(y_n)&-U_m(y_n)\\
0&-\chi^+_m(x_{n+1})D(y_n)&0&-\chi^+_m(x_{n+1})U_0(y_n)
\end{bmatrix}
\begin{bmatrix} P_{m-1}(x_{n})\\P_m(x_{n})\\N_{m-1}(x_{n})\\N_m(x_{n})\end{bmatrix}
\end{split}\end{equation}

which can be condensed as a \emph{scalar} fourth-order difference
equation for, say, $P_m$
\cite{BangRama,Ronv4}
whose solutions are products of solutions of second order differential or difference eq:
\cite[p.446]{IsWimp}, \cite[\S 4]{RonvFactor}


\medskip
Finally, when $B(y)\equiv 0$, %

$-U_m(y_n)P_{m-1}(x_{n+1})+ B_m(y_n) P_m(x_{n+1})= \chi^+_m(x_{n+1})V_0(y_n)P_{m-1}(x_{n})$,

$ -\chi_m(x_n) U_0(y_n) P_{m-1}(x_{n+1})=V_m(y_n)P_{m-1}(x_{n})-B_m(y_n) P_m(x_{n})$

\medskip

More on the system of difference equations of first order, in $\mathcal{D}-\mathcal{M}$ parlance, 

$-U_m(\mathcal{M}P_{m-1}+\Delta x_n \mathcal{D}P_{m-1}/2)+ B_m (\mathcal{M}P_{m}+\Delta x_n \mathcal{D}P_{m}/2)= \chi^+_m V_0(\mathcal{M}P_{m-1}-\Delta x_n \mathcal{D}P_{m-1}/2)$,

$ -\chi_mU_0 (\mathcal{M}P_{m-1}+\Delta x_n \mathcal{D}P_{m-1}/2)=V_m(\mathcal{M}P_{m-1}-\Delta x_n \mathcal{D}P_{m-1}/2)-B_m(\mathcal{M}P_{m}-\Delta x_n \mathcal{D}P_{m}/2) $

$ \Delta x_n\begin{bmatrix}-U_m+ \chi^+_m V& B_m \\ -\chi_mU_0 +V_m& -B_m \end{bmatrix}\begin{bmatrix} \mathcal{D}P_{m-1} \\  \mathcal{D}P_{m}\end{bmatrix}=2\begin{bmatrix}U_m+\chi^+_m V & -B_m \\ \chi_mU +V_m& -B_m \end{bmatrix}\begin{bmatrix} \mathcal{M}P_{m-1} \\  \mathcal{M}P_{m}\end{bmatrix}
$

or, isolating  $\mathcal{D}P_{m-1}$ and $\mathcal{D}P_{m}$, 
 using
$U_m=C_{m}/2+ A_{m}/\Delta x_n, V_m=C_{m}/2- A_{m}/\Delta x_n$,

{\small \begin{equation}\label{equa2}
(A_m+\Delta x_n(\chi_mU_0-\chi^+_m V_0)/2)\begin{bmatrix} \mathcal{D}P_{m-1} \\  \mathcal{D}P_{m}\end{bmatrix}=\begin{bmatrix}-C_m-\chi_mU_0+\chi^+_m V_0 & 2B_m \\  -2D_{m}    & C_m-\chi_mU_0-\chi^+_m V_0\end{bmatrix}\begin{bmatrix} \mathcal{M}P_{m-1} \\  \mathcal{M}P_{m}\end{bmatrix}.
\end{equation}   }

to be compared with (\ref{diffeqp1}-\ref{diffeqp2}) .  However, the calculations give
$2(\chi_m\chi^+_m  U_0V_0-U_mV_m)/B_m$
 at the
low-left matrix element, indeed,
\eqref{detABCD} is $U_mV_m-B_mD_m=\dfrac{\Upsilon^2_0\dotsb \Upsilon^2_{m-1}  F(x_0,y)\dotsb F(x_{m-1},y)}{Y_2^m}U_0V_0$, as $B_0=0$, and
$\chi_m\chi^+_m=\Upsilon_0^2\dotsb \Upsilon^2_{m-1}  (x_{n}-x_0) (x_{n+1}-x_0)\dotsb(x_{n}-x_{m-1})(x_{n+1}-x_{m-1})$.

\smallskip
Scalar \emph{second} order difference equation:
 put $m\to m+1, n\to n-1$, eliminate $P_{m+1}(x_{n})$ from

$P_{m+1}(x_{n}) =\dfrac{U_{m+1}(y_{n-1})P_{m}(x_{n})+ \chi^+_{m+1}(y_{n})V_0(y_{n-1})P_{m}(x_{n-1})                                                  }{B_{m+1}(y_{n-1}) }
     $

$=\dfrac{ \chi_{m+1}(y_{n}) U_0(y_{n}) P_{m}(x_{n+1})+V_{m+1}(y_{n})P_{m}(x_{n})}{B_{m+1}(y_{n})  }$

to have

$ \chi_{m+1}(y_{n})  U_0(y_{n}) B_{m+1}(y_{n-1}) P_{m}(x_{n+1})
     +   [  V_{m+1}(y_{n}) B_{m+1}(y_{n-1})
- U_{m+1}(y_{n-1})B_{m+1}(y_{n})  ]P_m(x_n)
    \\   - \chi^+_{m+1}(y_{n})V_0(y_{n-1}) B_{m+1}(y_{n})   P_{m}(x_{n-1})=0
$,

which is a scalar linear second order difference equation.   \qed

See \S \ref{secondor} for more on second order difference equations.





$\chi_m$ and $\chi_m^+(y_n)$, as well as $U_m$ and $V_m$,  are conjugate algebraic functions. For instance,  at $m=1$:

$\chi_1, \chi_1^+(y_n)=(Y_2(y_n))^{d_1} \dfrac{-Y_1(y_n)-2x_0Y_2(y_n)\mp\sqrt{Q(y_n)}}{2  (y_{n}-y_{-1})Y_2(y_n)}$.

$\chi_m(y_n)+\chi_m^+(y_n)$ is a rational function of $y_n$ with a numerator vanishing at
$n=0,\dotsc,m-2$, and poles when $x=\infty$, so, at zeros of $Y_2$, what remains has the form

$\chi_m(y)+\chi_m^+(y)=(Y_2(y))^{d_m-m} \omega_m (y-\eta_m)/(y-y_{-1})$,

the product is simpler!

  $\chi_m(y_n)\chi_m^+(y_n)=(Y_2(y_n))^{2d_m} {\displaystyle
\prod_{j=0}^{m-1}\dfrac{   (x_{n}-x_j) (x_{n+1}-x_j)  = X_2(x_j)(y_n-y_{j-1})(y_n-y_j)/Y_2(y_n) }
         {   (y_{n}-y_{j-1})^2   }   }\\
             = (Y_2(y))^{2d_m-2m}X_2(x_0)\dotsb X_2(x_{m-1})
            (y-y_{m-1})/(y-y_{-1})   $.

Ratios $ \dfrac{   (x-x_0)\dotsb(x-x_{m-1})   }         {   (y-y_{-1})\dotsb(y-y_{m-2})   }$ entering $\chi_m(y)$, where $F(x,y)=0$,and similar ones for $\chi^+_m(y)$ simplify nicely for simple grids, where $(x-x_j)/(y-y_{j-1})$ appears as a function of $j$ over the same function of $j-1$:

{\small
For instance, with $F(x,y)=(y-x-a)(y-x-a-h)\Rightarrow x_n=x_0+nh, y_n=a+x_{n+1},\\ \chi_m(y)= (x-x_0)(x-x_0-h)\dotsb(x-x_0-(m-1)h)/((y-y_{-1}=a+x+h-a-x_0)\dotsb(x-x_0-(m-2)h)) \\=(x-x_0-(m-1)h)/(x-x_0+h)=(y-a-x_0-mh)/(y-a-x_0),\\
\chi_m^+(y)= (x+h-x_0)(x+h-x_1)\dotsb(x+h-x_{m-1})/((y-y_{-1}=a+x+h-a-x_0)\dotsb(x+h-x_{m-1})=1$.

With $F(x,y)=(y-x)(y-qx), x_n=q^n, y_n=q^{n+1}, (x-x_j)/(y-y_{j-1})= (y/q-q^j)/(y-q^{j})$, and the products $\chi_m(y)=q^{-m}(y-q^{m})/(y-1),   \chi_m^+(y)=1$.


More difficult:  $x_n=\sin(nh+\theta), y_n=\sin((n+1/2)h+\theta), {\text{\ let\ }} x=\sin s, y=\sin(s+h/2)$ gives $(x-x_j)/(y-y_{j-1})=[ \sin((s-jh-\theta)/2)\cos((s+jh+\theta)/2)]/
[\sin((s-(j-1)h-\theta)/2)\cos((s+jh+\theta)/2)]$ and  \\
$\chi_m(y)=\sin((s-(m-1)h-\theta)/2)/\sin((s+h-\theta)/2)$ remains. Also,
$\chi_m^+(y)=\cos((s+mh+\theta)/2)/\cos((s+\theta)/2),  \\
\chi_m(y)+\chi_m^+(y)=2\dfrac{ \sin((s-(m-1)h-\theta)/2) \cos((s+\theta)/2)+\sin((s+h-\theta)/2)  \cos((s+mh+\theta)/2)}{2\sin((s+h-\theta)/2) \cos((s+\theta)/2)=\sin(s+h/2)- \sin(\theta-h/2)}\\
=\dfrac{ \sin(s-(m-1)h/2) -\sin(\theta+(m-1)h/2)+\sin(s+(m+1)h/2) - \sin(\theta+(m-1)h/2)}{y-y_{-1}}
\\=2\dfrac{ y \cos(mh/2) - \sin(\theta+(m-1)h/2)}{y-y_{-1}}
$.



A true elliptic case
$x_n=\sn(nh+\theta), y_n=\sn((n+1/2)h+\theta), x=\sn s, y=\sn(s+h/2)$ gives $(x-x_j)/(y-y_{j-1})=
\dfrac{ \sn((s-jh-\theta)/2)\cn\dn((s+jh+\theta)/2)}
{1-k^2 \sn^2((s-jh-\theta)/2)\sn^2((s+jh+\theta)/2)}
\dfrac{1-k^2\sn^2((s-(j-1)h-\theta)/2)\sn^2((s+jh+\theta)/2) }
{\sn((s-(j-1)h-\theta)/2)\cn\dn((s+jh+\theta)/2)}$ and the product
$\dfrac{\sn((s-(m-1)h-\theta)/2)}
{\sn((s+h-\theta)/2)  }$ times a product  of $1-k^2\sn^2\sn^2$ factors. No obvious simplification happens.

}  

\section{The adjoint operator.\label{adjointop}}\chead{\thesection  \ \ \ \   Adjoint op}

\subsection{Definitions\label{defadj}} \   \\

Let $g$ be defined on at least  a sequence $\{y_n\}$, then, the
adjoint operator to $\mathcal{D}$ is defined as

\begin{equation}\label{adj}(\mathcal{D}^\dagger g)(x_n)= -\dfrac{g(y_n)-g(y_{n-1})}{y_n-y_{n-1}}\end{equation}
at $x=x_n$.

One has

 $ \sum_0^{N-1} f(x_n)(\mathcal{D}^\dagger g)(x_n) (y_n-y_{n-1})
=\sum_0^{N-1} (\mathcal{D}f)(y_n)\  g(y_n)(x_{n+1}-x_n)   -f(x_0)g(y_{-1}) +f(x_N)g(y_{N-1})$.

Indeed,  the left-handed side sum is rearranged as

$\sum_0^{N-1} (f(x_{n+1}-f(x_n)) g(y_n)\ \ \   f(x_0)g(y_{-1})-f(x_N)g(y_{N-1})\\
=\sum_0^{N-1} (\mathcal{D}f)(y_n)\  g(y_n)(x_{n+1}-x_n)\ \ \  f(x_0)g(y_{-1}) -f(x_N)g(y_{N-1})$,

showing a discrete integration by parts, also called Abel's summation rule, as noted in \eqref{Abelsum} in \S \ref{diffsum}.

The difference elements $y_n-y_{n-1}$ at $x_n$, and $x_{n+1}-x_n$ at $y_n$
are not as crazy as they seem:

$d_qx= c(q^n-q^{n+1})$ in Jackson \cite[\S 1.15 ]{Koe}.

Milne-Thomson \cite[\S 12.6]{MilneT} and N\"orlund \cite[p. 288]{Norlund} give
the adjoint equation of $\alpha_n v_{n+d}+\beta_n v_{n+d-1}+\dotsb +\omega_n v_n=0$ to be
$\alpha_n u_n+\beta_{n+1}u_{n+1}+\dotsb +\omega_{n+d}u_{n+d}=0$ (equation for
a multiplier). Rearranging $\sum u_n(\alpha_n v_{n+d}+\dotsb +\omega_n v_n)$
yields $\sum v_n(\alpha_{n-d} u_{n-d}+\dotsb +\omega_n u_n)$, showing the
 adjoint operator. In particular,
the adjoint of $\Delta$ is $-\nabla$.

\medskip

We now come to a simple fraction case,
\begin{equation} \label{simpleratd}
\left(\mathcal{D}^\dagger \dfrac{1}{y-A}\right)(x)
=\dfrac{X_2(x)}{Y_2(A)(x-x_\kappa)(x-x_{\kappa+1})}
\end{equation}

where $\kappa$ is such that $A=y_\kappa$.

Indeed,

  $\dfrac{1/(y_n-A)-1/(y_{n-1}-A)}{y_n-y_{n-1}}=
\dfrac{-1}{(y_n-A)(y_{n-1}-A)}=\dfrac{-X_2(x_n)}{F(x_n,A=y_\kappa)=Y_2(A)(x-x_\kappa)(x-x_{\kappa+1})}$.

\medskip

Naturally enough, we define

\begin{equation}\label{adjM}(\mathcal{M}^\dagger g)(x_n)= \dfrac{g(y_n)+g(y_{n-1})}{2}
\ \ \ \text{at\  } x=x_n.\end{equation}

Note that the operator $\mathcal{M}^\dagger$ just defined in \eqref{adjM} is NOT
the adjoint of $\mathcal{M}$ with respect to the bilinear form associated to
$\mathcal{D}$ and $\mathcal{D}^\dagger$: if $\mathscr{M}$ is the adjoint to $\mathcal{M}$
with respect to the same bilinear form, one must have
$\langle f(x), (\mathscr{M}g)(x) \rangle = \sum_n (y_n-y_{n-1})f(x_n)(\mathscr{M}g)(x_n)
=\langle (\mathcal{M}f)(y),g(y)\rangle =
\sum_n (x_{n+1}-x_n)(f(x_n)+f(x_{n+1})g(y_n)/2$ rearranged as
$\sum  (y_n-y_{n-1})f(x_n) \dfrac{ (x_{n+1}-x_n)g(y_n)+(x_n-x_{n-1})g(y_{n-1}) }{2(y_n-y_{n-1})}$,

so, the "true" adjoint of $\mathcal{M}$ is $(\mathscr{M}g)(x_n)=\dfrac{ (x_{n+1}-x_n)g(y_n)+(x_n-x_{n-1})g(y_{n-1}) }{2(y_n-y_{n-1})}$, will not be used here.

\medskip

\subsection{Adjoint operators applied to products} \   \\

We now transpose \S \ref{Dprod}

\textbf{Theorem.\label{Daprod}}
\textsl{Let $\{(x_n,y_n)\}, \{(x'_n,y'_n)\}$ be two elliptic lattices on the same
$F-$curve, then}
\begin{equation} \label{daY}\begin{split}
 \mathcal{D}^\dagger\dfrac{(y-y_r)\dotsb(y-y_{r+m-1}) }{(y-y'_s)\dotsb (y-y'_{s+m-1})}
&=
     C^\dagger_{m,r,s} X_2(x) \dfrac{(x-x_{r+1})\dotsb(x-x_{r+m-1}) }{(x-x'_{s})\dotsb (x-x'_{s+m})}, \\
 \mathcal{M}^\dagger\dfrac{(y-y_r)\dotsb(y-y_{r+m-1}) }{(y-y'_s)\dotsb (y-y'_{s+m-1})}
&=
     D_{m,r,s}^\dagger(x)  \dfrac{(x-x_{r+1})\dotsb(x-x_{r+m-1}) }{(x-x'_{s})\dotsb (x-x'_{s+m})},
\end{split}\end{equation}
\textsl{where $C_{m,r,s}^\dagger$ is a constant, and $D_{m,r,s}^\dagger$ is a second degree polynomial.}

Indeed, let $Q_m(y)=\dfrac{(y-y_r)\dotsb(y-y_{r+m-1}) }{(y-y'_s)\dotsb (y-y'_{k+m-1})}$.
the operators $\mathcal{D}^\dagger$ and $\mathcal{M}^\dagger$ applied to each simple fraction
$1/(y-y'_j)$ yield fractions divided by  $(x-x'_{j})(x-x'_{j+1})$ by
\eqref{simpleratd},
the results $\mathcal{D}^\dagger Q_m$ and $\mathcal{M}^\dagger Q_m$ are therefore   rational functions of degree $m+1$ and poles
$x'_{s},\dotsc, x'_{s+m}$.

We now look at some $x=x_n$ involving the values of $Q_m$ at $y_{n-1}$
and $y_{n}$. We see that
$y_{n-1}-y_r$ and $y_{n}-y_{r+1}$   both vanish at $x=x_{r+1}$, etc.
up to  $y_{n-1}-y_{r+m-2}$ and $y_{n}-y_{r+m-1}$ both vanishing
at $x=x_{r+m-1}$,
  so that $Q_m(y_{n-1})$ and $Q_m(y_{n})$ both vanish at
 $y_{n-1}=y_r,\dotsc, y_{r+m-2}$ and
 $y_{n}=y_{r+1},\dotsc, y_{r+m-1}$, corresponding to $x=x_{r+1},\dotsc, x_{r+m-1}$ .

The first part of \eqref{daY} then follows from \eqref{simpleratd}.

For $\mathcal{M}^\dagger$, we still have a rational function of degree
$m+1$ with the same $m-1$ zeros, whence the factor $D_{m,r,s}^\dagger(x)$.

\smallskip

\noindent Examples $\mathcal{D}^\dagger\dfrac{y-y_r }{y-y'_s}  =\mathcal{D}^\dagger\dfrac{y'_s-y_r }{y-y'_s} $ at $x_n$

$=  \dfrac{y'_s-y_r}{(y_{n-1}-y'_s)(y_{n}-y'_s)}= \dfrac{(y'_s-y_r)X_2(x_n)}
{F(x_n,y'_s)=Y_2(y'_s)(x_n-x'_{s})(x_n-x'_{s+1})}$, so,

\begin{equation}\label{C1a}C_{1,r,s}^\dagger=\dfrac{y'_s-y_r}{Y_2(y'_s)}
  =\dfrac{  (x_r-x'_s)(x_r-x'_{s+1}) }{X_2(x_r)(y'_s-y_{r-1}) }.
\end{equation}

$2\mathcal{M}^\dagger\dfrac{y-y_r }{y-y'_s}$ at $x_n$ is $
\dfrac{y_n-y_r }{y_n-y'_s}+\dfrac{y_{n-1}-y_r }{y_{n-1}-y'_s}
=\dfrac{2y_ny_{n-1} -(y_r+y'_s)(y_n+y_{n-1})+2y_ry'_s}{(y_{n}-y'_s)
( y_{n-1}-y'_s)=F(x_n,y'_s)/X_2(x_n)}  \\=
\dfrac{2X_0(x_n)+(y_r+y'_s)X_1(x_n)+2y_ry'_sX_2(x_n)}{Y_2(y'_s)(x_n-x'_{s})(x_n-x'_{s+1})}$, and
\begin{equation}\label{D1a}D_{1,r,s}^\dagger(x)=\dfrac{X_0(x)+(y_r+y'_s)X_1(x)/2+y_ry'_sX_2(x)}{Y_2(y'_s)}\end{equation}


The constant $C_{m,r,s}^\dagger$ is found through particular values of $x_n$,
involving $y_{n-1}-y_r$ and $y_{n}-y_r$  the latter only  vanishing at $x=x_r$, so

\begin{subequations}

\begin{equation}\label{Cxm1a}
\begin{split}
C_{m,r,s}^\dagger &= \dfrac{(x_{r}-x'_s)\dotsb (x_{r}-x'_{s+m})  }
                                                  {X_2(x_{r}) (x_{r}-x_{r+1})\dotsb(x_{r}-x_{r+m-1})} \\
     &\times \dfrac{     (y_{r-1}-y_r)\dotsb(y_{r-1}-y_{r+m-1}) }
                {       (y_r-y_{r-1})     (y_{r-1}-y'_s)\dotsb (y_{r-1}-y'_{s+m-1})} \\
&D_{m,r,s}^\dagger(x_r) = (y_r-y_{r-1})C_{m,r,s}^\dagger X_2(x_r)/2,
\end{split}
\end{equation}

and  $y_{n-1}-y_{r+m-1}$ and $y_{n}-y_{r+m-1}$ the first only  vanishing
at $n=r+m$:

\begin{equation}\label{Cxn1a}
\begin{split}
  C_{m,r,s}^\dagger&=-\dfrac{(x_{r+m}-x'_s)\dotsb (x_{r+m}-x'_{s+m})  }
                                    {X_2(x_{r+m}) (x_{r+m}-x_{r+1})\dotsb(x_{r+m}-x_{r+m-1}) }\\
  &\times     \dfrac{    (y_{r+m}-y_r)\dotsb(y_{r+m}-y_{r+m-1}) }
       {   (y_{r+m}-y'_s)\dotsb (y_{r+m}-y'_{s+m-1})} \\
 &D_{m,r,s}^\dagger(x_{r+m})= -(y_{r+m}-y_{r+m-1}) C_{m,r,s}^\dagger X_2(x_{r+m})/2,
\end{split}
\end{equation}
for $m=1,2,\dotsc$


Through residues at $x=x'_s$, involving $1/(y_n-y'_s)=\infty$, while $1/(y_{n-1}-y'_s)=1/(y'_{s-1}-y'_s)$ remains bounded,
replacing $(x-x'_s)/(y-y'_s)$ by $dx/dy=-(y'_s-y'_{s-1})X_2(x'_s)/((x'_s-x'_{s+1})Y_2(y'_s))$
 of the tangent at $(x'_{s+1},y'_s)$ from \eqref{tang}

$C_{m,r,s}^\dagger =- \dfrac{(y'_s-y_r)\dotsb (y'_s-y_{r+m-1})(x'_s-x'_{s+1})\dotsb(x'_s-x'_{s+m})}
          {(y'_s-y'_{s+1})\dotsb (y'_s-y'_{s+m-1})(y'_s-y'_{s-1})(x'_s-x_{r+1})\dotsb (x'_s-x_{r+m-1})X_2(x'_s)} \; \dfrac{dx}{dy}$
\begin{equation}\label{Cxp0a}\begin{split}
C_{m,r,s}^\dagger &= \dfrac{(y'_s-y_r)\dotsb (y'_s-y_{r+m-1})(x'_s-x'_{s+2})\dotsb(x'_s-x'_{s+m})  }
          {(y'_s-y'_{s+1})\dotsb (y'_s-y'_{s+m-1})(x'_s-x_{r+1})\dotsb (x'_s-x_{r+m-1})Y_2(y'_s) }\\
D_{m,r,s}^\dagger (x'_s)& = C_{m,r,s}^\dagger  (y'_s-y'_{s-1})  X_2(x'_s)/2,
 \end{split}\end{equation}
and finally at  $x=x'_{s+m}$, involving $1/(y_{n-1}-y'_{s+m-1})=\infty$,  $C_{m,r,s}^\dagger     = \dfrac{(y'_{s+m-1}-y_r)\dotsb (y'_{s+m-1}-y_{r+m-1})}{(x'_{s+m}-x'_{s})\dotsb(x'_{s+m}-x'_{s+m-1})}  $ \\  $\times
        \dfrac  {(y'_{s+m-1}-y'_{s})\dotsb (y'_{s+m-1}-y'_{s+m-1})(y'_{s+m}-y'_{s+m-1})}{(x'_{s+m}-x_{r+1})\dotsb (x'_{s+m}-x_{r+m-1})X_2(x'_{s+m})} \; \dfrac{dx}{dy}$

  \begin{equation}\label{Cxpna}\begin{split}
C_{m,r,s}^\dagger&=-\dfrac{(y'_{s+m-1}-y_r)\dotsb (y'_{s+m-1}-y_{r+m-1})(x'_{s+m}-x'_{s})\dotsb(x'_{s+m}-x'_{s+m-2})}
          {(y'_{s+m-1}-y'_{s})\dotsb (y'_{s+m-1}-y'_{s+m-1})(x'_{s+m}-x_{r+1})\dotsb (x'_{s+m}-x_{r+m-1})Y_2(y'_{s+m-1})} \\     
&D_{m,r,s}^\dagger (x'_{s+m}) = C_{m,r,s}^\dagger  (y'_{s+m}-y'_{s+m-1})  X_2(x'_{s+m})/2,
\end{split}
\end{equation}

\end{subequations}


remark from (\ref{Cxm1a}-\ref{Cxp0a}) $\dfrac{C_{m+1,r,s}^\dagger}{C_{m,r,s}^\dagger}=
\dfrac{ (x_{r}-x'_{s+m+1})(y_{r-1}-y_{r+m}) }
          {(x_{r}-x_{r+m}) (y_{r-1}-y'_{s+m})}=
\dfrac{ (x'_s-x'_{s+m+1})(y'_s-y_{r+m})}
          { (x'_s-x_{r+m})(y'_s-y'_{s+m})
}$ transposing \eqref{slopes} .

\section{Second order  difference operators.\label{secondor}}
\chead{\thesection  \ \ \ \   Second order}

\subsection{Definitions and examples}  \    \\


A linear difference operator of the second order must involve
$f(x_n)$ and the two neighboring values $f(x_{n\pm 1})$

$\alpha_n f(x_{n+1}) +\beta_n f(x_n)+\gamma_n f(x_{n-1})$

the operator $\mathcal{S}$  extending the second derivative operator must
vanish on the constants,

$\alpha_n+\beta_n+\gamma_n=0$, so

$\mathcal{S}f(x_n)= \alpha_n (f(x_{n+1})- f(x_n))+\gamma_n (f(x_{n-1})-f(x_n))$
remains. Should the operator vanish on $x$,

$\alpha_n (x_{n+1}- x_n)+\gamma_n (x_{n-1}-x_n)=0$, whence the
operator

$\mathcal{D}^\dagger\mathcal{D}f(x_n)= -\left(\dfrac{f(x_{n+1})- f(x_n)}{x_{n+1}- x_n}
-\dfrac{f(x_n)-f(x_{n-1})}{x_n-x_{n-1}}\right)/(y_{n}-y_{n-1})$,  see also \S \ref{MarMedem}.

But $f(x)=x$ is a special case of a rational function of first degree,we
consider a family of operators, each of them vanishing on
$(\alpha' x+\beta')/(\gamma' x+\delta')$ for some $(\alpha',\dotsc, \delta')$, then

$\alpha_n \dfrac{x_{n+1}- x_n}{(x_{n+1}-t)(x_n-t)}+\gamma_n
\dfrac{x_{n-1}-x_n}{(x_n-t)(x_{n-1}-t)}=0$, where $t=-\delta'/\gamma'$,  so, using
$(x_{n+1}-t)(x_n-t)= F(t,y_n)/Y_2(y_n)= X_2(x_\kappa)(y_n-y_\kappa)(y_n-y_{\kappa-1})/Y_2(y_n)$,
where $t=x_\kappa$ for some $\kappa$,

$\alpha_n= \dfrac{(y_n-y_\kappa)(y_n-y_{\kappa-1})}{Y_2(y_n)(x_{n+1}-x_n)}$,
$\gamma_n= \dfrac{(y_{n-1}-y_\kappa)(y_{n-1}-y_{\kappa-1})}{Y_2(y_{n-1})(x_{n}-x_{n-1})}$, up to multiplication by constants. We will use

\begin{equation*}
\mathcal{S}f(x)= \mathcal{D}^\dagger\left( \dfrac{(y-y_\kappa)(y-y_{\kappa-1})}
 {Y_2(y)}  \right) \mathcal{D} f(x).
\end{equation*}

Applying this $\mathcal{S}$ to the product $P_m(x)=\dfrac{(x-x_0)\dotsb(x-x_{m-1})}
 {(x-x'_0)\dotsb(x-x'_{m-1})}$ yields interesting results. From \eqref{dY}
$\mathcal{S}P_m(x)= \mathcal{D}^\dagger \dfrac{(y-y_\kappa)(y-y_{\kappa-1})}
 {Y_2(y)}   C_{m,0,0} Y_2(y) \dfrac{(y-y_0)\dotsb(y-y_{m-2}) }{(y-y'_{-1})\dotsb (y-y'_{m-1})}$. As seen before, a big simplification occurs when
$y_\kappa=y'_0$, then,

$\mathcal{S}P_m(x)= \mathcal{D}^\dagger
   C_{m,0,0} \dfrac{(y-y_0)\dotsb(y-y_{m-2}) }{(y-y'_{1})\dotsb (y-y'_{m-1})}
$

Finally,

\begin{equation} \label{secondo}
\mathcal{D}^\dagger\left( \dfrac{(y-y'_0)(y-y'_{-1})}
 {Y_2(y)}  \right) \mathcal{D} \dfrac{(x-x_0)\dotsb(x-x_{m-1})}
 {(x-x'_0)\dotsb(x-x'_{m-1})}
= C_{m,0,0}
     C^\dagger_{m-1,0,1} X_2(x) \dfrac{(x-x_1)\dotsb(x-x_{m-2}) }{(x-x'_{1})\dotsb (x-x'_{m})},
\end{equation}

and we have table \ref{tableDDMM}, from \eqref{daY}

\begin{table}[htbp]
\begin{center}\begin{tabular}{|c|c|c|c|}\hline
 & $f(x_{n-1})\times$ & $f(x_n)\times$ & $f(x_{n+1})\times$\\ \hline
$(\mathcal{D}^\dagger w(y)\mathcal{D}f)(x_n)$ &
$-\dfrac{w(y_{n-1})}{\Delta x_{n-1}\  \Delta y_{n-1}}$
 &  $\dfrac{w(y_{n-1})}{\Delta x_{n-1}\ \Delta y_{n-1}}+\dfrac{w(y_{n})}{\Delta x_{n}\ \Delta y_{n-1}}$
& $-\dfrac{w(y_{n})}{\Delta x_{n}\ \Delta y_{n-1}}$        \\
   &  \    &    &   \\
$(\mathcal{M}^\dagger w(y)\mathcal{D}f)(x_n)$ &
$-\dfrac{w(y_{n-1})}{2\Delta x_{n-1}}$
 &  $\dfrac{w(y_{n-1})}{2\Delta x_{n-1}}-\dfrac{w(y_{n})}{2\Delta x_{n}}$
& $\dfrac{w(y_{n})}{2\Delta x_{n}}$        \\
&  \  $\vphantom{\dfrac{A}{B}  }$ &    &   \\
  $(\mathcal{D}^\dagger w(y)\mathcal{M}f)(x_n)$ &  $\dfrac{w(y_{n-1}}{2\Delta y_{n-1}}$ &
          $-\dfrac{ w(y_n)-w(y_{n-1})}{2\Delta y_{n-1}}$ & $-\dfrac{w(y_{n}}{2\Delta y_{n-1}}$ \\
&  \  $\vphantom{\dfrac{A}{B}  }$ &    &   \\
$(\mathcal{M}^\dagger w(y)\mathcal{M}f)(x_n)$ &  $w(y_{n-1})/4$ &
$(w(y_{n-1})+w(y_{n}))/4$ &$w(y_{n})/4$ \\
 \hline
\end{tabular}
\caption{Second order operators applied to $f$ at $x_n$ seen as
linear combinations of $f(x_{n-1}), f(x_n)$, and $f(x_{n+1})$
(and where $\Delta x_n=x_{n+1}-x_n, \Delta y_n=y_{n+1}-y_n$).\label{tableDDMM}} \end{center}\end{table}

\subsection{Adjoint operators} \   \\

We defined in \S \ref{defadj} the adjoint of $\mathcal{D}$ with respect to the
bilinear form $\sum_k f(x_k)g(x_k)\Delta y_{k-1}$. Keeping this bilinear form,
the adjoint of a general linear second order operator
$\mathcal{A}f(x_n)=\alpha_n f(x_{n+1})+\beta_n f(x_n)+\gamma_n f(x_{n-1})$ is such that
$\sum_n \Delta y_{n-1} (\alpha_n f(x_{n+1})+\beta_n f(x_n)+\gamma_n f(x_{n-1}))
  g(x_n)$ is reordered as
$\sum_n f(x_n) [ g(x_{n-1})\alpha_{n-1}\Delta y_{n-2} +(x_n)\beta_n\Delta y_{n-1}
   +g(x_{n+1})\gamma_{n+1}\Delta y_{n}]$
so, the adjoint is

 $\mathcal{A}^\dagger g(x_n)=\dfrac{g(x_{n-1})\alpha_{n-1}\Delta y_{n-2} +g(x_n)\beta_n\Delta y_{n-1}
   +g(x_{n+1})\gamma_{n+1}\Delta y_{n}}{\Delta y_{n-1}}$.

For instance, the adjoint of $\mathcal{D}^\dagger w(y)\mathcal{D}$ is, from
table~\ref{tableDDMM},

$\dfrac{-g(x_{n-1})w(y_{n-1})/\Delta x_{n-1} +g(x_n)(w(y_{n-1})/\Delta x_{-1}
    +w(y_n)/\Delta x_n)-g(x_{n+1})w(y_n)/\Delta x_{n}}{\Delta y_{n-1}}$
at $g(x_n)$, and we recover the same operator!

\textbf{Definition.} Let $f$ and $g$ be defined on at least one lattice $\{x_n\}$. The
operator $\mathcal{A}$ is selfadjoint
with respect to the bilinear form $\langle .,. \rangle$ if
$\langle f,\mathcal{A}g\rangle$ is symmetric in $f$ and $g$.

with respect to the bilinear form $\sum_n \omega_n f(x_n) g(x_n)$,
one must look at
$\sum_n (y_n-y_{n-1})f(x_n)(\mathcal{A}g)(x_n)$

With the matrix formulation
$\langle f,\mathcal{A}g\rangle \\ = [ \dotsb , f(x_{n-1}, f(x_n),f(x_{n+1}), \dotsb ]
  \begin{bmatrix} \ddots &    &    &    &    \\
                        & \omega_{n-1}\beta_{n-1} &\omega_{n-1}\alpha_{n-1} &   &  \\
                     &\omega_n\gamma_n & \omega_n\beta_n   & \omega_n\alpha_n &   \\
                  &   & \omega_{n+1}\gamma_{n+1} & \omega_{n+1}\beta_{n+1} &   \\
                    &   &   &    &  \ddots  \end{bmatrix}
\begin{bmatrix}  \vdots \\  g(x_{n-1}) \\  g(x_n) \\  g(x_{n+1}) \\  \vdots\end{bmatrix}
 $
indeed \\   symmetric  $f$ and $g$ if the matrix is symmetric, $  \omega_n\alpha_n=         \omega_{n+1}\gamma_{n+1}$

We also check the selfadjointness of  $\mathcal{D}^\dagger (w(y)\mathcal{D} )$
on the matrix form $  \alpha_n\Delta y_{n-1}=         \gamma_{n+1}\Delta y_{n}
= w(y_n)/\Delta x_n$ from table \ref{tableDDMM}.

\medskip

For operators involving $\mathcal{M}^\dagger$, which is not the true adjoint of $\mathcal{M}$,
it is better to check \emph{by the book}: let $\mathcal{L}=a(x)(\mathcal{M}^\dagger
+b(x)\mathcal{D}^\dagger)(\mathcal{M}+c(y)\mathcal{D})$,
then, with $h(y)=(\mathcal{M}g)(y)+c(y)(\mathcal{D}g)(y)$,

$\langle f,\mathcal{L}g\rangle= \sum_n \omega_n f(x_n)a(x_n)(( \mathcal{M}^\dagger
+b(x)\mathcal{D}^\dagger)h)(x_n)\\
= \sum_n \omega_n f(x_n)a(x_n)\left[ \left( \dfrac{1}{2}-\dfrac{b(x_n)}{y_n-y_{n-1}}\right)h(y_n)
+\left( \dfrac{1}{2}+\dfrac{b(x_n)}{y_n-y_{n-1}}\right)h(y_{n-1})\right]\\
= \sum_n h(y_n) \left[ \omega_{n+1} f(x_{n+1})a(x_{n+1})
    \left( \dfrac{1}{2}+\dfrac{b(x_{n+1})}{y_{n+1}-y_{n}}\right)
   +\omega_n f(x_n)a(x_n)\left( \dfrac{1}{2}-\dfrac{b(x_n)}{y_n-y_{n-1}}\right)\right]\\
=\sum_n [\mathcal{M}g)(y_n)+c(y_n)(\mathcal{D}g)(y_n)]
    \left[ (\mathcal{M} a\omega f)(y_n) + (x_{n+1}-x_n)
            \left(\mathcal{D}\dfrac{ \omega f ab }{y_n-y_{n-1}}\right)(y_n)\right]
$


 \subsection{Hypergeometric difference equation and expansion} \    \\

An elliptic hypergeometric expansion is a (formal) sum $\sum_n R_1 R_2\dots R_n$,
where $R_n$ is a rational function $R(x_n,y_n,x'_n,y'_n,\dotsc,x_n^{(r)},y_n^{(r)})$,
and where $(x_n,y_n),\dotsc, (x_n^{(r)},y_n^{(r)})$ are sequences of points on the
same biquadratic curve $F(x,y)=0$.


\subsection{Some classical hypergeometric expansions\label{someclassical}} \   \\

\subsubsection{The expansion of ${}_2F_1$ from the differential equation}  \   \\

Let $\mu(x)f''(x)+\nu(x)f'(x)-\lambda f(x)=0$ where $\mu, \nu,\lambda$ are polynomials of
degrees $2,1,0$.

Put a (formal) Taylor expansion $f(x)=\sum_0^\infty c_n(x-a)^n$ in the equation
to get

$[\mu(a)+\mu'(a)(x-a)+\mu''(a)(x-a)^2/2]\sum_0^\infty n(n-1)c_n (x-a)^{n-2}
+[\nu(a)+\nu'(a)(x-a)]\sum_0^\infty nc_n (x-a)^{n-1}
-\lambda\sum_0^\infty c_n (x-a)^{n}=0$, so, a three-term recurrence relation
linking $c_n$, $c_{n+1}$ and $c_{n+2}$

$\mu(a)(n+2)(n+1)c_{n+2} +\mu'(a)(n+1)nc_{n+1}+\mu''(a)n(n-1)c_n/2
+\nu(a)(n+1)c_{n+1}+\nu'(a)nc_n -\lambda c_n=0$, $n=0,1,\dotsc$ follows.

A big simplification occurs if $x=a$ is a root of $\mu(x)=0$ (singular
point), then

$c_{n+1} =\dfrac{  \lambda -n\nu'(a) -n(n-1)\mu''(a)/2}{(n+1)(n\mu'(a)+\nu(a))} c_n$, so,

$c_n=\dfrac{ \Gamma(\alpha+n)\Gamma(\beta+n) }{n! \Gamma(\gamma+n)}(-\mu''(a)/(2\mu'(a)))^n c_0$,
where $-\alpha$ and $-\beta$ are the $n-$roots of  $\lambda -n\nu'(a) -n(n-1)\mu''(a)/2=0$, and where $\gamma=\nu(a)/\mu'(a)$,

leading to a hypergeometric  \cite[10.3, 10.7]{WW}

$f(x)= {}_2F_1(\alpha,\beta;\gamma;-\mu'(a)(x-a)/(2\mu'(a))$

\medskip

Another technique:

Derivate to get $\mu(x)f'''(x)+(\mu'(x)+\nu(x))f''(x)+(\nu'-\lambda)f'(x)=0$ has the same form
with same degrees!

The equation for $f^{(n)}(x)$ has the same $\mu(x)$, and  $\nu_n(x)= n\mu'(x)+\nu(x)$,
and the new constant coefficient $\lambda_{n+1}= \lambda_n-\nu'_n =\lambda_n-n\mu''-\nu'$, or $\lambda_n=-n(n-1)\mu''/2-n\nu'+\lambda $

If $x=a$ is a singular point $(\mu(a)=0)$  with, say, $f(a)=1$, we have
$f^{(n+1)}(a)/f^{(n)}(a)= -\lambda_n/\nu_n(a) =[n(n-1)\mu''/2+n\nu'-\lambda ]/(n\mu'(a)+\nu(a)) $ leading to the same
Taylor expansion about $a$.

Darboux-B\"acklund \cite{DarbouxBack}


\subsubsection{Discrete equation on the lattice $\{x_n=x_0+nh\}$. }  \  \\

$\mu(x)\Delta\nabla f(x)/h^2 +\nu(x)(\Delta+\nabla)f(x)/(2h)-\lambda f(x)=0$ \cite{NS,NSU}

with degrees $2,1,0$ as before.

with the operators $\Delta$ and $\nabla$ of table~\ref{difftable} p.~\pageref{difftable}

note that $\Delta\nabla f(x)=\nabla\Delta f(x)=f(x+h)-2f(x)+f(x-h)$

we now put $f(x)=\sum_0^\infty c_n (x-a)(x-a-h)\dotsb (x-a-(n-1)h)$,

then 
$\Delta f(x)= \sum nhc_n (x-a)\dotsb (x-a-(n-2)h)$,
$\nabla f(x)= \sum nhc_n (x-a-h)\dotsb (x-a-(n-1)h)$,
$\Delta\nabla f(x)=\sum n(n-1)h^2c_n (x-a-h)\dotsb (x-a-(n-2)h)$

$0=\sum \{ n(n-1)\mu(x) +n\nu(x)(x-a-(n-1)h)/2 \}c_n(x-a-h)\dotsb(x-a-(n-2)h)
+\sum\{n\nu(x)/2-\lambda(x-a-(n-1)h)\}c_n(x-a)\dotsb(x-a-(n-2)h)$

we find a combination of products $(x-a-h)(x-a-2h)\dotsb (x-a-mh)$, $m=0,1,\dotsc$, each product involving $c_{m+2}, c_{m+1}$ and $c_m$.

A big simplification now occurs if $x=a$ is a root of $\mu(x)-h\nu(x)/2=0$, then,
let $\mu(x)-h\nu(x)/2=\mu''(x-a)(x-b)/2$ and we have

$0=\sum n(n-1)\mu'' c_n (x-a)\dotsb (x-a-(n-2)h)(x-b=x-a-(n-1)h+a-b+(n-1)h)/2
+\sum n c_n (x-a)\dotsb(x-a-(n-2)h)(\nu(x)=\nu(a+(n-1)h)+\nu'(x-a-(n-1)h))
-\lambda \sum c_n (x-a)\dotsb (x-a-(n-1)h)$ and

$n(n-1)\mu''c_n/2 +(n+1)n\mu''(a-b+nh)c_{n+1}/2
+n\nu' c_n +(n+1)\nu(a+nh)c_{n+1} -\lambda c_n=0$
follows, or

$c_{n+1}=\dfrac{ \lambda -n\nu' -n(n-1)\mu''/2 }
  { (n+1)[ n(a-b+nh)\mu''/2 +\nu(a+nh)]} c_n$

gives ${}_3F_2$  \cite[\S 2.2.1, 2.7.3.3]{NSU}  \cite[\S 18.20]{NIST}  \cite[chap. 5 \& 9]{Koe}

 with the $3-$term recurrence rel.

$(\mu(x)/h^2-\nu(x)/(2h))f(x-h)-(2\mu(x)/h^2+\lambda)f(x)
+(\mu(x)/h^2+\nu(x)/(2h))f(x+h)=0$,
we see that the condition $\mu(a)-h\nu(a)/2=0$ shows that $x=a$ is a singular point,
that $f(a+h)$ is not a linear combination of $f(a-h)$ and $f(a)$
but depends only on $f(a): f(a+h)=\dfrac{ \lambda h^2+2\mu(a)}{\mu(a)+h\nu(a)/2}f(a)=\dfrac{ \lambda h +\nu(a) }{\nu(a)}f(a)$

\medskip

other way: apply $\Delta$ and see that we have an equation of the same form
for $\Delta f$

$\Delta(\mu\Delta\nabla)f/h^2+\Delta(\nu(\Delta+\nabla))f/(2h)-\lambda\Delta f=0$

we need $\Delta(uv)(x)=u(x+h)v(x+h)-u(x)v(x)=v(x+h)\Delta u(x) +u(x)\Delta v(x)$

$\mu(x+h)\Delta^2\nabla f/h^2  +(\Delta \mu)
\Delta\nabla f/h^2
+\nu(x+h)(\Delta^2+\Delta\nabla)f/(2h)
+(\Delta \nu)(\Delta+\nabla)f/(2h)-\lambda\Delta f=0$

$\{\mu(x+h)\Delta\nabla /h^2  +(\Delta \mu)
\nabla /h^2
+\nu(x+h)(\Delta+\nabla)/(2h)
+(\Delta \nu)(2-\nabla)/(2h)-\lambda\} \Delta f=0$
using $\nabla-\Delta= -\Delta\nabla$,

$\{\mu(x+h)\Delta\nabla /h^2  +(\Delta (\mu-h\nu/2))
(\Delta+\nabla-\Delta\nabla) /(2h^2)
+\nu(x+h)(\Delta+\nabla)/(2h)
+(\Delta \nu)/h-\lambda\} \Delta f=0$

using $\nabla=(\nabla+\Delta)/2-\Delta\nabla/2$

so, $\mu_{\text{new}}= \mu(x+h)-\Delta (\mu-h\nu/2)/2
=\mu+\Delta (\mu+h\nu/2)/2
$,

$\nu_{\text{new}}= \Delta (\mu-h\nu/2)/h+\nu(x+h)=\nu(x+h)+\Delta \mu''(x-a)(x-b)/(2h)=\nu(x+h)+\mu''(x-(a+b-h)/2)$

$\lambda_{\text{new}} =\lambda -\Delta \nu/h=\lambda-\nu'$

note that $\mu-h\nu/2 $ is unchanged: $x=a$ remains a singular point of the new equation. After $n$ steps, $\nu_{n+1}(x)=\nu_n(x+h)+\mu''(x+nh-(a+b-h)/2)$, so,
$\nu_n(x)=\nu(x+nh)+n\mu''(x-(a+b-h)/2)+n(n-1)\mu''h/2=\nu(x+nh)+n\mu''(x-(a+b-nh)/2)$,

$\lambda_{n+1}= \lambda_n -\nu'_n=\lambda_n -\nu' -n \mu''$, so
$\lambda_n= \lambda-n\nu' -n(n-1)\mu''/2$, and

$\dfrac{\Delta^{n+1}f(a)}{\Delta^n f(a)}= \dfrac{ \lambda_n}{\nu_n(a)} =
\dfrac{   \lambda-n\nu' -n(n-1)\mu''/2 }{\nu(a+nh)+n\mu''(a-b+nh)/2}$ is confirmed.

\subsubsection{$q-$Sturm-Liouville difference equation.\label{MarMedem}} We follow  Marcell\'an and   Medem \cite{MarETNA}, and consider $\mathcal{L}=\mu \Theta\Theta^*  +\nu \Theta^*  - \lambda$, with
operators $\Theta f(x)=(f(qx)-f(x))/(qx-x),  \Theta^* f(x)=(f(x/q)-f(x))/(x/q-x)$, and $\mu, \nu$ are polynomials of degrees $\leqslant 2, 1$.

When $\mathcal{L}f$ is written as a combination of $f(qx), f(x), f(x/q)$, we have
$\mu(x)f(qx)/(x-qx)^2$ and $ q [ \mu(x)+(x-qx)\nu(x)  ]     f(x/q)/(x-qx)^2$,
and the zeros of $\mu$ and $\mu(x)+(x-qx)\nu(x))$, singular points of  $\mathcal{L}$, are involved in the support
of the relevant orthogonality measure \cite[eq. (2.6), table 3.1]{MarETNA}.

\subsubsection{ANSUW difference equation.}  \   \\

Askey and Wilson \cite[eq. (5.16)]{AW319} find the equation $D_q[w_1'x)D_qp_n(x)]+\lambda_n w_2(x)p_n(x)=0$
for their orthogonal polynomials $p_n$, where $D_q f(x)= \dfrac{ f(q^{1/2}x)-f(q^{-1/2}x) }
{q^{1/2}x-q^{-1/2}x }$.
\\
Nikiforov et al. \cite[eq. (3.1.
5)]{NSU} find $\tilde{\sigma}[x(s)]\dfrac{\Delta}{\Delta x(s-1/2)}
\left[ \dfrac{\nabla f(x(s)) }{\nabla x(s)}\right]
+\dfrac{\tilde{\tau}[x(s)]}{2}
\left[ \dfrac{\Delta f(x(s)) }{\Delta x(s)}+
 \dfrac{\nabla f(x(s)) }{\nabla x(s)}\right] \\  +\lambda f(x(s))=0$
on the NSU lattice $x(s)$ of table~\ref{difftable}.

We recognize the $\mathcal{D}, \mathcal{D}^\dagger$, and $\mathcal{M}$ operators
of the present sudy.

See also \cite{bang,kenfack,koorn2022,verdestar,newtbasis} for interesting constructions and  discussions.


\subsection{\label{elldefL}Second order elliptic difference equation}  \   \\

We introduce a set of second order equations allowing hypergeometric expansions   \cite{ellbib}.

\textbf{Theorem.} \textsl{ Let $\{x_n,y_n\}$ and $\{x'_n,y'_n\}$,
$n\in\mathbb{Z}$,  be two
elliptic sequences on  the same biquadratic curve $F(x,y)=X_0(x)+X_1(x)y+X_2(x)y^2=Y_0(y)+xY_1(y)+x^2Y_2(y)=0$, and
let $\mathcal{L}$ be the second order difference operator  }
\begin{multline}\label{defL}
   \mathcal{L}=\dfrac{1}{x-x'_1}\left( \dfrac{\mu(x)}{X_2(x)} \mathcal{D}^\dagger +\nu(x) \mathcal{M}^\dagger \right)
   \left( \dfrac{(y-y'_{-1})(y-y'_0)} {Y_2(y)}  \mathcal{D}\right) \\
-\lambda \left(\mathcal{M}^\dagger + \dfrac{y'_0-y'_{-1}}{2} \mathcal{D}^\dagger\right)
   \left( \mathcal{M} -\dfrac{R(y)}{Y_2(y)}
          \mathcal{D}\right),
\end{multline}
 \textsl{where $\mu$ and $\nu$ are polynomials of degrees 3 and 1,
 $\lambda$ is a constant, and $R$ is the polynomial }\\
 $R(y)= \dfrac{ (y-y_{-1})Y_2(y'_{-1}) (x'_0-x'_{-1})+(y-y'_{-1})Y_2(y_{-1}) (x_0-x_{-1})}{2(y'_{-1}-y_{-1}) }$. \\
\textsl{Moreover, }

\begin{equation}\label{munux0}\begin{split}
2\mu(x_0)  -(y_0-y_{-1} )X_2(x_0)\nu(x_0)&=0,  \\
2\mu(x'_1)  +(y'_1-y'_{0} )X_2(x'_1)\nu(x'_1)&=0.
\end{split}\end{equation}

\textsl{The equation $\mathcal{L}\phi=0$ has then a solution
with an interpolating expansion}
\begin{equation}\label{sumf}
\phi(x)=\sum_{n=0}^\infty c_n \dfrac{ (x-x_0)\dotsb (x-x_{n-1}) }
                        {(x-x'_0)\dotsb (x-x'_{n-1})}
\end{equation}
\textsl{where }$c_0=f(x_0)$ \textsl{may be chosen at will, and with}

$c_n= c_0 (x_{n-1}-x'_{n}) \dfrac{\zeta_0\dotsb \zeta_{n-1}}{\alpha_0\dotsb \alpha_{n-1}}, n=1,2,\dotsc$

\text{with}


\begin{multline}\label{alphan}\alpha_n=
  \left( -\dfrac{\mu(x_n)}{(y_n-y_{n-1})X_2(x_n)}           +\dfrac{\nu(x_n)}{2}\right)      \dfrac{                (y_{n}-y'_{-1})(y_{n}-y'_0)}{ (x_{n+1}-x_{n})   (x_n-x'_1)Y_2(y_{n}) }\\
-\dfrac{\lambda}{2} \left[  \left(1-\dfrac{y'_0-y'_{-1} }{y_n-y_{n-1}}\right)\left(\dfrac{1}{2}-\dfrac{R(y_{n}) }{(x_{n+1}-x_{n})Y_2(y_{n})}\right)\right]    ,
\end{multline}



The formula for  
$\zeta_n$ will be given later on in the proof,  all the $\zeta_n$s are first degree polynomials in $\lambda$,
$R$ is here the $R_{0,0}$ of \eqref{Dkinterp}.

Note that $R$ achieves the linear interpolation of $R(y_{-1})=-(x_0-x_{-1})Y_2(y_{-1})/2, \\
  R(y'_{-1})=(x'_0-x'_{-1})Y_2(y'_{-1})/2$,

 all the $\zeta_n$s are first degree polynomials in $\lambda$.

The first part of $\mathcal{L}$ generalizes the terms
$a(x)d^2/dx^2+b(x)d/dx$ of the original differential hypergeometric
operator. The second part corresponds to the last constant term.
Why not finish simply with $-\lambda$ times the identity operator?
This obvious choice does not lead to a workable ratio $c_{n+1}/c_n$
(the author tried). One hopes that a less awkward construction
will be found someday.

We start the proof by 
 writing $\mathcal{L}\phi=0$ in \eqref{defL} as a recurrence relation
at $x_n$, according to \eqref{adj} and \eqref{adjM},

\begin{equation}\label{recurL}
 \alpha_n \phi(x_{n+1}) +\beta_n \phi(x_n)+\gamma_n \phi(x_{n-1})=0, \end{equation}

{\small

$(   \mathcal{L}\phi)(x_n)=\dfrac{1}{x-x'_1}\left[
            \left( \dfrac{\mu(x_n)}{(y_n-y_{n-1})X_2(x_n)}           +\dfrac{\nu(x_n)}{2}\right)          \dfrac{            (y_{n-1}-y'_{-1})(y_{n-1}-y'_0)(\mathcal{D}\phi)(y_{n-1})}{ Y_2(y_{n-1})}
                   \right.    \\  \left.
          \ \ \ \ \ \ +   \left( \dfrac{-\mu(x_n)}{(y_n-y_{n-1})X_2(x_n)}           +\dfrac{\nu(x_n)}{2}\right)             \dfrac{         (y_{n}-y'_{-1})(y_{n}-y'_0)(\mathcal{D}\phi)(y_{n})}{ Y_2(y_{n})}
        \right]$

$-\dfrac{\lambda}{2} \left[  \left(1+\dfrac{y'_0-y'_{-1} }{y_n-y_{n-1}}\right)\left(\mathcal{M}\phi-\dfrac{R \mathcal{D}\phi}{Y_2}\right) (y_{n-1}) 
+\left(1-\dfrac{y'_0-y'_{-1} }{y_n-y_{n-1}}\right)\left( \mathcal{M}\phi -\dfrac{R\mathcal{D}\phi}{  Y_2}\right)   (y_{n})   \right]
$,

and we extract the coefficients of $\phi(x_{n+1})$ and $\phi(x_{n-1})$ from $(\mathcal{D}\phi)   (y_{n})$ and $ ( \mathcal{M}\phi)   (y_{n})$, giving
 $\alpha_n$ of \eqref{alphan},

\begin{multline}\label{gamman}\gamma_n =
  - \left( \dfrac{\mu(x_n)}{(y_n-y_{n-1})X_2(x_n)}           +\dfrac{\nu(x_n)}{2}\right)      \dfrac{                (y_{n-1}-y'_{-1})(y_{n-1}-y'_0)}{ (x_n-x_{n-1})   (x_n-x'_1)Y_2(y_{n-1}) }\\
-\dfrac{\lambda}{2} \left[  \left(1+\dfrac{y'_0-y'_{-1} }{y_n-y_{n-1}}\right)\left(\dfrac{1}{2}+\dfrac{R(y_{n-1}) }{(x_n-x_{n-1})Y_2(y_{n-1})}\right)\right]
\end{multline},

and \begin{equation}\label{abg} \alpha_n+\beta_n+\gamma_n=\mathcal{L}1= -\lambda.
\end{equation}




In matrix form
\begin{equation}\label{Lmatrix}\begin{split}
   \begin{bmatrix}\beta_0 & \alpha_0 &          &  \\
                    \gamma_1& \beta_1  & \alpha_1 &  \\
                            &\gamma_2  & \beta_2  & \alpha_2 \\
                            &          & \ddots   &  \end{bmatrix}
    \begin{bmatrix} \phi(x_0) \\  \phi(x_1)  \\  \phi(x_2) \\ \vdots \end{bmatrix}
  &=\begin{bmatrix} 0 \\  0  \\  0 \\ \vdots \end{bmatrix}\\
\text{or\ \  }  \mathsf{A}\mathbf{\phi}+\lambda\mathsf{B}\mathbf{\phi}&=0,
\end{split}\end{equation}
 where $\mathbf{\phi}$ is the column vector
$\phi(x_0),\phi(x_1),\dotsc$. Also, $\gamma_0=0$, see later on.

 One has $\phi(x_1)=-\beta_0\phi(x_0)/\alpha_0,
\phi(x_2)=(-\gamma_1\phi(x_0)-\beta_1 \phi(x_1))/\alpha_1=(\beta_0\beta_1-\alpha_0\gamma_1)
\phi(x_0)/(\alpha_0\alpha_1)$,etc.: each $\phi(x_n)$ is a polynomial in $\mu,\nu,\lambda$
divided by $\alpha_0\dotsb \alpha_{n-1}$.

Each $\phi(x_n)=\dfrac{ -\beta_{n-1}(\lambda)\phi(x_{n-1}) -\gamma_{n-1}(\lambda)\phi(x_{n-2})}{\alpha_{n-1}(\lambda)} $ depends on the two earlier values of $\phi$, but
as the factors $\dfrac{\mu(x_n) } {(y_n-y_{n-1})X_2(x_n) }
        -\dfrac{\nu(x_n)}{2}$ and
$\dfrac{1}{2}+\dfrac{R(y_{n-1})}{(x_{n}-x_{n-1})Y_2(y_{n-1})}$ in $\gamma_n(\lambda)$ do vanish
at $n=0$, from the first equation of \eqref{munux0} and the value of $R(y_{-1})$,
 $\gamma_0(\lambda)\equiv 0$,   $x_0$ is a singular point
of the recurrence relation, $\phi(x_1)$ can not be chosen at will but is related to $\phi(x_0)$ by
$\phi(x_1)=-\dfrac{\phi(x_0)\beta_0(\lambda)}{\alpha_0(\lambda)}= \phi(x_0)\dfrac{\lambda+\alpha_0(\lambda)}{\alpha_0(\lambda)}$,
$\phi(x_2)=-\dfrac{\beta_1\phi(x_1)+\gamma_1\phi(x_0)}{\alpha_1}
= \phi(x_0)\dfrac{ \beta_1\beta_0-\gamma_1\alpha_0}
       {\alpha_0\alpha_1}
=\phi(x_0)\dfrac{- \beta_1 \lambda+\lambda\alpha_0+\alpha_1\alpha_0}  {\alpha_0\alpha_1}$. Each $\phi(x_n)$ is therefore a rational
function of degree $n$ in $\lambda$, with denominator $\alpha_0(\lambda)\dotsb\alpha_{n-1}(\lambda)$, and so is $c_m$, a divided difference
involving $\phi(x_0),\dotsc, \phi(x_{m})$, from \eqref{cdivdiff}:
$\dfrac{c_1}{c_0}= \dfrac{x_1-x'_0}{x_1-x_0}\left[\dfrac{\phi(x_1)}{\phi(x_0)}-1\right] =\dfrac{(x_1-x'_0)\lambda}{(x_1-x_0)\alpha_0}$,
we shall also need the next step

$\dfrac{c_2}{c_1}=\dfrac{    \dfrac{(x_2-x'_0)(x_2-x'_1)}{(x_2-x_0)(x_2-x_1)}\left( \dfrac{\phi(x_2)}{\phi(x_0)}-1-\dfrac{c_1(x_2-x_0)}{c_0(x_2-x'_0)}\right)  }{ c_1/c_0}  \\   =
 \dfrac{(x_2-x'_0)(x_2-x'_1)}{(x_2-x_0)(x_2-x_1)}  \left((x_1-x_0)
   \dfrac{- \beta_1 \lambda+\lambda\alpha_0}  {(x_1-x'_0)\lambda\alpha_1}-\dfrac{x_2-x_0}{x_2-x'_0}\right)
= \dfrac{x_2-x'_1}{x_2-x_1}\left( \dfrac{(x_1-x_0)(x_2-x'_0)(\alpha_0-\beta_1)}{(x_1-x'_0)(x_2-x_0)\alpha_1}-1\right)
$

whence $\zeta_0=\dfrac{\lambda(x_1-x'_0)}{(x_0-x'_1)(x_1-x_0)}$ and $\zeta_1=\dfrac{(x_0-x'_1)(x_2-x'_1)}{(x_1-x'_2)(x_2-x_1)}\left( \dfrac{(x_1-x_0)(x_2-x'_0)(\alpha_0-\beta_1)}{(x_1-x'_0)(x_2-x_0)}-\alpha_1\right)
$.

\bigskip

We proceed with the proof by  entering \eqref{sumf} in \eqref{defL}, here is the intermediate step:

\subsubsection{\label{lemmaP}}\textbf{Lemma.} \textsl{let $P_m(x)=  \dfrac{(x-x_0)\dotsb(x-x_{m-1}) }{(x-x'_0)\dotsb (x-x'_{m-1})}$ and
$\widetilde{P}_m(x)= \dfrac{(x-x_0)\dotsb(x-x_{m-1}) }{(x-x'_1)\dotsb (x-x'_{m})}$.  Then there are rational functions $\mathfrak{r}_m$ and $\mathfrak{s}_m$ of
fixed degrees in $x$ and $y$ such that}
\begin{equation}\label{LPk}
\mathcal{L}P_0=-\lambda,\ \ \ \ \mathcal{L}P_m=\mathfrak{r}_m \widetilde{P}_m
  +\mathfrak{s}_m \widetilde{P}_{m-1}, \ \  m=1,2,\dotsc,
\end{equation} \textsl{with the operator $\mathcal{L}$ of \eqref{defL}
 and the conditions \eqref{munux0}.  }

See \cite{koorn2022,verdestar,newtbasis} on introducing interpolatory-type expansions in differential and difference equations.

We first apply $\mathcal{D}$ and $\mathcal{M} - \dfrac{R(y)}{Y_2(y)}  \mathcal{D}$
to $P_m$, finding from \eqref{dY}, \eqref{secondo} and \eqref{Dkinterp0} at $r=s=0$,

$C_{m,0,0}Y_2(y)\dfrac{(y-y_0)\dotsb(y-y_{m-2}) }{(y-y'_{-1})\dotsb (y-y'_{m-1})}$ and
$ \underbrace{[D_{m,0,0}(y) - C_{m,0,0}R(y)]}_{ \displaystyle
        \rho_{m,0,0} C_{m,0,0} (y-y_{-1})(y-y'_{-1})  }
\dfrac{(y-y_0)\dotsb(y-y_{m-2}) }{(y-y'_{-1})\dotsb (y-y'_{m-1})}.
$

Then,
$\mathcal{L}P_m(x)  \\  =\dfrac{1}{x-x'_1}\left( \dfrac{\mu(x)}{X_2(x)} \mathcal{D}^\dagger +\nu(x) \mathcal{M}^\dagger \right)
   \left( \dfrac{(y-y'_{-1})(y-y'_{0}) } {Y_2(y)}
 C_{m,0,0}Y_2(y)\dfrac{(y-y_0)\dotsb(y-y_{m-2}) }{(y-y'_{-1})\dotsb (y-y'_{m-1})}\right) \\
  -\lambda \left(\mathcal{M}^\dagger
+\dfrac{y'_0-y'_{-1}}{2X_2(x)} \mathcal{D}^\dagger\right)
          \rho_{m,0,0} C_{m,0,0}
\dfrac{(y-y_{-1})\dotsb(y-y_{m-2}) }{(y-y'_{0})\dotsb (y-y'_{m-1})}
$

$=   {C}_{m,0,0}
\dfrac{\mu(x)C_{m-1,0,1}^\dagger
+\nu(x) D_{m-1,0,1}^\dagger(x)}{x-x'_1}\ \ \dfrac{(x-x_1)\dotsb(x-x_{m-2}) }{(x-x'_1)\dotsb (x-x'_{m})}\\
 -\lambda \rho_{m,0,0} C_{m,0,0}
[ D_{m,-1,0}^\dagger(x) -  C_{m,-1,0}^\dagger (y'_0-y'_{-1}) X_2(x)/2]
  \dfrac{(x-x_{0})\dotsb(x-x_{m-2}) }{(x-x'_{0})\dotsb (x-x'_{m})}$, or

\begin{equation}\label{Sm}
\mathcal{L}P_m(x) = S_m(x)\dfrac{(x-x_{0})\dotsb(x-x_{m-2}) }{(x-x'_{1})\dotsb (x-x'_{m})}
=S_m(x)\dfrac{\widetilde{P}_{m-1}(x)}{x-x'_m},
\end{equation}
where $S_m$ is a polynomial of first degree, thanks to the conditions \eqref{munux0}, indeed,

$\mu(x)C_{m-1,0,1}^\dagger+\nu(x) D_{m-1,0,1}^\dagger(x)=
   \nu(x)[ \pm C_{m-1,0,1}^\dagger X_2(x)\Delta y/2+D_{m-1,0,1}^\dagger(x)]=0$
at $x=x_0$ and $x=x'_1$ from \eqref{munux0}, \eqref{Cxm1a} and \eqref{Cxp0a} (when $m>1$), and

$D_{m,-1,0}^\dagger(x) -  C_{m,-1,0}^\dagger (y'_0-y'_{-1}) X_2(x)/2=0$ at $x=x'_0$ from \eqref{Cxp0a}.


The use of the $\mathcal{D},\dotsc,\mathcal{M}^\dagger$ formalism
allowed to exhibit a rational function of limited complexity for $\mathcal{L}P_m/P_m$, which is easiest to study through the recurrence
relation \eqref{recurL}:

\begin{equation}\label{Sm2}
S_m(x_n)=\dfrac{(\mathcal{L}P_m)(x_n)}{\widetilde{P}_{m-1}(x_n)/(x_n-x'_m)}=
\dfrac{\ \ \alpha_n P_m(x_{n+1}) +\beta_n P_m(x_n)+\gamma_n P_m(x_{n-1})\ \  }{  \frac{(x_n-x_{0})\dotsb(x_n-x_{m-2}) }{(x_n-x'_{1})\dotsb (x_n-x'_{m})}=\frac{P_{m-1,0,2}(x_n)}{x_n-x'_1}   }
\end{equation}
from \eqref{Sm} and \eqref{PP}.


Consider  the linear interpolation $S_m(x)=\dfrac{ (x-x'_m)S_m(x_{m-1})-(x-x_{m-1})S_m(x'_m)
 }{x_{m-1}-x'_m}$ ,  \eqref{Sm} becomes

$\mathcal{L}P_m(x)
=\dfrac{x-x'_0}{x_{m-1}-x'_m}  \left[  S_m(x_{m-1})\dfrac{P_m(x)}{x-x_{m-1}} -  S_m(x'_{m})\dfrac{P_m(x)}{x-x'_{m}}\right]$,

$=\dfrac{S_m(x_{m-1})}{x_{m-1}-x'_m}
\widetilde{P}_{m-1}(x)-\dfrac{S_m(x'_m)}{x_{m-1}-x'_m}
\widetilde{P}_m(x)
$,
achieving the proof of the lemma \ref{lemmaP}.  \qed


We proceed with the expansion

$\displaystyle 0=\mathcal{L}(\sum_0^\infty c_m P_m)= -\lambda c_0+\sum_1^\infty c_m\left(\dfrac{S_m(x_{m-1})}{x_{m-1}-x'_m}
\widetilde{P}_{m-1}(x)-\dfrac{S_m(x'_m)}{x_{m-1}-x'_m}
\widetilde{P}_m(x) \right)$

$ -\lambda c_0 +c_1\dfrac{S_1(x_{0})}{x_{0}-x'_1}=0,\  \ \  \dfrac{S_{m+1}(x_{m})}{x_{m}-x'_{m+1}}c_{m+1}
-\dfrac{S_m(x'_m)}{x_{m-1}-x'_m}c_m=0, m=1,2,\dotsc$

$c_m=c_0 (x_{m-1}-x'_m)\dfrac{ S_0(x'_0)S_1(x'_1)\dotsb S_{m-1}(x'_{m-1})}
      {S_1(x_0)\dotsb S_m(x_{m-1}) }$

\medskip

Particular values $S_m(x_{m-1}), S_m(x'_{m})$ are found from \eqref{Sm2}:

at $x_n=x_{m-1}$, only $P_m(x_{n+1})\neq 0$, and we work $P_{m-1,0,2}(x_{m-1})$ in \eqref{Sm2} with \eqref{PP}:
\begin{equation}\label{Smm1}
S_m(x_{m-1})=\dfrac{\ \ \alpha_{m-1}\frac{(x_{m}-x_{0})\dotsb(x_m-x_{m-1}) }{(x_m-x'_{0})\dotsb (x_m-x'_{m-1})}}
  {   \frac{X_m}{x_{m-1}-x'_1}   \frac{(x_m-x_0)\dotsb(x_m-x_{m-2})}{(x_m-x'_2)\dotsb(x_m-x'_m)}
}=
\dfrac{ \alpha_{m-1}(x_m-x_{m-1}) (x_m-x'_{m})(x_{m-1}-x'_1)}{(x_m-x'_{0})(x_m-x'_{1})X_m},
\end{equation}
with $X_m=\dfrac{R_{0,2}(y_{m-1})   +\rho_{m-1,0,2}(y_{m-1}-y_{-1})(y_{m-1}-y'_{1})-(x_{m}-x_{m-1})Y_2(y_{m-1})/2}
{R_{0,2}(y_{m-1})  +\rho_{m-1,0,2}(y_{m-1}-y_{-1})(y_{m-1}-y'_{1})+(x_{m}-x_{m-1})Y_2(y_{m-1})/2}$,

and when $x_n\to x'_m$,  $P_m(x_{n})$ and $P_m(x_{n+1})$ remain bounded, only the coefficient of $\gamma_n$ and the denominator of \eqref{Sm2}  $\to\infty$, as they are divided by $x_{n-1}-x'_{m-1}$ and $x_n-x'_{m}$.  We again work $P_{m-1,0,2}(x_{n})$  with \eqref{PP}:
\begin{equation}\label{Smp}
S_m(x'_{m})=\lim_{x_n\to x'_m}\dfrac{\ \ \gamma_{n}\frac{(x_{n-1}-x_{0})\dotsb(x_{n-1}-x_{m-1}) }{(x_{n-1}-x'_{0})\dotsb (x_{n-1}-x'_{m-1})}}
  {   \frac{X'_m}{x_{n}-x'_1}   \frac{(x_{n-1}-x_0)\dotsb(x_{n-1}-x_{m-2})}{(x_{n-1}-x'_2)\dotsb(x_{n-1}-x'_m)}
}=
\dfrac{ \gamma'_{m}(x'_{m-1}-x_{m-1})    (x'_{m-1}-x'_{m})   (x'_{m}-x'_1)}{(x'_{m-1}-x'_{0})(x'_{m-1}-x'_{1})X'_m},
\end{equation}
with

 \begin{multline}\label{gammapn}
\gamma'_{n}=
  \left( -\dfrac{\mu(x'_{n})}{(y'_{n}-y'_{n-1})X_2(x'_{n})}           -\dfrac{\nu(x'_{n})}{2}\right)      \dfrac{                (y'_{n-1}-y'_{-1})(y'_{n-1}-y'_0)}{ (x'_{n}-x'_{n-1})   (x'_{n}-x'_1)Y_2(y'_{n-1}) }\\
-\dfrac{\lambda}{2} \left[  \left(1+\dfrac{y'_0-y'_{-1} }{y'_{n}-y'_{n-1}}\right)\left(\dfrac{1}{2}+\dfrac{R(y'_{n-1}) }{(x'_{n}-x'_{n-1})Y_2(y'_{n-1})}\right)\right]    ,
\end{multline}

  $X'_m=\dfrac{R_{0,2}(y'_{m-1})   +\rho_{m-1,0,2}(y'_{m-1}-y_{-1})(y'_{m-1}-y'_{1})+(x'_{m}-x'_{m-1})Y_2(y'_{m-1})/2}
{R_{0,2}(y'_{m-1})  +\rho_{m-1,0,2}(y'_{m-1}-y_{-1})(y'_{m-1}-y'_{1})-(x'_{m}-x'_{m-1})Y_2(y'_{m-1})/2}$.

So, $ \zeta_m=\alpha_m \dfrac{S_{m}(x'_{m})}
      { S_{m+1}(x_{m}) }= \gamma'_{m}Z_m$
with

$Z_m=\dfrac{ (x'_{m-1}-x_{m-1})    (x'_{m-1}-x'_{m})   (x'_{m}-x'_1)    (x_{m+1}-x'_{0})(x_{m+1}-x'_{1})X_{m+1}}
{(x'_{m-1}-x'_{0})(x'_{m-1}-x'_{1})    (x_{m+1}-x_{m}) (x_{m+1}-x'_{m+1})(x_{m}-x'_1) X'_m }, m=3,4,\dotsc$

when $m=2$, an annoying $0/0$ undetermination appears, $\lim_{x_n\to x'_2} Z_n $ contains the ratio
$(x'_{n-1}-x'_1)$ over the denominator of   $X'_n$.

We get $S_2$ by \eqref{Sm2}
$S_2(x_n)=(x_n-x'_{1})\dotsb (x_n-x'_{2})  \\    \times
\dfrac{\ \ \alpha_n\dfrac{(x_{n+1}-x_0) (x_{n+1}-x_1)}{(x_{n+1}-x'_0)(x_{n+1}-x'_1)}+\beta_n \dfrac{(x_{n}-x_0) (x_{n}-x_1)}{(x_{n}-x'_0)(x_{n}-x'_1)}+\gamma_n \dfrac{(x_{n-1}-x_0) (x_{n-1}-x_1)}{(x_{n-1}-x'_0)(x_{n-1}-x'_1)}\ \  }{  x_n-x_{0} }$ evaluated at $x_1$ and $x_2$.

\qed

{\small
\begin{table}[htbp]
\begin{center}\begin{tabular}{|c|c|}\hline
    $m$ &  $\zeta_m/\alpha_m$        \\
\hline
 0 & -0.1667290    $\lambda$/( $\lambda$ - 0.05884553) \\
 1 & 0.5579333    ($\lambda$-0.07890984)/( $\lambda$ - 0.03433201) \\
 2 & 0.5302872    ($\lambda$-0.06204682)/( $\lambda$ - 0.03131567) \\
 3 & 0.2925826    ($\lambda$-0.05612122)/( $\lambda$ - 0.03569202) \\
 4 & -4.245072    ($\lambda$-0.06095382)/( $\lambda$ - 0.04168631) \\
 5 & 0.03446440    ($\lambda$+2.026877)/( $\lambda$ - 0.04464025) \\
 6 & 1.149683    ($\lambda$-0.02042389)/( $\lambda$ - 0.04062707) \\
 7 & 0.9825896    ($\lambda$-0.03707535)/( $\lambda$ - 0.01107954) \\
 8 & 25.31608    ($\lambda$-0.4342644)/( $\lambda$ - 0.1923718) \\
 9 & 0.3151368    ($\lambda$+0.1705402)/( $\lambda$ - 0.06416716) \\
 10 & 0.5542091    ($\lambda$-0.08160555)/( $\lambda$ - 0.03582462) \\
 11 & 0.5397617    ($\lambda$-0.06348380)/( $\lambda$ - 0.03111010) \\
 12 & 0.3501657    ($\lambda$-0.05635686)/( $\lambda$ - 0.03497262) \\
 13 & -2.088893    ($\lambda$-0.05933810)/( $\lambda$ - 0.04102451) \\
 14 & -0.6157532    ($\lambda$-0.1993437)/( $\lambda$ - 0.04456265) \\
 15 & 1.110735    ($\lambda$-0.01791228)/( $\lambda$ - 0.04173657) \\
 16 & 1.088867    ($\lambda$-0.03383604)/( $\lambda$ - 0.01885579) \\
 17 & 19.45439    ($\lambda$-0.2384268)/( $\lambda$ - 0.3065949) \\
 18 & 0.2071817    ($\lambda$+0.8567446)/( $\lambda$ - 0.07013919) \\
 19 & 0.5482573    ($\lambda$-0.08417506)/( $\lambda$ - 0.03768646) \\
 20 & 0.5471349    ($\lambda$-0.06510092)/( $\lambda$ - 0.03101770) \\
 21 & 0.3952189    ($\lambda$-0.05672590)/( $\lambda$ - 0.03428834) \\
 22 & -1.103745    ($\lambda$-0.05815415)/( $\lambda$ - 0.04032115)   \\   \hline
\end{tabular}
\caption{Ratio of coefficients of the expansion $\sum c_m\dfrac{(x-x_0)\dotsb(x-x_{m-1})}{(x-x'_0)\dotsb(x-x'_{m-1})}$ of a solution of $\mathcal{L}\phi=0$, on $x_0, x_1,x'_0,x'_1\dotsc$ given in \S \ref{labo}, and polynomials $\mu,\nu$ satisfying \eqref{munux0}. Some values $\phi(x_0), \phi(x_1),\dotsc$ are computed with (\ref{recurL}-\ref{Lmatrix}), keeping $\lambda$  literal. The coefficients $c_m$ are then obtained through a conventional divided difference method, and the ratios,
actually $c_1/((x_0-x'_1)c_0),     (x_{m-1}-x'_m)c_{m+1}/((x_m-x'_{m+1} )c_m)$, are
given, showing indeed that exact divisions produce rational functions of first degree.  }     \end{center}\label{secondotable}\end{table}
}

\bigskip

See here short sketches of possible further developments:

\subsubsection{Equation for first differences\label{eqfirst}}  \   \\
From \eqref{recurL} and \eqref{abg},

$\alpha_n(\phi(x_{n+1})-\phi(x_{n})) -\lambda \phi(x_n) -\gamma_n(\phi(x_{n})-\phi(x_{n-1}))=0$,

subtract from $n+1$:

$\alpha_{n+1}(\phi(x_{n+2})-\phi(x_{n+1}) -( \alpha_{n}+\lambda +\gamma_{n+1})    (\phi(x_{n+1})-\phi(x_{n})) +\gamma_{n}(\phi(x_{n})-\phi(x_{n-1})) =0$, a second order linear difference equation for $\mathcal{D}\phi(y)$,
similar to eq. (8.41) in Tsujimoto et al.\cite{Tsujim2024}. also eq. (4) of \cite{Spi2024}


\subsubsection{Biorthogonality of eigenfunctions\label{self2}}\     \\

$ \begin{bmatrix}\beta_0 & \alpha_0 &          &  \\
                    \gamma_1& \beta_1  & \alpha_1 &  \\
                            &\gamma_2  & \beta_2  & \alpha_2 \\
                            &          & \ddots   &  \end{bmatrix}
    \begin{bmatrix} \phi(x_0) \\  \phi(x_1)  \\  \phi(x_2) \\ \vdots \end{bmatrix}=0,$

Recall from \eqref{Lmatrix} that $\mathsf{A}\mathbf{\phi}+\lambda\mathsf{B}\mathbf{\phi} =0$, where $\mathbf{\phi}$ is the column
vector of $\phi(x_0), f(x_1), $ etc.

Consider now a row vector $\mathsf{g}=[g_0, g_1,\dotsc]$ such that $\mathsf{g}(\mathsf{A}+\lambda\mathsf{B})=0$, i.e.,
$g_1=-\dfrac{\beta_0}{\gamma_1}g_0, g_2=\dfrac{\beta_0\beta_1-\alpha_0\gamma_1}{\gamma_1\gamma_2}g_0,\dotsc$

Is $g_n$ the value of some function at $x_n$?

When $\lambda$ is such that $\zeta_r=0$, let $\lambda=\lambda_r$, the equations
\eqref{recurL} are solved by $\phi(x)=\phi_r(x)$, a rational function of
degree $r$.  
So, for any $\mathsf{g}$, $\mathsf{gA}\mathbf{\phi}_r+\lambda_r \mathsf{gB}\mathbf{\phi}_r=0$.



Let $\phi_r$ and $\phi_s$ be associated to $\lambda_r$ and $\lambda_s$:
$g_1^{(s)}=-\dfrac{\beta_0(\lambda_s)}{\gamma_1(\lambda_s)}$, etc.

$
g^{(s)}(\mathsf{A}+\lambda_r \mathsf{B})\phi_r =0, g^{(s)}(\mathsf{A}+\lambda_s \mathsf{B})\phi_r=0$,  so, $g^{(s)}\mathsf{B}f_r=0$,
if $\lambda_r\neq\kappa_s$, which is the biorthogonality relation
of left and right eigenvectors (Wilkinson \cite[eq. (3.6)]{WilkAlg}).

\subsubsection{}

From \eqref{LPk},

$ 0=\mathcal{L}\phi_r=\sum_0^r c_m(\lambda_r)\mathcal{L}P_m$

$=
    \begin{bmatrix}c_0(\lambda_r),\dotsc, c_r(\lambda_r), c_{r+1}(\lambda_r)=0, 0, \dotsc\end{bmatrix}
 \begin{bmatrix} \mathfrak{r}_0(\lambda_r)              &             &         &              \\
                  \mathfrak{s}_1(\lambda_r)  &\mathfrak{r}_1(\lambda_r)&         &              \\
                                             &\ddots                   &\ddots   &    \\
                                            &                          &\mathfrak{s}_r(\lambda_r)&\mathfrak{r}_r(\lambda_r)=0  \\
&&&\ddots \end{bmatrix}\begin{bmatrix}\widetilde{P}_0 \\  \widetilde{P}_1 \\ \vdots \\
\widetilde{P}_r \\ \vdots \end{bmatrix}
$

could be an interesting representation of $\mathcal{L}$.

\subsubsection{Riccati.}

A  scalar difference Riccati equation relates two lattice values through a first
degree rational equation.  This certainly holds for ratios of values of a solution
of a second order difference equation, in the form \eqref{recurL} and \eqref{Riccati02}

$r(y_n)=\dfrac{\phi(x_{n+1})}{\phi(x_n)}=-\beta_n/\alpha_n -\dfrac{\gamma_n/\alpha_n}{
   \ r(y_{n-1})=\dfrac{\phi(x_n)}{\phi(x_{n-1})}\   }$

\bigskip


\subsubsection{Recurrence relation of eigenfuctions.A very special case.}  \   \\

$\lambda=0$, $\nu(x)\equiv 0$, then $x_{0}$ and $x'_{1}$ are zeros of $\mu$, let
$\mu(x)=(x-x_{0})(x-x'_{1})(x-\xi)$.

 $\zeta_m=0$ if $\xi=x'_m$

\eqref{eqfirst} and $\lambda=0$,
$-\alpha_n(\phi(x_n)-\phi(x_{n-1}))  +\gamma_n(\phi(x_{n+1})-\phi(x_n))=0$,

$\phi_m(x_{n+1})-\phi_m(x_n)= \dfrac{\alpha_1(\xi_m)\dotsb\alpha_n(\xi_m) }{\gamma_1(\xi_m)\dotsb \gamma_n(\xi_m)}(\phi_m(x_1)-\phi_m(x_0))$

Emphasizing the dependence in $\xi$:
$\dfrac{c_{n+1}}{c_n}=\dfrac{\mathcal{C}_{n+1}}{\mathcal{C}_n}
\dfrac{x'_n-\xi}{x_n-\xi}$,

$\phi_m(x)=\sum \mathcal{C}_n \dfrac{(x'_m-x'_0)\dotsb(x'_m-x'_{n-1})}{(x'_m-x_0)\dotsb(x'_m-x_{n-1})}
\dfrac{(x-x_0)\dotsb(x-x_{n-1}) }{(x-x'_0)\dotsb (x-x'_{n-1})  }$.

\section{ Equivalence with other formulations: theta functions\label{theta}}  \   \\
\chead{\thesection  \ \ \ \   $\theta$ functions}

Elliptic hypergeometric functions are usually defined as
series of products of theta functions \cite{ellbib,Rosengren16,SpZ1}.

\begin{multline}\label{theta1}\theta(u)= 2\sum_0^\infty (-1)^m p^{(m+1/2)^2} \sin((2m+1)u)
 =2p^{1/4}\sin u \prod_{m=1}^\infty (1-2p^{2m}\cos 2u +p^{4m})(1-p^{2m}) \\
= C \sin u\prod_{m=1}^\infty (e^{2iu}-p^{2m})(e^{-2iu}-p^{2m})
= (-iC/2)(e^{iu}-e^{-iu})\prod_{m=1}^\infty (1-p^{2m} e^{-2iu})\prod_{m=1}^\infty (1-p^{2m}e^{2iu}),
\end{multline}
where $p$ is related to the modulus,
$C=2p^{1/4} \prod_{m=1}^\infty (1-p^{2m})$, and $\theta$ is actually Jacobi's $\theta_1$
function\footnote{The parameter $p$ is often written $q$ in the theta
literature, but we will need $q$ associated to the step.} \cite[chap. 20]{NIST}.

Other theta functions will not be used here, excepting $\theta_2(u)=\theta(u+\pi/2)
=\theta(\pi/2-u)$.

With $|p|<1$, the zeros of the entire function $\theta$ are the integer multiples
of $\pi$ and the roots of $\cos 2u =(p^{2m}+p^{-2m})/2= (e^{2m\pi i \tau}+e^{-2im\pi\tau})/2$ which are $u=\pm m\pi \tau+$ integer multiples  $\pi$, where $p=e^{i\pi\tau}$, and
Im $\tau >0$.   \cite[16.27- 16.38]{AbrS}  \cite[chap. 20]{NIST}  \cite[\S 13.19]{Erd2} (with $u=\pi v$)  \cite[Chap. XX]{WW}

The function $\theta$ is not doubly periodic, as

$\theta(u+\pi\tau)= \theta(u-i\log p)
= (-iC/2)(pe^{iu}-e^{-iu}/p)\prod_{m=1}^\infty(1-p^{2m-2} e^{-2iu})\prod_{m=1}^\infty (1-p^{2m+2}  e^{2iu})  \\  = (-iC/2)\dfrac{(pe^{iu}-e^{-iu}/p)(1-e^{-2iu})}{1-p^2 e^{2iu}}\prod_{m=1}^\infty (1-p^{2m}e^{-2iu})(1-p^{2m}e^{2iu})=-\theta(u)/(pe^{2iu})
$

\cite[22.2.4,5,6]{NIST}  $\sn(2Ku/\pi)= \dfrac{ \theta_3(0)\theta_1(u)}{\theta_2(0)\theta_4(u)}, \cn(2Ku/\pi)= \dfrac{ \theta_4(0)\theta_2(u)}{\theta_2(0)\theta_4(u)}, \dn(2Ku/\pi)= \dfrac{ \theta_4(0)\theta_3(u)}{\theta_3(0)\theta_4(u)},$

\cite[20.2.14,12,13]{NIST} $= i \dfrac{ \theta_3(0)\theta_1(u)}{\theta_2(0)e^{iu+i\pi\tau/4}\theta_1(u+\pi\tau/2)}= i \dfrac{ \theta_4(0)\theta_1(u+\pi/2)}{\theta_2(0)e^{iu+i\pi\tau/4}\theta_1(u+\pi\tau/2)}
= i \dfrac{ \theta_4(0)\theta_1(u+\pi/2+\pi\tau/2)}{\theta_3(0)\theta_1(u+\pi\tau/2)}$

Rosengren  uses $z$ with $\cos u=(z+z^{-1})/2$ \cite{Rosengren16}.

The theta functions realize a kind of the fundamental theorem of
algebra for elliptic functions, they take into account the zeros
and poles of such functions.
 \cite[p. 474]{WW}  an elliptic function of the variable $u$, with given zeros and poles in a fundamental parallelogram,
must be a product of ratios $\dfrac{ \theta(\pi(u- \text{zero})/\text{period})}{\theta(\pi(u- \text{pole})/\text{period})}$, where "period" is a period of the elliptic function,
and where $\tau= -i \log p/\pi$ is the ratio of two periods.

As seen on \S , $x_n$ is the simplest possible elliptic function of $u=nh$, it is a bivalent
function in any fundamental parallelogram. As we always encounter differences of
the form $x_n-x_m$, one of the zeros is of course $n=m$, we have a factor
$\theta(\pi(n-m)h/(4K))=\theta(\pi(n-m)/n_{\text{period}})$.
As $\sn(2K-u)=\sn(u)$ \cite[16.8]{AbrS} etc. $\sn(nh+g)=\sn(mh+g)$ at $nh+g=2K-mh-g$,
the other factor is $\theta(\pi((n+m)h+2g-2K)/(4K))=\theta(\pi(n+m)/n_{\text{period}}+\pi g/(2K)-\pi/2)$.

We take again the example of \S \ref{essential} , interpolate $n$ to see that $x_n=x_1$ occurs
at $n\approx 3.2$ which should be $2(K-g)/h-1=-5.895$ from $4K=6.9713, g=-0.1275$, and $h=-0.7641$, from which we subtract $n_{\text{period}}=4K/h=-9.1234$ to have 3.228.

Remark also that the extremum of $\xi_n$ (as well as $x_n$) occurs at
$n$ such that $nh+g=\pm K$, giving $n=2.114)$.  The values of $n$ such that
$x_n$ takes a given value are symmetrically placed with respect to 2.114.

{\small
\begin{verbatim}
  n       -1       0      1       2   2.114    3      3.228    4    4.228     5       6
x_n    0.6614      0  -0.6955 -0.9968  -1   -0.8067  -0.695 -0.1752   0    0.5380  0.9309
xi_n   0.5834 -0.1270 -0.7557 -0.9975  -1   -0.8470  -0.756 -0.2957 -0.127 0.4412  0.9117
=sn(nh+g)
\end{verbatim}
}

The poles of $x_n$ are the values of $n$ such that $\xi_n=\sn(nh+g)=\alpha/\gamma$,
from $\xi=(\alpha x+\beta)/(1+\gamma x)$,  so,  $n=(\text{arcsn}\,(\alpha/\gamma)-g)/h$ and $n= (2K-g-\text{arcsn}\,(\alpha/\gamma))/h$, one finds $x_n-x_m$ to be a constant
times
 $\dfrac{ \theta\left(\pi\dfrac{(n-m)h}{4K}\right)
\theta\left(\pi\dfrac{(n+m)h+2g}{4K}-\dfrac{\pi}{2}\right)   }
 { \theta\left(\pi\dfrac{nh-\mathfrak{a}+g}{4K}\right)
\theta\left(\pi\dfrac{nh+\mathfrak{a}+g}{4K}-\dfrac{\pi}{2}\right)  },
$
where $g=$ arcsn$(\xi_0=(\alpha x_0+\beta)/(1+\gamma x_0))$, and $\sn(\mathfrak{a})=\alpha/\gamma)$.  See also \cite[eq. 4.2]{SpZ2007}.

Let "$\mathfrak{t}$" the notation for $\mathfrak{t}(x,y)=\theta\left(\dfrac{\pi(x-y)}{4K}\right)\theta_2\left(\dfrac{\pi(x+y)}{4K}\right)$.

The constant  is independent of $n$ (it is always good to know what constants depend on),
but still depends on $m$, and
 what about the often encountered $x_n-x_{n-1}$? Here is a formula with a constant
independent of $n$ and $m$:

\begin{subequations}
\begin{equation}\label{xtheta}
 x_n-x_m =     C \dfrac{ \mathfrak{t}(nh+g,mh+g)   }
 { \mathfrak{t}(mh+g,\mathfrak{a})
\mathfrak{t}(nh+g,\mathfrak{a})  },
\end{equation}
where $x_n=\dfrac{\beta-\sn(nh+g)}{\gamma\,\sn(nh+g)-\alpha}, \sn(\mathfrak{a})=\alpha/\gamma$.

We get $C$ by looking at particular values of $m$ and $n$,
numerical checks giving indeed the same value $-6.55287 -2.71261 i$
for various $m$ and $n$ in  our example, where  $\mathfrak{a}=$ arcsn$(\alpha/\gamma)   = iK'+$ arcsn $(1/(k\alpha/\gamma))=
-0.2177+2.0106 i $
\cite[16.8 ]{AbrS}

The periods of the sn function are $4K$ and $2iK'$, so\footnote{Many formulas relating
elliptic functions to theta functions use $\tau=iK'/K$ instead of half
of this value. The author spent some painful weeks chasing mistakes.}  $\tau=i K'/(2K)$.
Here, $K=1.7428, K'=2.0106, \tau=0.5768 i, p=0.1633$.

The most convenient choice is
 $m\to n=-g/h$:

$\dfrac{dx_n}{dn}$ at $n=-g/h$ is $(\beta\gamma-\alpha)h/\alpha^2$ and must be
  $C \dfrac{ (\pi h/(4K))\theta'(0)\theta(-\pi/2)   }
 { \theta^2(-\pi\mathfrak{a}/(4K))
\theta^2(\pi\mathfrak{a})/(4K)-\pi/2) }$, so

\begin{equation}\label{Cxtheta}
C= \dfrac{4K(\alpha-\beta\gamma)  \theta^2(\pi\mathfrak{a}/(4K))
\theta^2(\pi\mathfrak{a}/(4K)-\pi/2)  }
      {\pi\alpha^2 \theta'(0)\theta(\pi/2) }.
\end{equation}

One has $\theta'(0)=2p^{1/4}\prod_1^\infty (1-p^{2m})^3,
\theta_2(0)=\theta(\pi/2)=2p^{1/4}\prod_1^\infty(1-p^{2m})(1+p^{2m})^2$.

\end{subequations}

\begin{equation}\label{xxptheta}
\dfrac{x-x_n}{x-x'_n}=\dfrac{
   \mathfrak{t}(\kappa h+g,nh+g)
\mathfrak{t}(nh'+g',\mathfrak{a})
} {
\mathfrak{t}(\kappa h+g,nh'+g')
\mathfrak{t}(nh+g,\mathfrak{a})  },
\end{equation}
where $x=x_\kappa= \dfrac{\beta-\sn(\kappa h+g)}{\gamma\,\sn(\kappa h+g)-\alpha}, \sn(\mathfrak{a})=\alpha/\gamma,
x'_n=\dfrac{\beta-\sn(nh'+g')}{\gamma\,\sn(nh'+g')-\alpha}$, and
where either $h'=h$ or $h'=-h$.

 For the $y$s, $y_n= \dfrac{\beta_y-\sn((n+1/2)h+g)}{\gamma_y\,\sn((n+1/2) h+g)-\alpha_y}, \sn(\mathfrak{b})=\alpha_y/\gamma_y$.

\begin{equation}\label{yyptheta}
\dfrac{y-y_n}{y-y'_n}=\dfrac{
   \mathfrak{t}(\kappa_y h+g,(n+1/2)h+g)
\mathfrak{t}((n+1/2)h'+g',\mathfrak{b})
} {
\mathfrak{t}(\kappa_y h+g,(n+1/2)h'+g')
\mathfrak{t}((n+1/2)h+g,\mathfrak{b})  },
\end{equation}

Here, $\alpha_y=1.4500, \beta_y=-0.1059, \gamma_y=-0.1396, \mathfrak{b}=-0.1643+2.0106 i,
h'=-h=0.7641$.

Check of \eqref{slopes} 

$\dfrac{C_{n+1,r,s}}{C_{n,r,s}}=
 \dfrac{ (y_{r-1}-y'_{s+n})(x_{r-1}-x_{r+n})}{(y_{r-1}-y_{r+n-1})(x_{r-1}-x'_{s+n})}=
\dfrac{ (x'_s-x_{r+n})(y'_{s-1}-y'_{s+n})}
          { (x'_s-x'_{s+n})(y'_{s-1}-y_{r+n-1})  }$  both found equal to

{\small
$\dfrac{\mathfrak{t}((r-1/2)h+g,(s+n+1/2)h'+g')
 \theta\left(\pi\dfrac{-(n+1)h}{4K}\right)
\mathfrak{t}((s+n)h'+g',\mathfrak{a})\mathfrak{t}((r+n-1/2)h+g,\mathfrak{b})
}
 { \mathfrak{t}((r-1)h+g,(s+n)h'+g')
\theta\left(\pi\dfrac{-nh}{4K}\right)
\mathfrak{t}((s+n+1/2)h'+g',\mathfrak{b} )
 \mathfrak{t}((r+n)h+g,\mathfrak{a}) }
$
}   

The elliptic hypergeometric expansion is a series of
products $[\theta(\pi h a_j)\dotsb \theta(\pi h (a_j+m-1)]^{\epsilon_j}$,
where $\epsilon_j=1$ or $-1$.

For the 'elliptic logarithm', of \S \ref{elllogthm}, $f(x)=\displaystyle \sum_0^\infty c_n\dfrac{(x-x_0)\dotsb(x-x_{n-1}) }{(x-x'_0)\dotsb (x-x'_{n-1})}$, where

$c_n=$
$\dfrac{y'_{n-1}-y_{n-1}}{C_{n,0,0} }
   \dfrac{(y_\rho-y'_{-1})\dotsb (y_\rho-y'_{n-2})  }{(y_\rho-y_0)\dotsb (y_\rho-y_{n-1})}$,
where $y_\rho=A$ or $c$ or $d$. The ratio of two term is

$\dfrac{y'_{n}-y_{n}}{y'_{n-1}-y_{n-1}}\, \dfrac{C_{{n},0,0} }{C_{n+1,0,0} }\,  \dfrac
    {y_\rho-y'_{n-1}}{y_\rho-y_{n}}\,   \dfrac{x_\kappa-x_{n}}{x_\kappa-x'_{n}}$,
$n=1,2,\dotsc$
reduces (!) to

$  \dfrac{ \mathfrak{t}((n+1/2)h'+g',(n+1/2)h+g)
 } { \mathfrak{t}((n-1/2)h'+g',(n-1/2)h+g)
  },
$

$\times\dfrac{ \mathfrak{t}((-1)h+g,(n)h'+g')
\theta\left(\pi\dfrac{-nh}{4K}\right)
}{\mathfrak{t}((-1/2)h+g,(n+1/2)h'+g')
 \theta\left(\pi\dfrac{-(n+1)h}{4K}\right)
}
$

$\times\dfrac{
   \mathfrak{t}((\rho+1/2) h+g,(n-1/2)h'+g')
} {
\mathfrak{t}((\rho+1/2) h+g,(n+1/2)h+g)
 }\;\dfrac{
   \mathfrak{t}(\kappa h+g,nh+g)    
} {
\mathfrak{t}(\kappa h+g,nh'+g')  
 },
$

from \eqref{xtheta}

Note that $\mathfrak{a}$ and $\mathfrak{b}$ have disappeared, we could have
used directly the canonical variables

$\xi_n=\sn(nh+g),  \eta_n=\sn((n+1/2)h+g)$ etc.

We have a rational expression of degree 9 in the thetas, so ${}_9E_8$,
considering that the \\  $\theta(\pi(n+1)/(4K))$ in the denominator comes from the 'factorial'
in the hypergeometric expansion.

the full term is $  \mathfrak{t}((n-1/2)h'+g',(n-1/2)h+g)
$ times the product from $m=0$ to $m=n-1$ of

\noindent{\small
$\dfrac{ \mathfrak{t}((-1)h+g,mh'+g')
\theta\left(\pi\dfrac{-mh}{4K}\right)
}{\mathfrak{t}((-1/2)h+g,(m+1/2)h'+g')
 \theta\left(\pi\dfrac{-(m+1)h}{4K}\right)
}
\dfrac{
   \mathfrak{t}((\rho+1/2) h+g,(m-1/2)h'+g')
} {
\mathfrak{t}((\rho+1/2) h+g,(m+1/2)h+g)
 }\;\dfrac{
   \mathfrak{t}(\kappa h+g,mh+g)    
} {
\mathfrak{t}(\kappa h+g,mh'+g')  
 },
$}

the isolated factor

 $\mathfrak{t}((n-1/2)h'+g',(n-1/2)h+g)  =
 \theta( (n-1/2)(h-h')+g-g')\theta(\pi/2-(n-1/2)(h+h')-g-g')$,
where $h'=h$ or $h'=-h$, so that only one  $\theta$ factor actually depends on $n$,
is an indication of possible \textit{very well  poised} character, according to
\cite[eq. 2.3]{SpZ2},
 \cite[eq. 7.18]{SpZ2007}.

\medskip


\section{Scholium.\label{scholium}}  \   \\
\chead{\thesection  \ \ \ \   scholium}

We saw in \S  \ref{elldefL} a second order elliptic\footnote{and any special case: differential, Hahn, $q-$theory, etc. \S\ref{someclassical}.}  difference equation solved by a hypergeometric expansion

\[   \phi(x)=\sum_0^\infty c_n \dfrac{(x-x_0)\dotsb(x-x_{n-1})}{(x-x'_0)\dotsb(x-x'_{n-1})},  \]
where $c_{n+1}/c_n$ is a rational function of $x_n, y_n, x'_n, y'_n$ and
some of their neighbors $x_{n\pm 1}$ etc.

For some particular values of a parameter $\lambda$, the expansion of $\phi=\phi_m$ is reduced to its $m$ first terms, so that $\phi_m$
is a rational function of poles $x'_0,\dotsc,x'_{m-1}$.

Do the $\phi_m$s enter some orthogonality or biorthogonality scheme? Things are not clear.

Do we have a recurrence relation connecting $\phi_m,  \phi_{m+1},  \phi_{m-1}$?  I found none.

\medskip

On the other hand,   we found in \S\ref{freepol} and \S\ref{singlex}  polynomials $\mathsf{A}_m$ satisfying recurrence relations, as well as biorthogonality conditions
\[   \sum_j   w_j \dfrac{\mathsf{A}_m(x'_j)}{(x'_j-x_0)(x'_j-x_2)\dotsb(x'_j-x_{2m})}\dfrac{\mathsf{B}_n(x'_j)}{(x'_j-x_1)(x'_j-x_3)\dotsb(x'_j-x_{2n+1})} =0, n\neq m,  \]
with the weights $w_j$ satisfying the discrete  Pearson's relations \eqref{Pear2}
\[\dfrac{w_{j+1}}{(y'_{j+1}-y'_j)X_2(x'_{j+1}) }
 \left( \dfrac{A(y'_j)}{x'_{j+1}-x'_j}-\dfrac{C(y'_j)}{2}\right) =
\dfrac{w_{j}}{(y'_{j}-y'_{j-1})X_2(x'_{j}) }
 \left( \dfrac{A(y'_j)}{x'_{j+1}-x'_j}+\dfrac{C(y'_j)}{2}\right)\],
with given polynomial $A$ and $C$.

We established for these $\mathsf{A}_m$s the existence of  second order elliptic  difference equations, as extension of Laguerre's theory
\S\ref{diff24}.

\medskip

Do these two families of  second order elliptic  difference equations have a non empty intersection? The matter will perhaps be settled in the future.

\section{Acknowledgments.}
\chead{\thesection   \ \ \ \ Acknowledgments.}

Many thanks  to the organizers of the
  3\`emes Journ\'ees Approximation
May 15-16, 2008, organized by B. Beckerman and C. Brezinski,
Universit\'e de Lille 1, France,

and the
workshop "Elliptic integrable systems, isomonodromy problems, and hypergeometric functions",
Hausdorff Center for Mathematics, Bonn, Germany, July 2008, organized by  V.~Spiridonov {\it et al.}

Many thanks to L. Haine, V.~Spiridonov and A.~Zhedanov

\section{References }

\chead{\thesection   \ \ \ \ References.}


\begin{thebibliography}{133}
\thispagestyle{fancy}

\bibitem{Abel} N.H. Abel, Recherches sur la s\'erie $1+ \frac{m}{1}x +\frac{m(m-1)}{1.2}x^2+\frac{m(m-1)m-2)}{1.2.3}x^3+ \dotsb$, \textsl{  J. reine angew. Math.} \textbf{1} 1826, 311-339
= \OE uvres, 2nd ed. Christiania 1881, \textbf{1} 219-250.

\bibitem{Abelrho} N.H. Abel, Sur l'int\'egration de la formule $\frac{\rho dx}{\sqrt{R}}$, $R$ et $\rho$ \'etant des fonctions enti\`eres, \textsl{J. reine angew. Math.} \textbf{1}, 1826  = \OE uvres, 2nd ed. Christiania 1881, \textbf{1} 104-144.

\bibitem{AbrS} M. Abramowitz,I.A. Stegun,
\textsl{Handbook of Mathematical Functions with Formulas, Graphs, and Mathematical Tables},
National Bureau of Standards Applied Mathematics Series, \textbf{55},
1964 $=$ Dover, New York, 1965 etc.
\href{http://members.fortunecity.com/aands/}
{\tt http://members.fortunecity.com/aands/},
\href{http://www.convertit.com/Go/ConvertIt/Reference/AMS55.ASP}
{\tt http://www.convertit.com/Go/ConvertIt/Reference/AMS55.ASP}\\
\href{http://www.math.sfu.ca/~cbm/aands/}{\tt http://www.math.sfu.ca/\~{}cbm/aands/}



\bibitem{Aitken} A.C. Aitken, Determinants And Matrices,
 Oliver and Boyd,
1944

\href{https://archive.org/download/in.ernet.dli.2015.205508}{\tt https://archive.org/download/in.ernet.dli.2015.205508}


\bibitem{Akhell} N.I.~Akhiezer, {\sl Elements of the Theory of Elliptic
      Functions\/}, 2nd ed., ``Nauka'', Moscow, 1970 (in Russian) =
      Translations of Math.\ Monographs {\bf 79\/}, Amer.\ Math.\
      Soc.\ , Providence, 1990.




\bibitem{AAR} G.E. Andrews , R. Askey , and R. Roy, \textsl{Special Functions},
   Cambridge University Press,  Encyclopedia of Mathematics and its Applications.  Volume \textbf{71}, Cambridge. 1999.


\bibitem{Appell} P. Appell and E. Goursat,
\textsl{Th\'eorie des fonctions alg\'ebriques et de leurs int\'egrales:
 \'etude des fonctions analytiques sur une surface de Riemann},
Gauthier-Villars 1895,

\href{http://www.archive.org/details/theoriefonctions00apperich}
{\texttt{http://www.archive.org/details/theoriefonctions00apperich}}

\bibitem{AW319} R.~Askey and J.~Wilson, Some basic hypergeometric
      orthogonal polynomials that generalize Jacobi polynomials,
      {\sl Memoirs Amer. Math. Soc.} {\bf 54} no. 319, 1985.

\bibitem{Baker}G.A.~Baker, Jr., {\sl Essentials of Pad\'e Approximants},
  Ac. Press, new York, 1975.



\bibitem{bang}
Gaspard Bangerezako,
\textsl{An Introduction to
$q-$Difference Equations},
Universit\'e catholique de Louvain,  Institut de
math\'ematique pure et appliqu\'ee, s\'eminaire de math\'ematique,
rapport $\text{n}^\circ$ 354,  Nov. 2008
\href{https://uclouvain.be/cps/ucl/doc/math/documents/RAPSEM354.pdf}{\tt https://uclouvain.be/cps/ucl/doc/math/documents/RAPSEM354.pdf}





\bibitem{BangRama} G. Bangerezako, The fourth order difference equation for the Laguerre-Hahn polynomials orthogonal on special non-uniform lattices, \textsl{The Ramanujan Journal : an international journal devoted to areas of mathematics influenced by Ramanujan}  Vol. \textbf{5}, p. 167-181 (2001).

\bibitem{Baxter} Rodney Baxter, \textsl{Exactly Solved Models in Statistical Mechanics}, Academic Press, 1982

You can view this book online in PDF format (17Mb) or in DjVu format (2.3Mb). The file is published with the permission of Academic Press and the author; it can be downloaded for private non-commercial use only.

\href{https://physics.anu.edu.au/theophys/baxter_book.php}{\tt https://physics.anu.edu.au/theophys/baxter\_book.php}


\bibitem{biorthKiran}
Kiran Kumar Behera, A. Swaminathan,  Biorthogonal rational functions of $R_{II}-$type, \textsl{Proc. Amer. Math. Soc.} {\bf 147} (2019), no.~7, 3061--3073.


 \bibitem{Belm92} S. Belmehdi, On semi-classical linear functionals of class $s=1$. Classification and integral representations,
\textsl{Indagationes Mathematic\ae}
Volume \textbf{3}, Issue 3, 1992, Pages 253-275.

\bibitem{BelRonv} S. Belmehdi, A. Ronveaux,  Laguerre-Freud's Equations for the Recurrence Coefficients of Semi-classical Orthogonal Polynomials,
     \textsl{Journal of Approximation Theory}, Volume \textbf{76}, Issue 3, March 1994, Pages 351-368.

\bibitem{BenMenahem} A. Ben-Menahem, \textsl{Historical Encyclopedia of Natural and Mathematical Sciences} vol. \textbf{1}, Springer 2009.


\bibitem{bliss} G.A. Bliss, Algebraic Functions, Amer. Math. Soc.
Colloq. Pub. \textbf{XVI}, 1933 = Dover  1966,

\href{https://archive.org/download/in.ernet.dli.2015.218603}{\tt https://archive.org/download/in.ernet.dli.2015.218603}



\bibitem{Boole} G. Boole, \textsl{
A treatise on the calculus of finite differences}. Reprint of the 1860 ed. 
Cambridge Library Collection - Mathematics. Cambridge: Cambridge University Press 
 (2009).

\bibitem{BourbakiHist} Nicolas Bourbaki, \'El\'ements d'histoire des math\'ematiques,
   Hermann 1974 $=$ Elements of the history of mathematics, Springer 1994.

\bibitem{Boyer}Carl B. Boyer, Uta C. Merzbach - \textsl{A history of mathematics} (1991, Wiley)

\bibitem{Branq} A. Branquinho, A. Foulqui\'e Moreno, A. Paiva, M.N. Rebocho,
Second-order differential equations in the Laguerre-Hahn class, \textsl{Applied Numerical Mathematics}, \textbf{94} (2015), 16-32.

\bibitem{BrezinskiH}  C. Brezinski, \textsl{  History of Continued Fractions and Pad\'e Approximants},
 Springer-Verlag, 1991.



\bibitem{Brezbio}  C.\ Brezinski,
{\sl Biorthogonality and its Applications to Numerical Analysis\/},
M.\ Dekker, 1992.




\bibitem{Brez2020} Claude Brezinski , Michela Redivo-Zaglia, \textsl{Extrapolation and Rational Approximation, The Works of the Main Contributors}, Springer 2020.

\bibitem{deBruin}  M. G. de Bruin,  Pad\'e approximation and its applications, book review,  \textsl{Acta Appl Math} \textbf{5}, 200-203 (1986).

\bibitem{BursZhed1} V. P. Burskii and A. S. Zhedanov,
The Dirichlet and the Poncelet problems,
Reports of RIAM Symposium No.16ME-S1
{\sl Physics and Mathematical Structures of Nonlinear Waves}
Proceedings of a symposium held at Chikushi Campus, Kyushu Universiy,
Kasuga, Fukuoka, Japan, November 15 - 17, 2004.

\href{http://www.riam.kyushu-u.ac.jp/fluid/meeting/16ME-S1/papers/Article_No_24.pdf}
{\tt http://www.riam.kyushu-u.ac.jp/fluid/meeting/16ME-S1/papers/Article\_No\_24.pdf}

\bibitem{BursZhed} V. P. Burskii  and A. S. Zhedanov,
Dirichlet and Neumann Problems for String Equation,
Poncelet Problem and Pell-Abel Equation,
{\sl Symmetry, Integrability and Geometry: Methods and Applications} Vol. {\bf 2} (2006), Paper 041, 5 pages
arXiv:math.AP/0604278 v1 12 April 2006
{\it SIGMA} {\bf 2} (2006), 041, 5 pages,

\href{http://www.emis.de/journals/SIGMA/2006/Paper041/}{http://www.emis.de/journals/SIGMA/2006/Paper041/}

\bibitem{BursZhed09}  Burskii V.P.  and  Zhedanov A.S.,
    On Dirichlet, Poncelet and Abel problems, arXiv:0903.2531 [math.CA], 47 pages, 2009.


\bibitem{Ce1}
    O. Salazar Celis, Adaptive Thiele interpolation, \textsl{ACM Communications in Computer Algebra}, \textbf{56}(3): 125--132, 2022.

\bibitem{Ce2}   O. Salazar Celis, Numerical continued fraction interpolation, \textsl{Ukrainian Mathematical Journal}, \textbf{76} 635--648, 2024.

\bibitem{chi}
   T.S.~Chihara,
   \textsl{An Introduction to Orthogonal Polynomials},
   Gordon and Breach, New York 1978.

\bibitem{AbelChu1} Wenchang Chu, Cangzhi Jia, Abel's method on summation by parts for elliptic hypergeometric series, \textsl{Communications in Contemporary Mathematics} Vol. \textbf{11}, 2009,  pp. 337-353.


\bibitem{AbelChu2} Wenchang Chu, Abel's method on summation by parts and balanced q-series identities,
   \textsl{Applicable Analysis and Discrete Mathematics} \textbf{4}, 2010, 54-65.

\bibitem{jjqell}
Ciet, M., Quisquater, J. J., \& Sica, F. , Compact elliptic curve representations, \textsl{Journal of Mathematical Cryptology}, \textbf{5}, 2011,  89-100, \href{https://doi.org/10.1515/JMC.2011.007}{\tt https://doi.org/10.1515/JMC.2011.007}


\bibitem{ClarksonFreud}
P A Clarkson and K Jordaan, "Properties of generalized Freud polynomials",  \textsl{Journal of Approximation Theory}, \textbf{225} (2018) pp 148-175.

\bibitem{pari} H. Cohen,  PARI/GP is a widely used computer algebra system designed for fast computations in number theory (factorizations, algebraic number theory, elliptic curves, modular forms, L functions...), but also contains a large number of other useful functions to compute with mathematical entities such as matrices, polynomials, power series, algebraic numbers etc., and a lot of transcendental functions. PARI is also available as a C library to allow for faster computations.
Originally developed by Henri Cohen and his co-workers (Universit\'e Bordeaux I, France), PARI is now under the GPL and maintained by Karim Belabas with the help of many volunteer contributors.
\href{https://pari.math.u-bordeaux.fr/}{\tt https://pari.math.u-bordeaux.fr/}

\bibitem{Coop} 
S. Clement Cooper, William B. Jones,  Arne Magnus,
General t-fraction expansions for ratios of hypergeometric functions,
\textsl{Applied Numerical Mathematics}, \textbf{4} (1988), 241-251.

\bibitem{Cretney} Rosanna Cretney, The origins of Euler's early work on continued fractions, \textsl{Historia Mathematica}
Vol. \textbf{41}, Issue 2, May 2014, Pages 139-156.

\bibitem{Cuyt1}Annie A.M. Cuyt, Lisa Jacobsen, Brigitte M. Verdonk, Instability and modification of Thiele interpolating continued fractions,
\textsl{Applied Numerical Mathematics},
Volume \textbf{4}, Issues 2–4, June 1988, Pages 253-262.


\bibitem{Cuyt} A. Cuyt, V. Petersen, B. Verdonk, H. Waadeland, W. B. Jones,
{\sl Handbook of Continued Fractions for Special Functions}, Springer 2008,
softcover reprint Dec.  2010.

\bibitem{Davis} P.J.~Davis, \textsl{Interpolation and Approximation},
  Blaisdell, Waltham, 1963 = Dover, New York, 1975.


\bibitem{Deaux} Roland Deaux, \textsl{Introduction to the geometry of complex numbers}. Transl. from the rev. French edition by Howard Eves,
Frederick Ungar Publishing Co. 208 p. (1957)=
 Dover, 2008.

\bibitem{Devroye} L. Devroye, Type design, typography, typefaces and fonts: An encyclopedic treatment of type design, typefaces and fonts. This site is also known as on snot and fonts.
\href{http://luc.devroye.org/fonts.html}{\tt http://luc.devroye.org/fonts.html}

\bibitem{Erd2} A. Erd\'elyi {\it et al.}, \textsl{Higher transcendental functions. Vol. II. Based on notes left by Harry Bateman},
McGraw-Hill Book Company  1953 .

\bibitem{Erdel} T. Erd\'elyi,  A.P. Magnus, P. Nevai, ``Generalized Jacobi
         weights, Christoffel functions and Jacobi polynomials,''
        \textsl{SIAM J. Math. An.} \textbf{25} (1994),     602-614, 1461.


\bibitem{ACm18} R. W. Floyd, Algorithm 18: Rational Interpolation by Continued Fractions,
\textsl{Communications of the ACM}, Volume \textbf{3}, Page 508
https://doi.org/10.1145/367390.3639419 (1960).


\bibitem{RonvFactor} M. Foupouagnigni, W. Koepf,  A. Ronveaux,
Factorization of fourth-order differential equations for perturbed classical orthogonal polynomials
\textsl{Journal of Computational and Applied Mathematics} Vol. \textbf{162}  2004,  299-326.

\bibitem{GammelNuttall}
Gammel, J.L., Nuttall, J.  Note on generalized jacobi polynomials. In: Chudnovsky, D.V., Chudnovsky, G.V. (eds) \textsl{The Riemann Problem, Complete Integrability and Arithmetic Applications. Lecture Notes in Mathematics}, vol \textbf{925}. Springer, Berlin, Heidelberg. 1982  https://doi.org/10.1007/BFb0093514

\bibitem{discri} I. M. Gelfand, M. M. Kapranov, and A. V. Zelevinsky,
\textsl{Discriminants, Resultants, and Multidimensional Determinants},
Birkh\"auser 2008.

\bibitem{Gon5} A.A.~Gonchar, G.~L\'opez, On Markov's theorem for
  multipoint Pad\'e approximation,
   \textsl{Math. USSR Sbornik} \textbf{34} (1978) 449-459.

\bibitem{DiscreteCalculus}
Leo J. Grady, Jonathan R. Polimeni,  \textsl{Discrete Calculus. Applied Analysis on Graphs for Computational Science}, Springer 2010 (softcover  2014),  \href{https://doi.org/10.1007/978-1-84996-290-2}{\tt https://doi.org/10.1007/978-1-84996-290-2}

\bibitem{GravesMorris}
P. R. Graves-Morris,
 Practical, Reliable, Rational Interpolation,
\textsl{IMA Journal of Applied Mathematics}, Volume \textbf{25}, Issue 3, March 1980, Pages 267-286, \href{https://doi.org/10.1093/imamat/25.3.267}{\tt https://doi.org/10.1093/imamat/25.3.267}.

\bibitem{GrH} F.A. Gr\"unbaum and L. Haine,
   On a $q$-analogue of Gauss equation and some $q$-Riccati equations, {\it in}
    Ismail, Mourad E. H. (ed.) et al., \textsl{Special functions, $q$-series
   and related topics}. Providence, RI: American Mathematical Society,
   \textsl{Fields Inst. Commun}. \textbf{14}, 77-81 (1997).

\bibitem{Gutchn} M.H. Gutknecht,
The rational interpolation problem revisited,
\textsl{Rocky Mountain J. Math.}, Volume \textbf{21} (1991) no. 4, p. 263-280.

\bibitem{Hahn1949} W. Hahn,  \"Uber Orthogonalpolynome, die $q$-Differenzengleichungen gen\"ugen, \textsl{Mathematische Nachrichten}, vol. \textbf{2}, 4-34;
Berichtigungen, \textsl{ibid.} \textbf{2}, 379-379.

\bibitem{Hahn1952} W. Hahn,  \"Uber lineare Differentialgleichungen, deren L\"osungen einer Rekursionsformel gen\"ugen, II. \textsl{Math. Nachr.} \textbf{7} (1952), 85-104.

\bibitem{HahnE} W.Hahn, On differential equations for orthogonal
      polynomials, {\sl Funk. Ekvacioj\/}, {\bf 21\/} (1978) 1-9.

\bibitem{Hahnmonh}  W. Hahn, \"Uber Differentialgleichungen f\"ur Orthogonalpolynome,
\textsl{Monatshefte f\"ur Mathematik},
December 1983, Volume \textbf{95}, Issue 4, pp 269-274.


\bibitem{Haine} L. Haine, P. Iliev,
The bispectral property of a $q-$deformation of the Schur polynomials and the $q-$KdV hierarchy,
\textsl{J. Phys. A: Math. Gen.} \textbf{30} (1997), 7217-7227.

\bibitem{Haine2026} L. Haine, Krall-type orthogonal polynomials and integrable isomonodromic deformations,
    arXiv:2603.01132 [math.CA]
    https://doi.org/10.48550/arXiv.2603.01132

\bibitem{Hal} G.H.~Halphen, \textsl{Trait\'e des fonctions elliptiques
  et de leurs applications.  \textbf{II},  Applications \`a la m\'ecanique, \`a la physique, \`a la g\'eod\'esie, \`a la g\'eom\'etrie et au calcul int\'egral. } Gauthier-Villars, Paris, 1888.

\href{http://catalogue.bnf.fr/ark:/12148/cb30571604z}{\tt http://catalogue.bnf.fr/ark:/12148/cb30571604z}

\href{https://archive.org/details/traitedesfonctio02halprich}{\tt https://archive.org/details/traitedesfonctio02halprich}

\bibitem{HaydNex} R Haydock and C M M Nex, A general terminator for the recursion method,
  \textsl{J. Phys. C: Solid State Phys.} \textbf{18} (1985) 2235-2248.

\bibitem{HazeP} M. Hazewinkel, editor: \textsl{Encyclopaedia of Mathematics},
Springer-Verlag, 2002, \href{http://eom.springer.de}{\textsl{http://eom.springer.de}},
  entry ``Classical orthogonal polynomials'',
\href{http://eom.springer.de/C/c022420.htm}{\texttt{http://eom.springer.de/C/c022420.htm}}


\bibitem{HvR2}  E. Hendriksen, H. van Rossum, Semi-classical orthogonal
     polynomials,
                                       pp. 354-361 {\it in\/}
      {\sl Polyn\^omes Orthogonaux et Applications, Proceedings,
      Bar-le-Duc 1984\/}, (C.Brezinski \& al., editors),
      {\sl Lecture Notes Math.\/} {\bf 1171\/}, Springer, Berlin 1985.





\bibitem{Hild1931}
Emanuel Henry Hildebrandt,
Systems of Polynomials Connected with the Charlier Expansions and the Pearson Differential and Difference Equations,
  \textsl{Ann. Math. Statist.}
    Volume \textbf{2}, Number 4 (1931), 379-439.


\bibitem{Houzel} C. Houzel, The Work of Niels Henrik Abel, p. 21-177   in O.A. Laudal, R. Piene, editors, \textsl{The Legacy of Niels Henrik Abel}, The Abel Bicentennial, Oslo, 2002, Springer 2004.

 \bibitem{iser1} A.~Iserles,  Rational interpolation to exp(-x) with application to certain stiff
     systems,  \textsl{SIAM J. Numer. Anal.} \textbf{18} (1981), 1-12.



\bibitem{Ism2005} M.E.H. Ismail,  {\sl Classical and quantum orthogonal polynomials in one variable,  with two chapters by Walter Van Assche},  Cambridge University Press,
Encyclopedia of Mathematics and its Applications.  Volume \textbf{98}, Cambridge. 2005.


\bibitem{IsMa} M.E.H. Ismail and D.R. Masson,
Generalized orthogonality and continued fractions.
    {\sl J. Approx. Theory} \textbf{83}, 1-40 (1995).

\bibitem{IsWimp} M.E.H. Ismail, J. Wimp,
On differential equations for orthogonal polynomials,
\textsl{Methods and Applications of Analysis}
\textbf{5} (1998) 439-452


\bibitem{Jackson}
 F. H. Jackson, $q-$Difference Equations,
\textsl{American Journal of Mathematics}, Volume \textbf{32}, 1910, 305-314.  
\href{https://archive.org/download/jstor-2370183}{\tt https://archive.org/download/jstor-2370183}

\bibitem{Jbeli}
S. Jbeli,  Description of the symmetric $H_q$-Laguerre-Hahn orthogonal q-polynomials of class one. \textsl{Period Math Hung} (2024). \href{https://doi.org/10.1007/s10998-024-00574-5}{https://doi.org/10.1007/s10998-024-00574-5}

\bibitem{JonesT1} William B. Jones,
Olav Nj\aa stad, W.J. Thron,
Two-point Pad\'e expansions for a family of
analytic functions, \textsl{Journal of Computational and Applied Mathematics} \textbf{ 9}   (1983) 105-123.


\bibitem{JonesT} William B. Jones , W. J. Thron
Continued Fractions
Analytic Theory and Applications
Cambridge University Press
 1984


\bibitem{kenfack}  Maurice Kenfack, Nangho and Kerstin Jordaan,
 Structure Relations of Classical Orthogonal Polynomials in the Quadratic and q-Quadratic Variable,
https://www.emis.de/journals/SIGMA/2018/126/
\textsl{SIGMA} \textbf{14} (2018), 126, 26 pages

\bibitem{Khov}  A. N. Khovanskii, \textsl{ The Applications Of Continued Fractions And Their Generalizations To Problems In Approximation Theory},
    Noordhoff, Ltd.  1963.

\href{https://archive.org/details/the-applications-of-continued-fractions-and-their-generalizations-to-problems-in}{\tt https://archive.org/details/the-applications-of-continued-fractions-and-their-generalizations-to-problems-in}

\bibitem{Kline} M. Kline, \textsl{Mathematical Thought from Ancient to Modern Times},
Oxford U.P. 1972.

\bibitem{Knuth}Donald E. Knuth, Art of Computer Programming, Volume 2: Seminumerical Algorithms, 3rd edition,
 Addison-Wesley Professional, 
1998


\bibitem{Koe98} Roelof Koekoek and Ren\'e  F. Swarttouw: The Askey-scheme of hypergeometric orthogonal polynomials and its q-analogue. Delft University of Technology, Faculty of Information Technology and Systems, Department of Technical Mathematics and Informatics, Report no. 98-17, 1998.
https://tudelft.roelofkoekoek.nl/askey.html

\bibitem{Koe} Roelof Koekoek,  Peter A. Lesky, Ren\'e F. Swarttouw,
\textsl{Hypergeometric
Orthogonal Polynomials
and Their $q-$ Analogues},
Springer 2010.    See also the online version in \href{http://aw.twi.tudelft.nl/~koekoek/askey.html}{\texttt{http://aw.twi.tudelft.nl/\~{}koekoek/askey.html}}

\bibitem{Kru} Sergey Khrushchev, Orthogonal Polynomials and Continued Fractions.
From Euler's Point of View
Cambridge U.P.,
 Encyclopedia of Mathematics and its Applications (No. 122), 2008.

\bibitem{Ko} T.H.~Koornwinder,
  {\sl Compact quantum groups and $q$-special functions},
  in {\sl Representations of Lie groups and quantum groups},
  V. Baldoni \& M. A. Picardello (eds.),
  Pitman Research Notes in Mathematics Series 311,
  Longman Scientific \& Technical, 1994,
pp. 46--128.
  See also the part  ``q-Special functions, a tutorial'', in\\
    \href{ftp://unvie6.un.or.at/siam/opsf_new/koornwinder/koornwinder3/koornwinder3.tex}{\texttt{ftp://unvie6.un.or.at/siam/opsf\_new/koornwinder/koornwinder3/koornwinder3.tex}}
   and \emph{ibid.}\href{ftp://unvie6.un.or.at/siam/opsf_new/koornwinder/wilsonqhahn.ps}
{\texttt{/wilsonqhahn.ps}}

\bibitem{Koornstru1}Tom H. Koornwinder, Lowering and raising operators for some special orthogonal polynomials, \textsl{Contemporary Math.} \textbf{ 417},  Amer. Math. Soc., 2006, pp. 227-238.


\bibitem{Koornstru2}Tom H. Koornwinder, The structure relation for Askey-Wilson polynomials, \textsl{J. Comput. Appl. Math.}  (2007), 214-226.


\bibitem{koorn2022} Tom H. Koornwinder
Charting the q-Askey scheme. II. The q-Zhedanov scheme
arXiv:2209.07995v1 [math.CA] 16 Sep 2022

\bibitem{SpiriConvergence2023}
Krotkov, Danil I.; Spiridonov, Vyacheslav P.
Infinite elliptic hypergeometric series: convergence and difference equations. (English. Russian original) Zbl 1543.30007
Sb. Math. 214, No. 12, 1751-1778 (2023); translation from Mat. Sb. 214, No. 12, 106-134 (2023).


\bibitem{Lag1885}  E.~Laguerre, Sur la r\'eduction en fractions continues d'une
      fraction qui satisfait \`a une \'equation diff\'erentielle lin\'eaire
      du premier ordre dont les coefficients sont rationnels,
      {\sl J. Math. Pures Appl. (S\'erie 4)\/} {\bf 1\/} (1885), 135-165 =
      pp. 685-711 {\it in\/}
      {\sl \OE uvres\/}, Vol.II, Chelsea, New-York 1972.

see
\href{http://math-tenere.ujf-grenoble.fr/JMPA/feuilleter.php?id=JMPA_1885_4_1}{\texttt{http://math-tenere.ujf-grenoble.fr/JMPA/feuilleter.php?id=JMPA\_1885\_4\_1}}

\href{http://portail.mathdoc.fr/cgi-bin/oetoc?id=OE_LAGUERRE__2}{\texttt{http://portail.mathdoc.fr/cgi-bin/oetoc?id=OE\_LAGUERRE\_\_2}}


\href{https://archive.org/download/oeuvresdelaguer02fragoog}
{\tt https://archive.org/download/oeuvresdelaguer02fragoog}

\bibitem{PellLasserre} J.B. Lasserre, Pell's equation, sum-of-squares and equilibrium measures of a compact set. \textsl{Comptes Rendus. Mathématique}, 2023, \textbf{361} (1), pp.935-952.


\bibitem{DarbouxBack}
S.B. Leble (originator), Darboux transformation, \textsl{Encyclopedia of Mathematics} - ISBN 1402006098. \href{http://encyclopediaofmath.org/index.php?title=Darboux_transformation&oldid=50455}{\tt
http://encyclopediaofmath.org/index.php?title=Darboux\_transformation\&oldid=50455}



\bibitem{Lo1}   G.~L\'opez Lagomasino,  Survey of multipoint Pad\'e
  approximation to Markov type meromorphic functions and asymptotic
  properties of the orthogonal polynomials generated by them,
  pp.~309-316 {\it in} \textsl{Polyn\^omes orthogonaux et applications,
  Proceedings, Bar-le-Duc 1984} (C.~Brezinski {\it et al.}, editors),
  \textsl{Springer Lecture Notes Math.} \textbf{1171}, 1985.

\bibitem{Lo2}   G.~L\'opez Lagomasino,  On the asymptotic of the ratio of
      orthogonal polynomials and convergence of multipoint Pad\'e
  approximants, \textsl{Math. USSR Sb.} \textbf{56} (1987) 207-219.



\bibitem{Maillard} M. Maillard and S. Boukraa, Modular invariance in lattice statistical mechanics,
\textsl{Annales de la Fondation Louis de Broglie},
Volume \textbf{26} num\'ero sp\'ecial, 2001,
287-328.
\href{http://www.ensmp.fr/aflb/AFLB-26j/aflb26jp287.pdf}{\texttt{http://www.ensmp.fr/aflb/AFLB-26j/aflb26jp287.pdf}}

\href{https://fondationlouisdebroglie.org/AFLB-26j/table26j.htm}{ \tt https://fondationlouisdebroglie.org/AFLB-26j/table26j.htm}



\bibitem{Mag84}  A.P.Magnus, Riccati acceleration of Jacobi continued fractions
and Laguerre-Hahn orthogonal polynomials,
pp 213-230 {\it in\/}:
H. Werner and H.J. B\"unger, editors: {\sl Pad\'e Approximation and its Applications,
Bad Honnef 1983\/},
Lecture Notes in Mathematics {\bf 1071\/},
Springer Verlag , Berlin , 1984.

\bibitem{Segovia} A.P.~Magnus, Associated Askey-Wilson polynomials as
   Laguerre-Hahn orthogonal polynomials, pp.~261-278 {\it in\/}:
   M.~Alfaro {\it et al.\/}, editors: {\sl Orthogonal Polynomials and
   their Applications, Proceedings, Segovia 1986}, Lecture
   Notes in Mathematics {\bf 1329\/}, Springer Verlag, Berlin, 1988.

\bibitem{magsnul} A.P.~Magnus,
    Special non uniform lattice ($snul$) orthogonal polynomials
    on discrete dense sets of points,
    \textsl{J. Comp. Appl. Math.} \textbf{65},  253-265 (1995).

\bibitem{MagnusCanterb} A.P. Magnus, Freud's equations for orthogonal polynomials
as discrete Painlev\'e
    equations, pp. 228-243 in
  \textsl{Symmetries and Integrability of Difference Equations},
        Edited by Peter A. Clarkson \& Frank W. Nijhoff,
   Cambridge U.P., \textsl{Lond. Math. Soc. Lect. Note Ser.} \textbf{255}, 1999.


\bibitem{MagLuminy} Alphonse P. Magnus,  Rational interpolation to solutions of Riccati difference equations on elliptic lattices, {\sl Journal of Computational and Applied Mathematics}, {\bf 233} (Issue 3, 1 Dec. 2009,  9th OPSFA Conference
) 793-801.


\bibitem{magsigma2009} A.P.~Magnus, Elliptic Hypergeometric Solutions
to Elliptic Difference Equations, SIGMA 5 (2009), 038, 12 pages,
\href{http://www.emis.de/journals/SIGMA/2009/038/}{\tt
\ http://www.emis.de/journals/SIGMA/2009/038/}

\bibitem{graphene} A.P.Magnus, Gaussian integration formulas for logarithmic
weights and application to 2-dimensional solid-state lattices,
\textsl{Journal of Approximation Theory}, \textbf{228}  (2018)  21-57.

\bibitem{MNR} Alphonse P. Magnus,
 Fran\c cois Ndayiragije,
Andr\'e Ronveaux, About families of orthogonal polynomials satisfying Heun's differential equation, \textsl{Journal of Approximation Theory}, \textbf{263} (2021) 105522,
\href{https://doi.org/10.1016/j.jat.2020.105522}{\tt https://doi.org/10.1016/j.jat.2020.105522}

\bibitem{Manas1} M. Ma\~nas, I. Fern\'andez-Irisarri, O.F. Gonz\'alez-Hern\'andez,  Pearson equations for discrete orthogonal polynomials: I. Generalized hypergeometric functions and Toda equations. \textsl{Stud Appl Math.} 2022; \textbf{148}, 1141–1179. https://doi.org/10.1111/sapm.12471

\bibitem{Manas3} Manuel Ma\~nas. Pearson equations for discrete orthogonal polynomials: III. Christoffel and Geronimus transformations. \textsl{Revista de la Real Academia de Ciencias Exactas, Físicas y Naturales. Serie A. Matemáticas} \textbf{116} (2022), article number 128 . doi:10.1007/s13398-022-01296-4.

\bibitem{MarETNA} F. Marcell\'an and  J.C. Medem,  $q-$Classical orthogonal polynomials:
a very classical approach,  \textsl{Electronic Transactions on Numerical Analysis},.
\textbf{ 9} , 1999,  pp. 112-127.\\
\href{http://www.maths.tcd.ie/EMIS/journals/ETNA/vol.9.1999/pp112-127.dir/pp112-127.pdf}{\texttt{http://www.maths.tcd.ie/EMIS/journals/ETNA/vol.9.1999/pp112-127.dir/pp112-127.pdf}}

\bibitem{Marstruc}Francisco Marcell\'an,  Ridha Sfaxi,
Second structure relation for semiclassical orthogonal polynomials,
\textsl{Journal of Computational and Applied Mathematics}, Vol. \textbf{ 200}, Issue 2, 15 March 2007, Pages 537-554.

\bibitem{Nuttall25} Andrei Mart\' inez-Finkelshtein, E. A. Rakhmanov, S. P. Suetin:
Heine, Hilbert, Pad\'e, Riemann, and Stieltjes: John Nuttall's work 25 years later,
\textsl{Contemporary Mathematics}, 2012, Vol. \textbf{578}, Pages 165-193,
DOI: https://doi.org/10.1090/conm/578/11474

\bibitem{Masson}
D. R. Masson, The Last of the Hypergeometric
Continued Fractions , {\it in} \textsl{Mathematical Analysis, Wavelets,
and Signal Processing
An International Conference on Mathematical Analysis
and Signal Processing
January 3-9, 1994
Cairo University, Cairo, Egypt},
M. E. H. Ismail \& {\it al.}, editors, \textsl{Contemporary Mathematics}
Volume \textbf{190}, 1995,
p. 287--294.

\bibitem{Mboustru} D. Mbouna,  A Structure Relation for Some Orthogonal Polynomials. \textsl{Mediterr. J. Math.} \textbf{20}, 237 (2023). https://doi.org/10.1007/s00009-023-02432-z

\bibitem{McCabe} John McCabe,
Some remarks on a result of Laguerre concerning continued fraction solutions of first order linear differential equations
    pp 373-379 {\it in}
   C.Brezinski \textit{et al.}, editors:
\textsl{Polyn\^omes Orthogonaux et Applications
Bar-le-Duc 1984},
Lecture Notes in Mathematics \textbf{1171}, Springer, Berlin 1985.

\bibitem{McCabeM} J.H. McCabe, J.A. Murphy,
Continued Fractions which Correspond to Power Series Expansions at Two Points,
\textsl{IMA Journal of Applied Mathematics}, Volume \textbf{17}, Issue 2, April 1976, Pages 233-247.


\bibitem{MCa} J. Meinguet, On the solubility of the Cauchy interpolation problem,
A. Talbot (Ed.), \textsl{Proceedings of the University of Lancaster Symposium on Approximation Theory and Its Applications}, Academic P. (1970), pp. 137-164.

\bibitem{MilneT}  L. M. Milne-Thomson,  \textsl{The Calculus Of Finite Differences}, Macmillan And Company., Limited,  1933.\\
\href{http://www.archive.org/details/calculusoffinite032017mbp}{\texttt{http://www.archive.org/details/calculusoffinite032017mbp}}

\bibitem{Nevai1984} P. Nevai, Asymptotics for orthogonal polynomials associated with $\exp(-x^4)$, \textsl{ SIAM J. Math. Anal.}, \textbf{15} (1984), 1177-1187.


\bibitem{NS} A.F.~Nikiforov and S.K.~Suslov, Classical orthogonal
   polynomials of a discrete variable on nonuniform lattices,
   \textsl{Letters in Math.\ Phys.} \textbf{11} (1986) 27-34.

\bibitem{NSU} A.F.~Nikiforov, S.K.~Suslov, and V.B.~Uvarov, {\sl Classical
    Orthogonal Polynomials of a Discrete Variable\/}, Springer, Berlin,
    1991.


\bibitem{Norlund} N.E. N\"orlund,
Vorlesungen \"uber Differenzenrechnung. (German)
Die Grundlehren der mathematischen Wissenschaften in Einzeldarstellungen, Bd. 13. Berlin: J. Springer.  (1924).


\bibitem{Nuttall}  J Nuttall, Asymptotics of diagonal Hermite-Pad\'e polynomials
    Journal of Approximation Theory, Volume 42, Issue 4, December 1984, Pages 299-386


\bibitem{NuttallS}   J Nuttall, S. R Singh
Orthogonal polynomials and Pad\'e approximants associated with a system of arcs
    Journal of Approximation Theory, Volume 21, Issue 1, September 1977, Pages 1-42





\bibitem{NIST} F. W. J. Olver, D. W. Lozier, R. F. Boisvert, C. W. Clark, editors:
  \textsl{NIST Handbook of Mathematical Functions},
NIST (National Institute of Standards and Technology) and Cambridge University Press,
 2010.
Also digital library of mathematical functions \href{http://dlmf.nist.gov}{\tt http://dlmf.nist.gov/}

\bibitem{Perron}  O.Perron, {\sl Die Lehre von den Kettenbr\"uchen\/},
 Teubner, Leipzig, 1913,

\href{https://archive.org/details/dielehrevondenk00perrgoog}{\tt https://archive.org/details/dielehrevondenk00perrgoog};
  $2^{\text{nd}}$ edition, 1929 = Chelsea,  3rd ed. 1954, 1957.

\bibitem{bispectrencyc}
Emma Previato (originator), Bispectrality, \textsl{Encyclopedia of Mathematics}, \\
\href{https://encyclopediaofmath.org/wiki/Bispectrality}{\tt https://encyclopediaofmath.org/wiki/Bispectrality}


\bibitem{Rahman81}Mizan Rahman, Families of biorthogonal rational functions in a discrete variable,  \textsl{SIAM J. Math. Anal.} \textbf{12} (1981), 355-367.

\bibitem{RahmanSuslov} Mizan Rahman and Sergei K. Suslov, Classical biorthogonal rational functions, pp. 131-146 {in}   Andrei A. Gonchar, Edward B. Saff, editors: \textsl{Methods of Approximation Theory in Complex Analysis and Mathematical Physics    Leningrad, May 13-24, 1991},  Lecture Notes in Math. , \textbf{1550}, Springer, 1993.

\bibitem{Rainsg}
Eric M. Rains, The (noncommutative) geometry of difference equations,
arXiv:2504.16187 [math.AG]
https://doi.org/10.48550/arXiv.2504.16187



\bibitem{Rebo} M.N. Rebocho,  On Laguerre-Hahn orthogonal polynomials on the real
line, \textsl{Random Matrices: Theory and Applications} Vol. \textbf{09}, No. 01, 2040001 (2020) DOI:  10.1142/S2010326320400018

\bibitem{Rebo2025}
 Maria das Neves Rebocho,  A note on Laguerre-Hahn orthogonal polynomials: from 1984 to 2024,
  \textsl{Numerical Algorithms},
  2025,
  https://doi.org/10.1007/s11075-025-02145-2


\bibitem{Ronv} A.~Ronveaux, private communication, October 2006.

\bibitem{Ronv4} A. Ronveaux, Fourth order differential equations and orthogonal polynomials of the Laguerre-Hahn class, in: C. Brezinski, L. Gori and A. Ronveaux, Eds., Orthogonal Polynomials and their Applications, \textsl{IMACS Ann. Comput. Appl. Math.}, \textbf{9}, Baltzer, Basel (1991), pp. 379-385.


\bibitem{ellbib} Hjalmar Rosengren. Bibliography of Elliptic Hypergeometric Functions

\href{http://www.math.chalmers.se/~hjalmar/bibliography.html}{\tt http://www.math.chalmers.se/\~{}hjalmar/bibliography.html}


\bibitem{RosengrenRahman} Hjalmar Rosengren,  Rahman's biorthogonal functions and superconformal indices, \textsl{Constr. Approx.} \textbf{47} (2018), 529-552.

\bibitem{Rosengren16}H. Rosengren, Elliptic hypergeometric functions, in H. S. Cohl and M. E. H. Ismail (eds.), \textsl{Lectures on Orthogonal Polynomials and Special Functions}, Cambridge University Press, 2020, pp. 213-279, =
arXiv:1608.06161 [math.CA]
\href{https://doi.org/10.48550/arXiv.1608.06161}{\tt  https://doi.org/10.48550/arXiv.1608.06161}

\bibitem{Shohat} J. Shohat,   A differential equation for orthogonal polynomials,
\textsl{Duke Math. J.}  \textbf{5}(2), 401-417 (June 1939).  Articl DOI: 10.1215/S0012-7094-39-00534-X

\bibitem{Silverman}
Joseph H. Silverman, \textsl{The Arithmetic of Elliptic Curves},  2nd edition,
Springer-Verlag New York 2009.


\bibitem{spir}
V. P. Spiridonov: classical elliptic hypergeometric
functions and their applications,
Lectures delivered at the workshop ``Elliptic integrable systems'' (RIMS, Kyoto, November 8-11, 2004).
253-287 \href{http://www.math.kobe-u.ac.jp/publications/rlm18/17.pdf}
{\tt http://www.math.kobe-u.ac.jp/publications/rlm18/17.pdf}



\bibitem{Spi2024} Vyacheslav P. Spiridonov,
"An extended scheme of classical special functions" , 
 arXiv:2403.16264v1 [math.CA] 24 Mar 2024.

 \bibitem{SpZ1} V.P Spiridonov and A.S.   Zhedanov,
Generalized eigenvalue problem and a new family of rational functions biorthogonal on elliptic grids, {\it in}
Bustoz, Joaquin (ed.) et {\it al}.,
{\sl Special functions 2000: current perspective and future directions. Proceedings of the NATO Advanced Study Institute, Tempe, AZ, USA, May 29--June 9, 2000}, Dordrecht: Kluwer Academic Publishers. NATO Sci. Ser. II, Math. Phys. Chem. 30, 365-388 (2001).

\bibitem{SpZ2}
  V.P. Spiridonov and A.S.   Zhedanov,
Classical biorthogonal rational functions on elliptic grids.
{\sl C. R. Math. Acad. Sci., Soc. R. Can.} \textbf{22}, No.2, 70-76 (2000).

\bibitem{SpZComm2000}
Spiridonov, V., Zhedanov, A. Spectral Transformation Chains and Some New Biorthogonal Rational Functions. \textsl{Comm Math Phys}
\textbf{210}, 49-83 (2000). \href{https://doi.org/10.1007/s002200050772}{\tt https://doi.org/10.1007/s002200050772}

\bibitem{SpZ2007}
V. P. Spiridonov  and A. S. Zhedanov: Elliptic grids, rational functions, and the Pad\'e interpolation
\textsl{The Ramanujan Journal} \textbf{13}, Numbers 1-3, June, 2007,
p.~285--310.

\bibitem{Steffensen} J. F.
Steffensen, Interpolation, Williams \& Wilkins, Baltimore, 1927.


\bibitem{Sz} G. Szeg\H o, \textsl{Orthogonal Polynomials}, Amer. Math. Soc., 1939,
$4^{\text{th}}$ edition, 1975.



\bibitem{Thiele} T. N. Thiele, \textsl{Interpolationsrechnung}, Teubner (1909).


\bibitem{TrueType} TrueType,   \href{https://en.wikipedia.org/wiki/TrueType}{\tt https://en.wikipedia.org/wiki/TrueType}

\bibitem{Tsujim}
Satoshi Tsujimoto, Luc Vinet,  Alexei Zhedanov,   "An algebraic description of the bispectrality of the biorthogonal rational functions of Hahn type", \textsl{Proceedings of the American Mathematical Society}, \textbf{149} (November 2020), pp.715-728.

\bibitem{Tsujim2024}
Satoshi Tsujimoto, Luc Vinet,  Alexei Zhedanov,   "Meta algebras and biorthogonal rational functions:
The Hahn case'',  \textsl{Forum of Mathematics, Sigma} (2025), Vol. \textbf{13:e139} 1-27
doi:10.1017/fms.2025.10040

\bibitem{VAssche2018}W. Van Assche, \textsl{ Orthogonal polynomials and Painlev\'e equations}, Cambridge University Press,  \textsl{ Australian Mathematical Society Lecture Series:} \textbf{27},  2018.

\bibitem{VdP}A. J. Van Der Poorten and C. S. Swart,
Recurrence Relations for Elliptic Sequences: Every Somos 4 is a
Somos k, \textsl{Bull. Lond. Math. Soc.} \textbf{38} (2006)
546-554.


\bibitem{verdestar} Luis Verde-Star
 A unified construction of all the hypergeometric and basic hypergeometric families of orthogonal polynomial sequences
Linear Algebra and its Applications 627 (2021) 242-274

\bibitem{newtbasis}Luc Vinet and Alexei Zhedanov,
Hypergeometric Orthogonal Polynomials
with respect to Newtonian Bases,
\textsl{Symmetry, Integrability and Geometry: Methods and Applications SIGMA }  \textbf{12}  (2016), 048, 14 pages

\bibitem{Wall}H. S. Wall,  {\sl Analytic Theory of Continued Fractions},
New York,
Van Nostrand, 1948
(= Chelsea Publishing, 1973).

\bibitem{Walsh1932}  J.L. Walsh, On Interpolation and Approximation by Rational Functions with Preassigned Poles,
\textsl{ Transactions of the American Mathematical Society}, Vol. \textbf{34}, No. 1 (Jan., 1932), pp. 22-74.

\bibitem{Walsh1935}  J.L. Walsh, \textsl{ Interpolation and
Approximation by  Rational Functions in the Complex
Domain}, Amer. Math. Soc. Colloquium Pub., vol. \textbf{XX},
New York City, 1935.

\bibitem{We1}
Werner, H., (1979), A reliable method for rational interpolation,  in \textsl{Pad\'e approximation and its applications}, ed. L. Wuytack, Springer-Verlag, 257-277,
\textsl{Lecture Notes in Mathematics}, vol \textbf{765}.

\bibitem{We2}
Werner, H., (1980). Ein Algorithmus zur rationalen Interpolation,  in \textsl{Numerische Methoden der Approximationtheorie}, Band \textbf{5}, eds. L. Collatz, G. Meinardus and H. Werner, Birkh\"auser, 319-337,
excerpts of The Conference on Numerical Methods of Approximation, March 18-24, 1979 at the Mathematical Research Institute Oberwolfach = ISNM - International Series of Numerical Mathematics - Vol. \textbf{52}.


\bibitem{WW}  E. T. Whittaker and G. N. Watson, \textsl{A Course of Modern Analysis },
4th ed, Cambridge U.P. 1927;  5th ed., with V.H. Moll,  2021.




\bibitem{WilkAlg} J.H. Wilkinson, \textsl{The Algebraic Eigenvalue Problem}, Oxford U.P. 1965.

\bibitem{Zannier} Umberto Zannier, Hyperelliptic Continued Fractions and Generalized Jacobians,
\textsl{American Journal of Mathematics},
Vol. \textbf{141}, No. 1 (2019), pp. 1-40 = math arXiv:1602.00934.

\bibitem{zhed99}  A.S.  Zhedanov,
Biorthogonal Rational Functions and the Generalized Eigenvalue
Problem, {\sl J.  Approx. Theory}, Volume {\bf 101}, Issue 2,
December 1999, Pages 303-329.

\bibitem{ZhedPad} A.S. Zhedanov,
Pad\'e interpolation table and biorthogonal rational functions. In: {\sl Proceedings of RIMS Workshop on Elliptic Integrable Systems. Kyoto, November 8-11} (2004) to be published

\href{http://www.math.kobe-u.ac.jp/publications/rlm18/20.pdf}{\texttt{
http://www.math.kobe-u.ac.jp/publications/rlm18/20.pdf}}


\bibitem{Zou}
Le Zou, Liang-Tu Song, Xiao-Feng Wang, Qian-Jing Huang, Yan-Ping Chen, Chao Tang \& Chen Zhang, Univariate Thiele Type Continued Fractions Rational Interpolation with Parameters, {\it in}: Huang, DS., Huang, ZK., Hussain, A. (eds) \textsl{Intelligent Computing Methodologies. ICIC 2019. Lecture Notes in Computer Science  (LNAI)}, vol \textbf{11645} (2019)  pp 399-410, Springer, Cham. \href{https://doi.org/10.1007/978-3-030-26766-7_37}{\tt https://doi.org/10.1007/978-3-030-26766-7\_37}.


\bibitem{x}
\href{
https://math.stackexchange.com/questions/1603172/digamma-function-in-expectation}{\tt https://math.stackexchange.com/questions/1603172/digamma-function-in-expectation}



\end{thebibliography}
\end{document}